\documentclass[a4paper,11pt]{amsart} 
\usepackage{amscd}             
\usepackage{xypic}  %commu
\usepackage{amssymb}
\usepackage{amsthm}
\usepackage{epsf} 
\usepackage[T1]{fontenc}            
\makeindex             % used for the subject index
\newtheorem{thm}{Theorem}[section]   
\newtheorem{cor}[thm]{Corollary}
\newtheorem{exa}[thm]{Example}
\newtheorem{Lemma}[thm]{Lemma}
\newtheorem{prop}[thm]{Proposition}

\newtheorem{defn}[thm]{Definition}

\newtheorem{rem}[thm]{Remark}
\newtheorem{conj}[thm]{Conjecture}
\renewcommand{\proofname}{Proof}

\def\Hilb{\operatorname{Hilb}}
\def\Sing{\operatorname{Sing}}
\def\ker{\operatorname{ker}}

\def\min{\operatorname{min}}
\def\im{\operatorname{im}}
\def\rank{\operatorname{rank}}
\def\max{\operatorname{max}}
\def\length{\operatorname{length}}
\def\c1{\operatorname{c_1}}
\def\c2{\operatorname{c_2}}
\def\Cliff{\operatorname{Cliff}}
\def\gon{\operatorname{gon}}

\def\Grass{\operatorname{Grass}}

\def\Sym{\operatorname{Sym}}
\def\CC{{\mathbf C}}
\def\ZZ{{\mathbf Z}}

\def\QQ{{\mathbf Q}}
\def\PP{{\mathbf P}}

\def\A{{\mathcal A}}
\def\B{{\mathcal B}}
\def\C{{\mathcal C}}
\def\L{{\mathcal L}}
\def\D{{\mathcal D}}
\def\R{{\mathcal R}}

\def\O{{\mathcal O}}
\def\I{{\mathcal J}}
\def\Z{{\mathcal Z}}
\def\E{{\mathcal E}}
\def\T{{\mathcal T}}
\def\H{{\mathcal H}}
\def\F{{\mathcal F}}
\def\K{{\mathcal K}}
\def\M{{\mathcal M}}
\def\x{\times}                   % product (fiber)
\def\v{^{\vee}}                  % dual
\def\iso{\simeq}                 
\def\eqv{\equiv}
\def\sub{\subseteq}

\def\sup{\supseteq}

\def\+{\oplus}                   % direct sum
\def\*{\otimes}                  % tensor product
\def\hpil{\longrightarrow}       % ----->
\def\khpil{\rightarrow}

\def\Aut{\operatorname{Aut}}

\def\Ext{\operatorname{Ext}}

\def\Pic{\operatorname{Pic}}

\def\disc{\operatorname{disc}}
\def\det{\operatorname{det}}
\def\hs{\hspace{.1in}}

\begin{document}

  \title{$K3$ projective models in Scrolls}
  \author{Trygve Johnsen and Andreas Leopold Knutsen}  
  
  \address{Dept. of Mathematics\\ 
    University of Bergen\\ Johs. Brunsgt 12\\ N-5008 Bergen\\ Norway}

\email{johnsen@mi.uib.no, andreask@mi.uib.no}

  \begin{abstract}
     We study the projective models of complex $K3$ surfaces polarized by a line bundle $L$ such that all smooth curves in $|L|$  
    have non-general Clifford index. Such models are in a natural way 
    contained in rational normal scrolls. 

   We use this study to classify and describe all projective models of
   $K3$ surfaces of genus $g \leq 10$. 
  \end{abstract}

\maketitle
\tableofcontents
\eject

\section{Introduction} \label{intro}

A pair $(S,L)$ of a $K3$ surface $S$ 
and a base point free line bundle $L$ with $L^2=2g-2$ will be called a polarized 
$K3$ surface of genus $g$. The sections of $L$ give
a map $\phi_L$ of $S$ to $\PP^g$, and the image is called a projective model
of $S$. When $\phi_L$ is birational, the 
image is a surface of degree $2g-2$ in $\PP^g$. 

The following is well known:
For $g=3$ the birational projective model is a quartic surface and for $g=4$ a complete intersection of a quadric and a cubic hypersurface.
For $g=5$ the general model is a complete intersection of three hyperquadrics.
For $6 \leq g \leq 10$ and $g=12$ it is shown by Mukai in \cite{Mu1} and \cite{Mu2} that 
the general projective model
is a complete intersection in certain
homogeneous varieties contained in projective spaces of larger dimension than
$g$. 
%More precisely, in \cite{Mu2} a polarized $K3$ surface $(S,L)$ is
%said to 
%be {\it Brill-Noether (BN) general} if for all non-trivial 
%decompositions $L \sim M+N$ one has
%\[h^0(M) h^0(N) < h^0(L)=g+1.\] (This is for instance satisfied if
%any smooth curve $C \in |L|$ is Brill-Noether general.)
%
%This is clearly satisfied for the general $K3$ surface.
%It is moreover shown that this is a necessary and suffucient
%condition for the surface 
%to have a general projective model as above, for $2 \leq g \leq 10$
%and $g=12$ 
%(for $g=2$ and $g=3,4$ with $\phi_L$ birational, all surfaces are
%automatically BN general).
%Also note that a BN general $K3$ surface does 
%not necessarily have Picard number 
%one. 
%
%We will study a slightly different notion of generality for polarized $K3$ surfaces.

We will study the projective models of a particular kind of polarized
$K3$ surfaces. 
Recall the result of Green and Lazarsfeld \cite{gl}, which states
that if $L$ is a base point free line bundle on a $K3$ surface $S$, then
$\Cliff C$ is constant for all smooth irreducible $C \in |L|$, and if
$\Cliff C < \lfloor \frac{g-1}{2} \rfloor $, then there exists a line
bundle $M$ on $S$ such that $M_C := M \* \O _C$ computes the Clifford index of
$C$ for all smooth irreducible $C \in |L|$.

It therefore makes sense to define the Clifford index $\Cliff L$ of
a base point free line bundle, or the Clifford index $\Cliff _L (S)$
of a polarized $K3$ surface $(S,L)$ as the Clifford index of the
smooth curves in $|L|$. 

By the {\it Existence Theorem} \ref{exthm} below we have: For any pair of integers $(g,c)$ such that $g \geq 2$ and 
$0 \leq c \leq \lfloor  \frac{g-1}{2} \rfloor $, there exists an 
$18$-dimensional family of polarized $K3$ surfaces of genus $g$ and of 
Clifford index $c$.

In this paper we study and classify 
projective models of non-Clifford general polarized $K3$ surfaces (i.e. with $c < \lfloor \frac{g-1}{2} \rfloor $) of any genus
larger than two. 

The central point is that by the result of Green and Lazarsfeld, there
exists in these cases a divisor class $D$ on $S$ computing the Clifford
index of $L$. We will show that $0 \leq D^2 \leq c+2$. We can choose such a divisor that is base point free and such that the general member of $|D|$ is a smooth curve.
Such a divisor (class) will be called a {\it free Clifford divisor} for $L$.
(The definition only depends on the class of $D$.)

The images of the members of $|D|$ by  $\phi_L$
span sublinear spaces inside $\PP^g$. Each subpencil $\{ D _{\lambda}
\}$ within the complete
linear system $|D|$ then gives rise to a pencil of sublinear spaces. For each
fixed pencil the union of these spaces will be a rational normal scroll $\T$. 
We will investigate under which conditions this scroll is smooth.
These scrolls are the natural
ambient spaces for non-Clifford general $K3$ surfaces. 

In the cases $c=1$ and $2$ with $D^2=0$, the description of the
projective models is particularly nice, since they are then complete
intersections in their corresponding scrolls. 

If $\T$ is smooth, we will be able to find a resolution (up to certain
invariants) of 
$\phi_L(S)$ inside $\T$  for arbitrary $c$. 

If $\T$ is singular, we take the blow up $f: \tilde S \khpil S$
at the $D^2$ base points of the
pencil $\{ D _{\lambda}\}$ and show that the projective model $\phi_H
(\tilde S)$ of $\tilde S$ by the base point free line bundle $H :=
f^*L+f^*D-E$, where $E$ is the exceptional divisor, is contained in a
smooth rational normal scroll $\T_0$ which is a desingularization of $\T$. 
We find a resolution of $\phi_H (\tilde S)$
inside $\T_0$ and investigate the possibility of pushing
down this resolution to a resolution of $\phi_L(S)$ inside $\T$.

We also give a description of those projective models for $g \leq 10$
that are Clifford general, but still not general in the sense of Mukai
(i.e. they are not complete intersections in homogeneous spaces).
These models are also contained in scrolls, and can be analysed in a
similar manner. This can also be done for $g=12$, which we leave to
the reader.

Together with Mukai's results this will then give a complete picture
of the birational projective models for $g\leq 10$ (and $g=12$).
For $g=11$ and $g \geq 13$ 
our description of non-Clifford general projective models is not supplemented
by any description of general projective models at all. We hope, however, that
our description of the non-general models may have some interest in
themselves. 

%An additional motivation for classifying projective models for low $g$ is
%of course the well known fact that up to and including the case $g=9$ the 
%general canonical curve appears as a hyperplane section with some projective
%model of a $K3$ surface.

The paper is organized as follows. 

In Section \ref{sins} we recall some basic facts about
rational normal scrolls, and how to obtain such scrolls from surfaces
with pencils on them.
In Section \ref{cliff} we define the essential concept of a free Clifford 
divisor (Definition \ref{free}). In Section \ref{exi} we recall two important 
results from \cite{kn3},  the above mentioned  {\it Existence Theorem}, and 
a similar result concerning existence of curves with a 
prescribed gonality on $K3$ surfaces.

In Section \ref{sing} we study in detail the singular loci
of the projective model $\phi_L(S)$ and the scroll $\T$ in which we
choose to view this model as contained. We show (Theorem \ref{mainsing}) 
that we 
can always
find a free Clifford divisor $D$ such that the singular locus of
$\T$ is ``spanned'' by the images of the base points of the pencil $\{
D _{\lambda} \}$ and the contractions of smooth rational curves across
the members of the pencil. A free Clifford divisor with this extra property will be called a {\it perfect Clifford divisor} (Definition \ref{perfect}).
 We also 
include a study of the
projective model if $c=0$ (the hyperelliptic case), which is due to 
Saint-Donat \cite{S-D}. Some of the longer proofs of the results in this section are postponed until Section \ref{tedious}.

In Section \ref{resol} we study the resolution of
$\phi_L(S)$ inside its scroll $\T$ when $\T$ is smooth. In this case a general
hyperplane section of $\T$ is a scroll formed in a similar way from a
pencil computing the gonality on a canonical curve $C$ of genus g (the gonality is $c+2$).
Such scrolls were studied in \cite{Sc}, and our results (Lemma \ref{Betti0} and
Proposition \ref{resolv}) for $K3$ surfaces
in smooth scrolls are quite parallel to those of \cite{Sc}.

In Section \ref{singscrolls} we treat the case when the scroll $\T$ is 
singular. The approach is to study  the blow up $f: \tilde S \khpil S$
at the $D^2$ base points of the
pencil $\{D _{\lambda} \}$ and the projective model $S'':=\varphi_H
(\tilde S)$ of $\tilde S$ by the base point free line bundle $H :=
f^*L+f^*D-E$, where $E$ is the exceptional divisor. The pencil
$|f^*D-E|$ defines a smooth rational normal scroll $\T_0$ that
contains $S''$ and is a desingularization of $\T$.

We use Koszul cohomology and 
techniques inspired by Green and Lazarsfeld to compute some Betti-numbers of the 
$\varphi_L(D_{\lambda})$ and we obtain that they all have the same Betti-numbers for low values of $D^2$ and this is a necessary and sufficient condition for ``lifting'' the resolutions of the fibers to one of the surface 
$S''$ in $\T_0$. 
 We prove that $S''$ is normal, 
and use this to give more details about the resolution. We give conditions 
under which we can push down the resolution to one of $\varphi_L(S)$ in $\T$.
Here we use results from \cite{Sc}. We end the section by
investigating some examples for low genera.

In Section \ref{lowind} we study the set of projective models in smooth scrolls
for $c=1$, $2$ and $3$ ($< \lfloor \frac{g-1}{2} \rfloor$). 
We study the sets of projective models in 
$(c+2)$-dimensional scrolls of given types. Since the scroll type is
dependent on which rational curves that exist on $S$, and therefore
on the Picard lattice, it is natural that the dimension of the set of
models in question in a scroll as described, is dependent on the
scroll type. We study this interplay, and obtain a fairly clear picture
for $c=1$ and $2$. For $c=3$ we study a Pfaffian map of the resolution of 
$\phi_L(S)$ in the scroll. In Remark \ref{exctypes}  we predict the 
dimension of the set of projective $K3$ models inside a fixed smooth scroll 
of a given type, for arbitrary $c < \lfloor \frac{g-1}{2} \rfloor$. We state 
the special case $c=3$ as Conjecture \ref{c3dim}.  

In Section \ref{bncl} we study the projective models for $g \leq 10$
that are Clifford general but not general in the sense of Mukai
(i.e. they are not complete intersections in homogeneous spaces).
By the concrete description in \cite{Mu2} of such surfaces, it follows
that their projective models are also contained in scrolls. We analyse
them in a similar manner.

We conclude by giving a complete list and descripton of {\it all}
birational projective models of $K3$ surfaces for $g \leq 10$ 
(including both the general one in the sense of Mukai, and the remaining ones,
that we give a detailed classification of here)

\subsection{Acknowledgements} 
The authors thank Kristian Ranestad for suggesting the problem. 
We are also grateful to M.~Coppens, S.~Ischii, S.~Lekaus, R.~Piene, J.~Stevens, 
S.~A.~Strømme, 
B.~Toen and J.~E.~Vatne for useful conversations.

This paper was
written while the authors were visitors at the Department of Mathematics,
University of Utah, and at the Max-Planck-Institut f\"ur Mathematik, Bonn.
We thank both institutions for their hospitality.

The second author was supported by a grant from the Research Council of Norway.

\subsection{Notation and conventions}
\label{prel}

We use standard notation from algebraic geometry, as in \cite{Ha}.

The ground field is the field of complex numbers. All surfaces are 
reduced and irreducible \textit{algebraic surfaces}. 

By a $K3$ surface\index{$K3$ surface} is meant a smooth surface $S$ with trivial 
canonical bundle and such that $H^1 (\O_S)=0$. In particular $h^2 (\O_S)= 1$ 
and $\chi (\O_S)= 2$. \label{K3}

By a \textit{curve} is always meant a \textit{reduced} and
\textit{irreducible curve} (possibly singular). The \textit{adjunction formula}\index{adjunction formula}
for a curve $C$ on a surface $S$ reads $\O_C(C+K_S) \iso \omega_C$, where $\omega_C$ is the dualising sheaf of $C$, which is just the canonical bundle when $C$ is smooth. In particular, the arithmetic genus 
$p_a$ of $C$ is given by $C.(C+K_S)=2p_a-2$.

On a smooth 
surface we use line bundles and divisors, as well as the multiplicative and additive notation,  with little or no distinction. We denote by $\Pic S$\index{$\Pic S$} the 
\textit{Picard group}\index{Picard group} of $S$, i.e. the group of linear equivalence classes of line bundles on $S$. The 
\textit{Hodge index theorem}\index{Hodge index theorem} yields that if $H \in \Pic S$ with $H^2 >0$, then $D^2H^2 \leq (D.H)^2$ for any $D \in \Pic S$, with equality if and only if $(D.H)H \eqv H^2D$.

Linear equivalence of divisors is denoted by $\sim$, and numerical equivalence by $\eqv$. Note that on a $K3$ surface $S$ linear and numerical equivalence is the same, so that $\Pic S$ is torsion free. The usual intersection product of line bundles (or divisors) on surfaces therefore makes the Picard group of a $K3$ surface into a lattice, the \textit{Picard lattice}\index{Picard lattice} of $S$, which we also denote by  $\Pic S$.

For two divisors or line bundles $M$ and $N$ on a surface, we use the notation $M \geq N$ to mean $h^0(M-N) >0$ and $M > N$, if in addition $M-N$ is nontrivial.

If $L$ is any line bundle on a smooth surface, $L$ is said to be \textit{numerically effective}, or simply \textit{nef}\index{nef}, if $L.C \geq 0$ for all curves $C$ on $S$. In this case $L$ is said to be \textit{big}\index{big} if $L^2 >0$.

If $\F$ is any coherent sheaf on a variety $V$, we shall denote by 
$h^i (\F)$ the complex dimension of $H^i (V, \F)$, and by $\chi (\F)$\index{$\chi$} 
the \textit{Euler characteristic}\index{Euler characteristic} $\sum (-1)^i h^i (\F)$. In particular, if $D$ is any divisor on a normal surface $S$, the \textit{Riemann-Roch formula}\index{Riemann-Roch} for $D$ is 
$\chi (\O_S (D))= \frac{1}{2}D.(D-K_S)+ \chi (\O_S)$. Moreover, if $D$ is effective and nonzero and $\L$ is any line bundle on $D$, the Riemann-Roch formula\index{Riemann-Roch} for $\L$ on $D$ is 
$\chi (\L)= \deg \L +1-p_a(D) = \deg \L -\frac{1}{2}D.(D+K_S)$.

We will make use of the following results of Saint-Donat on line bundles on 
$K3$ surfaces. The first result will be used repeatedly, without further 
mention. 

\begin{prop}
{\rm \cite[Cor. 3.2]{S-D}} A complete linear system on a $K3$ surface
has no base points outside of its fixed components.
\end{prop}

\begin{prop} \label{sd1}
  {\rm \cite[Prop. 2.6]{S-D}} Let $L$ be a line bundle on a $K3$
  surface $S$ such that $|L| \not = \emptyset$ and such that $|L|$ has no
  fixed components. Then either
  \begin{itemize}
  \item [(i)] $L^2 >0$ and the general member of $|L|$ is a smooth
    curve of genus $L^2/2 +1$. In this case $h^1 (L) =0$, or
  \item [(ii)] $L^2=0$, then $L \iso \O _S (kE)$, where $k$ is an integer
    $\geq 1$ and $E$ is a smooth curve of arithmetic genus 1. In
    this case $h^1 (L) = k-1$ and every member of $|L|$ can be written
    as a sum $E_1 +  \cdots  + E_k$, where $E_i \in |E|$ for $i=1, \ldots , k$. 
  \end{itemize}
\end{prop}

\begin{Lemma} \label{sd2}
  {\rm \cite[2.7]{S-D}} Let $L$ be a nef line bundle on a $K3$
  surface $S$. Then $|L|$ is not base point free if and only if there
  exist smooth irreducible curves $E$, 
    $\Gamma$ and an integer $k \geq 2$ such that 
    \[ L \sim kE+ \Gamma, \hspace{.2in} E^2=0, \hspace{.2in}  \Gamma^2=-2,
    \hspace{.2in} E.\Gamma=1. \] 
  In this case, every member of $|L|$ is of the form
  $E_1+  \cdots  +E_k+\Gamma$, where $E_i \in |E|$ for all $i$. 
\end{Lemma}

To show the existence of $K3$ surfaces with certain divisors on it, a very 
useful result is the following by Nikulin \cite{nik} (the formulation we use is due to Morrison):

\begin{prop} \label{morrison}
  {\rm \cite[Cor. 2.9(i)]{Morrison}} Let $\rho \leq 10$ be an integer. Then every even lattice of signature 
  $(1, \rho-1)$ occurs as the Picard group of some algebraic $K3$ surface.
\end{prop}

Consider now the group generated by the Picard-Lefschetz reflections
\begin{eqnarray*}
  \phi_{\Gamma}: \Pic S & \hpil & \Pic S \\
   D & \mapsto & D+ (D.\Gamma)\Gamma
\end{eqnarray*}
where $\Gamma \in \Pic S$ satisfies $\Gamma ^2=-2$. Note that a reflection 
leaves the intersections between divisors invariant.

The following result will also be useful for our purposes:

\begin{prop} \label{plr}
  {\rm \cite[VIII, Prop. 3.9]{BPV}} A fundamental domain for this action,
  restricted to the positive cone, is the big-and-nef cone of $S$.
\end{prop}

This means that given a certain Picard lattice, we can perform
Picard-Lefschetz reflections on it, and thus assume that some
particular 
line bundle chosen (with positive self-intersection) is nef.

We will need some results about higher order embeddings of $K3$ surfaces from \cite{kn1}, which we here recall:

\begin{prop} \label{kn1prop}
  Let $k \geq 0$ be an integer and $L$ a big and nef line bundle on a $K3$ surface with $L^2 \geq 4k$. Assume $\Z$ is a $0$-dimensional subscheme of $S$ of length $h^0(\O_{\Z}) = k+1$ such that
\[ H^0 (L) \hpil H^0 (L  \* \O _{\Z}) \]
is not onto, and for any proper subscheme $\Z'$ of $\Z$, the map 
\[ H^0 (L) \hpil H^0 (L  \* \O _{\Z'}) \]
is onto.  

Then there exists an effective
  divisor $D$ passing through $\Z$ satisfying $D^2 \geq -2$, $h^1(D)=0$ and the numerical conditions
\begin{eqnarray}
\nonumber & 2D^2 \stackrel{(i)}{\leq} L.D \leq D^2 + k+1 \stackrel{(ii)}{\leq} 2k+2   \\
\label{eq:kva0} & \mbox{ with equality in (i) if and only if } L \sim 2D \mbox{ and } L^2 \leq 4k+4,   \\
& \mbox{and
\nonumber equality in (ii) if and only if } L \sim 2D \mbox{ and } L^2= 4k+4.     
\end{eqnarray}
Furthermore, either $L-2D \geq 0$, or $L^2=4k$ and $h^0(\O_S(L-D) \* \O_{\Z}) >0$.

Finally, if $L^2=4k+4$ and $L \sim 2D$, then $h^0(\O_S(D) \* \I_{\Z})=2$, and if
$L^2=4k+2$ and $D^2=k$, then $L \sim 2D+ \Gamma$, for a smooth rational curve $\Gamma$ satisfying $\Gamma.D=1$
and $\Gamma \cap \Z \neq \emptyset$.
\end{prop}

\begin{proof}
  All the statements are implicitly contained in \cite{kn1}, but we will go through the main steps in the proof for the sake of the reader.

Under the above hypotheses, it follows from the first part of the proof of
  \cite[Thm. 2.1]{BS2} or from \cite[(1.12)]{tyu} that there exists a rank 2 vector bundle 
$E$ on $S$ fitting into the exact sequence 
\begin{equation} \label{eq:kva1}
 0 \hpil \O_S \hpil E \hpil L \* \I _{\Z} \hpil 0, 
\end{equation}
and such that the coboundary map
\[ \delta : H^1 (L \* \I _ {\Z}) \hpil H^2 (\O_S) \iso \CC, \]
is an isomorphism. In particular $H^1(E)=H^2(E)=0$ and we also have $\det E = L$ and $c_2(E)=\deg \Z=k+1$.

Secondly, since $L^2 \geq 4k$, one computes by Riemann-Roch
\[ \chi (E \* E^{*} ) =  c_1 (E)^2 - 4 c_2 (E) + 4 \chi (\O_S) \geq 4, \]
whence $h^0 (E \* E^{*} ) \geq 2$ by Serre duality. This means that $E$ has nontrivial endomorphisms, and by standard arguments, as for instance in \cite[Lemma 4.4]{DM}, there are line bundles $M$ and  $N$ on $S$ and a zero-dimensional subscheme $A
  \subset S$ such that $E$ fits in an exact sequence 
\begin{equation} \label{eq:kva2}
 0 \hpil N \hpil E \hpil M \* \I _ {A} \hpil 0 
\end{equation}  
and either $N \geq M$, or $A= \emptyset$ and the sequence splits. In the latter case, we can and will assume by symmetry that $N.L \geq M.L$ (which is automatically fulfilled in the first case, by the nefness of $L$). 

Combining (\ref{eq:kva1}) and (\ref{eq:kva2}) we find
\begin{equation}
  \label{eq:kva3}
  \det E = L \sim M+N \hs \mbox{ and } \hs c_2(E)=\deg \Z=k+1=M.N+ \deg A.
\end{equation}
It follows that $N.L \geq \frac{1}{2}(N.L+M.L) \geq \frac{1}{2}L^2 >0$ and $N^2 =N.L-M.N \geq \frac{1}{2}L^2-M.N \geq 2k-(k+1) > -2$, so $N >0$ by Riemann-Roch. It folloms that $h^0(N \v)=0$, so tensoring 
(\ref{eq:kva1}) and (\ref{eq:kva2}) by $N \v$ and taking cohomology, we get 
$h^0(M \* \I_{\Z}) \geq h^0(E \* N \v) >0$, whence $M >0$ as well and there is an effective divisor $D \in |M|$ such that $D \sup \Z$, as stated.

Moreover, since $h^1(E)=h^2(E)=h^2(N)=0$, we get from (\ref{eq:kva2}) that $h^1(D)=h^1(M)=0$,
whence $D^2 \geq -2$ by Riemann-Roch.

From (\ref{eq:kva3}) we have $M.N \leq k+1$, in other words $L.D \leq D^2+k+1$, and by the Hodge index theorem 
\[ 2(D.L)D^2 \leq D^2L^2 \leq (D.L)^2. \]
Hence $2D^2 \leq D.L$, with equality if and only if $2D.L=L^2$ and $L \sim 2D$, in which case we have $L^2 =4D^2 \leq 4(k+1)$. It also follows that $D^2 \leq D.L-D^2=D.M \leq k+1$, with equalities if and only if $L \sim 2D$ and $D^2=k+1$, so that $L^2=4k+4$. This establishes the numerical criteria.

Now we want to show that, possibly after interchanging $M$ and $N$, either $L-2D \geq 0$, or $L^2=4k$ and $h^0(N \* \I_{\Z})>0$. So assume 
that $h^0(L-2D)=0$. Then the sequence (\ref{eq:kva2}) splits, and by arguing as above with $N$ and $M$ interchanged, we find 
$h^0(N \* \I_{\Z})>0$. Since $(L-2D)^2=L^2-4M.N \geq 4k-4(k+1)=-4$ and $(L-2D).L \geq 0$, we see by Riemann-Roch and the Hodge index theorem that we must have $(L-2D)^2=-2$ or $-4$. If $(L-2D)^2=-2$, Riemann-Roch yields that $2D-L >0$. By the nefness of $L$ we must have $(L-2D).L =0$, whence $M.L=N.L$ and $M^2=N^2$, and we get the desired result after interchanging $M$ and $N$. If $(L-2D)^2=-4$,
then $L^2=4k$, as stated.

We now prove the two last assertions.

If $L^2=4k+4$  and $L \sim 2D$, then from (\ref{eq:kva3}) we get $A= \emptyset$, so by tensoring 
(\ref{eq:kva1}) and (\ref{eq:kva2}) by $M \v$ and taking cohomology, we get 
$h^0(M \* \I_{\Z}) \geq h^0(E \* M \v) =2$, as stated.

If $L^2=4k+2$ and $D^2=k$, then clearly $L \not \sim 2D$, so by the numerical conditions above, we get $2k < L.D \leq 2k+1 <2k+2$, whence $L.D=2k+1$. Moreover, we find that $L \sim 2D + \Delta$, for a $\Delta >0$ satisfying $\Delta^2=-2$ and $\Delta.D=1$. Since $\Delta.L=0$, we have that $h^0(\Delta)=1$ and  $\Delta$ is supported only on smooth rational curves, and there has to exist a smooth rational curve $\Gamma$ with $\Gamma.D >0$. Since $L$ is big and nef, we get by Riemann-Roch
\begin{eqnarray*}
  h^0(L) & = & \frac{1}{2}L^2+2 = \frac{1}{2}(2D+\Delta)^2+2= \frac{1}{2}(4k+4-2)+2 \\
         & \leq & \frac{1}{2}(2D+\Gamma)^2+2 = h^0(2D+\Gamma),
\end{eqnarray*}
and since $L$ is not of the particular form in Lemma \ref{sd2} above, $L$ is base point free, so we must have $L \sim 2D +\Gamma$. So $N \sim D+\Gamma$, and $\Gamma$ is fixed in $N$. Since $N^2=D^2$ and $h^0(N)=h^0(D)$, it follows by Riemann-Roch that $h^1(N)=h^1(D)=0$. Moreover, we see from 
(\ref{eq:kva3}) that $A = \emptyset$, and by tensoring (\ref{eq:kva1}) and (\ref{eq:kva2}) with $N \v$ and $M \v$ respectively, using $H^1(\Gamma)=H^1(N)=0$, we get $h^0(M \* \I_{\Z})=1$ and 
$h^0(N \* \I_{\Z})=2$, respectively. This means that we can choose two distinct elements $N_1$ 
  and $N_2$ in $|N|$ both containing $\Z$ (scheme-theoretically). But since 
  $\Gamma$ is a base component of $|N|$, we must have 
  $N_1 =D_1 + \Gamma$ and  $N_2 =D_2 + \Gamma$, for two distinct elements 
  $D_1$  and $D_2$ of $|D|$. If $\Z$ does not meet $\Gamma$, 
  we would have both $D_1$  and $D_2$ containing $\Z$ (scheme-theoretically). 
  But this contradicts the fact that $h^0 (\O_X (D) \* \I _ {\Z}) =1$. So $\Z$ meets $\Gamma$ and we are done.
\end{proof}

\begin{rem} \label{kn1rem}
  {\rm If we replace the assumptions that $L$ be big and nef with $L^2 >0$ and $h^1(L)=0$, one can check from the proof of
  \cite[Thm. 2.1]{BS2} or from \cite[(1.12)]{tyu} that we still have a rank 2 vector bundle 
$E$ on $S$ fitting into an exact sequence as in (\ref{eq:kva1}). Moreover, if the stronger condition $L^2 >4k+4$ is fulfilled, then $c_1(E)^2 > 4c_2(E)$, and we can use Bogomolov's theorem (see \cite{bog} or \cite{re}) to find an exact sequnce as (\ref{eq:kva2})
with the properties that $(N-M)^2 >0$ and $(N-M).H >0$ for any ample line bundle on $S$. These two numerical conditions yield with Riemann-Roch that $N >M$ and it follows almost automatically that $h^0(N \v)=0$, so as in the proof of Proposition \ref{kn1prop} we find that $h^0(M \* \I_{\Z}) >0$, $h^1(M)=0$ and $M^2 \geq -2$. Furthermore, 
(\ref{eq:kva3}) still holds.

Since $L$ is not necessarily nef, we cannot assume that $N.L \geq M.L$, so we do not get the numerical conditions as in Proposition \ref{kn1prop}.

To sum up, under the assumptions of Proposition \ref{kn1prop}, with $L$ being big and nef replaced by $h^1(L)=0$, and $L^2 \geq 4k$ replaced by $L^2 >4k+4$, we get the weaker result there is a nontrivial effective decomposition $L \sim D+N$ such that $N >D$, $N.D \leq k+1$, $h^1(D)=0$, $D^2 \geq -2$ and $D$ passes through $\Z$. }
\end{rem}

\section{Surfaces in Scrolls}
\label{sins} % Always give a unique label
% use \chaptermark{}
% to alter or adjust the chapter heading in the running head

In the beginning of this section we briefly review some basic facts that can 
be found in \cite{Sc}.

\begin{defn} \label{ratscro}
Let $\E=\O_{\PP^1}(e_1) \+  \cdots \+\O_{\PP^1}(e_d)$, with $e_1 \geq  \ldots  \geq e_d
\geq 0$ and $f=e_1+ \cdots +e_d \geq 2.$ Consider the linear system 
$\L=\O_{\PP(\E)}(1)$ on the corresponding $\PP^{d-1}$-bundle $\PP(\E)$\index{$\PP(\E)$} over 
$\PP^1$. We map $\PP(\E)$ into $\PP^r$ with the complete linear system 
$H^0 (\L)$, where $r=f+d-1$. The image $T$ is by definition a rational normal 
scroll\index{rational normal scroll} of type $\bf{e}$$=(e_1, \ldots ,e_d)$. The image is smooth, and 
isomorphic to $\PP(\E)$, if and only if $e_d \geq 1$.  
\end{defn}

\begin{rem}
{\rm Some authors, like in \cite{PS}, use the term} rational normal scroll 
{\rm only if $e_d \geq 1$ (so that $T$ is smooth), but for our purposes it 
will be more convenient to use the more liberal definition above.  The 
definition of rational normal scrolls goes back at least to C. Segre, see 
\cite{Se1} and \cite{Se2}}. 
\end{rem}

\begin{defn} \label{balscro}
Let $\T$ be a rational normal scroll of type $(e_1, \ldots ,e_d)$.
We say that $\T$ is a scroll of maximally balanced type\index{maximally balanced type} if $e_1-e_d \le 1$.
\end{defn}

Let $L$ be a base point free and big line bundle on a smooth surface $S$. We denote by 
$\varphi_L$ the natural morphism
\[ \varphi_L: S \hpil \PP ^{h^0(L) -1} := \PP ^g \]
defined by the complete linear system $|L|$.
  
Assume that $L$ can be decomposed as
\begin{equation}
  \label{eq:decomposition}
  L \sim D+F, \hs \mbox{ with } \hs h^0 (D)\geq 2 \mbox{ and } h^0 (F) \geq 2.
\end{equation}
Choose a $2$-dimensional subspace $W \sub H^0 (S,D)$, which then defines a pencil\index{$\{ D_{\lambda} \}$}
\[ \{ D_{\lambda} \} _{ \lambda \in \PP^1} \sub |D|. \]

Each $\varphi_L (D_{\lambda})$ will span a 
$(h^0 (L) - h^0 (L-D) -1)$-dimensional subspace of $\PP ^g$, which is 
called the linear span of $\varphi_L (D_{\lambda})$ and denoted by 
$\overline{D_{\lambda}}$\index{$\overline{D_{\lambda}}$}. The variety swept out by 
these linear spaces,
\[ T = \cup _{\lambda \in \PP^1} \overline{D_{\lambda}} \sub \PP ^g, \]
is a rational normal scroll:

\begin{prop} \label{normscroll}
{\rm \cite{Sc}} The multiplication map
\[ W \* H^0 (S,F) \hpil H^0 (S,L) \]
yields a $2 \times h^0 (F)$ matrix with linear entries whose $2 \times 2$ 
minors vanish on $\varphi_L (S)$. The variety $T$ defined by these minors 
contains $\varphi_L (S)$ and is a rational normal scroll of degree $f :=h^0 (F)$
and dimension $d := h^0(L)-h^0(L-D)$. In particular $d+f=g+1$.
\end{prop}

Decomposing the pencil $\{ D_{\lambda} \}$ into its moving part 
$\{D'_{\lambda} \}$ and fixed part $\Delta$,
\[ D_{\lambda} \sim D'_{\lambda} + \Delta, \]
the type\index{scroll type} $(e_1,  \ldots  ,e_d)$ of the scroll $T$ is given by
\index{$e_i$}\begin{equation} \label{eq:etype}
  e_i = \# \{ j \hs | \hs d_j \geq i \}-1 ,
\end{equation}
where\index{$d_i$} \label{duinv}
\begin{eqnarray*}
  d =  d_0 & := & h^0(L)-h^0(L-D),                         \\[.05in]
       d_1 & := & h^0(L-D)-h^0(L-2D'-\Delta),              \\
       \vdots   &    &                                          \\ 
       d_i & := & h^0(L-iD'-\Delta)-h^0(L-(i+1)D'-\Delta), \\
       \vdots   &    &                                          \\ 
\end{eqnarray*}

\begin{rem}   \label{convex}
{\rm The $d_i$ form a non-increasing sequence. This follows essentially as in
the proof of Exercise B-4 in \cite{acgh}, using the socalled ``base-point-free
pencil trick''.
 }
\end{rem}

Conversely, if $\varphi_L (S)$ is contained in a scroll $T$ of degree $f$, the 
ruling of $T$ will cut out on $S$ a pencil of divisors (possibly with base 
points) $\{ D_{\lambda} \} \sub |D|$ with $h^0 (L-D) = f \geq 2$, whence 
inducing a decomposition as in (\ref{eq:decomposition}).

For any decomposition $L \sim D+F$, with $h^0(D) \geq 2$ and $h^0(F) \geq 2$, denote by 
$c$ the integer $D.F-2$. We may assume $D.L \leq F.L$, or equivalently 
$D^2 \leq F^2$. Then we have by the Hodge index theorem that $D$ satisfies the numerical conditions below:
     \begin{eqnarray*}
& 2D^2 \stackrel{(i)}{\leq} L.D =   D^2 + c + 2 \stackrel{(ii)}{\leq} 2c+4   \\
 & \mbox{ with equality in (i) or (ii) if and only if } L \eqv 2D \mbox{ and } L^2= 4c+8. 
     \end{eqnarray*}

Indeed, the condition $D.L \leq F.L$ can be rephrased as $2D.L \leq
L^2$, and by the Hodge index theorem $2D^2(D.L) \leq D^2 L^2 \leq (D.L)^2$, with equalities if and only if $L \eqv 2D$.

If the set\index{$\A(L)$}
\[ \A (L) := \{ D \in \Pic S \hs | \hs h^0(D) \geq 2 \mbox{   and   } h^0(L-D) \geq 2 \} \]
is nonempty, define the integer $\mu (L)$\index{$\mu (L)$} as 
\begin{eqnarray*}
  \mu (L)  & := & \min \{D.F-2 \hs | \hs L \sim D+F \mbox{   and   } D, F \in 
          \A (L) \}                        \\
      & = &  \min \{D.L-D^2-2 \hs | \hs D \in \A (L) \}   \\ 
\end{eqnarray*}

and set\index{$\A^0 (L)$}
\[ \A^0 (L) := \{ D \in \A (L) \hs | \hs D.(L-D) = \mu (L)+2\} \]

For $K3$ surfaces we have the following result:

\begin{prop} \label{mu1}
  Let $L$ be a base point free and big line bundle on a $K3$ surface $S$ such that 
  $\A (L) \not = \emptyset$. Then $\mu (L) \geq 0$ and any divisor $D$ 
  in $\A^0 (L)$ will have the following properties:

  \begin{itemize}
  \item [(i)]  the (possibly empty) base divisor $\Delta$ of $D$ satisfies 
              $L.\Delta=0$,
  \item [(ii)]  $h^1(D) =0$.
\end{itemize}
\end{prop}

\begin{proof}
  The first statement follows from the fact that any member of the complete 
  linear system of a base point free and big line bundle on a $K3$ surface is 
  numerically 
  $2$-connected (see \cite[(3.9.6)]{S-D}, or \cite[Thm. 1.1]{kn1} for a 
  more general statement).

  We first show (i).

  If $D$ is nef but not base point free, then by Lemma \ref{sd2}, 
  $D \sim kE+ \Gamma$, for an integer $k \geq 2$ and divisors $E$ and 
  $\Gamma$ satisfying $E^2=0$, $\Gamma^2=-2$ and $E.\Gamma=1$. Since $L$ is 
  base point free, we must have $E.L \geq 2$ (see \cite{S-D} or \cite[Thm. 1.1]{kn1}), so 
  $D.L-D^2 = (kE+\Gamma).L -(2k-2) \geq kE.L -2(k-1) \geq E.L +2(k-1) - 2 (k-1)
  = E.L$, which implies $E.L =2$, $\mu (L)=0$, and as asserted $\Gamma.L=0$. 
  
  If $D$ is not nef, there exists a smooth 
  rational curve $\Gamma$ such that $\Gamma.D <0$. Letting $D':=D-\Gamma$ and 
  we have $D' \in \A(L)$ and $D'.(L-D')= D.(L-D) -L.\Gamma+ 
  2\Gamma.D +2 \leq D.(L-D)$, whence 
  $L.\Gamma=0$, $\Gamma.D=-1$, ${D'}^2 = D^2$, $L.D'=L.D$ and 
  $D'(L-D')= D.(L-D) = c+2$. Continuing 
  inductively, we get that $\Delta.L =0$, as desired.

  Since $\Delta.L =0$ and $(D-\Delta).L -(D-\Delta)^2 \geq  D.L-D^2$, we must 
  have $D^2 \geq (D-\Delta)^2 \geq 0$. 

  We now prove (ii).
  
  If $h^1 (D) \not =0$, there exists by Ramanujam's lemma an effective 
  decomposition $D \sim D_1 +D_2$ such that $D_1.D_2 \leq 0$. By the Hodge index theorem (and the fact that $D^2 \geq 0$)
  we can assume $D_1 ^2 \geq 0$ and $D_2 ^2 \leq 0$, with equalities 
  occurring simultaneously. The divisor $D_1$ is in $\A(L)$, and writing 
  $F:=L-D$\index{$F$} we get $D_1.(F+D_2)= D.F + D_1.D_2 - D_2.F \geq F.D$, whence 
  $D_2.F \leq D_1.D_2 \leq 0$. But $L$ is nef, so 
  $D_2.L = D_2.D + D_2.F \geq 0$, which implies $D_2.F = D_1.D_2 = D_1^2 
  = D_2^2 = 0$. Now the same argument works for $D_1$, so $D_1.F=0$ and we 
  get the contradiction $D.F = (D_1+D_2).F = 0$.
\hspace{0.07cm} $\square$ 
\end{proof} 

Writing $L \sim D + F$, the above result is of course symmetric in $D$ and $F$.
It turns out that we can choose one of them to have an additional property. 
More precisely, we have :

\begin{prop} \label{mu2}
  Let $L$ be a base point free line bundle on a $K3$ surface $S$ such that 
  $\A (L) \not = \emptyset$. We can find a divisor $D$ in 
  $\A^0 (L)$ such that either $|D|$ or 
  $|L-D|$ (but not necessarily both at the same time) is base point free 
  and its general member is smooth and irreducible. If $L$ is ample, then 
  for any divisor $D$ in $\A^0 (L)$  the above 
  conditions will be satisfied for both $|D|$ and $|L-D|$.
\end{prop}

\begin{proof}
  Let $D \in \A^0 (L)$. Denote its base locus by $\Delta$ and assume it is 
  not zero. Then $L.\Delta=0$ 
  by the previous proposition. 
  
  If $D$ is nef but not base point free, then $D \sim kE+ \Gamma$ as above 
  and the smooth curve $E$ will satisfy the desired conditions.

  If $D$ is not nef, there exists a smooth 
  rational curve $\Gamma$ such that $\Gamma.D <0$. Letting $D':=D-\Gamma$, 
  we can argue inductively as above until we reach a divisor which is base 
  point free or of the form $kE+ \Gamma$.

  This procedure can of course not be performed on both $D$ and $L-D$ 
  simultaneously, but if $L$ is ample, they are both automatically base point 
  free.

  The fact that the general member of $|D|$ (or $|L-D|$) is a smooth curve 
  now follows from Proposition \ref{sd1}, since $h^1 (D)= 0$.
\hspace{0.07cm} $\square$
\end{proof} 

\begin{rem}  \label{mu2rem}
  {\rm Note that by the proofs of the two previous propositions, if $D$ is 
  not nef, then ${D'}^2= D^2$. This means that given a divisor $D \in \A (L)$, 
  we can find a divisor $D_0 \in \A (L)$ satisfying the additional conditions 
  in Proposition \ref{mu2} and such that $D_0 ^2 \leq D^2$. }
\end{rem}

\section{The Clifford index of smooth curves in $|L|$ and the
  definition of the scrolls $\T(c,D,\{ D _{\lambda} \})$}
\label{cliff}

We briefly recall the definition and some properties of gonality and Clifford 
index of curves.

Let $C$ be a smooth irreducible curve of genus $g \geq 2$.
We denote by $g^r_d$\index{$g^r_d$} a linear system of dimension $r$ and degree $d$ and say 
that $C$ is $k$-{\it gonal} (and that $k$ is its {\it gonality}\index{gonality}) if $C$
posesses a $g^1_k$ but no $g^1_{k-1}$. In particular, we call a $2$-gonal 
curve {\it hyperelliptic}\index{hyperelliptic curve} and a $3$-gonal curve {\it trigonal}\index{trigonal curve}. We denote 
by $\gon C$ the gonality of $C$. Note that if 
$C$ is $k$-gonal, all $g^1_k$'s must necessarily be base point free and 
complete.

If $A$ is a line bundle
on $C$, then the {\it Clifford index}\index{Clifford index of a line bundle on a curve} of $A$ (introduced by H.~H.~Martens in \cite{Mar}) is the integer
\[ \Cliff A = \deg A - 2(h^0 (A) -1). \]
If $g \ge 4$, then the {\it Clifford index of $C$}\index{Clifford index of a curve} itself is defined as 
\[ \Cliff C = \min \{ \Cliff A \hspace{.05in} | \hspace{.05in} h^0 (A) \geq 2, h^1 (A) \geq 2 \}. \]
 Clifford's theorem then states that 
$\Cliff C \geq 0$ with equality if
and only if $C$ is hyperelliptic and $\Cliff C =1$ if
and only if $C$ is trigonal or a smooth plane quintic.

At the other extreme, we get from Brill-Noether theory (cf.
\cite[V]{acgh}) that the gonality of $C$ satisfies $\gon C \leq \lfloor 
  \frac{g+3}{2} \rfloor $, whence $\Cliff C \leq \lfloor 
  \frac{g-1}{2} \rfloor $. For the general curve of genus
  $g$, we have $\Cliff C = \lfloor \frac{g-1}{2} \rfloor $.

We say that a line bundle $A$ on $C$ {\it contributes to the Clifford
  index of $C$}\index{contribute to the Clifford
  index} if both $h^0 (A) \geq 2$ and $ h^1 (A) \geq 2$ and that it {\it computes 
the Clifford index of $C$}\index{compute
the Clifford index} if in addition $\Cliff C = \Cliff A$. 

Note that $\Cliff A = \Cliff (\omega _C \* A^{-1})$.

The {\it Clifford dimension}\index{Clifford dimension} of $C$ is defined as
\[ \min \{ h^0 (A) -1 \hs | \hs A  
\mbox{ computes the Clifford index of } C \} . \]

A line bundle $A$ which achieves the minimum and computes the Clifford index, is said to {\it compute}\index{compute
the Clifford dimension} the Clifford dimension. A curve of Clifford index $c$ is $(c+2)$-gonal if and only if it has Clifford dimension $1$. For a general curve $C$, we have  $\gon C = c+2$.

Following \cite{gl} we give ad hoc definitions of $\Cliff C$ for $C$ of genus
$2$ or $3$: We set
$\Cliff C =0$ for $C$ of genus $2$ or hyperelliptic of genus $3$, and
$\Cliff C =1$ for $C$ non-hyperelliptic of genus $3$. This convention will be used throughout the paper, with no further mention. 

\begin{Lemma}[Coppens-Martens \cite{cm}] \label{gon1}
 The gonality $k$ of a smooth irreducible projective curve $C$ of genus $g \ge 2$ satisfies 
\[ \Cliff C +2 \leq k \leq \Cliff C +3.\]
\end{Lemma}

The curves satisfying $\gon C = \Cliff C +3$ are conjectured to be very rare and called {\it exceptional}\index{exceptional curve} (cf. \cite[(4.1)]{ma}).

Recall the following result of Green and Lazarsfeld already mentioned in the introduction:

\begin{thm}[Green-Lazarsfeld \cite{gl}] \label{grelaz}
Let $L$ be a base point free line bundle on a $K3$ surface $S$ with $L^2 >0$. Then
$\Cliff C$ is constant for all smooth irreducible $C \in |L|$, and if
$\Cliff C < \lfloor \frac{g-1}{2} \rfloor $, then there exists a line
bundle $M$ on $S$ such that $M_C := M \* \O _C$ computes the Clifford index of
$C$ for all smooth irreducible $C \in |L|$.
\end{thm}

Note that since $(L-M) \*
\O _C \iso \omega _C \* {M_C}^{-1}$, the result is symmetric in $M$ and
$L-M$.

With Theorem \ref{grelaz} in mind we make the following definition:

\begin{defn} \label{cloff}
  Let $L$ be a base point free and big line bundle on a $K3$ surface. We define the 
  {\rm Clifford index of}\index{Clifford index of a line bundle on a $K3$ surface} $L$ to be the Clifford index of all the smooth 
  curves in $|L|$ and denote it by $\Cliff L$.

  Similarly, if $(S,L)$ is a polarized $K3$ surface we will often call 
  $\Cliff L$ the {\rm Clifford index of}\index{Clifford index of a polarized $K3$ surface} $S$ and denote it by $\Cliff _L (S)$. 
\end{defn}

\begin{defn} \label{cloffnongen}
  A polarized $K3$ surface $(S,L)$ of genus $g$ is called {\rm Clifford general} if 
  $\Cliff L < \lfloor \frac{g-1}{2} \rfloor$.
\end{defn}

It turns out that we can choose the line bundle $M$ appearing in Theorem 
\ref{grelaz} above so that it satisfies certain
properties. We will need the following result in the sequel.

\begin{Lemma} \label{chooseD}
  {\rm \cite[Lemma 8.3]{kn1}} Let $L$ be a  base point free line bundle on a $K3$ surface $S$ with
  $L^2 = 2g-2 \geq 2$ and $\Cliff L=c$.

  If $c < \lfloor \frac{g-1}{2} \rfloor $, then there exists a smooth
  curve $D$ on $S$ satisfying $0 \leq D^2 \leq c+2$, $2D^2 \leq D.L$
  (either of the latter two inequalities being an equality if and only if 
  $L \sim 2D$) and
\[ \Cliff C = \Cliff (\O _S (D) \* \O_C) = D.L - D^2 -2 \] 
  for any smooth curve $C \in |L|$. 
\end{Lemma}

It is also known (see e.g. \cite{ma}) that $D$ satisfies 
$h^0 (D \*\O_C) = h^0 (D)$ and 
$h^0 ((L-D) \*\O_C) = h^0 (L-D)= h^1 (D \* \O_C)$ for any smooth curve 
$C \in |L|$.

From the results in the previous section, it is also clear that 
\[ \Cliff C = \min \{ \mu (L), \lfloor \frac{g-1}{2} \rfloor \} \]
for any smooth $C \in |L|$.

Summarizing the results of the previous section, and using Propositions \ref{mu1} and \ref{mu2}, we have that if $L$ is a base point free line bundle on a $K3$ surface $S$, of {\it sectional genus}
\[ g= g(L) = \frac{1}{2}L^2 +1,\]
and the smooth curves in $|L|$ have Clifford index
\[ c < \lfloor \frac{g-1}{2} \rfloor =  \lfloor \frac{L^2}{4} \rfloor, \]
(which in particular implies 
$L^2 \geq 4c+4$),
then there exists a divisor (class) $D$ on $S$ with the following properties (with \label{f}
$F:= L-D$):
\begin{itemize}
\item[(C1)] \hs $c= D.L - D^2 -2 =D.F-2$,
\item[(C2)] \hs $D.L \leq F.L$ (eqv. $D^2 \leq F^2$) and if equality occurs, 
            then either $L \sim 2D$ or  $h^0 (2D-L)=0$,
\item[(C3)] \hs $h^1 (D) = h^1 (F) = 0$.
\end{itemize}
A divisor (class) $D$ with the properties (C1) and (C2) above will 
be called \label{cli} a {\it Clifford divisor}\index{Clifford divisor} for $L$. This 
means in other words that $D$ and $L-D$ compute the Clifford index of all 
smooth curves in $|L|$. The property (C2) can be considered an ordering of 
$D$ and $L-D$. 
Any Clifford divisor will automatically fulfill 
property (C3) and the (possibly empty) base loci $\Delta'$ of $|D|$ and 
\label{Delt} $\Delta$ of $|L-D|$ will satisfy $L.\Delta = L.\Delta' =0$ by Proposition 
\ref{mu1}.

By Proposition \ref{mu2} we can find a Clifford divisor $D$ satisfying the 
properties
\begin{itemize}
\item[(C4)] \hs the (possibly empty) base divisor $\Delta$ of $F$ 
            satisfies $L.\Delta=0$,
\item[(C5)] \hs $|D|$ is base point free and its general member is a smooth curve,
\end{itemize}
\begin{defn} \label{free}
A divisor $D$ satisfying all properties (C1)-(C5) 
will be called a {\it free Clifford divisor}\index{free Clifford divisor} for $L$. 
\end{defn}
Since this definition only depends on the class of $D$, we will by abuse of
notation never distinguish between $D$ and its divisor class.
Hopefully, this will not cause any confusion. From now on, for a free Clifford divisor $D$, the term  $\Delta$\index{$\Delta$} will always denote the base divisor of $F=L-D$.  

Note that any Clifford divisor $D$ will satisfy the numerical 
conditions: 
\begin{eqnarray*}
& 2D^2 \stackrel{(i)}{\leq} L.D =   D^2 + c + 2 \stackrel{(ii)}{\leq} 2c+4   \\
 (*) & \mbox{ with equality in (i) or (ii) if and only if } L \sim 2D \mbox{ and } L^2= 4c+8. 
     \end{eqnarray*}
In particular, 
\begin{equation} \label{eq:c+2}
    D^2 \leq c+2, \mbox{with equality if and only if 
    $L \sim 2D$ and $L^2= 4c+8$}, %\label{eq:c+2}
\end{equation}
and by the Hodge index theorem
\begin{equation}  \label{eq:index}
    D^2 L^2 \leq (L.D)^2 = (D^2+c+2)^2.
 \end{equation}

The special limit case $D^2=c+2$ (where by (\ref{eq:c+2}) we necessarily have 
$L \sim 2D$ and $L^2= 4c+8$) will henceforth be denoted by (Q).\label{q} \index{(Q)} 

We will now take a closer look at two particular kinds of free Clifford divisors, namely:
\begin{itemize}
\item[(a)] $D^2 =c+1$, or
\item[(b)] $D^2=c$, $L \sim 2D + \Delta$, with $\Delta>0$. 
\end{itemize}

It turns out that these free Clifford divisors are of a particular form.

\begin{prop} \label{E0-E2}
   Let $L$ be a base point free and big line bundle of Clifford index 
   $c < \lfloor \frac{g-1}{2} \rfloor$ on a $K3$ surface, and 
   let $D$ be a free Clifford divisor.

   If $D$ is as in (a) above, then $L^2=4c+6$ and 
\begin{itemize}
\item [(E0)]  \index{(E0)} \hs $L \sim 2D + \Gamma$, where $\Gamma$ is a smooth rational 
             curve satisfying $\Gamma.D =1$.
\end{itemize}
   
   If $D$ is as in (b) above, then $L^2=4c+4$ and (with all $\Gamma _i$ 
denoting  smooth rational curves) either 
\begin{itemize}
\item [(E1)] \index{(E1)} \hs $L \sim 2D + \Gamma_1 + \Gamma_2 $, $D^2=c$, 
           $D.\Gamma_1 = D.\Gamma_2 = 1$, $\Gamma_1.\Gamma_2=0$, or
\item [(E2)] \index{(E2)} \hs $L \sim 2D + 2\Gamma_0 + 2\Gamma_1 + \cdots + 2\Gamma_N + \Gamma_{N+1}
           + \Gamma_{N+2}$, $D^2=c$, and the following configuration: 
\end{itemize}

\[ \xymatrix{
{D} \ar@{-}[r] & {\Gamma_0} \ar@{--}[r] & {\Gamma_N} \ar@{-}[d] \ar@{-}[r] & 
           {\Gamma_{N+1}} \\
& & {\Gamma_{N+2}}   &  
} \]

Furthermore, there can only exist  Clifford divisors of one of the three types
(E0)-(E2) (on the same surface), and all such are linearly equivalent.
   
\end{prop}

\begin{proof}
   If $D$ is as in case (a), it follows that $D$ is of the desired form in the same way as in the proof of the last statement of Proposition \ref{kn1prop}.

 Now assume $D$ is as in case (b). Then we have 
  $L^2=2D.L+ \Delta.L= 2(D^2+c+2)=4c+4$. This gives 
  $\Delta ^2= (L-2D)^2 = -4$ and 
  $D.\Delta= \frac{1}{2}(L.\Delta-\Delta ^2) =2$. By 
  Proposition \ref{Gset} below, any smooth rational curve $\Gamma$ component
  of $\Delta$ such
  that $\Gamma.D >0$, satisfies $\Gamma.D =1$, and any two such curves are 
  disjoint.  So we have to distinguish between two cases.

  If there exist two distinct rational curves $\Gamma _1$ and $\Gamma _2$ in 
  $\Delta$ such that $\Gamma _1.D = \Gamma _2.D = 1$ and 
  $\Gamma _1.\Gamma _2=0$, write 
  $L \sim 2D + \Gamma _1 + \Gamma _2 + \Delta '$, for some $\Delta ' \geq 0$.
  Then $0=\Gamma _i.L =2-2 + \Gamma _i.\Delta '$ gives 
  $\Delta '.\Gamma _i=0$, for $i=1,2$. Clearly $D.\Delta '=0$, so 
\[ (2D + \Gamma _1 + \Gamma _2).\Delta '=0, \]
  whence $\Delta '=0$, since $L$ is numerically $2$-connected, and we are in 
  case (E1).

  If there exists a rational curve $\Gamma$ occurring with multiplicity $2$ 
  in $\Delta$ such that $\Gamma.D = 1$, write 
  $L \sim 2D + 2\Gamma + \Delta '$, for some $\Delta ' \geq 0$. Since 
  $0= \Gamma.L =2-4 + \Delta '.\Gamma$, we get $\Delta '.\Gamma = 2$. 
  Iterating the process, we get case (E2).
  
  Assume now $D$ is given and $B$ is any free Clifford divisor as in (E0)-(E2). We want to show that $B \sim D$.

If $D$ is of type (E0), we must have $L^2 = 4c+6$, so $B$ 
must also be of type (E0), which means that 
$L \sim 2B + \Gamma_0$, where $\Gamma$ is a smooth rational curve such that $B.\Gamma_0=1$.

  Since
\[ 0 = \Gamma_0.L = 2D.\Gamma_0 + \Gamma.\Gamma_0,\]
  and $D$ is nef, we get the two possibilities
\[     \mbox{ (i) }   D.\Gamma_0 = \Gamma.\Gamma_0=0 \hs  \mbox{ or } 
   \hs \mbox{ (ii) }  D.\Gamma_0 = 1, \Gamma=\Gamma_0.   \]

  If (i) were to happen, we would get
\[ 4D.B = (L-\Gamma)(L-\Gamma_0) = L^2 = 4c+6,\]
  which is clearly impossible.

  Hence $\Gamma=\Gamma_0$ and $D \sim B$.

If $D$ is of type (E1) or (E2), we have $L^2=4c+4$, so $B$ must also be of type (E1) or (E2) and will therefore 
  satisfy either 
\begin{itemize}
  \item [(1)] $L \sim 2B+ \Gamma _1'  + 
             \Gamma _2' $, or
  \item [(2)] $L \sim 2B + 2\Gamma' _0  + \cdots  + 2\Gamma' _N  + \Gamma' _{N+1}
             + \Gamma' _{N+2} $,
  \end{itemize}
  where the $\Gamma' _i$ are smooth rational curves with 
   configurations as in the cases (E1) and (E2). 

   Assume now that $D$ is of type (E1). The proof if $D$ is of type (E2) is 
   similar.
   
   If (2) holds, we get from
\[ 0= \Gamma _i' .L = 2D.\Gamma _i' + 
   \Gamma _1.\Gamma _i' + \Gamma _2.\Gamma _i' ,\]
   and the fact that $D$ is nef, that 
\[ \Gamma _j.\Gamma _i' =0 \hs \mbox{  or  } \hs \Gamma _j=\Gamma _i',
    \hs i=0, \ldots ,N+2, \hs j=1,2.\]

   This gives 
\[ 4D.B = (L-\Gamma _1-\Gamma _2)(L-2\Gamma _0 ' - \cdots  - 
             2\Gamma _N ' + \Gamma _{N+1} ' 
             - \Gamma _{N+2} ') \leq L^2 =4(c+1),\]
   whence $D.B \leq c+1$ and $(D-B)^2 \geq -2$, and by 
   Riemann-Roch, if $D \not \sim B$, either $D-B$ or $B-D$ is effective. The 
   argument below is symmetric in those two cases, so assume 
   $D \sim B+ \Sigma$, for $\Sigma$ effective and $\Sigma ^2 \geq -2$. 
   Then $\Sigma.L =0$ and $\Sigma ^2 = -2$ by the Hodge index theorem. Furthermore,
\[ L \sim 2D + \Gamma _1 + \Gamma _2 
\sim 2B + \Gamma _1 + \Gamma _2 + 2\Sigma
\sim 2B + 2\Gamma _0' + \cdots  + 2\Gamma _N ' + 
      \Gamma _{N+1} ' + \Gamma _{N+2} ',\]
   whence 
\[ 2\Sigma \sim 2\Gamma _0 ' + \cdots  + 2\Gamma _N' + 
  \Gamma _{N+1} ' + \Gamma _{N+2} '-\Gamma _1 - \Gamma _2,\] 
   and $2\Sigma.B=2-(\Gamma _1 + \Gamma _2).B \leq 2$ (since $B$ is nef). By
\[ 0= \Sigma.L=2B.\Sigma + (\Gamma _1 + \Gamma _2).\Sigma + 
      2\Sigma ^2,\]
   we get $(\Gamma _1 + \Gamma _2).\Sigma \geq 2$, and
\[ 2\Sigma.D = (L-\Gamma _1 - \Gamma _2).\Sigma \leq -2,\]
contradicting the nefness of $D$.
   So we are in case (1) above and again from
\[ 0= \Gamma _i '.L = 2D.\Gamma _i '+ 
   \Gamma _1.\Gamma _i '+ \Gamma _2.\Gamma _i ',\]
   and the fact that $D$ is nef, we get the three possibilities:
   \begin{itemize}
\item [(i)] \hs  $D.\Gamma _1 '=1$, $\Gamma _1 ' = \Gamma _1$, 
            $D.\Gamma _2'  = \Gamma _2.\Gamma _2 ' =0$,
\item [(ii)]  \hs $D.\Gamma _i ' = \Gamma _1.\Gamma _i ' = 
            \Gamma _2.\Gamma _i '=0$, \hs i=1,2,
\item [(iii)]  \hs $D.\Gamma _i '=1$, $\Gamma _i ' = \Gamma _i$, 
             $i=1,2$.
      \end{itemize}
   
   In case (i) we get the absurdity
   $4D.B=(L-\Gamma _1 - \Gamma _2).(L-\Gamma _1 ' - 
                   \Gamma _2 ')=4(c+1)-2$.

   In case (ii) we get  $4D.B=4(c+1)$, whence $D.B=c+1$. 
   We calculate $(D-B)^2 \geq -2$, and by 
   Riemann-Roch, if $D \not \sim B$, either $D-B$ or $B-D$ is effective. 
   Writing $D \sim B+ \Sigma$, for $\Sigma$ effective and $\Sigma ^2 =-2$, 
   we get the same contradiction as above.

   So we are in case (iii) and $B \sim D$.
\hspace{0.07cm} $\square$ \end{proof}

The following proposition describes the case (E0) further.
           
\begin{prop} \label{onlyex}
  Let $L$ be a base point free and big line bundle on a $K3$ surface and let
  $c$ be the Clifford index of all smooth curves in $|L|$. Then the
  following conditions are equivalent:
           \begin{itemize} 
\item [(i)  ] \hs all smooth curves in $|L|$ are exceptional (i.e. have gonality $c+3$), 
\item [(ii) ] \hs there is a free Clifford divisor of type (E0), 
\item [(iii)] \hs all free Clifford divisors are linearly
             equivalent and of type (E0). \end{itemize}
                                 
Furthermore, if any of these conditions are satisfied, then all the smooth
curves in $|L|$ have Clifford dimension $r = h^0(D)-1= \frac{1}{2}(c+3)$
and $D_C$ computes the Clifford dimension of all smooth $C \in |L|$.
\end{prop}
            
\begin{proof}
  The equivalence between (i) and (iii) follows from the proof of
  \cite[Prop. 8.6]{kn1}. We will however go through the whole proof for the sake of the reader. 

 We first prove that (i) implies (iii).

 Assume, to get a contradiction, that all smooth curves in $|L|$ have gonality $c+3$,
 and that there is a free Clifford divisor $D$ which is not of type (E0). We claim that in this case, the line bundle $F_{D'}$ is base point free, for any smooth $D' \in |D|$. This clearly holds in the case (Q), so we can assume that $D^2 \leq c$, whence for any smooth $D' \in |D|$ we have
$\deg F_{D'} = c+2 \geq D_0 ^2 +2 = 2g(D_0)$, so $F_{D'}$ is base point free. By \cite[Lemma 2.2]{cp}, 
we have that there exists a smooth curve in $|L|$ of gonality 
$F.D = c+2$, a contradiction.

Next we prove that (iii) implies (i).

By Lemma \ref{gon1}, we have $c+2 \leq \gon C \leq c+3$, for any smooth curve $C \in |L|$.
Assume, to get a contradiction, that there is a smooth curve  $C \in |L|$ with $\gon C=c+2$.
Let $|B|$ be a $g^1_{c+2}$ on $C$ and pick any $Z \in |A|$ lying outside the finitely many rational curves $\Gamma'$ on $S$ satisfying $\Gamma'.L \leq 2c+4$ (we can find such a $Z$ since $B$ is base point free).
By Riemann-Roch, one easily computes $h^1(\O_C(Z))=h^0(\omega_C(-Z))=g-1-c=h^0(L)-c-2$.
From the short exact sequence
\[ 0 \hpil \O_S \hpil L \* \I_Z \hpil \omega_C(-Z) \hpil 0 \]
we then find $h^0(L \* \I_Z)=h^0(L)-c-1$. In particular, the restriction map
$H^0 (L) \khpil H^0 (L  \* \O _Z)$
is not surjective. One easily sees that for any proper subscheme $Z'$ of $Z$, the map 
$H^0 (L) \hpil H^0 (L  \* \O _{Z'})$
is surjective, since otherwise $h^0(\O_C(Z'))=2$, and $\gon C \leq \deg Z' < \deg Z$, a contradiction. So $Z$ satisfies the conditions in Proposition \ref{kn1prop} and there exists an effective divisor $D_0$ passing through $\Z$ satisfying the conditions in Proposition \ref{kn1prop}.
Since $D_0.L \leq 2c+4$ and $Z$ by assumption does not meet any smooth rational curve $\Gamma'$ on $S$ satisfying $\Gamma'.L \leq 2c+4$, we must have $D_0^2 \geq 0$. let $F_0:= L-D_0$. Since $L^2=4c+6$, we have that $h^0(F_0) \geq h^0(D_0) \geq 2$, so $D_0 \in \A(L)$, and since $D_0.F_0 \leq c+2$, we must have $D_0.F_0 = c+2$ and $D_0 \in A_0(L)$, i.e. $D_0$ is a Clifford divisor for $L$.
By our assumptions, the moving part of $|D_0|$ is of type (E0), whence either $D_0 \sim D$ or $D_0 \sim D+ \Gamma$, but the latter is ruled out since $L-2D_0 \geq 0$. So $D_0 \sim D$, and by the last statement in Proposition \ref{kn1prop}, we have that $Z$ meets $\Gamma$, but this is a contradiction on our choice of $Z$.

It is clear that (iii) implies (ii) and we now show that (ii) implies (iii) by showing that  
if $D$ is a free Clifford
divisor of type (E0) and $D'$ is any other free Clifford divisor, then $D'
\sim D$.
       
Let $B :=D-D'$. Define $R':=L-2D'$ and note that $R' \sim 2B+\Gamma$.
We have \begin{eqnarray*}
 c+2 = D.(L-D) & = & (D'+B).(D'+B+ \Gamma) \\
                 & = & {D'}^2 + (2B+\Gamma).D'+ B.(B+\Gamma) =
c+2+B.(B+\Gamma), \end{eqnarray*}
  whence $B^2+B.\Gamma=0$. Combined with $\Gamma.D = \Gamma.(D'+B)=1$, and
  since $\Gamma.D'=0$ or $1$ by Lemma \ref{propertiesR}(c), we get \[
\Gamma.D'=1, B^2=0, B.\Gamma=0, \]
  whence \[ B.R' = B.(2B +\Gamma)=0 \hs \mbox{and} \hs R'.D'=2B.D' + 1. \]
  This gives \[ L.B = 2D'.B+ R'.B= 2D'.B= R'.D'-1.\]
  But this implies \[ B.L - B^2 - 2 = R'.D' -3 < {D'}^2 + R'.D' -2 = c, \]
  whence we must have $B \sim 0$, as desired.
     
It remains to prove the last statement. If $D^2=2$, then $h^0(D)=3$, and
clearly all the smooth curves in $|L|$ have Clifford dimension $2$, so
there is nothing more to prove. We therefore can assume $D^2 \geq 4$.
   
We first show that $D$ cannot be decomposed into two moving classes, i.e.
that we cannot have $D \sim D_1 + D_2$, with $h^0(D_i) \geq 2$ for
$i=1,2$.
    
Indeed, if this were the case, then since $D.\Gamma=1$, we can assume that
$D_2.\Gamma \geq 1$, whence the contradiction \begin{eqnarray*}
  D_1.L-D_1^2 & = & (D-D_2).L - (D-D_2)^2=D.L-D^2-D_2.L+2D.D_2-D_2^2 \\
              & \leq & D.L-D^2-D_2(L-2D) = c+2 -D_2.\Gamma < c+2.
\end{eqnarray*}
   
It follows that $D_C$ is very ample for any smooth $C \in |L|$. Indeed, if
$Z$ is a length two scheme that $|D|$ fails to separate, then by the
results in \cite{kn1} (see Proposition \ref{kn1prop}) and the fact that $D^2 \geq 4$, we have that $Z$
is contained in a divisor $B$ satisfying $B^2=-2$, $B.D=0$, or $B^2=0$,
$B.D=2$, or $B^2=2$, $D \sim 2D$. One easily sees that the two last cases
would induce a decomposition $D \sim B + (D-B)$ into two moving classes,
which we have just seen is impossible. So $B^2=-2$ and $B.D=0$. Now
$(D-B)^2 \geq 2$ and since \[ (D-B).L- (D-B)^2 = D.L-D^2-\Gamma.L+2, \] we
must have $B.L \leq 1$, by the condition (iii). This means that none of
the smooth curves in $|L|$ contain $Z$, whence $D_C$ is very ample for any
smooth $C \in |L|$, as claimed.
    
Since $h^0 (D-C) =h^1(D-C)=0$, we have $r:=h^0(D)-1=h^0(D_C)-1$, which
means that $|D|$ embeds $C$ as a smooth curve of genus $g=4r-2$ and degree
$d:=g-1$ in $\PP^r$. To show that the Clifford dimension of $C$ is $r$, it
suffices by \cite[Thm. 3.6 (Recognition Theorem)]{elms} to show that $C$
(embedded by $|D|$) is not contained in any quadric of rank $\leq 4$.
  
But if this were the case, the two rulings would induce a decomposition of
$D$ into two moving classes, which is impossible by the above.
  
So $C$ has Clifford dimension $r$. 
\hspace{0.07cm} $\square$ \end{proof}
   
  \begin{rem}
  {\rm This result can be seen as a generalization of \cite[Rem.
7.13]{S-D} and \cite[Thm. 4.3]{elms}. In \cite[Thm. 4.3]{elms} the authors
prove essentially the same as above, but with the hypotheses that
 $\Pic S \iso \ZZ D \+ \ZZ \Gamma$.
               
Moreover, note that for $r \geq 3$, given any of the equivalent conditions
in Proposition \ref{onlyex}, all the smooth curves in $|L|$ satisfy the
conjecture in \cite[p. 175]{elms}. Indeed, one immediately sees that it
satisfies condition (1) in that conjecture, and in \cite{elms} it is also
shown that any curve satisfying condition (1) also satisfies the remaining
conditions (2)-(4) in that conjecture.} 
\end{rem}
            
The following result shows that, except for one particular case, any free Clifford divisor 
is itself Clifford general. We will need this result later.

\begin{prop} \label{cliffgenD}
  Let $L$ be a base point free and big line bundle of Clifford index 
   $c < \lfloor \frac{g-1}{2} \rfloor$ on a $K3$ surface, and 
   let $D$ be a free Clifford divisor with $D^2 \geq 2$.

  Then $D$ is Clifford general (i.e. all the smooth curves in $|D|$ have Clifford index
  $\lfloor \frac{g(D)-1}{2} \rfloor$), except for the case (Q), with $c=2$ (in particular $D^2=4$ and $L \sim 2D$), when there exists a smooth elliptic curve $E$ such that $E.D=2$
  In this case $D$ is hyperelliptic.
\end{prop}

\begin{proof}
  Assume $D$ is not Clifford general. Then we can assume $D^2 \geq 4$, and there is an effective decomposition $D \sim A + B$, with $h^0(A) \geq 2$, $h^0(B) \geq 2$ and $A.B =\Cliff D +2 \leq \lfloor \frac{g(D)-3}{2} \rfloor +2 \leq  \lfloor \frac{1}{4} D^2 \rfloor +1$.

By symmetry, we can assume $B.F \geq \frac{1}{2}D.F=\frac{1}{2}(c+2)$. Moreover, we must have   $A.L-A^2 \geq c+2$. Hence
\begin{eqnarray*}
 c+2 & \leq & A.L-A^2 = A.(B+F) = D.F-B.F+A.B \\
     &   =  & c+2-B.F +A.B \leq  c+2  -\frac{1}{2}(c+2) + \frac{1}{4} D^2  +1,
 \end{eqnarray*}
so $D^2 \geq 2c$. Combining with (\ref{eq:c+2}) we get $c=2$, $D^2=4$ and $L \sim 2D$, as asserted.

Since $D$ is hyperelliptic, there either exists a smooth curve $E$
  satisfying $E^2=0$ and $E.D=2$, or a smooth curve $B$ such that
  $B^2=2$ and $D \sim 2B$. However, in the second we get the contradiction
\[ D.L -D^2 -2 = 6 > 4 = B.L - B^2 -2, \]
 so we must be in the first case.
\hspace{0.07cm} $\square$ \end{proof}

Now we return to the theory of scrolls. Let $D$ be a free Clifford
divisor.

If $D^2=0$, then $|D| = \{ D_{\lambda} \} _{ \lambda \in \PP^1}$ is a pencil, 
which defines in a natural way a scroll containing $\varphi_L (S)$.

If $D^2 >0$, then $\dim |D| = \frac{1}{2} D^2 +1 >1$, and we choose a 
subpencil $\{ D_{\lambda} \} _{ \lambda \in \PP^1} \sub |D|$ as follows: Pick 
any two smooth members $D_1$ and $D_2 \in |D|$ intersecting in $D^2$ distinct 
points and such that none of these points belong to the union of the finite 
set of curves
\begin{equation}   \label{eq:degree}
  \{ \Gamma \hs | \hs \Gamma \mbox{ is a smooth rational curve, } 
               \Gamma.L \leq c+2 \}.
\end{equation}

Then \index{$ \{ D_{\lambda} \}$}
\[ \{ D_{\lambda} \} _{ \lambda \in \PP^1} := \mbox{the pencil generated by  
                                              $D_1$ and  $D_2$}.  \]

The pencil $\{ D_{\lambda} \}$ will be without fixed components 
(but with $D^2$ base points) and define in a natural way a scroll containing 
$\varphi_L (S)$ and of type determined as in equation (\ref{eq:etype}) by the 
integers
\begin{equation}  \label{dtype}
  d_i = h^0 (L-iD) - h^0(L-(i+1)D), \hs i \geq 0.
\end{equation}

Since all choices of subpencils of $|D|$ will give scrolls of the same type, 
and scrolls of the same type are isomorphic, the scrolls arising are
up to isomorphism only dependent on $D$. We denote these scrolls by
\label{tau} $\T = \T(c,D,\{ D _{\lambda} \})$.\index{$\T = \T(c,D,\{ D _{\lambda} \})$} If $h^0(D)=2$, we sometimes write only $\T(c,D)$.

Since $h^1 (L)=0$, we get by the conditions (C1) and (C5) that $\T$ has 
dimension
\begin{equation} \label{eq:dim}
  \dim \T = d_0 = h^0 (L) - h^0 (F) = c+2+ \frac{1}{2}D^2,
\end{equation}
and degree
\begin{equation} \label{eq:deg}
  \deg \T = h^0(F) = g-c-1-\frac{1}{2}D^2.
\end{equation}
Furthermore, each $D_{\lambda} \in \{ D_{\lambda} \}$ has linear span\index{$\overline{D_{\lambda}}$} 
\begin{equation} \label{eq:span}
  \overline{D_{\lambda}} = \PP ^{c+1+\frac{1}{2}D^2}.  
\end{equation}

\begin{rem} \label{storskrue} 
{\rm If $h^0(D)=r+1 \geq 3$, then $|D|$ is parametrized by a $\PP^r$.
For each $D_{\lambda}$ in $|D|$ we may take the linear span 
$  \overline{D_{\lambda}} = \PP ^{c+1+\frac{1}{2}D^2}$. Taking the union of
{\it all} these linear spaces, and not only of those corresponding to a
subpencil of $|D|$, we obtain some sort of ``ruled'' variety, which perhaps
is a more natural ambient variety for $S':=\varphi_L (S)$\index{$S'$} than the scrolls 
described above (since it is independent of a choice of pencil). Such a 
variety is an image of a $\PP ^{c+1+\frac{1}{2}D^2}$-bundle over $\PP^r$. 
The main reason why we choose to study the scrolls described above rather 
than these big ``ruled'' varieties, is that we know too little about the 
latter ones to be able to use them constructively. By using the scrolls 
above we are able to utilize the results in \cite{Sc} and in many cases 
find the resolutions of $\O_{S'}$ as an $\O_{\T}$-module. A detailed 
explanation will be given in Section \ref{singscrolls}. }
\end{rem}

\section{Two existence theorems}
\label{exi}

Given integers $g \geq 2$ and $0 \leq c \leq \lfloor \frac{g-1}{2} \rfloor$, 
one may ask whether there actually exists a pair $(S,L)$, where $S$ is a $K3$ 
surface, $L^2 = 2g-2$ and all smooth curves in $|L|$ have Clifford index $c$.

Theorem \ref{exthm} below gives a positive answer to this question. Theorem \ref{gon} below answers the same kind of question concerning the possible gonalities of a curve on a $K3$ surface.

The results in this section were first given in \cite{kn3}.
We also include the material here, to obtain a complete exposition.  

\begin{thm} \label{exthm}
  Let $g$ and $c$ be integers such that $g \geq 3$ and $0 \leq c \leq \lfloor 
  \frac{g-1}{2} \rfloor $. Then there exists a polarized $K3$ surface of genus
   $g$ and Clifford index $c$.
\end{thm}

The theorem is an immediate consequence of the following

\begin{prop} \label{exprop}
  Let $d$ and $g$ be integers such that $g \geq 3$ and $2 \leq d \leq
  \lfloor \frac{g-1}{2} \rfloor   +2$. Then there exists a $K3$
  surface $S$  with $\Pic S = \ZZ L \+ \ZZ E$, where $L^2 = 2(g-1)$, $E.L = d$ and $E^2 = 0$.
  Moreover $L$ is base point free, and
  \[ c := \Cliff L = d-2 \le \lfloor \frac{g-1}{2} \rfloor \]
  Furthermore, $E$ is the only Clifford divisor for $L$ (modulo
  equivalence class) if $d < \lfloor \frac{g-1}{2} +2 \rfloor$.
\end{prop}

To prove this proposition, we first need the following basic existence result:

\begin{Lemma} \label{help}
  Let $g \geq 3$ and $d \geq 2$ be integers. Then there exists a $K3$
  surface $S$ with $\Pic S = \ZZ L \+ \ZZ E$, such that $L$ is base point free and $E$ is a smooth curve, $L^2 = 2(g-1)$, 
  $E.L = d$ and $E^2 = 0$.
\end{Lemma}

\begin{proof}
  By Propositions \ref{morrison} and \ref{plr}, we can find a $K3$
  surface $S$ with $\Pic S = \ZZ L \+ \ZZ E$, with intersection matrix 
\[  \left[ 
  \begin{array}{cc}
  L^2  &  L.E   \\ 
  E.L  &  E^2 
    \end{array} \right]  = 
    \left[
  \begin{array}{cc}
  2(g-1)  &  d   \\ 
  d   &  0      
    \end{array} \right]     \] 
  and such that $L$ is nef. If $L$ is not base point free, there exists by Proposition \ref{sd2} a curve 
  $B$ such that $B^2=0$ and $B.L =1$. An easy calculation shows that this is 
  impossible. By \cite[Proposition 4.4]{kn2}, we have that  $|E|$ contains a smooth curve. 
\hspace{0.07cm} $\square$
\end{proof}

\renewcommand{\proofname}{Proof of Proposition {\rm \ref{exprop}}}
  
  \begin{proof} 
  Let $S$, $L$ and $E$ be as in Lemma \ref{help}, with 
  $d \leq \lfloor \frac{g-1}{2} \rfloor +2$. Note that 
  since $E$ is irreducible, we have $h^1(E) = 0$. By the cohomology of the
  short exact sequence
\[ 0 \hpil \O_S (E-L) \hpil \O_S (E) \hpil \O_{C  } (E) \hpil 0, \]
  where $C$ is any smooth curve in $|L|$, we find that $h^0 (\O_{C  } (E)) \geq h^0(E) =2$ and 
  $h^1 (\O_{C  } (E))= h^0(L-E) \geq 2$, so $\O_{C  } (E)$ contributes to the 
Clifford index of
$C  $ and 
\[ c \leq \Cliff \O_{C  } (E) \leq E.L - E^2 - 2 = d-2 < \lfloor
  \frac{g-1}{2} \rfloor. \]
If $c = \lfloor \frac{g-1}{2} \rfloor $, then we are finished.
If $c < \lfloor \frac{g-1}{2} \rfloor $, then
there has to exist an effective divisor $D$ on $S$ satisfying
\[ c = \Cliff \O_{C} (D) = D.L - D^2 - 2. \]

Since both $D$ and $L-D$ must be effective and $E$ is nef, we must have
\[ D.E \geq 0 \mbox{     and     } (L-D).E \geq 0. \]
Writing $D \sim xL + yE$ this is equivalent to
\[ dx  \geq 0 \mbox{     and     }  d(1-x)  \geq 0, \]
which gives $x=0$ or $1$. These two cases give, respectively, $D=yE$ or
$L-D =-yE$. Since $h^1(D)=h^1(L-D)=0$ by (C3), we must have $y=1$ and $D \sim E$. This shows that $c= E.L - E^2 - 2 = d-2$ and that there are no other Clifford divisors but $E$ (modulo equivalence class). 
\hspace{0.07cm} $\square$ \end{proof}

\renewcommand{\proofname}{Proof}

This concludes the proof of Theorem \ref{exthm}. 

The proof of this theorem also gives the following result, which is of its 
own interest:

\begin{thm} \label{gon}
  Let $g$ and $k$ be integers such that $g \geq 2$ and $2 \leq k \leq \lfloor 
  \frac{g+3}{2} \rfloor $. Then there exists a $K3$ surface containing a smooth
  curve of genus $g$ and gonality $k$.
\end{thm}

  The surfaces constructed in Proposition \ref{exprop} all have the property 
  that the only free Clifford divisor (modulo equivalence class) is
  a smooth elliptic curve $E$. One could 
  also perform the same construction with lattices of the form
\[  \left[ 
    \begin{array}{cc}
  L^2  &  L.D   \\ 
  D.L  &  D^2 
    \end{array} \right]     \] 
with $D^2 >0$ (and satisfying the constraints given by equations 
(\ref{eq:c+2}) and (\ref{eq:index})), but for each pair $(g,c)$ there might be 
values of $D^2$ that 
cannot occur.  We will in sections \ref{bncl} and \ref{Conc} perform more such constructions, also with lattices of higher ranks. See Proposition \ref{dimfam} for  a result
concerning low values of $c$.

\section{The singular locus of the surface $S'$ and the scroll $\T$}
\label{sing}

We start this section by describing the image \label{norms} $S':=\varphi_L(S)$\index{$S'$} by the complete
linear system $|L|$ on the $K3$ surface $S$.

\begin{prop} \label{normsing}
Let $L$ be a base point free and big line bundle on a $K3$ surface $S$, and denote 
by $\varphi_L$ the corresponding morphism and by $c$ the Clifford index of the 
 smooth curves in $|L|$.
             \begin{itemize}
  \item [(i)] If $c=0$, $\varphi_L$ is $2:1$ onto a surface of degree 
             $\frac{1}{2}L^2$,
  \item [(ii)] If $c >0$, then $\varphi_L$ is birational onto a surface of 
             degree 
             $L^2$ (in fact it is an isomorphism outside of finitely many contracted smooth rational curves), and $S':=\varphi_L(S)$ is normal and has only
             rational double points as singularities. In particular 
             $K_{S'} \iso \O_{S'}$, and  $p_a(S')=1$.
             \end{itemize}
\end{prop}

\begin{proof}
These are  well-known results due to Saint-Donat \cite{S-D} (see 
also \cite{kn1} for further discussions).
\hspace{0.07cm} $\square$ \end{proof}

Let $D$ be a free Clifford divisor and $ \{ D_{\lambda} \}$ a subpencil of $|D|$ chosen as described in the previous section.

Define the subset $\D$\index{$\D$} of the pencil $ \{ D_{\lambda} \}$ by \label{tjukkd} 
\[ \D := \{ D_{\lambda} \in \{ D_{\lambda} \} \hs | \hs \varphi_L \mbox{  does not contract any component of  } D_{\lambda} \}. \]

We then have

\begin{Lemma} \label{embfib}
  If $c>0$, then $L_{D_{\lambda}}$ is very ample for all $D_{\lambda} \in \D$.
\end{Lemma}

\begin{proof}
  By \cite[Thm. 3.1]{CF} it is sufficient to show that for any effective 
  subdivisor $A$ 
  of ${D_{\lambda}}$ we have $L.A \geq A^2 +3$.

  If $A^2 \geq 0$, then we have $L.A \geq A^2 +c+2 \geq A^2 +3$ 
  (which actually holds for any divisor $A$ on $S$). If $A^2 \leq -2$, then 
  $L.A \geq 1 \geq A^2 +3$, unless $L.A=0$, which proves the lemma.
\hspace{0.07cm} $\square$ \end{proof}

In the rest of this section we focus on the singular locus of the 
rational normal scroll $\T = \T(c,D,\{ D _{\lambda} \})$.

It is well-known that the singular locus of a rational normal scroll of 
type $(e_1, \ldots , e_d)$ is a projective space of dimension $r-1$, where 
\begin{equation}  \label{eq:r}
  r := \# \{e_i \hs | \hs e_i=0 \}. 
\end{equation}
From equation 
(\ref{eq:etype}) we have
\[ r = d_0 - d_1=h^0(L)+h^0(L-2D)-2h^0(L-D). \]

By property (C2), when $L \not \sim 2D$, we have for \label{r}
$R:=L-2D$\index{$R$} by 
Riemann-Roch
\[ h^0 (R) = \frac{1}{2} R^2 +2 + h^1 (R). \]
Note that we have $h^0 (R) >0$ if $L^2 \geq 4c+6$, and $h^0 (R) = 0$ if and 
   only if $L^2 = 4c+4$ and $h^1 (R) =0$. 

If $L \sim 2D$, we have $D^2 =c+2$ and $h^0 (L-2D)= h^0 (\O_S)=1$. 
In general, using Riemann-Roch and the fact that $H^1(L)=H^1(L-D)=0$, we get the following expression for $r$:

\begin{equation} \label{eq:sing}
  r = \left\{ \begin{array}{ll}
             D^2 + h^1 (L-2D) & \mbox{ if $ L \not \sim 2D$ (equiv. 
                                       $D^2 \not = c+2$)}, \\
             D^2-1            & \mbox{ if $ L \sim 2D$ (equiv. 
                                       $D^2 = c+2$)}
          \end{array}
    \right .     
\end{equation}

The next results will show that the term $D^2$ (or $D^2-1$)
can be interpreted geometrically as follows: The pencil 
$ \{ D_{\lambda} \}$ has $n=D^2$ distinct base points, denote their images by 
$\varphi_L$ by $x_1, \ldots ,x_n$. The linear spaces 
$ \overline{ D_{\lambda} }$ that sweep out the scroll $\T$ will intersect in 
the linear space spanned by these points, which we denote by $<x_1, \ldots ,x_n>$. 
This is a $\PP^{n-1}$ when $L \not \sim 2D$ and a $\PP^{n-2}$ when $L \sim 2D$.

Define the set
\index{$\R_{L,D}$} \begin{equation}   \label{eq:contra} 
  \R_{L,D} := \{ \Gamma \hs | \hs \Gamma \mbox{ is a smooth rational curve, } 
               \Gamma.L=0 \mbox{ and } \Gamma.D >0 \}.
\end{equation} 
The members of 
$ \varphi_L (\{ D_{\lambda} \})$ will intersect in the points 
$\{ \varphi_L (\Gamma) \} _{\Gamma \in \R_{L,D}}$ in addition to the images of 
the $D^2$ base points of $ \{ D_{\lambda} \}$. If these extra points pose new 
independent conditions, they will contribute to the singular locus of $\T$.
We will show below that among all free Clifford divisors, we can choose one such that 
the term $h^1 (L-2D)$ will correspond exactly to the singularities of the scroll arising 
from the contractions of the curves in $\R_{L,D}$.

The contraction of smooth rational curves $\Gamma$ which are not in 
$\R_{L,D}$, will occur in some fiber. Indeed, since $D.\Gamma=0$ one 
calculates $h^0 (D-\Gamma) = h^0 (D)-1$, whence $\Gamma$ will be a component of 
a unique reducible member of $ \{ D_{\lambda} \}$. Clearly, such contractions 
which occur in some fiber, and not transversally to the fibers, will not 
influence the singularities of $\T$. 

The proofs of the next three propositions are rather long and tedious, and will 
therefore be postponed until the next section.

\begin{prop} \label{Gset}
  Let $D$ be a free Clifford divisor for $L$ and $\Gamma$ a curve in 
  $\R_{L,D}$.
  
  Then $D.\Gamma=1$, $F.\Gamma=-1$ and $\Gamma$ is contained in the base 
  locus $\Delta$\index{$\Delta$} of $F$. As a consequence, $ \Delta.D = \# \R_{L,D}$, where 
  the elements are counted with the multiplicity they have in $\Delta$.
  
  Furthermore, if $\gamma$ is any reduced and connected effective divisor 
  such that $\gamma.L=0$ and $\gamma.D >0$, then $D.\gamma=1$.

  In particular, the curves in $\R_{L,D}$ are disjoint.
\end{prop}

We defined the cases (E0)-(E2) in Section \ref{cliff}. We also need
to define the following two cases for $c=0$:
\begin{itemize} \label{mere}
\item [(E3)] \index{(E3)}\hs $L \sim 3D + 2\Gamma_0 + \Gamma_1 $, $\Gamma_0$ and $\Gamma_1$ 
           are smooth rational curves, $c=D^2=0$, $L^2=6$,
           $D.\Gamma_0 = 1$, $D.\Gamma_1 = 0$, $\Gamma_0.\Gamma_1=1$.
\item [(E4)] \index{(E4)}\hs $L \sim 4D + 2\Gamma$, $\Gamma$ is a smooth rational curve, 
           $c=D^2=0$, $L^2=8$, $D.\Gamma= 1$. 
\end{itemize}

Note that in all cases (E0)-(E4) we have $h^1 (L-2D)= \Delta.D -1$.
More precisely we have:
           \begin{itemize}
\item [(E0)]  \hs $\Delta=\Gamma$, $\Delta.D =1$, $h^1 (L-2D)=0$,
\item [(E1)]  \hs $\Delta=\Gamma_1 + \Gamma_2$, $\Delta.D =2$, $h^1 (L-2D)=1$,
\item [(E2)]  \hs $\Delta=2\Gamma_0 + 2\Gamma_1 + \cdots + 2\Gamma_N + \Gamma_{N+1}
           + \Gamma_{N+2}$, $\Delta.D =2$, $h^1 (L-2D)=1$, 
\item [(E3)]  \hs $\Delta=2\Gamma_0 + \Gamma_1$, $\Delta.D =2$, $h^1 (L-2D)=1$,
\item [(E4)]  \hs $\Delta=2\Gamma$, $\Delta.D =2$, $h^1 (L-2D)=1$.
           \end{itemize}

\begin{rem} \label{perf1}
  {\rm If $D$ is any free Clifford divisor not of type (E0)-(E4), then it
  will follow from equation (\ref{ferrara}) in
  Section \ref{tedious} that $h^1 (L-2D) \geq \Delta.D$.}
\end{rem}

\begin{prop}  \label{h1delta}
   Among all free Clifford divisors for $L$ there is one, call it $D$,
   with the following property (denoting by $\Delta$ the base locus of
   $F:=L-D$):

   If $D$ is not of type (E0)-(E4), then
\[ h^1 (L-2D)= \Delta.D.\]
   \end{prop}

We will also need the following:

\begin{prop}  \label{bound}
  We have for $D$ a free Clifford divisor 
\[ h^1(L-2D) \leq \frac{1}{2}c+1-D^2, \]
  except possibly for the case $L^2 \leq 4c+6$ and $\Delta =0$, the cases {\rm(E0)-(E2)} above, 
  and the case
  \begin{equation}  \label{eq:bound}
     L^2=4c+4, D.\Delta=1, \Delta ^2 = -2.
  \end{equation}
In this latter case, $D^2 <c$. 
\end{prop}

We now study the singular locus  $V$\index{$V$} of the scroll $\T$. By equation
(\ref{eq:sing}) we know its dimension $r-1$, and in the following 
results we 
will see which points in $\varphi_L (S)$ that span $V$ and how $\varphi_L (S)$ 
intersects $V$. We will divide the treatment into the two cases $c=0$ and 
$c >0$. We recall from Proposition \ref{normsing} that these two cases are 
naturally different.

We will now treat the case $c>0$. Since we choose the base points of the pencil
$\{ D_{\lambda} \}$ to be distinct and to lie outside of the finitely
many curves in 
$\R_{L,D}$, the images by $\varphi_L$ of these points will be $n=D^2$ distinct 
points in $\varphi_L (S)$, denote them by $x_1, \ldots ,x_n$\index{$x_i$}, and their preimages by $p_1, \ldots ,p_n$\index{$p_i$}. Let $m=D.\Delta$ and let
\begin{equation} \label{eq:att}
 \R_{L,D} = \{ \Gamma_1, \ldots , \Gamma_t \}, 
\end{equation}
and define
\begin{equation}\label{eq:btt}
 m_i:= \mbox{  multiplicity of  } \Gamma_i \mbox{  in  } \Delta. 
\end{equation}\index{$m_i$}
Then $m = \sum_{i=1}^t m_i$.
Denote by $y_1, \ldots ,y_t$ the images (distinct from $x_1, \ldots ,x_n$) of the 
contractions of the curves in $\R_{L,D}$, and by 
$q_{1,\lambda}, \ldots ,q_{t,\lambda}$\index{$q_{i,\lambda}$} their corresponding preimages in each fiber. So $q_{i,\lambda}= \Gamma_i \cap D_{\lambda}$.

In the cases (E0)-(E2) of Proposition \ref{h1delta}, we use the following 
notation:

\begin{itemize}
\item [(E0)] \hs $y= \varphi_L (\Gamma)$,
\item [(E1)] \hs $y_1= \varphi_L (\Gamma_1)$, $y_2= \varphi_L (\Gamma_2)$,
\item [(E2)] \hs $y_0= \varphi_L (\Gamma_0)$.
\end{itemize}

We will denote by $q_{\lambda}$, $q_{1,\lambda}$, $q_{2,\lambda}$ and
$q_{0,\lambda}$ their respective preimages in the fiber $D_{\lambda}$.
 
Also, recall from p. \pageref{q} that we denote the special case  $L \sim 2D$ by (Q).

For each $D_{\lambda} \in \D$, we can identify $D_{\lambda}$ with its image $D_{\lambda}' := \varphi_L (D_{\lambda})$\index{$D_{\lambda}$} on $S'$ by Lemma \ref{embfib}. Moreover, we clearly have that the multiplicities of the points 
$p_1, \ldots ,p_n,q_{1,\lambda}, \ldots ,q_{t,\lambda}$ on each $D_{\lambda}$ is one, hence these points are all smooth points of $D_{\lambda}$, and consequently all 
$x_1, \ldots ,x_n,y_1, \ldots ,y_t$ are smooth points of $D_{\lambda}'$.

For any $D_{\lambda} \in \D$, we define $Z_{\lambda}$ to be the 
zero-dimensional
subscheme of length $n+m$ of $D_{\lambda}$ defined by
\index{$Z_{\lambda}$}\begin{equation} \label{eq:ctt} 
 Z_{\lambda}:= p_1 + \cdots +p_n + m_1 q_{1,\lambda}+ \cdots +m_t q_{t,\lambda}. 
\end{equation}
In particular
\begin{equation} \label{eq:dtt}
 \O _{D_{\lambda}} (Z_{\lambda}) \iso 
\O _{D_{\lambda}} (D + \sum_{i=1}^t m_i \Gamma_i). 
\end{equation}
This zero-dimensional scheme can, by the isomorphism between $D_{\lambda}$ 
and $D_{\lambda}'$, be identified with the following zero-dimensional
subscheme of $D_{\lambda}'$, which we by abuse of notation denote by the same name:
\[ Z_{\lambda}= x_1 + \cdots +x_n + m_1 y_{1,\lambda}+ \cdots +m_t y_{t,\lambda}. \]
Note that in the case (Q) all the $Z_{\lambda}$ are equal to $p_1 + \cdots +p_n$ and will be denoted by $Z$.

In the special cases (Q), (E0)-(E2) we will also define the following 
zero-dimensional subschemes of $Z_{\lambda}$ (which we again will identify to their corresponding subschemes of $D_{\lambda}'$):
 
\begin{itemize}
\item [(Q)] \hs \index{$Z^{i}$}$Z^{i}:= p_1 + \cdots +\hat{p_i}+ \cdots +p_n$,
\item [(E0)] \hs  \index{$Z^{i,\lambda}$}$Z_{0,\lambda}:= p_1 + \cdots +p_n$,
\item [(E1)] \hs  $Z_{i,\lambda}:= p_1 + \cdots +p_n+q_{i,\lambda}$, $i=1,2$, \label{mark2}
\item [(E2)] \hs  $Z_{0,\lambda}:= p_1 + \cdots +p_n+q_{0,\lambda}$.
\end{itemize}

By $<Z>$\index{$<Z>$} we will mean the linear span of a zero-dimensional scheme $Z$ on 
$S'$.

The following is the main result of this section:

\begin{thm} \label{mainsing}
  Assume $c>0$. Among all free Clifford divisors for $L$ there is one, call it $D$, satisfying the property in Proposition \ref{h1delta} and with the 
following three additional properties:
\begin{itemize}
\item[(a)] If $D$ is not of type (Q), (E0), (E1) or (E2), then for all 
  $D_{\lambda} \in \D$ we have
\[ V := \Sing \T = <Z_{\lambda}> \iso \PP ^{n+m-1}, \]
  and if $D$ is of one of the particular types above, then:
  \begin{itemize}
  \item [(Q)]  \hs $V = <Z^{i}> = <Z> \iso \PP ^{n-2}$, all $i$.
  \item [(E0)] \hs $V= <Z_{0,\lambda}> = <Z_{\lambda}> \iso \PP ^{n-1}$,              
  \item [(E1)] \hs $V= <Z_{1,\lambda}> = <Z_{2,\lambda}> = <Z_{\lambda}>  \iso \PP ^n$, 
  \item [(E2)] \hs $V= <Z_{0,\lambda}> = <Z_{\lambda}> \iso \PP ^n$. 
  \end{itemize}
\item[(b)] $V$ does not intersect $S'$ (set-theoretically) outside the points in the support of $Z_{\lambda}$.
\item[(c)] For any irreducible $D_{\lambda}$, we have
\[ V \cap D_{\lambda} = Z_{\lambda}. \] 
\end{itemize}
\end{thm}

In the theorem above, the following convention is used: 
$\PP ^{-1} = \emptyset$ (which happens if and only if 
$n=m=0$ and implies that the scroll is smooth).

\begin{rem} \label{perf2}
  If $D$ is any free Clifford divisor, we have $V \sup <Z_{\lambda}>
  \iso \PP ^{n+m-1}$,   except in the cases (Q), (E0)-(E2), where the
  property (a) is automatically fulfilled.

  If $D$ is not of type (E1) or (E2), the properties (b) and (c)
  automatically hold. If $D$ is of 
  type (E1) or (E2), then it might be that $V$ intersects $S'$ outside of the support of 
  $<Z_{\lambda}>$.
\end{rem}

The proof of Theorem \ref{mainsing} will be divided in the general
case and in the special cases (Q), (E0)-(E2). We will only prove the
two first properties. The last one will be left to the reader.

In this section, we give the proofs for the general case and the cases (Q) and (E0). The proof of the case (E1) is postponed until the next section, and the proof of the case (E2) is similar and therefore left to 
the reader.

We will write $\lambda \in \D$ for a $\lambda$ such that 
$D_{\lambda} \in \D$.

\renewcommand{\proofname}{Proof of Theorem \ref{mainsing} in the general case}

\begin{proof}
  Let $s:=n+m$. To prove that $<Z_{\lambda}> \iso \PP ^{n+m-1}$ it suffices to prove that the natural map
\[ H^0(L) \hpil H^0 (L \* \O_{Z_{\lambda}}) \]
  is surjective for all $\lambda \in \D$.

  So assume this map is not surjective for some $\lambda$. Then there exists a subscheme $Z' \sub Z_{\lambda}$ of length $s' \leq s$, for 
  some integer $s' \geq 2$ (since $L$ is base point free), such that the 
  map $H^0(L) \khpil H^0 (L \* \O_{Z'})$
  is not surjective, but such that the map $H^0(L) \khpil H^0 (L \* \O_{Z''})$ 
  is surjective for all proper subschemes $Z'' \subsetneqq Z'$.
  We now use Propositions \ref{h1delta} and \ref{bound}.

  If $\Delta=0$ and $L^2 \leq 4c+6$, we have $n=D^2 \leq c$ and $m=0$, so 
\[ L^2 \geq 4(c+1) \geq 4(n+1) = 4 (s+1).\]
  If we are in the case given by (\ref{eq:bound}), we have
\[ L^2 \geq 4(c+1) \geq 4(n+2) = 4 (s+1). \] 
  In all other cases, we have $s = m +n \leq \lfloor \frac{1}{2}c \rfloor +1$ by Propositions \ref{h1delta} and  
  \ref{bound}, so 
\[ L^2 \geq 4(c+1) \geq 4(\lfloor \frac{1}{2}c \rfloor +2) \geq 4(s+1). \]

  Therefore, by Proposition \ref{kn1prop}, there exists an effective divisor $B$ passing 
  through $Z'$ and satisfying $B^2 \geq -2$, $h^1(B)=0$ and the numerical conditions
\[ 2B^2 < B.L \leq B^2 +s' < 2s'.\]

  If $B^2 \geq 0$, then $B$ would induce a Clifford index 
  $c_B \leq s'-2 \leq n+m-2$ on the smooth curves in $|L|$. If $\Delta=0$ and $L^2 \leq 4c+6$, we 
  get the contradiction $c_B \leq n-2 \leq c-2$. If we are in the case given by 
  (\ref{eq:bound}), we get the contradiction $c_B \leq n-1 \leq c-2$. 
  Finally, in all other cases, we have 
  $c_B \leq n+m-2 < \lfloor \frac{1}{2}c \rfloor $, again a contradiction.

  Hence $B^2 =-2$ and $B$ is 
  supported on a union of smooth rational curves. Furthermore, 
  $B.L \leq s'-2$ 
  and $B.D \geq s'$ (the last inequality follows since $D_{\lambda}$ passes 
  through $Z_{\lambda}$).

  We now consider the effective decomposition 
\[ L \sim (D+B) + (F-B). \]

  Firstly note that $L.(D+B) \leq n+s'+c$ and $(D+B)^2 \geq n+2s'-2$, whence 
  $(F-B)^2=(L-D-B)^2 \geq 2c-n+2 \geq c+2 > 0$, so that $h^0 (F-B) \geq 2$.
  
  Secondly, $L.(D+B) - (D+B)^2 -2 \leq c-s' <c$, a contradiction.
  
 For the second statement, it suffices to show that there is no point 
   $x_0 \in S'- \{x_1, \ldots ,x_n,y_1, \ldots ,y_t \}$ such that 
$S'$ has an 
  $(s+1)$-secant $(s-1)$-plane through $Z_{\lambda}$ and $x_0$ for all 
  $\lambda$.

  Assume, to get a contradiction, that there is such a point $x_0$. Choose 
  any preimage $p_0$ of $x_0$, and denote by $X_{\lambda}$ the zero-dimensional
  scheme defined as the union of $Z_{\lambda}$ and $p_0$. 
  Fix any ${\lambda}$ such that $D_{\lambda}$ is irreducible.

  In these terms we have that the natural map
\[ H^0(L) \hpil H^0 (L \* \O_{X_{\lambda}}) \]
  is not surjective.

  Then there exists a subscheme $X' \sub X_{\lambda}$ of length $s'+1 \leq s+1$, for 
  some integer $s' \geq 1$, such that the map $H^0(L) \khpil H^0 (L \* \O_{X'})$
  is not surjective, but such that the map $H^0(L) \khpil H^0 (L \* \O_{X''})$ 
  is surjective for all proper subschemes $X'' \subsetneqq X'$.

  Since $L^2 \geq 4 (s+1)$ by the above, there exists by Proposition \ref{kn1prop}
  again an effective divisor 
  $B$ passing 
  through $X'$ and satisfying $B^2 \geq -2$, $h^1(B)=0$ and the numerical conditions
\[ 2B^2 \leq  B.L \leq B^2 +s'+1 \leq 2s'+2.\]

  As above, if $B^2 \geq 0$, we would get a contradiction on the Clifford index $c$.
  Hence $B^2 =-2$ and $B$ is 
  supported on a union of smooth rational curves. Furthermore, 
  $B.L \leq s'-1$ and $B.D \geq s'$ (the last inequality follows since 
  $D_{\lambda}$ is irreducible).

  As above, the effective decomposition 
\[ L \sim (D+B) + (F-B) \]
  induces a Clifford index $<c$ on the smooth curves in $|L|$, unless
  $s'=1$, $B.L=0$ and $B.D =1$. This means that $p_0$ lies in some divisor 
  which is contracted to one of the points $y_1, \ldots ,y_t$. Hence $x_0$ is one 
  of these points, a contradiction.
\hspace{0.07cm} $\square$ \end{proof}

\renewcommand{\proofname}{Proof of Theorem \ref{mainsing} in the case 
$L\sim 2D$}

 \begin{proof}
  It suffices to prove that if there is a point 
  $x_0 \in S' - \{x_1, \ldots , \hat{x_i}, \ldots ,x_n \}$ for some $i$, such that
  $S'$ has an $n$-secant $(n-2)$-plane through $x_0$ and $Z^i$, then $x_0=x_i$.

  Choose any preimage $p_0$ of $x_0$, and denote by $X_i$ the zero-dimensional
  scheme defined by $p_0$ and $Z^i$. We will show that if the natural map
\[ H^0(L) \hpil H^0 (L \* \O_{X_i}) \]
  is not surjective, then $p_0=p_i$.

  Let $X' \sub X_i$ be a subscheme of length $n' \leq n$, for 
  some integer $n' \geq 2$, such that the map 
  $H^0(L) \khpil H^0 (L \* \O_{X'})$
  is not surjective, but such that the map 
  $H^0(L) \khpil H^0 (L \* \O_{X''})$ 
  is surjective for all proper subschemes $X'' \subsetneqq X'$.
 
  By assumption, we have $n=D^2 = c+2$ and $L^2 = 4c+8 =4n$. Hence, by 
  Proposition \ref{kn1prop}, there exists an effective divisor $B$ passing 
  through $X'$ and satisfying $B^2 \geq -2$ and the numerical conditions
\[ 2B^2 \stackrel{(a)}{\leq} L.B \leq   B^2 + n' \stackrel{(b)}{\leq} 2n', \]    
  with equality in (a) or (b) implying $L \sim 2B$. 

  Since $B$ passes through $X'$, we have $B.D \geq n'-1$, whence 
  $B.L \geq 2n'-2$.
  From the inequalities above, we get $B^2 \geq n'-2 \geq 0$, so we have 
  $n'=n$ and $B.L=B^2+n$, since otherwise $B$ would induce a Clifford 
  index $<n-2 =c$ 
  on the smooth members of $|L|$. This leaves us with the two possibilities: 
\[ (i) \hs B^2=n  \mbox{ and } L \sim 2B, \hs \mbox{   or   } \hs (ii) \hs B^2=n-1.  \]     
  
  But in the second case, by Proposition \ref{kn1prop}, we have $L \sim 2B + \Gamma$, for $\Gamma$ a 
  smooth rational curve, which is impossible, since $L \sim 2D$. So we are in 
  case (i), and $B \in |D|$. By the last assertion in Proposition \ref{kn1prop} we have
  $h^0 (B \* \I _{X'}) = h^0 (E-B) = 2$,
  so there is a pencil $P$ of divisors in $|D|$ passing through $X'$. 
  
  We claim that any divisor $D_0 \in |D|$ passing 
  through $n-1$ of the points $p_1, \ldots ,p_n$, will also pass through the last 
  one. Indeed, by the surjectivity of the map 
  $H^0 (D) \khpil H^0 (\O _{D_0}(D))$, we reduce to the same statement for 
  $\O _{D_0}(D)$. By Riemann-Roch, this is equivalent to 
  $h^0 (\O _{D_0}(Z)) \geq 2$ and $\O _{D_0}(Z)$ base point free, which are 
  both satisfied since $\O _{D_0}(Z) \iso \O _{D_0}(D)$, and $\O _S(D)$ is 
  base point free.

  Since $Z$ contains the points $p_0, \ldots , \hat{p_i}, \ldots ,p_n$, we 
  have that all the members in $P$ contain all the points 
  $p_0, \ldots ,p_n$. Therefore, $P$ is the pencil $ \{ D_{\lambda} \}$, whose 
  general member is smooth and irreducible. Since all the members intersect 
  in $n$ points, we have $p_0 = p_i$, as asserted.
\hspace{0.07cm} $\square$ \end{proof}

\renewcommand{\proofname}{Proof of Theorem \ref{mainsing} in the case (E0)}

\begin{proof}
  We first prove that $<Z_{0,\lambda}> \iso \PP^{n-1}$ for all $\lambda$. If 
  this were not true, the natural map
\[ H^0(L) \hpil H^0 (L \* \O_{Z_{0,\lambda}}) \]
  would not be surjective for some $\lambda$.

  As usual let $Z' \sub Z_{0,\lambda}$ be a subscheme of length $n' \leq n$, 
  for some integer $n' \geq 2$, such that the map 
  $H^0(L) \khpil H^0 (L \* \O_{Z'})$
  is not surjective, but such that the map $H^0(L) \khpil H^0 (L \* \O_{Z''})$ 
  is surjective for all proper subschemes $Z'' \subsetneqq Z'$.

  Since $L^2=4n+2=4(n-1)+6$, we get by Proposition \ref{kn1prop} that there exists an 
  effective divisor $B$ passing through $Z'$ such that  
  $B^2 \geq -2$, $h^1(B)=0$ and 
\[ 2B^2 < B.L \leq B^2 +n' < 2n'.\]

  If $B^2 \geq 0$, we would get that $B$ induces a Clifford index 
  $c_B \leq n'-2 \leq c-1$ on the smooth curves in $|L|$, a contradiction.

  So $B^2 =-2$, and $B$ is necessarily supported on a union of smooth rational curves, since $h^1(B)=0$. But $B.L \leq n'-2 \leq n-2 =c-1$ and $Z'$ consists of base points of $\{ D_{\lambda} \}$. This means that $B$ passes through some of these base points, which contadicts the fact that we have chosen these base points to lie outside of smooth rational curves of degree $\leq c+2$ with respect to $L$.

  So $<Z_{0,\lambda}> \iso \PP^{n-1}$, and by 
  equation (\ref{eq:sing}) and Proposition \ref{h1delta} we know that 
  $V \iso \PP^{n-1}$, so the point $y$ 
  does not pose any additional conditions.

  To prove the last assertion, assume to get a contradiction that there exists
  a point $x_0 \in S' -\{x_1, \ldots ,x_n,y \}$ such that $S'$ has an 
  $(n+1)$-secant $(n-1)$-plane through $x_0$ and $Z_{0,\lambda}$. Choose any 
  preimage $p_0$ of $x_0$ and denote by $X_{\lambda}$ the zero-dimensional 
  scheme defined by $p_0$ and $Z_{0,\lambda}$. We then have that
  the natural map
\[ H^0(L) \hpil H^0 (L \* \O_{X_{\lambda}}) \]
  is not surjective. Fix a $\lambda$.

  Again let $X' \sub X_{\lambda}$ be a subscheme of length $n'+1 \leq n+1$, 
  for some integer $n' \geq 1$, such that the map 
  $H^0(L) \khpil H^0 (L \* \O_{X'})$
  is not surjective, but such that the map $H^0(L) \khpil H^0 (L \* \O_{X''})$ 
  is surjective for all subschemes $X'' \sub X'$.

  Since $L^2 = 4n+2$ and $L$ is not divisible by assumption, we get by 
  Proposition \ref{kn1prop} again that there exists an effective divisor $B$ passing 
  through $X'$ satisfying 
  $B^2 \geq -2$, $h^1(B)=0$ and  the numerical conditions
\[ 2B^2 < B.L \leq B^2 +n'+1 < 2n'+2.\]

  If $B^2 =-2$ then, since $Z'$ has length $ \geq 2$, we must have that $B$ passes through some of the base points of  $\{ D_{\lambda} \}$, a contradiction as above.

  So $B^2 \geq 0$, $n=n'$ and $X'=X_{\lambda}$. Since $h^0(L-B) \geq h^0(B) \geq 2$ by Proposition \ref{kn1prop} again, we have that $B$ is a Clifford divisor, and by 
  Proposition \ref{onlyex}, 
  we have $D \sim B$. By the last statement in  Proposition \ref{kn1prop}, we have $\Gamma \cap X_{\lambda} \not = \emptyset$, whence 
  we conclude that $p_0 \in \Gamma$. This gives the desired contradiction 
  $x_0=y$. 
\hspace{0.07cm} $\square$ \end{proof}

It will be convenient to make the following definition:

\begin{defn} \label{perfect}
  A free Clifford divisor satisfying the properties described in Proposition 
\ref{h1delta} and Theorem \ref{mainsing} will be called a {\rm perfect} 
Clifford divisor\index{perfect Clifford divisor}.   
\end{defn}

In the next section we will prove Proposition \ref{h1delta} and
Theorem \ref{mainsing}, thus proving that we can find a perfect Clifford divisor.

The main advantage of  choosing a perfect Clifford divisor is that we
then get a nice description of the singular locus of $\T$ and how it
intersects $S'$ as in Theorem \ref{mainsing} above. This theorem
states that $\Sing \T$ is ``spanned'' by the images of the base points
of the chosen subpencil of $|D|$ and the contracted curves, and
moreover that it intersects $S'$ in only these points.  If $D$ is not
perfect, then $\Sing \T \supsetneqq  Z _{\lambda}$, as seen in Remarks
\ref{perf1} and \ref{perf2}. In Proposition \ref{non2-uple} below we
will see an example where this occurs.

It will also be practical, for classification purposes, to restrict the attention to perfect Clifford divisors, as we will do in Section \ref{Conc}.

Apart from this, any free Clifford divisor will be equally fit for our purposes.

\renewcommand{\proofname}{Proof}

We include an additional description of the case (Q):

\begin{prop} \label{2-uple}
  Assume $D$ is a free Clifford divisor of type (Q) and $c \geq 2$.
  Then $\varphi_L(S)$ is the $2$-uple embedding of $\varphi_D (S)$, except in the special case described in Proposition \ref{cliffgenD} (where
  $c=2$ and there exists a smooth elliptic curve $E$ such that $E.D=2$,
in which case $D$ is hyperelliptic).
\end{prop}

\begin{proof}
  By \cite[Thm. 6.1]{S-D} $\varphi_L$ is the $2$-uple embedding of
  $\varphi_D (S)$, when $D$ is not hyperelliptic.

  Conversely, if $D$ is hyperelliptic, then $\varphi_D$ is not
  birational, so $\varphi_L$ cannot be the $2$-uple embedding of
  $\varphi_D (S)$.

  Since we assume $c \geq 2$, we have $D^2 \geq 4$, and we can use Proposition \ref{cliffgenD}
  to conclude the proof.
\hspace{0.07cm} $\square$ \end{proof}

The special case appearing in the proposition will be thouroughly
described in
Proposition \ref{non2-uple} below.

If $c=0$, there exist two kinds of (free) Clifford divisors for $L$, 
namely:
\begin{itemize}
\item [1.] $D^2=0$, $D.L=2$ and
\item [2.] $D^2=2$, $L \sim 2D$.
\end{itemize}

In both these cases $\varphi_L (S)$ is $2:1$ on each fiber. 

In the case $c=0$ we have the following result:

\begin{prop} \label{singc=0}
   Assume $c=0$. Let $D$ be a free Clifford divisor for $L$. Then 
   $D^2=0$ and $V = \emptyset$ except in the following cases:
  \begin{itemize}
  \item [(Q)]  \hs $L \sim 2D$, $D^2=2$, $V = \{ x \}$, where $x$ is the 
               common image 
               of the two base points of the chosen pencil $\{ D_{\lambda} \}$,
  \item [(E1)] \hs $D^2=0$, $L \sim 2D + \Gamma_1 + \Gamma_2 $, 
               $V= \{ \varphi_L(\Gamma_1) \}=\{ \varphi_L(\Gamma_2) \} $, 
  \item [(E2)] \hs $D^2=0$, $L \sim 2D + 2\Gamma_0 + 2\Gamma_1 + \cdots + 
               2\Gamma_N + \Gamma_{N+1}+ \Gamma_{N+2}$, 
               $V= \{ \varphi_L(\Gamma_0) \}$, 
  \item [(E3)] \hs $D^2=0$, $L \sim 3D + 2\Gamma_0 + \Gamma_1 $, 
               $V= \{ \varphi_L(\Gamma_0) \}$, 
  \item [(E4)] \hs $D^2=0$, $L \sim 4D + 2\Gamma$, $V= \{ \varphi_L(\Gamma) \}$.
  \end{itemize}  
\end{prop}

\begin{proof}
  For $D^2=0$, this follows from the fact that except for the cases (E1)-(E4), 
  the base locus $\Delta$ of $L-D$ is zero, which is shown in the proof of
  \cite[Prop. 5.7]{S-D}. In the other case, it follows from the equation 
  (\ref{eq:r}).
\hspace{0.07cm} $\square$ \end{proof}

All these cases have been completely described in 
\cite[Prop. 5.6 and 5.7]{S-D}.

When $V =\emptyset$, then $\varphi_L (S)$ is a rational ruled surface.

The cases where there are contractions across the fibers, are the cases 
(E1)-(E4). In these cases $\varphi_L (S)$ is a cone.

In the case (Q), $\varphi_L (S)$ is the Veronese 
surface in $\PP^5$.

\section{Postponed proofs}
\label{tedious}

In this section we will give the proofs omitted in the previous section.

Throughout this section $L$ will be a base point free and big line bundle of non-general Clifford index $c$. In particular,
this implies $L^2 \geq 4c+4$.

Also we write $F:=L-D$\index{$F$}  and $R:=L-2D=F-D$\index{$R$}, and denote the 
(possibly zero) base 
divisor of 
$|F|$ by $\Delta$. Recall that $L.\Delta=0$ and that we have $h^0 (R) = 0$ 
if and only if $L^2 = 4c+4$ and $h^1 (R) =0$. In particular, $h^1 (R) >0$ 
implies that $R>0$.

Furthermore, we have

\begin{Lemma} \label{delta0}
  If $h^0 (R) = 0$, then $\Delta=0$.
\end{Lemma}

\begin{proof}
  We have $L^2 = 4c+4$, so we cannot be in the cases (Q) or (E1), whence 
  $D^2 \leq c$. Choose any smooth curve $D_0 \in |D|$ and let 
  $F_{D_0} := F \* \O _{D_0}$. Then $\deg F_{D_0} = c+2 \geq D^2+2 = 2 g(D_0)$,
  whence $F_{D_0}$ is base point free.

  We first will show that this implies that $F$ is nef.

  Taking cohomology of the short exact sequence
\[ 0 \hpil R \hpil F \hpil F_{D_0} \hpil 0,\]
  and of the same sequence tensored with $-\Delta$, we get the following 
  two exact sequences (using $h^0(R)=h^0(R-\Delta)=h^1(R)=0$)
\[ 
\xymatrix{
   0 \ar[r] & H^0 (F) \ar[r] \ar@{=}[d] &  H^0 (F_{D_0}) \ar[r]  & 0  
  \\ 
   0 \ar[r] & H^0 (F-\Delta) \ar[r]  & H^0 ((F-\Delta)_{D_0}).    & 
} \]
This gives $h^0 ((F-\Delta)_{D_0}) \geq h^0 (F_{D_0})$, whence $\Delta.D=0$, since $F_{D_0}$ is base point free. This means that for any smooth rational 
curve $\Gamma$ in the support of $\Delta$, we have $\Gamma.D=\Gamma.F=0$. 
Hence $F$ is nef.

By Lemma \ref{sd2} it now suffices to show that $F$ is not of the type 
$F \sim kE + \Gamma$, for $E$ a smooth elliptic curve and $\Gamma$ a smooth 
rational curve satisfying $E.\Gamma=1$ and an integer $k \geq 2$. But if this 
were the case, we would have $E.L = 2 + c/k$. If $c \not =0$, this would mean 
that $E$ induces a lower Clifford index than $c$ on the smooth curves in 
$|L|$, a contradiction. If $c=0$, we 
get $D.F=2$ and $D^2 =0$. But this would give $R^2 =(F-D)^2 \geq -2$ and by Riemann-Roch, we would then get the contradiction 
$h^0(F-D) \geq 1$.
\hspace{0.07cm} $\square$ \end{proof}

By this lemma, if $h^0(R)=0$, the Propositions \ref{Gset}, \ref{h1delta} and 
\ref{bound} will automatically be satisfied. So for the rest of this section, we will assume $R>0$.

Let  $F_0$\index{$F_0$} be the moving component of $|F|$. Since $R>0$, we
can write 
$F_0 \sim D+A$ for some divisor  $A \geq 0$\index{$A$}. Thus we have
\begin{equation} \label{eq:F}
  F \sim D + R \sim F_0 + \Delta \sim D + A + \Delta,
\end{equation}
and
\begin{equation} \label{eq:L}
  L \sim 2D + A + \Delta.
\end{equation}

We will first study the divisors above more closely.

\begin{Lemma} \label{F0smooth}
  Except for the cases (E3) and (E4), the general member of $|F_0|$ is smooth 
  and irreducible.
\end{Lemma}

\begin{proof}
  Since $|F_0|$ is base point free, by Proposition \ref{sd1} we only need to 
  show that $F_0 \not \sim kE$, for $E$ a smooth elliptic curve and an 
  integer $k \geq 2$.

  Assume, to get a contradiction, that $F_0 \sim kE$, then by (\ref{eq:F}), we 
  have $D \sim E$ and $A \sim (k-1)E$. Let 
  $d:= c+2 =E.L \geq 2$.

  Since $L \sim (k+1)E + \Delta$, we get $L^2 = d(k+1)$, so 
  $h^0(L) = d(k+1)/2 +2$, and 
  $h^0(F) = h^0(kE) = k+1$. On the other hand, by equation (\ref{eq:deg}), we 
  have $h^0(F) = h^0(L)-d$.

  Combining the last three equations, we get 
\[ k+1 = d(k-1)/2 +2, \]
  which is only possible if $d=2$, i.e. $c=0$. A case by case study as in the 
  proof of \cite[Prop. 5.7]{S-D} establishes the lemma in this latter case.
\hspace{0.07cm} $\square$ \end{proof}

We gather some basic properties of $R$.

\begin{Lemma} \label{propertiesR}
\begin{itemize}
  \item[(a)] If $R = R_1+R_2$ is an effective decomposition, then 
            $R_1.R_2 \geq 0$.
  \item[(b)] If $\gamma$ is an effective divisor satisfying $\gamma ^2=-2$ and 
            $\gamma.R <0$, then $\gamma.R =-1$ or $-2$.
  \item[(c)] If $\gamma$ is an effective divisor satisfying $\gamma ^2=-2$ and 
            $\gamma.L=0$, then either $\gamma.D=\gamma.F=\gamma.R=0$ or 
            $\gamma.D=1$, $\gamma.F=-1$ and $\gamma.R=-2$.
  \item[(d)] If $\Gamma$ is a smooth rational curve, then 
            $\Gamma \in \R _{L,D}$ if and only if $\Gamma.R=-2$ and 
            $\Gamma.L=0$.
\end{itemize}
\end{Lemma}

\begin{proof}
  To prove (a), one immediately sees that if $R_1.R_2 < 0$, then the effective 
  decomposition $L \sim (D+R_1) + (D+R_2)$ would induce a Clifford index $<c$.
  
  The other assertions are immediate consequences of (a).
\hspace{0.07cm} $\square$ \end{proof}

This concludes the proof of Proposition \ref{Gset}.

\begin{Lemma} \label{deltaA}
  Except for the cases (E3) and (E4), the following holds:

  $\Delta ^2 = -2D.\Delta$ and $\Delta.A=0$.
\end{Lemma}

\begin{proof}
  By Lemma \ref{F0smooth}, we have $h^1(F_0)=0$. From $0=\Delta.L = 2\Delta.D+
  \Delta.A + \Delta ^2$, we get
  \begin{equation} \label{eq:a}
    \Delta ^2 = -2\Delta.D - \Delta.A.
  \end{equation}
  Furthermore, we also have 
\[     h^0 (F_0) = h^0 (F) = \frac{1}{2}F_0^2+F_0.\Delta + 
     \frac{1}{2}\Delta ^2+2= h^0 (F_0) + (D+A).\Delta + \frac{1}{2}\Delta ^2,\]
  which implies
  \begin{equation} \label{eq:b}
    \Delta ^2 = -2\Delta.D -2 \Delta.A.
  \end{equation}
  Combining equations (\ref{eq:a}) and (\ref{eq:b}), we get $\Delta.A=0$ and 
  $\Delta ^2 = -2D.\Delta$.
\hspace{0.07cm} $\square$ \end{proof}

We have seen in Proposition \ref{onlyex}, that if there exists a free Clifford divisor of type (E0), then all free Clifford divisors are linearly equivalent and of type (E0).

We now take a closer look at the types (E1) and (E2).

\begin{prop} \label{C6ref}
   Let $L$ be a base point free and big line bundle of non-general Clifford index $c$ 
   on a $K3$ surface and let $D$ be a free Clifford divisor of type (E1) or (E2).

   If $D' \not \sim D$ is any other free Clifford divisor, then $B:= D-D' >0$ and 
\begin{equation} \label{eq:S}
  \Delta .D' = 0, \Delta .B = 2, B^2=-2.  
\end{equation}
\end{prop}

\begin{proof}
Let $R':=L-2D'$ as usual, and note that $R' \sim 2B+\Delta$.

Since ${R'}^2 = L^2 - 4(c+2) = \Delta ^2 = -4$, we get $B^2+B.\Delta =0$. 
 Combined with $\Delta .D = \Delta .(D'+B) = 2$, we get the two possibilities

\begin{itemize}
  \item[(a)] $\Delta .D' \geq 2$, $\Delta .B \leq 0$, $B^2 \geq 0$,
  \item[(b)] $\Delta .D' =0$, $\Delta .B =2$, $B^2 = -2$.
\end{itemize}

Using ${D'}^2 \leq c$, we calculate
\begin{eqnarray}
\label{eq:new2} B.L & = & \frac{1}{2}(R'-\Delta).L = \frac{1}{2}R'.L = \frac{1}{2}(L-2D').L \\
\nonumber  &  = & \frac{1}{2}(L^2-2({D'}^2+c+2)) \geq 2c+2-c-c-2 =0.  
\end{eqnarray}

In case (a) we then must have $B.L >0$ by the Hodge index theorem, so $B >0$ by Riemann-Roch. We get
\[ B.R' = B.(2B +\Delta )=B^2 \geq 0 \hs \mbox{and} 
   \hs R'.D'=2B.D' + \Delta .D' \geq 2B.D' +2, \]
 which gives
\[ L.B=2D'.B +R'.B = 2D'.B +B^2 \leq R'.D' +B^2 -2. \]
 But this implies
\[ B.L - B^2 -2 \leq R'.D' - 4 < {D'}^2 + R'.D' -2 =c, \]
 whence we must have $B \sim 0$. 

 So we must be in case (b), and by Riemann-Roch we have either $B >0$
 or $-B >0$. We see from (\ref{eq:new2}) that $B.L >0$ unless ${D'}^2=D^2=c$. But if the
 latter holds, since both $D'$ and $D$ are assumed to be free
 Clifford divisors (so that $h^1(D)=h^1(D')=0$), we have $h^0(D)=h^0(D')$, whence $D \sim D'$ and $B\sim 0$, a contradiction. Hence $D.L >0$, so $B >0$ and we are done.
\hspace{0.07cm} $\square$ \end{proof}

As seen below, we will distinguish between inclusions $D' < D$ as in
Proposition \ref{C6ref} with $\Delta' =0$ and $\Delta' \not =0$ (where
$\Delta'$ is the base divisor of $|L-D'|$).

By Propositions \ref{E0-E2} and \ref{onlyex} it is clear that we can choose a free Clifford divisor $D$ with the two additional properties (recall that $\Delta$ as usual denotes the base divisor of $|L-D|$):

\begin{itemize}
\item[(C6)] \hs If $D'$ is any other free Clifford divisor
            such that $D' >D$, then $\Delta \not =0$ and $D'$ is of type (E1) or (E2). 
\item[(C7)] \hs If $D$ is of type (E1) or (E2) above, and $D'$ is any other 
free Clifford divisor satisfying (C6), then $D' \sim D$.
\end{itemize}

Property (C6) is a {\it maximality condition}: it means that we choose a 
free Clifford divisor which is not contained in any other free Clifford 
divisor, unless possibly when $\Delta \not =0$ and it is contained in some free Clifford divisor of type (E1) or (E2). 

Property (C7) means that if we can, we will choose among all free Clifford 
divisors satisfying (C6), one that is not of type (E1) or (E2). 

It turns out, as we will show in this section,  that free Clifford divisors satisfying the additional properties (C6) and (C7) will be perfect, i.e. they will satisfy Propositions \ref{h1delta} and \ref{mainsing}.

Now assume $R= R_1 + R_2$ is an effective decomposition such that $R_1.R_2=0$.
Then $L \sim (D+R_1) + (D+R_2)$ is an effective decomposition satisfying
\[ (D+R_1).(D+R_2) = D^2 +D.( R_1 + R_2) = D.F = c+2, \]
so this decomposition induces the same Clifford index $c$. This means that 
either $D+R_1$ or $D+R_2$ is a Clifford divisor. This enables us to prove the following:

\begin{prop} \label{Rdec}
  Assume $D$ is not of type (E3) or (E4) and satisfies (C6) and (C7). Assume furthermore that there exists 
  an effective decomposition $R= R_1 + R_2$ such that $R_1.R_2=0$ and such 
  that $D+R_1$ is a Clifford divisor.
 
  Then either $\Delta \not =0$ and $D+R_1$ is of type (E1) or (E2), or there exists a smooth rational curve $\Gamma$ satisfying 
  either
\begin{itemize}
  \item[(I) ] \hs $\Gamma.D=\Gamma.F=\Gamma.L=0$, $\Gamma.R_1 =-1$, 
             $\Gamma.R_2 =1$, or
  \item[(II)] \hs $\Gamma.D=1$, $\Gamma.F=-1$, $\Gamma.L=0$, $\Gamma.R_1 =-2$, 
             $\Gamma.R_2 =0$.
\end{itemize}

\end{prop}

\begin{proof}
  Let $D_1 := D+R_1$ and $D_2:=D+R_2$. Since $D_1$ is a Clifford divisor 
  containing $D$, we have by condition (C6) that either $D_1$ is not a free
  Clifford divisor, or $\Delta \not =0$ and $D_1$ is of type (E1) or (E2).

  So we can assume $D_1$ is not a free
  Clifford divisor, which means that $D_1$ is not base 
  point free.

  If $D_1$ is nef, then by Lemma \ref{sd2} it is of the form
\[ D_1 \sim lE + \Gamma_0, \]
  for some smooth elliptic curve $E$ and smooth rational curve $\Gamma_0$ 
  satisfying $E.\Gamma_0=1$, and some integer $l \geq 2$.
  This gives
\[ R_1 \sim (l-1)E + \Gamma_0 \hs \mbox{and} \hs D \sim E. \]
  Write
\[ D_2 =D + R_2 \sim D + M  + B,\]
  where $B \geq 0$ is the base divisor of $D_2$ and $M \geq 0$. Note that 
  $M + B \sim R_2$.

  We have 
  \begin{eqnarray}
    \label{eq:new1} 0  =  R_1.R_2 & = & ((l-1)D + \Gamma_0).(M + B) \\
     \nonumber                             & = & (l-1)D.M + (l-1)D.B + \Gamma_0.M + \Gamma_0.B.
  \end{eqnarray}
  Also, we have an effective decomposition
\[ R \sim ((l-1)D + M + B) + \Gamma_0,\]
  such that, using (\ref{eq:new1}),  
  \begin{equation}
    \label{eq:new3}
((l-1)D + M + B). \Gamma_0 = l-1 + \Gamma_0.M + \Gamma_0.B = 
   (l-1)(1-D.M - D.B). 
  \end{equation}
  By \cite[Lemma 3.7]{S-D}, if $M \not =0$, either $M \sim kD$, with $D^2=0$, for some 
  integer $k \geq 1$, or $D.M \geq 2$. In this latter case, the latter product in 
 (\ref{eq:new3})
  would be negative, contradicting Lemma \ref{propertiesR}. So we must have 
  $M \sim kD$, for some integer $k \geq0$ and $D.B =0$ or $1$.

  So $R \sim R_1 + R_2 \sim (l-1)D + \Gamma_0 + M + B \sim (k+l-1)D +\Gamma_0
  + B$ and 
\[ c+2 = D.F = D.(D+R) = ((k+l)D + \Gamma_0 + B).D \leq 2, \]
  which gives $c=0$ and $B.D =1$. A short analysis as in part (b) of 
  the proof of \cite[Lemma 5.7.2]{S-D} shows that $D$ is then of type (E3) or (E4). 

  So $D_1$ is not nef, which means that there exists a smooth rational curve 
  $\Gamma$ such that $\Gamma.D_1 < 0$, whence $\Gamma$ is fixed in $|D_1|$ 
  and $\Gamma.L=0$, by Proposition \ref{mu1}. Combining 
  $\Gamma.D_1 = \Gamma.D + \Gamma.R_1 \leq -1$ and 
  $0 = \Gamma.L = 2\Gamma.D + \Gamma.R_1 + \Gamma.R_2$, we get
\begin{equation} \label{gamd3}
  1 - \Gamma.R_2 \leq \Gamma.D \leq -1 - \Gamma.R_1.
\end{equation}
  Furthermore, by Lemma \ref{propertiesR}(b), we have
\begin{equation} \label{gamd4}
  \Gamma.R = \Gamma.R_1 + \Gamma.R_2 \geq -2. 
\end{equation}

 If $R_1 = \Gamma$, we are done by Lemma \ref{propertiesR}(c), so we can assume 
 that $R_1-\Gamma >0$. Then by Lemma \ref{propertiesR}(a) we have 
 $(R_1-\Gamma).(R_2+\Gamma)= R_1.R_2 + \Gamma.R_1 - \Gamma.R_2 +2 \geq 0$, 
 which implies
\begin{equation} \label{gamd5}
 \Gamma.R_1 - \Gamma.R_2 \geq -2. 
\end{equation}
  Combining (\ref{gamd4}) and (\ref{gamd5}), we get
\begin{equation} \label{gamd6}
  -2 - \Gamma.R_1 \leq \Gamma.R_2 \leq 2 + \Gamma.R_1 \hs \mbox{and} \hs 
                       \Gamma.R_1 \geq -2.
\end{equation}
  Combining (\ref{gamd6}) with (\ref{gamd3}) and Lemma 
  \ref{propertiesR}(c), we end up with the two possibilities given by (I) and (II)
  above.
\hspace{0.07cm} $\square$ \end{proof}

We now need a basic lemma about $A$.

\begin{Lemma} \label{lemA}
 If $A=0$, then $D$ is of one of the types (E0)-(E2).

 If $A^2 \leq -2$, then one of the following holds:
\begin{itemize}
  \item[(a)] $A^2 = -4$, $\Delta=0$, $L^2=4c+4$,
  \item[(b)] $A^2 = -2$, $\Delta=0$, $L^2=4c+6$,
  \item[(c)] $A^2 = -2$, $\Delta ^2 =-2$, $D.\Delta=1$, $L^2=4c+4$.
\end{itemize}
Moreover, in case (c) we have $D^2 < c$.
\end{Lemma}
 
\begin{proof}
 If $A=0$, we must have $-4 \leq \Delta ^2 = R^2 \leq -2$, whence 
 $\Delta ^2 = -4$ or $-2$, $D. \Delta = 2$ or $1$ respectively (by 
 Lemma \ref{deltaA}), and 
 $L^2=4c+4$ or $4c+6$ respectively. An analysis as in Proposition 
 \ref{E0-E2} now gives that $D$ is as in one of the cases 
 (E0)-(E2).

 If $A^2 \leq -2$, we have by $R^2= A^2 + \Delta ^2 = L^2 - 4(c+2)$ (where 
 we have used Lemma \ref{deltaA}) that either $\Delta=0$ and we are in case (a) 
 or (b) above, or that $A^2 = -2$, $\Delta ^2 =-2$, $D.\Delta=1$ (by Lemma 
 \ref{deltaA}) and $L^2=4c+4$, i.e. case (c).

 In this latter case, we have
\[ c+2 = D.F = D^2+ D.A + D.\Delta = D^2+ D.A + 1, \]
 whence $D^2 = c+1 - D.A$. Since $D+A \sim F_0$ is base point free, we have
 $D.A \geq 2$ by \cite[(3.9.6)]{S-D}, whence $D^2 < c$.
\hspace{0.07cm} $\square$ \end{proof}

We can now prove Proposition \ref{bound}.

First note that the Proposition is true for the cases (E3) and (E4), so we will from now on assume that we are not in any of these two cases.

When we are not in the exceptional cases of the proposition (which are the cases (E0)-(E2) and the cases (a)-(c) of the last lemma), we have $A \not =0$ and $A^2 \geq 0$. In particular $h^0(A) \geq 2$. Moreover $h^0(L-A) \geq h^0(2D) \geq 3$. From the standard exact sequence for any $C \in |L|$
\[ 0 \hpil A-L \hpil A \hpil A_C \hpil  0, \]
we see that $A_C$ contributes to the Clifford index of $C$, and moreover that $h^0(A_C) \geq h^0(A)$.

We first claim that 
\begin{equation}
  \label{eq:nn}
h^1(A)= D^2-c-2+D.A+h^1(R).  
\end{equation}

Indeed, we have by Lemma \ref{F0smooth} that $h^1(F_0)=0$, whence 
\begin{eqnarray*}
  h^0(F) = h^0(F_0)=h^0(D+A) & = & \frac{1}{2}D^2 + D.A + \frac{1}{2}A^2+2 \\
                             & = & \frac{1}{2}D^2 + D.A + h^0(A)-h^1(A),
\end{eqnarray*}
which gives
\begin{eqnarray*}
  h^0(A) = h^0(R) & = & h^0(F)-L.D+ \frac{3}{2}D^2 +h^1(R) \\
  & = & (\frac{1}{2}D^2 + D.A + h^0(A)-h^1(A)) -L.D+ \frac{3}{2}D^2 +h^1(R) \\   & = & h^0(A) + D^2-c-2+D.A +h^1(R)-h^1(A),
\end{eqnarray*}
whence (\ref{eq:nn}) follows.

Now we get
\begin{eqnarray*}
  \Cliff A_C &  =   &  \deg A_C - 2(h^0(A_C)-1) \\
             & \leq &  L.A-2(\frac{1}{2}A^2+1+h^1(A)) \\
             &  =   &  L.A-A^2-2-2h^1(A) \\
             &  =   &  2D.A-2-2(D^2-c-2+D.A+h^1(R)) \\
             &  =   &  2(c+1-D^2-h^1(R)). 
\end{eqnarray*}

But since $A_C$ contributes to the Clifford index of $C$, we must have $\Cliff A_C \geq c$, whence Proposition \ref{bound} follows.
 
Before proving the next result, we will need the following easy lemma.

\begin{Lemma}  \label{Adec}
  Assume $D$ is not as in (E3) or (E4). If $A \sim A_1 + A_2$ is an effective 
  decomposition such that $A_1.A_2 \leq 0$, then
\[ A_1.A_2 = A_1.\Delta=  A_2.\Delta =0, \]
  and either $D+A_1$ or $D+A_2$ is a Clifford divisor.
\end{Lemma}

\begin{proof}
  By Lemma \ref{deltaA} we have $\Delta.A=0$, so we can assume (possibly after
  interchanging $A_1$ and $A_2$) that $A_1.\Delta \leq 0$ and 
  $A_2.\Delta \geq 0$. Then $R \sim A_1 + (A_2 + \Delta)$ is an effective 
  decomposition of $R$ such that
\[ A_1.(A_2 + \Delta) =   A_1.A_2 + A_1.\Delta \leq 0.\]
  By Lemma \ref{propertiesR}(a) we must have equality,  
  whence $A_1.A_2 = A_1.\Delta=  A_2.\Delta =0$.

  If $A_1.L > A_2.L$ (resp. $A_2.L > A_1.L$), then clearly $D + A_1$ (resp. 
  $D + A_2$) is a Clifford divisor by condition (C2). 

  If $A_1.L = A_2.L$, then $D+A_i$ is not a Clifford divisor if and only if
  $h^0 ((D+A_i)-(D+A_{3-i}+\Delta)) = h^0 (A_i - A_{3-i} - \Delta) >0$. 
  Clearly this condition cannot hold for both $i=1$ and $2$. So we are done. 
\hspace{0.07cm} $\square$ \end{proof}

The next result is the crucial one to prove Proposition 
\ref{h1delta}.

\begin{prop} \label{h1A}
 If $D$ satisfies (C6) and (C7), then $H^1 (A)=0$ except for the case (E4).
\end{prop}

\begin{proof}
 The result is trivial if $A =0$. So we will assume $A >0$. 
 Also, the result is fulfilled in the case (E3), so we can assume $D$ is not as in (E3) or (E4). In particular, we can use the Lemmas \ref{deltaA} and \ref{Adec}. 
 
 If $h^1 (-A) = h^1 (A) >0$, then $A$ cannot be numerically $1$-connected, 
 whence there exists a nontrivial effective decomposition $A \sim A_1 + A_2$ such that 
 $A_1.A_2 \leq 0$. By the previous lemma, we have 
 $A_1.A_2 = A_1.\Delta=  A_2.\Delta =0$, and (possibly after interchanging 
 $A_1$ and $A_2$) we can assume that $D+A_1$ is a Clifford divisor.

 Assume first that $D$ and $D+A_1$ are as in the special case where $\Delta \not =0$ and $D':=D+A_1$ is a free Clifford divisor of type (E1) or 
 (E2), so the base divisor $\Delta '$ of 
\[ F' := L-D' \sim D + A_2 + \Delta \sim D + A_1 + \Delta' \]
 satisfies ${\Delta '}^2 =-4$. Furthermore, by Proposition \ref{C6ref}, $\Delta '.D =0$, $A_1^2 =-2$ and $A_1.\Delta '=2$. Also, since $\Delta '.L =0$ and $D'L=F'.L$, we must have $A_1.L=A_2.L$. Also note that $A_1 \not \sim A_2$, since $A_1.A_2=0$.

Since $(A_1-A_2).L=0$, we must by the Hodge index theorem have 
$A^2= (A_1-A_2)^2 \leq -2$. By Lemma \ref{lemA} and the fact that $\Delta \not =0$, this gives us 
\[ A_1^2=-2, \hs A_2^2=0,  \hs \Delta ^2 =-2, \]
We then get from 
$L^2 = 4c+4 = (2D + A + \Delta).L = 2D.L + A.L = 2D^2 + 2c +4 + A.L$,
that
\[ A.L = 2(c-D^2). \]
Since $A_1.L=A_2.L$, we have 
\[ A_1.L=A_2.L= c-D^2. \]
So $A_2$ would induce a Clifford index $\leq L.A_2- A_2^2-2 = c-D^2-2 <c$ on the smooth curves in $|L|$, a contradiction.

So we can now use Proposition \ref{Rdec} and find a smooth rational curve satisfying one of the two conditions:
\begin{itemize}
  \item[(I) ] $\Gamma.D=0$, $\Gamma.A_1 =-1$, 
             $\Gamma.(A_2 + \Delta) =1$, 
  \item[(II)] $\Gamma.D=1$, $\Gamma.A_1 =-2$, 
             $\Gamma.(A_2 + \Delta) =0$.
\end{itemize}

In case (I) we get $\Gamma.A = \Gamma.F_0 - \Gamma.D \geq 0$, whence 
$\Gamma.A_2 \geq 1$. Since $\Gamma.A_1 =-1$, we have $A_1-\Gamma >0$, and we get an effective decomposition $A \sim (A_1-\Gamma) + (A_2 + \Gamma)$ such that
\[ (A_1-\Gamma).(A_2 + \Gamma) = A_1.A_2 - \Gamma.A_2 +  \Gamma.A_1 - \Gamma ^2  \leq 0, \]
so by Lemma \ref{Adec}, we must have $\Gamma.A_2=1$ and 
$(A_1-\Gamma).(A_2 + \Gamma) = 0$. Obviously, $D + A_1-\Gamma$ is a Clifford divisor, and we can now repeat the process with $A_1$ and $A_2$ replaced by 
$A_1-\Gamma$ and $A_2 + \Gamma$. This will eventually bring us in case (II) after a finite number of steps. 

So we can assume that $A_1$ and $A_2$ are as in case (II). Again, by $\Gamma.A = \Gamma.F_0 - \Gamma.D \geq -1$, we have $\Gamma.A_2 \geq 1$, whence $\Gamma \not = A_1$ and $A_1-\Gamma >0$. Since 
\[ (A_1-\Gamma).(A_2 + \Gamma) = A_1.A_2 - \Gamma.A_2 +  \Gamma.A_1 - \Gamma ^2  \leq -1, \]
we have a contradiction by Lemma \ref{Adec}.

This concludes the proof of the proposition.
\hspace{0.07cm} $\square$ \end{proof}

We can now prove Proposition \ref{h1delta}.

By Lemma \ref{deltaA}, we can assume $A.\Delta=0$ 
and $\Delta ^2 = -2D.\Delta$.

One easily sees that the base divisor of $R$ must contain $\Delta$, so $h^0(A) = h^0 (R) = h^0 (A + \Delta)$.

If $A >0$, we have 
\begin{eqnarray} 
 \label{ferrara} h^0 (A) =  h^0 (A + \Delta) & = & \frac{1}{2}A^2 +2 +  \frac{1}{2}\Delta ^2 + h^1 (R) \\
  \nonumber                           & = & h^0 (A) - h^1(A) + D.\Delta + h^1 (R), 
\end{eqnarray}
whence $h^1(R) =  D.\Delta+ h^1(A)$. If we choose $D$ such that it satisfies (C6) and (C7), then $h^1(A)=0$ by Proposition \ref{h1A}.

If $A=0$, then $R = \Delta$ and  $D$ is of one of the types (E0)-(E2) by Lemma \ref{lemA}. 
This concludes the proof of Proposition \ref{h1delta}.

Note that in the case $A=0$, we have 
\begin{equation} \label{montero}
1 = h^0 (R) = \frac{1}{2}R^2 +2 + h^1 (R) = -D.\Delta+2 + h^1 (R), 
\end{equation}
whence $h^1 (R) = D.\Delta-1$, as we have already noted.

We now give the proof of Theorem \ref{mainsing} in the case (E1), which was left out in the previous section. The proof in the case (E2) is similar, and therefore
left to the reader.

\renewcommand{\proofname}{Proof of Theorem \ref{mainsing} in the case (E1)}

\begin{proof}
  We first show that  $<Z_{i,\lambda}> \iso \PP ^n$ for $i=1,2$ and any 
  $\lambda$. (Recall the definition of $Z_{i,\lambda}$ on p.¨\pageref{mark2}. In particular, $\deg Z_{i,\lambda}=n+1$.) If this were not true, 
   the natural map
\[ H^0(L) \hpil H^0 (L \* \O_{Z_{i,\lambda}}) \]
  would not be surjective.

  Let $Z' \sub Z_{i,\lambda}$ be a subscheme of length $n'+1 \leq n+1$, for 
  some integer $n' \geq 1$, such that the map 
  $H^0(L) \khpil H^0 (L \* \O_{Z'})$
  is not surjective, but such that the map $H^0(L) \khpil H^0 (L \* \O_{Z''})$ 
  is surjective for all proper subschemes $Z'' \subsetneqq Z'$.

  Since $L^2=4c+4=4(n+1)$, we have by Proposition \ref{kn1prop} that there exists an 
  effective divisor $B$ passing through $Z'$ such that  
  $B^2 \geq -2$, $h^1(B)=0$ and 
\[ 2B^2 \leq B.L \leq B^2 +n'+1 \leq 2n'+2.\]

  If $B^2 \geq 0$, we would get that $B$ induces a Clifford index 
  $c_B \leq n'-1 \leq c-1$ on the smooth curves in $|L|$, a contradiction.

  So $B^2 =-2$, and $B$ is necessarily supported on a union of smooth rational curves, since $h^1(B)=0$. But $B.L \leq n'-1 \leq n-1 =c-1$ and $Z'$ has length $ \geq 2$, so $B$ passes through some of the base points of $\{ D_{\lambda} \}$. This contradicts the fact that we have chosen these base points to lie outside of smooth rational curves of degree $\leq c+2$ with respect to $L$.

  To prove the second assertion, we will show that if there is 
  a point $x_0 \in S'- \{x_1, \ldots ,x_n,y_1 \}$ such that 
  $S'$ has an $(n+2)$-secant 
  $n$-plane through $x_0$ and $Z_{1,\lambda}$, then $x_0 = y_2$. By symmetry,
  this will suffice. 
  
  As usual choose any 
  preimage $p_0$ of $x_0$ and denote by $X_{1,\lambda}$ the zero-dimensional 
  scheme defined by $p_0$ and $Z_{1,\lambda}$. We then have that
  the natural map
\[ H^0(L) \hpil H^0 (L \* \O_{X_{1,\lambda}}) \]
  is not surjective for any $\lambda$.
  
  As usual let Let $X'_{1,\lambda} \sub X_{1,\lambda}$ be a subscheme of 
  length $n'_{1,\lambda}+2 \leq n+2$, for 
  some integer $n'_{1,\lambda} \geq 0$, such that the map 
  $H^0(L) \khpil H^0 (L \* \O_{X'_{1,\lambda}})$
  is not surjective, but such that the map 
  $H^0(L) \khpil H^0 (L \* \O_{X''_{1,\lambda}})$ 
  is surjective for all proper subschemes 
  $X''_{1,\lambda} \subsetneqq X'_{1,\lambda}$.

  Since $n=D^2=c \geq 2$ and $L^2=4c+4=4(n+1)$, we again have by 
  Proposition \ref{kn1prop} that there for each $\lambda$ exists an effective divisor $B_{1,\lambda}$ passing 
  through $X_{1,\lambda}$ and satisfying $B_{1,\lambda}^2 \geq -2$, $h^1(B_{1,\lambda})=0$ and the 
numerical 
  conditions
\[ 2B_{1,\lambda}^2 \stackrel{(a)}{\leq} L.B_{1,\lambda} \leq   
        B_{1,\lambda}^2 + n'_{1,\lambda} +2 \stackrel{(b)}{\leq} 2n'_{1,\lambda}+4, \]    
  with equality in (a) or (b) implying $L \sim 2B_{1,\lambda}$. 

  Assume first that $B_{1,\lambda}^2 =-2$ for some $\lambda$. We then get the same contradiction on the choice of the base points of $\{ D_{\lambda} \}$, since $B_{1,\lambda}.L \leq n'_{1,\lambda} \leq n =c$.

So we must have $B_{1,\lambda}^2 \geq 0$ for all $\lambda$.
Then $n'_{1,\lambda}=n$, $X' _{1,\lambda}=X_{1,\lambda}$, $L.B_{1,\lambda}= B_{1,\lambda}^2 + n +2$, and $B_{1,\lambda}$ is a Clifford 
  divisor. The  moving part 
  $B_{1,\lambda}'$ of $|B_{1,\lambda}|$ is then a free Clifford divisor, so by condition (C7) we have that either 
$B_{1,\lambda}' \sim D$ or there exists a free Clifford divisor
$P_{1,\lambda}$ such that $B_{1,\lambda}' \leq P_{1,\lambda} < D$ with
the last inclusion as described in Proposition \ref{C6ref}, with the
additional property that $|L-P_{1,\lambda}|$ has no fixed divisor, by
the conditions (C6) and (C7).
  
We will show that this latter case cannot occur.

We have that $B_{1,\lambda}$ passes through $X_{1,\lambda} $. Now a (possible) base divisor in $|B_{i,\lambda}|$ cannot pass through any of the points $p_0, \ldots ,p_n$, since these points lie outside all the rational curves contracted by $L$. So we must have $B_{i,\lambda}'.D \geq n$.

In addition, by Proposition \ref{C6ref} we must have
\[ D \sim B_{1,\lambda}' + \gamma_{1,\lambda}, \]
for some $\gamma_{1,\lambda} > 0$ satisfying $\gamma_{1,\lambda} ^2 =-2$, and $B_{1,\lambda}'.\Delta=0$. Hence
\[ B_{1,\lambda}'.L = B_{1,\lambda}'.(L - \Delta) = 2B_{1,\lambda}'.D \geq 2n, \]
so that 
\[ {B_{1,\lambda}'} ^2 \geq n-2 = D^2-2. \]
Since $h^1 (B_{1,\lambda}') = h^1 (P_{1,\lambda}) = h^1 (D) =0$, we 
have $h^0(B_{1,\lambda}') \geq \frac{1}{2}D^2+1 \geq h^0(P_{1,\lambda})$, so 
$B_{1,\lambda}' \sim P_{1,\lambda}$.

Since $|L-P_{1,\lambda}|=|L-B_{1,\lambda}|$ has no fixed divisor, we have
$B_{1,\lambda}' \sim B_{1,\lambda}$, so 
$B_{1,\lambda} < D$ and $B_{1,\lambda}$ is a free Clifford divisor. Since $h^0 (B_{1,\lambda} \* \I_{X_{1,\lambda}}) >0$,
there must exist an element of $|D|$ of the form $B_{1,\lambda} + A_{1,\lambda}$ passing through $Z_{1,\lambda}$, for $A_{i,\lambda} >0$. But since there is only one element of $|D|$ passing through $p_1, \ldots ,p_n,q_{1,\lambda}$, which we called 
   $D_{\lambda}$ and which is smooth and irreducible, we have $B_{1,\lambda}=D_{\lambda}$, a contradiction.

So we must have $B_{1,\lambda}' \sim D$. 
  By Proposition \ref{kn1prop} either $L-B_{1,\lambda} \geq B_{1,\lambda}$, or 
  both $h^0 (B_{i,\lambda} \* \I_{X_{1,\lambda}} ) \not =0$ and $h^0 ((L-B_{1,\lambda}) \* \I_{X_{1,\lambda}}) \not =0$.
  This gives us the two possibilities:
          \begin{itemize}
   \item [1.] $B_{1,\lambda} \in |D|$, 
   \item [2.] $B_{1,\lambda} \in |D| + 
          \Gamma_{j(\lambda)}$, for $j(\lambda)=1$ or $2$, and there 
              exists an $F_{1,\lambda} \in |D| + \Gamma_{3-j(\lambda)}$ passing 
              through $X_{1,\lambda}$.
          \end{itemize}

In case 1., since there is only one 
   member of $|D|$ containing $p_1, \ldots ,p_n,q_{1,\lambda}$, which we called 
   $D_{\lambda}$, we have $B_{1,\lambda}=D_{\lambda}$. But this would mean 
   that $p_0 \in D_{\lambda}$ for all $\lambda$, a contradiction.

   In case 2. one easily sees that the only option is $p_0 \in \Gamma_2$, which means that $x_0=y_2$, as desired.
\hspace{0.07cm} $\square$ \end{proof}

\renewcommand{\proofname}{Proof}

The proof of Theorem \ref{mainsing} in the case (E2) is similar, and therefore
left to the reader.

Since we have seen that the crucial point in proving Propositon
\ref{h1delta} is to prove that $h^1(A)=0$, we get the following result
(by checking that the proof of Theorem \ref{mainsing} goes through):

\begin{Lemma} \label{pseudo-perfect}
  Let $D$ be a free Clifford divisor, not of type (E1) or (E2). 
If $h^1(A)=0$, then $D$ is perfect.
\end{Lemma}

\section{Projective models in smooth scrolls}
\label{resol}

Let $D$ be a free Clifford divisor on a non-Clifford general polarized $K3$ 
surface $S$. Assume that $\T=\T(c,D)=\T(c,D,\{ D _{\lambda} \})$ is 
smooth. This is equivalent to the conditions $D^2=0$ and 
$\R_{L,D}= \emptyset$ when $D$ is perfect. In any
case these two conditions are necessary to have $\T$ smooth, so $|D|$ has
projective dimension $1$ and the pencil $D _{\lambda}$ is uniquely determined.
We recall that $\varphi_L(S)$ is denoted by $S'$.

Since $\T$ is smooth, it 
can be identified with the $\PP ^1$-bundle $\PP (\E)$, where  
$\E = \O_{{\PP}^1}(e_1) \+ \O_{{\PP}^1}(e_2) \+  \cdots  \+ 
\O_{{\PP}^1}(e_{c+2})$, and  $(e_1,e_2, \ldots .,e_{c+2})$ is the type of the 
scroll.

We will construct a resolution of the structure sheaf 
$\O_{S'}$ as an $\O_{\T}$-module.

The contents in this section will be very similar to that in \cite{Sc},
where canonical curves of genus $g$ are treated. This is quite natural, since 
a general hyperplane section of $S'$ is indeed such a canonical curve.

The following are well-known facts about $\T$ in ${\PP}^g$ (see \cite{H} and 
\cite{E-H}):
\begin{itemize}
\item[(1)] $\deg \T  = g-c-1$.
\item[(2)] $\dim \T= c+2$.
\item[(3)] The Chow ring of $\T$ is $\ZZ [\H,\F]/(\F^2,\H^{c+3},\H^{c+2}\F, 
\H^{c+2} - (g-c-1)\H^{c+1}\F)$, where  $\H$\index{$\H$} is the
hyperplane section, and
$\F$\index{$\F$} is the class of the ruling.
\item[(4)] The canonical class of $\T$ is $-(c+2)\H + (g-c-3)\F$.
\item[(5)] The class of $S'$ in the Chow ring of $\T$ is $(c+2)\H^c + 
(c^2+3c-cg)\H^{c-1}\F$.
\end{itemize}

We will need the Betti-numbers of the $\varphi_L(D_{\lambda})$ in ${\PP}^{c+1}$. 
These can be found also when $\T$ is singular, and will be needed in this case later on.

\begin{Lemma} \label{Betti0}
Let $(S,L)$ be a polarized $K3$ surface of genus $g$ and of non-general 
Clifford index $c >0$. Let $D$ be a free Clifford divisor satisfying $D^2=0$. 
For $c \geq 2$, all the $\varphi_L(D_{\lambda})$ in ${\PP}^{c+1}$ have minimal resolutions
\begin{eqnarray*}
0 \hpil \O_{\PP^ {c+1}}(-(c+2)) \hpil \O_{\PP^{c+1}}(-c)^{\beta _{c-1}} \hpil \O_{\PP^{c+1}}(-(c-1))^
{\beta _{c-2}} \hpil  \\
\cdots \hpil \O_{\PP^{c+1}}(-3)^{\beta _2} \hpil \O_{\PP^{c+1}}(-2)^{\beta _1} 
\khpil  \O_{\PP^{c+1}} \khpil \O_ {\varphi_L(D_{\lambda})} \hpil 0,
\end{eqnarray*}
where \index{${\beta }_i$}\[{\beta }_i = i {c+1 \choose i+1} - {c \choose i-1}.\]

For $c=1$ all the $\varphi_L(D_{\lambda})$ in ${\PP}^2$ have the resolution
\[ 0 \hpil \O_{\PP ^2}(-3) \hpil \O_{\PP ^2} \hpil \O_{\varphi_L(D_{\lambda})} \hpil 0.\]
\end{Lemma}

\begin{proof}
 Pick any $D'_{\lambda} := \varphi_L(D_{\lambda})$. We will show in Proposition 
\ref{arnorm} below that any such $D'_{\lambda}$ is arithmetically normal, whence
 projectively
Cohen-Macaulay, since the $D_{\lambda}$ have pure dimension one. Then
its Betti-numbers (see p.~\pageref{Betti} below for the definition) are equal 
 to those of a general hyperplane section of it. It is sufficient that the 
 linear term defining the hyperplane is not a zero divisor in its coordinate 
 ring $R_{\lambda}$. This is 
essentially \cite[Theorem 27.1]{N}.

Now choose a sufficiently general hyperplane $H_{\lambda}$ in $\PP^g$ so that
$C_{\lambda} := H_{\lambda} \cap S'$ is a smooth canonical curve, $H_{\lambda}$ does not contain 
any of the linear spaces  $\overline{ D_{\lambda} }$, and the 
hyperplane section $A_{\lambda} := H_{\lambda} \cap D_{\lambda}$ is not a 
zero divisor of $R_{\lambda}$. 

We can identify $C_{\lambda}$ with an element in $|L|$, and by abuse of notation write
$\O_{C_{\lambda}} (A_{\lambda}) = \O_{C_{\lambda}} (D_{\lambda}) = \O_{C_{\lambda}} (D)$. This linear 
system is 
complete and base point free (in fact it is a pencil computing the gonality) 
of degree $c+2$ on $C_{\lambda}$. By \cite[Lemma p.119]{Sc} (where there is a misprint)
and \cite[Proposition 4.3]{Sc} the zero-dimensional scheme 
$A_{\lambda}$ then has the Betti-numbers 
${\beta}_{i,i+1}={\beta}_i = i {c+1 \choose i+1} - {c \choose i-1}$. 

In particular, these numbers are independent of $\lambda$.
\hspace{0.07cm} $\square$ \end{proof}

The following result is analogous to \cite[Corollary (4.4)]{Sc}.

\begin{prop} \label{resolv}
  Let $S$ be a polarized $K3$ surface of non-general Clifford-index $c>0$, 
whose associated scroll $\T$ as above is smooth.
                 \begin{itemize}
\item [(a)] $\O_{S'}$ has a unique $\O_{\T}$-resolution $F_*$\index{$F_*$} 
           (up to isomorphism). If $c=1$, the resolution is:
\[0 \hpil \O_{\T}(-3\H+(g-4)\F) \hpil \O_{\T} \hpil \O_{S'} \hpil 0,\]
If $c \geq 2$, the resolution is of the following type:
\begin{eqnarray*}
0  & \hpil & \O_{\T}(-(c+2)\H+(g-c-3)\F)  \hpil \+ _{k=1}^{\beta _{c-1}} 
\O_{\T}(-c\H+b_{c-1}^k\F)  \hpil   \\
& \cdots  & \hpil   \+ _{k=1}^{\beta _1} \O_{\T}(-2\H + b_1^k \F)  \hpil 
\O_{\T} \hpil \O_{S'} \hpil 0,
\end{eqnarray*}
where ${\beta}_i = i {c+1 \choose i+1} - {c \choose i-1}$.
\item [(b)] $F_*$ is self-dual:
           $\mathcal{H}om (F_*,\O_{\T}(-(c+2)\H+(g-c-3)\F)) \iso F_*$.
\item [(c)] If all $b_i^k \geq -1$, then an iterated mapping cone
           \[  [[\C ^{g-c-3}(-(c+2)) \hpil \+ _{k=1}^{\beta_{c-1}} 
           \C ^{b_c^k}
           (-c)] \ldots ] \hpil \C ^0 \]
           is a (not necessarily minimal) resolution of $\O_{S'}$ as an 
           $\O_{\PP^g}$-module.
\item [(d)] The $b_i^k$ satisfy the following polynomial equation in $n$ if 
$c \geq 2$:
\begin{eqnarray*}
{n+c+1 \choose c+1}(\frac{n(g-c-1)}{c+2}+1) - n^2(g-1)-2  = & & \\
\sum _{i=1}^{c-1}((-1)^{i+1}{n-i+c \choose c+1}(\frac{((n-i-1)
(g-c-1)+(c+2))\beta _i}{c+2} + \sum _{k=1}^{\beta _i} b_i^k) & + & \\ 
(-1)^{c+1}{n-1 \choose c+1}(\frac{(n-c-2)(g-c-1)}{c+2} + g-c-2). &&  
\end{eqnarray*}
\end{itemize}
\end{prop}
 
\begin{proof}  
We start by proving (a). 
We have  $\overline{D_{\lambda}} \iso {\PP}^{c+1}$ by
(\ref{eq:span}). The $\varphi_L(D_{\lambda})$ have Betti-numbers 
$\beta ^{\lambda}_{i,j}=
\dim(Tor_i^{R_{\lambda}}(R,k)_j)$, where $R$ is the homogeneous coordinate 
ring of ${\PP}^{c+1}$, and $R_{\lambda}$ the coordinate ring   
$R/I_{\lambda}$ of $\varphi_L(D_{\lambda})$. 
Following \cite{Sc}, for $c \geq 2$ it is enough to prove:

\begin{itemize}
\item[(1)]For fixed $i,j$ the $\beta ^{\lambda}_{i,j}$ are the same for all 
${\lambda}$. 
\item[(2)] If $c \geq 2$, then $\beta ^{\lambda}_{i,j}=0$, unless $j=i+1$ 
and $i \leq c-1$, or $(i,j)=(c,c+2)$. 
\item[(3)] The common value $\beta _{i,i+1}=\beta ^{\lambda}_{i,i+1}$ is 
${\beta}_i = i {c+1 \choose i+1} - {c \choose i-1}$ for $i \leq c-1$, and
$\beta _{c,c+2}=1$.
\end{itemize}

This follows immediately from the lemma above.

The easier case $c=1$ is dealt with in an analogous manner.

The proof of (b) is almost identical to that of \cite[Corollary 4.4(ii)]{Sc}.
In our case we have $\mathcal{E}xt_{\T}^i(\O_{S'}, \omega _{\T}) = \omega _{S'}$ if 
$i=c$, and zero otherwise, $\omega _{S'}= \O_{S'}$, and $\omega _{\T}= \O_{\T}(-(c+1)\H+(g-c-3)\F)$.

The proof of (c) is identical to that of \cite[Corollary 4.4(iii)]{Sc}.

Denote the term $i$ places to the left of $\O _{\T}$ in the resolution $F_*$
by $F_i$. The proof of (d) then follows from the identity 

\[ \chi (\O_{\T}(n\H)) - \chi (\O_{S'}(n\H)) = \sum_i(-1)^i \chi (F_i(n\H)). \] 

The contribution from the $F_c$-term is written out separately.
Moreover it is clear that for all large $n$, we have 
$\chi (F_i(n\H))=h^0(F_i(n\H)$, for all $i$, and  
$\chi (\O_{\T}(n\H))=h^0(\O_{\T}(n\H))$
since $\H$ is (very) ample on $\T$.
Then one uses the following well-known fact for $a \geq 0$:
\begin{equation} \label{dyttned}
  h^0(\PP (\E),\O _{\PP(\E)} (a\H+b\F)) = h^0({\PP}^1, \Sym ^a(\E ) \* \O_{{\PP}^1}(b)).
\end{equation}
\hspace{0.07cm} $\square$ \end{proof}

\begin{rem} \label{sumb}
{\rm Part (d) of the proposition only gives us the sums of the $b_i^k$ for 
each fixed $i$. The values $n=2,3, \ldots ,c$ give enough equations to determine 
these sums. The duality of part (c) gives $\beta _i=\beta _{c-i}$, for 
$i=1, \ldots ,c-1$, and $i \ne c/2$, and after a possible renumeration
of the $b_i^k$, for 
$k=1, \ldots ,\beta _{i}$, we also have $b_{c-i}^k=g-c-3-b_i^k$ for these $k$.
In particular this  enables us to identify the sums of the 
$b_i^k$ with those of the  $g-c-3-b_{c-i}^k$. 
To obtain more information about the individual $b_i^k$ a more refined 
study is necessary. }
\end{rem}

\section{Projective models in singular scrolls}
\label{singscrolls}

 Let $D$ be a free Clifford divisor on a non-Clifford general polarized $K3$ 
surface $(S,L)$. In this section we will make a  thourough study of
the case where the scroll $\T=\T(c,D,\{ D _{\lambda} \})$ is singular.

\begin{rem} \label{singprat}
 {\rm As seen above, the scroll  $\T=\T(c,D,\{ D _{\lambda} \})$ is
   singular
if $D^2 >0$ or the set $\R _{L,D}$ is non-empty. Moreover $\T$ 
is singular if and only if one of these two conditions holds, if $D$ is 
perfect.

We will always assume $c>0$, so that $\varphi_L: S \khpil S'$ is birational.

The type $(e_1, \ldots ,e_d)$ of the scroll, where $d= \frac{1}{2}D^2 +c+2$, is such that the last $r$ of the $e_i$ are zero, where $r$ is defined as in equation (\ref{eq:r}) and can be computed as in equation (\ref{eq:sing}).

As we have seen, when $D$ is perfect we have
\begin{equation} 
  r = \left\{ \begin{array}{lll}
      D^2 -1   & \mbox{ if $D$ is of type (Q)}, \\
      D^2 + D.\Delta-1   & \mbox{ if $D$ is of one of the types (E0)-(E2)}, \\
      D^2 + D.\Delta     & \mbox{ otherwise.}
          \end{array}
    \right .     
\end{equation}

We will however not assume that $D$ is perfect, unless explicitly
stated. }
\end{rem}

Let $n:=D^2$ and denote by $p_1, \ldots, p_n$ the $n$ base points of the pencil
$\{ D_{\lambda} \}$. Let 
\index{$\tilde S$}\[ 
\xymatrix{
    \tilde S \ar[r]^{f}   &  S 
} \]
be the blow up of $S$ at $p_1, \ldots, p_n$. Denote by $E_i$ the exceptional line over $p_i$ and let
\[ E := \sum_{i=1}^n E_i \]
denote the exceptional divisor. Define 
\index{$H$}\[ H := f^*L + f^*D -E. \]
The first observation is:

\begin{Lemma} \label{birS''}
  $H$ is generated by its global sections, $h^1(H)=0$ and $\varphi_H$ is birational; in fact $\varphi_H$ is an isomorphism outside of finitely many contracted smooth rational $(-2)$-curves.

Moreover, a smooth rational curve $\gamma$ is contracted by $H$ if and only if $\gamma = f^*\Gamma$, for some smooth rational curve $\Gamma$ on $S$ such that $\Gamma.L=\Gamma.D=0$.
\end{Lemma}

\begin{proof}
Since $H-E \sim (f^*L-E) + (f^*D-E)$ is clearly nef and $(H-E)^2 \geq 10$, we have
$h^1(H)=0$. Furthermore, since $|f^*D-E|$ is a base point free pencil and $f^*L$ is base point free, $H$ is base point free as well. 

The morphism given by $|f^*L|$ is clearly an isomorphism outside of the $n$ exceptional curves and the strict transforms of the finitely many smooth rational curves on $S$ which are contracted by $|L|$. By our choice of pencil (see \ref{eq:degree}), these curves do not intersect the $n$ blown up points.

Since $E_i.H =1$ for all $i$, every exceptional curve $E_i$ is mapped by $\varphi_H$ isomorphically to a line, so $\varphi_H$ is an isomorphism along the exceptional curves. Moreover if 
$\gamma = f^*\Gamma$ for some smooth rational curve $\Gamma$ on $S$ such that $\Gamma.L=0$ and 
$\Gamma.D >0$ then $\Gamma.D=1$ by Lemma \ref{propertiesR}(c), so $\gamma$ is mapped isomorphically to a line by $\varphi_H$ and $\varphi_H$ is an isomorphism along these curves as well.

Hence $\varphi_H$ is an isomorphism outside of finitely many contracted smooth rational $(-2)$-curves, which are precisely the ones of the form $f^*\Gamma$, for some smooth rational curve $\Gamma$ on $S$ such that $\Gamma.L=\Gamma.D=0$. 
\hspace{0.07cm} $\square$ \end{proof}

We have $h^0(H)= \frac{1}{2}H.(H-E)+2 =
\frac{1}{2}L^2+\frac{1}{2}D^2+c+4= g+
\frac{1}{2}D^2+c+2=g+d+1$. Denote by  $S''$\index{$S''$} the surface $\varphi_H(\tilde S)$ in $\PP^{g+ d}$.

One easily obtains $\deg S''=2g+2c+2+2D^2$.

\begin{prop} \label{partialblowup}
The surface $S''$ is normal, $p_a(S'')=1$, and $K_{S''} \iso \O_{S''}(E')$, 
where $E'$ is the sum of $D^2$ lines that are $(-1)$-curves on $S''$.
\end{prop}

\begin{proof}
The two last assertions are immediate consequences of $S''$ being normal, by 
\cite{Ar}.

Consider the blow-up $f: \tilde S \khpil S $ described above.

Denote by $\E _H$\index{$\E _H$} the set of irreducible curves $\tilde \Gamma$ on $\tilde S$ such that $\tilde \Gamma.H=0$. From the Hodge Index theorem it follows that such a curve has negative self-intersection. Moreover, by Lemma \ref{birS''}
\[ \tilde \Gamma = f^* \Gamma, \]
for some smooth rational curve $\Gamma$ on $S$ such that $\Gamma.L=\Gamma.D=0$.
Thus we can write
\[ \E _H = f^* (\E _L - \R _{L,D}). \]

Now let $\tilde \delta$ be the fundamental cycle of a connected component of 
$\E _H$, $p$ the image of $\tilde \delta$ on $S''$ and $U$ the inverse image of an affine open neighborhood of $p$. To prove the normality of $p$ it will be sufficient to prove the surjectivity of
\[ H^0 (U, \O_U (H-\tilde \delta)) \hpil H^0 (\tilde \delta, \O_{\tilde \delta}(H-\tilde \delta)), \]
hence of
\[ H^0 (\tilde S, \O_{\tilde S} (H-\tilde \delta)) \hpil H^0 (\tilde \delta, \O_{\tilde \delta}(H-\tilde \delta)). \]
To show the latter, it will suffice to show
\[ H^1 (\tilde S, \O_{\tilde S} (H-2\tilde \delta))=0. \]
By the degeneration of the Leray spectral sequence
\[ 0 \hpil H^1 (S, f_* (H-2\tilde \delta)) \hpil  
     H^1(\tilde S, H-2\tilde \delta) \hpil H^0(S, R^1 f_* (H-2\tilde \delta)) \]
it will suffice to show that 
\[ h^1 (S, f_* (H-2\tilde \delta)) = h^0 (S, R^1 f_* (H-2\tilde \delta)) =0.\]

Denote by $\delta$ the divisor on $S$ such that $f^* \delta = \tilde \delta$. Then $\delta$ is connected and $\delta^2 =-2$ ($\delta$ is in fact a fundamental cycle for a connected component of $\E _L$ minus a curve $\Gamma$ that is a 
tail of $\delta$ and is such that $\Gamma.D=1$. The fact that $\delta ^2=-2$
can be checked by inspection for each of the five platonic configurations \cite{Ar}). 
We then have
\[ f_* (H-2\tilde \delta) = (L+D-2\delta) \* \I _Z, \]
where $Z$ is the zero-dimensional scheme corresponding to the $n$ blown up points, and
\[ R^1 f_* (H-2\tilde \delta) = R^1 f_* (-E) \* (L+D-2\delta). \]
Since $f_* \O _E \iso \O _Z$, we have $R^1 f_* (-E)=0$, whence we are reduced to proving the vanishing of $H^1 ((L+D-2\delta) \* \I _Z)$. This will be proved in Lemma \ref{normallemma} below. 
\hspace{0.07cm} $\square$ \end{proof}

\begin{Lemma}  \label{normallemma}
 With the notation as above, $H^1 ((L+D-2\delta) \* \I _Z)=0$.
\end{Lemma}

\begin{proof}
 We will first need the following fact:
\[ h^1 (L+D-2\delta)=0. \]
 The proof for this is rather long and tedious, but does not involve any new 
 ideas and is similar in principle to the proof of \cite[Lemma 5.3.5]{cos}. 
 We therefore leave it to the reader.

 Note that if $D^2=0$, then $Z = \emptyset$, and we are done. So we will from now on assume that $n=D^2 >0$.

 Because of the vanishing of $H^1 (L+D-2\delta)$, the vanishing of $H^1 ((L+D-2\delta) \* \I _Z)$ is 
 equivalent to the surjectivity of the map
\[ H^0 (L+D-2\delta) \hpil H^0 ((L+D-2\delta) \* \O _Z). \]

Assume, to get a contradiction, that this map is not surjective. Let $Z' \sub Z$ be a subscheme of length $l+1 \leq n=D^2$, for some integer $l \geq 1$, such that $H^0 (L+D-2\delta) \hpil H^0 ((L+D-2\delta) \* \O _{Z'})$ is not surjective, but $H^0 (L+D-2\delta) \hpil H^0 ((L+D-2\delta) \* \O _{Z''})$ is for all proper subschemes $Z''$. 

Since $(L+D-2\delta)^2 > 4l + 4$ and $h^1 (L+D-2\delta)=0$, we get by Remark \ref{kn1rem} that there is an effective decomposition $L+D-2\delta \sim A+B$ such that $A > B$, $A.B \leq l+1$, $h^1 (B)=0$, $B^2 \geq -2$ and $B$ passes through $Z'$.

If $B^2 =-2$ (so that $B$ is necessarily supported on a union of smooth 
rational curves), then we use the fact that we have chosen $Z$ to lie outside of any rational curve $\Gamma$ such that $\Gamma.L \leq c+2$ by (\ref{eq:degree}) and 
\[ L.B \leq (L+D).B \leq l-1 + 2\delta.B \leq D^2-2+2\delta.B \leq c + 2\delta.B, \]
to conclude that we must have $\delta.B \geq 2$. Hence $(\delta+B)^2 \geq 0$.

This yields that we in all cases have
\[ h^0 (\delta+B) \geq 2. \]

We now want to show that also
\[ h^0 (\delta+A-D) \geq 2. \]

We can write 
\[ F \sim A +B + 2\delta -2D \sim (A+\delta-D) + (B+\delta-D) := F_1 + F_2. \] 
This is not necessarily an effective decomposition, but we have $F_1 > F_2$, since 
$A > B$.

We can easily calculate
\[ F_1.F_2 = A.B -D^2-c \leq -c <0, \]
and since $F^2 = {F_1}^2 + {F_2}^2 + 2F_1.F_2 \geq D^2 \geq 2$, we must have
${F_1}^2 \geq 2$ or ${F_2}^2 \geq 2$.

If ${F_1}^2 \geq 2$, then either $h^0 (F_1) \geq 2$ or $h^0 (-F_1) \geq 2$ by Riemann-Roch. Since $L.F_1=L.A-L.D > \frac{1}{2}(L^2+L.D) - L.D = 
\frac{1}{2}(c+2+F^2) >0$, we must have $h^0 (F_1) \geq 2$, and we are done.

If ${F_2}^2 \geq 2$, then either $h^0 (F_2) \geq 2$ or $h^0 (-F_2) \geq 2$. In the first case, we get $h^0 (F_1) \geq h^0 (F_2) \geq 2$. In the second, we get
$F_1 \sim F - F_2 > F$, since $-F_2$ is effective, whence $h^0 (F_1) \geq h^0 (F) \geq 2$ again.

So we have an effective decomposition of L as
\[ L \sim (B +\delta) + (A+ \delta-D), \]
such that both $h^0 (B +\delta)$ and $h^0 (A+ \delta-D) \geq 2$ and such that
\[ (B +\delta).(A+ \delta-D)= A.B -D.B +2 \leq l - D.B +3. \]
Since $l+1 \leq D^2 \leq c+2$, and $D.B \geq 2$, since $D$ is base point free and $h^0(B +\delta) \geq 2$, we must have $D.B=2$ and $l+1=n=D^2=c+2$. But since
$B$ passes through $Z$, we must have $D.B \geq n$, whence the contradiction $c=0$.

This concludes the proof of the lemma and hence of Proposition 
\ref{partialblowup}.
\hspace{0.07cm} $\square$ \end{proof}

Define the following line bundle on $\tilde S$:
\[ \tilde D := f^*D-E. \]
The members of $| \tilde D |$ are in one-to-one correspondence with the members of the pencil $ \{ D_{\lambda} \}$. One computes ${\tilde D}^2=0$, so 
$| \tilde D |$ is a pencil of disjoint members. Furthermore 
\[ h^0(H- \tilde D) = h^0 (f^*L) >2, \]
so $| \tilde D |$ defines a rational normal scroll  $\T_0$\index{$\T_0$} containing $S''$.

\begin{prop} \label{typet0}
  $\T _0$ has dimension $d$ and degree $g+1$ and is smooth of type 
  $(e_1+1, \ldots , e_d+1)$.
\end{prop}

\begin{proof}
  The two first assertions are easily checked.

We have to calculate the numbers 
$h^0(\tilde S, H-i\tilde D)= h^0(\tilde S, f^*(L-(i-1)D) + (i-1)E)$ for all $i \geq 0$.

 One easily sees that $(i-1)E$ is a fixed divisor in $|f^*(L-(i-1)D) + (i-1)E|$ for all $i \geq 1$, so we get for all $i \geq 1$:
 \begin{eqnarray} 
\label{eq:zidane} h^0(\tilde S, f^*(L-(i-1)D + (i-1)E) & = & h^0(\tilde S, f^*(L-(i-1)D) \\
 \nonumber                                           & = & h^0(S, L-(i-1)D).   
 \end{eqnarray}

We also have 
\begin{equation} \label{eq:inzaghi}
h^0(H)-h^0(H-\tilde D) =d.  
\end{equation}
Defining $d'_i:= h^0(\tilde S, H-i\tilde D)- h^0(\tilde S, H-(i+1)\tilde D)$, we get by combining (\ref{eq:zidane}) and (\ref{eq:inzaghi}) that
\[ d'_0 =d_0 \hs \mbox{ and } \hs d'_i =d_{i-1} \hs \mbox{ for } \hs i \geq 1. \]
It follows immediately that the type of $\T_0$ is as claimed.
\hspace{0.07cm} $\square$ \end{proof}

Since $\T_0$ is smooth, we have $\T_0 \iso \PP (\E)$, where 
$\E = \+ _{i=1}^d \O _{\PP^1} (e_i+1)$. Also, we have the maps
\[ 
\xymatrix{
   \PP (\E) \ar[r]^{j} \ar[d]^{\pi} & \T_0   \sub \PP ^{g+d}   
  \\ 
   \PP ^1  & 
} \]
where $j$ is an isomorphism. Then the Picard group of 
$\PP (\E)$ satisfies
\[ \Pic \PP (\E) \iso \ZZ \H_0 \+ \ZZ \F, \]
where $\H_0 := j^{*} \O _{\PP ^{g+d}} (1)$\index{$\H_0$} and  
$\F := \pi ^{*} \O _{\PP ^1} (1)$\index{$\F$}.

Furthermore, the Chow ring of $\PP (\E)$ is
\begin{equation} \label{chowring}
 \ZZ [\H_0,\F]/(\F^2,\H_0^{s+2},\H_0^{s+1}\F, 
\H_0^{s+1} - (g+1)\H_0^s\F),
 \end{equation}
where we set $s:=c+1+\frac{1}{2}D^2$.

Consider now the morphism $i$ given by the base point free line bundle 
$\H := \H_0 - \F$\index{$\H$}, where $\H_0 = \H+\F$:
\[ i : \PP (\E) \hpil \PP ^{g}. \]
One easily sees that $i$ maps $\PP (\E)$ onto a rational normal scroll of dimension $d$ and type $(e_1, \ldots , e_d)$, whence isomorphic to $\T$. So we can assume that $i$ maps $\PP (\E)$ onto $\T$. By abuse of notation we write
\[ i : \T _0 \hpil \T, \]
and this is a rational resolution\index{rational resolution} of singularities of $\T$ (in the sense that
$\T_0$ is smooth and $R^1i_{*}\O_{\T_0}=0$). Furthermore one easily sees that 
by construction $i$ restricts to a map
\[ g : S'' \hpil S' \]
which is a resolution of some singularities of $S'$ (precisely the singularities of $S'$ arising from the contractions of rational curves across the fibers in $S$, i.e. the curves in $\R_{L,D}$) and a blow up at the images of the base points of $\{ D _{\lambda} \}$. 

We get the following commutative diagram: \label{commutativepage} 
\[ 
\xymatrix{
    \PP (\E) \ar[r]^j & \T _0  \ar[r]^i  &  \T  \ar[r]  &\PP ^g  
  \\ 
       & S''  \ar[r]^g \ar[u]  &  S' \ar[u]  &
  \\  
       & \tilde S  \ar[u]^{\varphi _H}\ar[r]^{f}   &  S \ar[u]_{\varphi _L}. 
} \]
By construction, one has $ g \circ \varphi_H = \varphi _{f^*L}$.

\begin{prop} \label{ideal}
  Let $\I_{S''/\T_0}$ denote the ideal sheaf of $S''$ in $\T_0$ and
 $\I_{S'/\T}$ the ideal sheaf of $S'$ in $\T$.

  We have $\I_{S'/\T} = i_* \I_{S''/\T_0}$.
\end{prop}

\begin{proof}
  This follows since $i_* \O_{\T_0} =\O_{\T}$ and $i_* \O_{S''}=
  \O_{S}$. The latter fact is a consequence of $g$ being a birational
  map of normal surfaces.
\hspace{0.07cm} $\square$ \end{proof}

We recall that the Chow ring of $\T_0$ is given by (\ref{chowring}). Define $\H_{\T}$\index{$\H_{\T}$} and $\F_{\T}$\index{$\F_{\T}$} to be the push-down of cycles by $i$ of $\H$ and $\F$ respectively.

We have the following description:

\begin{prop} \label{Chowring}
  \begin{itemize}
\item [(a)] The class of $S''$ in the Chow ring of $\T_0$ is 
\[(D^2+c+2)\H_0 ^{d-2} + (c-cg-D^2(g-1))\H_0^{d-3}\F.\]
\item [(b)] The class of $S'$ in the Chow group of $\T$ is 
\[(D^2+c+2)(\H_{\T})^{d-2} + (D^2(d-1-g)-4-cg-c+cd+2d)(\H_{\T})^{d-3}\F_{\T}.\]
  \end{itemize}
\end{prop}

\begin{proof}
The class of $S''$ is of the type 
$m\H_0^{d-2}+n\H_0^{d-3}\F$, for two integers $m$ and $n$. To determine $m$
and $n$ one has the equations $S''\H_0^2 = \deg S'' = 2g+2c+2+2D^2$ and 
$S''\H_0\F=\deg( \varphi _L(D))=c+2+D^2$.

Statement (b) is an immediate consequence of $i$ being birational by using the cap product map $A^* (\PP^g) \* A_* (\T) \khpil A_* (\T)$.
\hspace{0.07cm} $\square$ \end{proof}

We would like to study the resolution of $S''$ in $\PP (\E) \iso \T_0$. 
We say that $S''$  has constant Betti-numbers $\beta_{ij}= \beta_{ij} (\lambda)$ over $\PP ^1$ if the one-dimensional schemes  obtained by intersecting $S''$ by the linear spaces $F_{\lambda}$ in the pencil of fibres of $\T_0$ have Betti-numbers 
in $\PP^{c+1+\frac{1}{2}D^2}$ that are independent of $\lambda$.
By \cite{Sc}, if $S''$  has constant Betti-numbers  over $\PP ^1$, we can 
(at least in principle) find a resolution of $\O_{S''}$ by free 
$\O_{\PP (\E)}$-modules which restricts to the minimal resolution of 
$\O _{S''_{\lambda}}$ on each fiber 
$\PP (\E) _{\lambda} \iso \PP ^{c+1+\frac{1}{2}D^2}$. 

Clearly, since the map $i$ is the identity on each fiber, the Betti-numbers of $S''_{\lambda}$ are the same as the Betti-numbers of $\varphi_L (D _{\lambda})$.

Recall that a projective scheme $V$ is called arithmetically normal\index{arithmetically normal} if the natural map
\[ S_kH^0(V,\O_V(1)) \hpil H^0(V,\O_V(k))\]
is surjective for all $k \geq 0$.

We start by showing that the $\varphi_L (D _{\lambda})$ are all arithmetically 
normal.

\begin{prop} \label{arnorm}
  All the $\varphi_L (D _{\lambda})$ are arithmetically normal in 
$\overline{D_{\lambda}} = \linebreak \PP ^{c+1+\frac{1}{2}D^2}$.
\end{prop}

\begin{proof}
  We can easily show that 
\begin{equation} \label{waw}
  h^1(\O _S (qL-D))=0  \mbox{ for all } q.
\end{equation}
Furthermore, by \cite[Thm. 6.1]{S-D}, we have that
\begin{equation} \label{wow}
  S_k H^0 (S,L) \hpil H^0 (S,kL) \mbox{ is surjective for all } k \geq 0.
\end{equation}
We have a commutative diagram
\[ 
\xymatrix{
    H^0 (\O_{\PP ^g}(q)) \ar[r] \ar[d]^{\alpha _1} & H^0 (\O_{\overline{D_{\lambda}}}(qL))  \ar[d]^{\alpha _3}    
  \\ 
    H^0 (\O_S (qL))  \ar[r]^{\alpha _2}   &  H^0 (\O_{D_{\lambda}}(qL)). 
} \]
Now $\alpha _2$ is surjective by (\ref{waw}) and $\alpha _1$ is surjective by (\ref{wow}). Hence $\alpha _3$ is surjective and $\varphi_L (D_{\lambda})$ is arithmetically normal.
\hspace{0.07cm} $\square$ \end{proof}

For each $\lambda \in \PP ^1$ define
\index{$B ^{\lambda}$}\index{$V ^{\lambda}$}\[ B ^{\lambda} := \+ _{q \in \ZZ} H^0 (D_{\lambda}, qL) \hs \mbox{   and   } 
   \hs V ^{\lambda} := H^0 (D_{\lambda}, L). \]

The symmetric algebra $S(V ^{\lambda})$\index{$S(V ^{\lambda})$} of $V ^{\lambda}$ satisfies
\[ S(V ^{\lambda}) \iso R_{\lambda}, \]
where $R_{\lambda}$\index{$R_{\lambda}$} is the homogeneous coordinate ring of 
$\PP (H^0((D_{\lambda}, L) \iso \PP ^{c+1+\frac{1}{2}D^2}$, and $B ^{\lambda}$ is a graded 
$R_{\lambda}$-module. Since all the $R_{\lambda}$ are isomorphic, we will sometimes suppress the $\lambda$, hoping to cause no confusion.

We have the {\it Koszul complex}\index{Koszul complex}
\[ 
\xymatrix{
  \cdots  \ar[r]  & \wedge ^{i+1} V ^{\lambda} \* B ^{\lambda} _{j-1} 
                                           \ar[r] ^{d^{\lambda} _{i+1,j-1}} &
\wedge ^i V ^{\lambda} \* B ^{\lambda} _j   \ar[r] ^{d^{\lambda} _{i,j}} &  \cdots 
} \]
with the {\it Koszul cohomology groups}\index{Koszul cohomology group} defined by
\[ \K^{\lambda} _{i,j} :=\K _{i,j} (B ^{\lambda}, V ^{\lambda}) := 
\frac{ \ker d^{\lambda} _{i,j}}{ \im d^{\lambda} _{i+1,j-1}}. \]

For each $\lambda$ we have a minimal free resolution\index{minimal free resolution} of $B ^{\lambda}$ as an 
$R_{\lambda}$-module:
\begin{eqnarray*}
  \cdots  & \hpil &  \+ _{j} R_{\lambda}(-j) ^{\beta ^{\lambda}_{i,j}} \hpil 
\cdots  \hpil \+ _{j} R_{\lambda}(-j) ^{\beta ^{\lambda}_{1,j}} \\
          & \hpil &  \+ _{i} R_{\lambda}(-j) ^{\beta ^{\lambda}_{0,j}} \hpil B ^{\lambda} \hpil 0,
\end{eqnarray*}
and the $\beta ^{\lambda}_{i,j}$ are the {\it (graded) Betti-numbers}\index{Betti-numbers} \label{Betti} for 
$\varphi_L(D_{\lambda})$ (since $\varphi_L(D_{\lambda})$ is arithmetically normal).

By the well-known {\it Syzygy Theorem}\index{Syzygy Theorem} \cite[Thm. (1.b.4)]{gre}, we have
\[ \beta ^{\lambda}_{i,i+j} = \dim \K ^{\lambda}_{i,j} \]
(where the dimension is as vector space over $\CC$).

\begin{exa} \label{zeroint}
{\rm As an example we look at the case where $D^2=0$ and $\T$ is
  singular 
(i.e. $D$ is not perfect or $\R_{L,D}$ is
non-empty). In this case the scroll $\T_0$ can be analyzed with the techniques
of Section \ref{resol}. Proposition \ref{partialblowup} gives that the
canonical sheaf on $S''$ is trivial. Lemma \ref{Betti0}
gives us the Betti-numbers of all 
the $D_{\lambda}$. Hence the analogue of Proposition  \ref{resolv}
goes through completely (we need the triviality of the canonical sheaf to
prove the analogue of part (b)) to give a resolution of $\O_{S''}$ as an
$\O_{\T_0}$-module. Set $g_0=g+d=g+c+2$. Since $\T_0$ has degree $g_0-c-1$, dimension
$c+2$, and spans $\PP^{g_0}$, we only need to replace $g$ by $g_0$ in Theorem 
\ref{resolv}. 
}
\end{exa}

Unfortunately, finding the Betti-numbers $\beta ^{\lambda}_{ij}$ for the 
$\varphi_L (D _{\lambda})$ when $D^2 >0$ is not as easy as in the case $D^2=0$.
In fact, we are not able to compute all of them, nor to show that
  they 
are constant over $\PP^1$, in general, but we will manage for the cases $D^2=2$ and $4$, which are the cases we need for the classification of projective models of genus $g \leq 10$.

By our choice the general element in the pencil $ \{ D_{\lambda} \}$ is smooth and irreducible, whence by 
Lemma \ref{embfib} also the general $\varphi_L (D _{\lambda})$ is a smooth irreducible curve. To compute its Betti-numbers in $\overline{D_{\lambda}} = \PP ^{c+1+\frac{1}{2}D^2}$, we can use several results of Green and Lazarsfeld, and it will turn out that these results are sufficient to determine its Betti-numbers uniquely
for $D^2 \leq 4$. However, there might be singular, reducible or even nonreduced elements in
the pencil $\{ \varphi_L (D _{\lambda}) \}$, and one then has to check that the results of Green and Lazarsfeld can still be applied to these cases. Roughly speaking, since the Betti-numbers do not change when taking general hyperplane sections (since all the $\varphi_L(D_{\lambda})$ are arithmetically normal whence projectively Cohen-Macaulay), we can avoid the isolated singularities, so the biggest problems arise from nonreduced fibers. It is therefore convenient to choose a pencil $ \{ D_{\lambda} \}$ with as few such cases as possible.
Also note that the existence of a reducible element in $ \{ D_{\lambda} \}$, will require the existence of
some effective divisors linearly independent of $L$ and $D$, so in the general case of every family we study, all elements in $ \{ D_{\lambda} \}$ will be reduced and irreducible.

It will be of use to us that we can choose a pencil $ \{ D_{\lambda} \}$ subject to the following additional condition when $D^2 >0$:

\begin{equation} \label{eq:redpencil}   
\mbox{Any member of $ \{ \varphi_L (D_{\lambda}) \}$ is 
one of the following:}  
\end{equation}

\begin{itemize}
\item A smooth irreducible curve of genus $p_a(D)$.
\item A singular irreducible curve of arithmetic genus $p_a(D)$ or $p_a(D)+1$ with exactly one node or one cusp.
\item $E_1+ E_2$, where $E_1$ and $E_2$ are distinct smooth elliptic curves intersecting in $\frac{1}{2}D^2$ points or in one point (the latter happening if and only if we are in the special case of Proposition \ref{cliffgenD}, where $D^2=4$, $L \sim 2D$ and $D$ is hyperelliptic).
\item $\overline{D} + \Omega$ with $\overline{D}$ a smooth  irreducible curve of genus 
$p_a(D)-1$ and $\Omega$ of degree $1$ or $2$
\item $\overline{D} + \Omega$ with $\overline{D}$ an  irreducible curve of genus $p_a(D)$ with exactly one node or cusp and $\Omega$ of degree $1$ or $2$.
\item[] (Note that $\Omega$ is either a conic, a union of two distinct lines, a double line or a line. In particular a nonreduced component of a member of $ \{ \varphi_L (D_{\lambda}) \}$ has to be a double line.)
\end{itemize}

\begin{Lemma} \label{redpencil}
  Let $D$ be a free Clifford divisor with $D^2 >0$. Then we can choose a pencil $ \{ D_{\lambda} \}$ such that (\ref{eq:redpencil}) is satisfied.
\end{Lemma}

\begin{proof}
Any irreducible element of $|D|$ is mapped isomorphically by $\varphi_L$ by Lemma \ref{embfib}.
Since the codimension of the set of irreducible elements in $|D|$ having more than one node or cusp as singularity is well-known to be $>1$, we can find a pencil so that all irreducible elements are mapped to irreducible curves which are either smooth of genus $p_a(D)$ or has at most one node or cusp and therefore have arithmetic genus $p_a(D)+1$.

Now we have to consider reducible elements of $|D|$ living in codimension one.

Assume that an element of $|D|$ has two components of arithmetic genus $\geq 1$. This means that $D \sim A +B$ with $h^0(A) \geq 2$ and $h^0(B) \geq 2$. A quick analysis as in the proof of 
Proposition \ref{cliffgenD} shows that $A^2=B^2=0$ (otherwise either $A$ or $B$ would induce a Clifford index $<c$ on $L$). So $D \sim E_1+E_2+\Sigma$ for $E_1$ and $E_1$ smooth elliptic curves and
an effective $\Sigma$ which is either zero or only supported on smooth rational curves. 
In the first case, since the general elements in both $|E_1|$ and $|E_2|$ are smooth elliptic curves, we can choose a pencil containing at most the union of two smooth elliptic curves $E_1$ and $E_2$.
Such a $D_0 =E_1+E_2$ is mapped isomorphically by $\varphi_L$ by Lemma \ref{embfib}.
In the second, as in the rest of the proof, we are reduced to studying the cases where
$B$ is an effective divisor on $S$ only supported on smooth rational curves 
such that $h^0(D \* \I _{B})=h^0(D-B)=h^0(D)-1 \geq 2$. 

By adding base divisors to $B$, we can assume that $|D-B|$ is base point free. Hence, by Proposition \ref{sd1} either $h^1(D-B)=0$, or $D-B \sim kE$, for an integer $k \geq 2$ and a smooth elliptic 
curve $E$.

In the first case we have $h^0(D-B)=\frac{1}{2}D^2-D.B+  \frac{1}{2}B^2+2=h^0(D)-D.B+  \frac{1}{2}B^2$.
If $B^2 >0$, then by the Hodge index theorem and the fact that $B^2 < D^2$ (since $h^0(B) < h^0(D)$) we 
get $(D.B)^2 \geq D^2B^2 > (B^2)^2$, so $D.B > B^2$, and in particular $D.B \geq 3$, whence
$h^0(D-B) < h^0(D)-\frac{1}{2}D.B \leq h^0(D)-2$, a contradiction. If $B^2=0$, then $B.D \geq 2$, since $D$ is base point free (by \cite[(3.9.6)]{S-D} or \cite[Thm. 1.1]{kn1}), so again 
$h^0(D-B) \leq h^0(D)-2$. 

So the only possibility remaining is $B^2 \leq -2$, and we see that $h^0(D-B)=h^0(D)-1$ if and only if 
$B^2=-2$ and $D.B=0$.  Since $h^0(D) \geq 3$, we have that $L \sim (D-B) + (F+B)$ is a decomposition into two moving classes with $(D-B).(F+B)=D.F+2-B.L=c+2-(B.L-2)$, so we must have $B.L \leq 2$.

This means that there is a codimension one subset of $|D|$ whose elements are of the form $D'+B$,
with $D'$ base point free with $p_a(D')=p_a(D)-1$, $h^1(D')=0$ and $B$ only supported on smooth rational curves
and satisfying $B^2=-2$, $B.D'=2$ and $B.L \leq 2$. Clearly, since the general element in $|D'|$ is a smooth irreducible curve, we can choose a pencil in $|D|$ such that elements of this form
are of the form $D'+B$ with $D'$ a smooth curve of genus $p_a(D)-1$. Now the contracted part $B_0$ of $B$ satisfies $B_0.D' \leq B.D'=2$, whence $D'$ is mapped by $\varphi_L$ to a curve with at worst one point of multiplicity two, i.e. either a node or a cusp. If $\varphi_L(D')$ is smooth then it has genus $p_a(D)-1$, if not it has arithmetic genus $p_a(D)$. The divisor $B$ is either zero or is mapped to a point or to an effective divisor $\Omega$ on $S'$ of degree $B.L \leq 2$, whence a line, a conic, a union of two distict lines, or a double line.

In the second case we have $h^0(D-B)=k+1=h^0(D)-1=\frac{1}{2}D^2+1$, whence $D^2=2k \geq 4$,
so $2k=D^2=(kE+B)^2=2kE.B+B^2 =2kE.D+B^2$. At the same time, by the base point freeness of $D$, we have 
$B.D=B.(kE+B)= kE.B+B^2=kE.D+B^2 \geq 0$ and $E.D \geq 2$, so the only possibility is
$E.B=E.D=2$, $B.D=0$ and $B^2 = -2k \leq -4$. In particular $D$ is hyperelliptic, so
$c=2$, $D^2=4$ and $L \sim 2D$ by Proposition \ref{cliffgenD}, which means that $k=2$.
Since $D.L=8$, and $D \sim 2E+B$, we must have $E.L=4$ and $B.L=0$. Now there is a codimension
one subset of $|D|$ whose elements are of the form $E_1+E_2+B$ where $E_1$ and $E_2$ are smooth elliptic curves in $|E|$. Since $B.E_1=B.E_2=1$ and $B$ is contracted by $\varphi_L$, the elements are mapped to a union of two smooth elliptic curves intersecting in one point.
\hspace{0.07cm} $\square$ \end{proof}

\begin{rem} \label{changeofcd}
  {\rm We see from the proof above that in the cases where there exists a reducible fibre $\varphi_L(D_{\lambda})$, then we are either in the case with $D \sim A+B$ into two moving classes or
$D \sim D'+B$ with $D'$ either irreducible or twice an elliptic pencil and $B$ supported on rational curves with $B^2=-2$ and $B.L \leq 2$. In the first case we find that $A$ and $B$ are Clifford divisors for $L$ and in the second that $D'$ is a free Clifford divisor.  In particular we see that we can always find a free Clifford divisor $D$ satisfying either $D^2=0$ or that $|D|$ contains a subpencil $ \{ D_{\lambda} \}$
such that all the members of $ \{ \varphi_L (D_{\lambda}) \}$ are 
irreducible. 

Such a $D$ need however not be perfect.}
\end{rem}

Note that the property (\ref{eq:redpencil}) also yields that the singular locus of 
any $\varphi_L (D _{\lambda})$ is either a finite number of points or at most a finite number of points and a double line, so it has an open set of regular points. The same applies for any $D_{\lambda}$.

Moreover, again by the property (\ref{eq:redpencil}), a general hyperplane section of any 
$\varphi_L (D _{\lambda})$ is of a scheme of length $L.D=D^2+c+2$ which either consists of 
distinct points (outside of $\Sing \varphi_L (D _{\lambda})$) or of a union of 
$L.D-2=D^2+c$ distinct points (outside of $\Sing \varphi_L (D _{\lambda})$) and a scheme of length two situated in one point,
namely the intersection with the double line, or equivalently, the image by $\varphi_L$ of the unique element 
in $|\O_C(2\Gamma)|$ for a general $C \in |L|$.

We will from now on always work with a pencil satisfying (\ref{eq:redpencil}).

We will need the following \index{GPS}{\it general position statement}:

\begin{Lemma} \label{gps}
  Assume $c >0$ and $D^2 >0$ and let $D' \sub \PP ^{c+1+\frac{1}{2}D^2}$ be any member of $\{ \varphi_L (D _{\lambda}) \}$. Then a general hyperplane section $Z$ is a scheme of length $L.D=D^2+c+2$ in general position, i.e. any subscheme of length $c+\frac{1}{2}D^2$ spans $\PP ^{c+1+\frac{1}{2}D^2}$.
\end{Lemma}

\begin{proof}
  A general hyperplane section $Z$ consists either of $D^2+c+2$ distinct points, or of 
$D^2+c+1$ distinct points, where one carries an additional tangent direction.

Set $r:=c+1+\frac{1}{2}D^2$, then $r \geq 3$. The proof now follows the lines of the proof of the 
well-known General Position Theorem on p.~109 in \cite{acgh}. We leave it to the reader to verify that the steps (i)-(iii) in that proof go through and that we can reduce to showing (correspondingly to the lemma on 
p.~109 in  \cite{acgh}) that a general hyperplane section of $D'$ contains no subscheme of length $3$ spanning only a $\PP^1$.

So assume there is a general hyperplane section $Z$ of $D'$ containing a subscheme $Z_0$ of length three spanning a $\PP^1$. Since we assume $Z$ is general, we can avoid it to touch the singular points of $S'$.
So we can consider $Z$ and $Z_0$ as subschemes of $S$ and we get that the natural map 
$H^0(L) \hpil H^0(L \* \O_{Z_0})$ fails to be surjective. Since $c \geq 1$, we must have $L^2 \geq 4c+4 \geq 8$, and we can use Proposition \ref{kn1prop} to conclude that there is an effective divisor $B$ passing though $Z_0$ satisfying either $B^2=-2$ and $B.L \leq 1$, $B^2=0$ and $B.L \leq 3$ or $B^2=2$ and $B.L \leq 5$.

In the first case we have $B.L=1$ and $B$ irreducible, since we assume that $Z$ lies outside the singular locus of $S'$.
So $Z_0$ lies on a line, and a general hyperplane will only meet this line in one point, a contradiction.

In the two other cases we see that $B$ induces the Clifford index one on $L$ and we must have 
$(B^2,B.L)=(0,3)$ or $(2,5)$. Since we assume $D^2 \geq 2$, we must have $D^2=c+1=2$, which means that we are in the case (E0), where $L \sim 2D+\Gamma$ for a smooth rational curve $\Gamma$ satisfying $\Gamma.D=1$.
Since $D$ is base point free, any other free Clifford divisor $D'$ must satisfy $D'.D \geq 2$, whence $D'.L
\geq 4$. Now the moving part of $B$ is a free Clifford divisor, whence we must have $(B^2,B.L)=(2,5)$.
It follows from Proposition \ref{kn1prop} that $Z_0 \cap \Gamma \neq \emptyset$, and since $\Gamma.L=0$,
it follows that $Z_0$ meets the singularities of $S'$, a contradiction.
\hspace{0.07cm} $\square$ \end{proof}

We will make use of the following lemma, which is well-known if $D_0$ is a smooth curve (see e.g. 
\cite[Lemma 3.1]{gl2} or \cite[Exc. K-2 p.~152]{acgh} for (a)):

\begin{Lemma} \label{indepcond}
  Let $D_0 \in |D|$. 
  \begin{itemize}
  \item[(a)] If $x_1, \ldots, x_n$ are $n:=h^0(L_{D_0})-2= \frac{1}{2}D^2+c$ distinct general points of $D_0$, outside of $\Sing D_0$,  then $L_{D_0}-x_1-\ldots - x_n$ is base point free, $h^0(L_{D_0}-x_1-\ldots - x_n)=2$ 
and $h^1(L_{D_0}-x_1-\cdots - x_n)=0$.  
  \item[(b)] If $x_1, \ldots, x_k$ are  $k \geq p_a(D)$ distinct general points of $D_0$,
outside of $\Sing D_0$, then $h^1(\O_{D_0}(x_1+ \cdots + x_k))=0$.
  \end{itemize}
\end{Lemma}

\begin{proof}
Since $n = L.D - \frac{1}{2}D^2-2 \leq L.D-3$ the statement (a) immediately follows from
the previous lemma.

As for (b), by Serre duality we have 
$h^1(\O_{D_0}(x_1+ \cdots + x_k))=h^0(\O_{D_0}(D) \linebreak (-x_1- \cdots - x_k))$. Denoting the ideal defined by the points $x_1, \ldots ,x_k$ by $Z$, we have an exact sequence
\begin{equation} \label{eq:negdeg}
  0 \hpil \O_S  \hpil \O_S(D) \* \I_Z  \hpil \O_{D_0}(D)(-Z)  \hpil   0,
\end{equation}
so $h^1(\O_{D_0}(x_1+ \cdots + x_k))=0$ if and only if $h^0(\O_S(D) \* \I_Z)=1$.
Clearly we can assume that $k=p_a(D)$. Then $h^0(\O_S(D) \* \I_Z)=1$ if and only if the $k$ points pose independent conditions on $|D|$. Proceeding inductively, we only have to show that for
$k'$ distinct points on $D_0$, with $1 \leq k' \leq k$, posing independent conditions on $|D|$, then a general point $p \in D_0$
away from $\Sing D_0$ poses one more additional condition. Let $\Sigma$ be the base divisor of 
$|D \* \I_{Z'}|$, where $Z'$ is the scheme defined by the $k'$ distinct points. Then we are done, unless
all the regular points of $D_0$ are contained in $\Sigma$. However, by the property (\ref{eq:redpencil}),
it would then follow that $h^0(\Sigma) \geq h^0(D)-1$. But then the moving part of $|D \* \I_{Z'}|$ has dimension zero, i.e. it consists only of $D_0$ itself, so $h^0(D \* \I_{Z'})=1$, and it follows that $k=k'$ and
we are done.
\hspace{0.07cm} $\square$ \end{proof}

We write $L_{\lambda}:= L_{D _{\lambda}}$.

We first define a vector bundle $\E_{\lambda}$ on every $D _{\lambda}$, as follows. If $B$ is an effective divisor on $S$ and $\A$ is any globally generated invertible sheaf on $B$, then the evaluation map $H^0(\A) \* \O_B \khpil \A$ is surjective, and the kernel is a vector bundle on $B$:
\begin{equation}
  \label{eq:vb}
  0 \hpil  \E_{\A} \hpil H^0(\A) \* _{\CC} \O_{B}   \hpil  \A \hpil  0. 
\end{equation}
Note that $\det \E_{\A} = \A \v$ and $\rank \E_{\A}=h^0(\A)$, so that $\E_{\A} = \A \v$ when $h^0(\A)=2$.

For every $\lambda$ we set $\E_{\lambda}:= \E_{L_{\lambda}}$.

Taking exterior powers in (\ref{eq:vb}) and twisting by suitable powers of $L$, we get for any $i \geq 0$ and any $j \geq 0$
\begin{equation}
  \label{eq:vbpk}
  0 \hpil  \wedge ^i \E_{\lambda} \*  L_{\lambda}^{\*j} \hpil 
\wedge ^i H^0(L_{\lambda}) \* _{\CC} L_{\lambda}^{\*j}   \hpil   
\wedge ^{i-1}  \E_{\lambda} \*  L_{\lambda}^{\*(j+1)}  \hpil  0. 
\end{equation}

Moreover, we get

\[ 
\xymatrix{
    & 0  \ar[d]  & & &
\\
& \wedge ^{i+1} \E_{\lambda} \*  L_{\lambda}^{\*(j-1)}  \ar[d]  & & &
\\
& \wedge ^{i+1} H^0(L_{\lambda}) \* _{\CC} L_{\lambda}^{\*(j-1)}  \ar[d] \ar[dr]^{f^{\lambda}_{i+1,j-1}} &  & 0 \ar[d] &
\\
0 \ar[r] & \wedge ^i \E_{\lambda} \*  L_{\lambda}^{\*j} \ar[r] \ar[d] & 
               \wedge ^i H^0(L_{\lambda}) \* _{\CC} L_{\lambda}^{\*j} \ar[r] \ar[dr] ^{f^{\lambda}_{i,j}} &  
              \wedge ^{i-1}  \E_{\lambda} \*  L_{\lambda}^{\*(j+1)} \ar[r]   \ar[d] &  0
\\
& 0  & &  \wedge ^{i-1} H^0(L_{\lambda}) \* _{\CC} L_{\lambda}^{\*(j+1)},     &
} \]
and we see that ${d^{\lambda} _{i,j}}= H^0(f^{\lambda}_{i,j})$ for all $i,j \geq 0$.

Chasing the diagram, and using that $h^1(L_{\lambda}^{\*k})=0$ for all $k \geq 1$, together with 
$h^0(L_{\lambda})=c+2+\frac{1}{2}D^2$, we easily get that the Koszul cohomology groups $\K^{\lambda} _{i,j}$ satisfy
 \begin{eqnarray}
   \label{eq:betti1}  \K^{\lambda} _{i,j} & = & 0 \mbox{  for all } j \geq 3. \\
   \label{eq:betti1'} \dim \K^{\lambda} _{i,2}& = & h^1(\wedge ^{i+1} \E_{\lambda} \*  L_{\lambda}).\\
   \label{eq:betti1''} \dim \K^{\lambda} _{i,1} & = &  h^1(\wedge ^{i+1} \E_{\lambda})- 
{c+2+\frac{1}{2}D^2 \choose i+1}(\frac{1}{2}D^2+1) \\
 \nonumber & & + h^1(\wedge ^i \E_{\lambda} \*  L_{\lambda})
 \end{eqnarray}

Of course we also have
\begin{equation}
  \label{eq:betti1'''}
  \K^{\lambda} _{i,j} =  0 \mbox{ for }  \hs i \geq h^0 (L _{D _{\lambda}}) -1 = 
                  c+1+ \frac{1}{2}D^2.
\end{equation}

We now want to show that $h^1(\wedge ^{i+1} \E_{\lambda})$ is independent of $\lambda$ for 
$i \leq h^0 (L _{D _{\lambda}}) -2 = c+ \frac{1}{2}D^2$.

By Lemma \ref{indepcond}(a), if $x_1, \ldots, x_n$ are $n:=h^0(L_{\lambda})-2= c+\frac{1}{2}D^2$ general distinct points of $D_{\lambda}$, outside of the singular points of $D_{\lambda}$, then $L_{\lambda}-x_1-\ldots -x_n$ is generated by its global sections and $h^1(L_{\lambda}-x_1-\ldots -x_n)= h^1(L_{\lambda})=0$, so from
\cite{glroma} or \cite[Lemma 1.4.1]{lazsamp} we have an exact sequence
\begin{equation}
  \label{eq:betti3}
  0 \hpil  \E_{L_{\lambda}-x_1-\ldots -x_n} \hpil  \E_{L_{\lambda}}\hpil \Sigma \hpil  0,
\end{equation}
where $\Sigma := \+_{i=1}^n \O_{D_{\lambda}}(-x_i)$. (We leave it to the reader to check that this also holds in our case when
$D_{\lambda}$ is singular or possibly reducible). Set $\B:= \O_{D_{\lambda}}(x_1 + \cdots + x_n)$. Since $h^0(L_{\lambda}-\B)=2$, we have $\E_{L_{\lambda}-x_1-\ldots -x_n} = \B-L_{\lambda}$.
Taking exterior products yields
\begin{equation}
  \label{eq:betti5}
  0 \hpil \wedge ^i \Sigma \*(\B-L_{\lambda}) \hpil  \wedge ^{i+1} \E_{L_{\lambda}} \hpil \wedge ^{i+1} \Sigma \hpil  0.
\end{equation}

The term on the right is a direct sum of ${n \choose i+1}$ line bundles of the form
$\O_{D_{\lambda}}(-x_{k_1}- \cdots -x_{k_{i+1}})$, whence for all $i \geq 0$ we have 
$h^0(\wedge ^{i+1} \Sigma)=0$ and by Riemann-Roch $h^1(\wedge ^{i+1} \Sigma)={n \choose i+1}(\frac{1}{2}D^2+1+i)$. 

The term on the left is a direct sum of ${n \choose i}$ line bundles of the form
$\O_{D_{\lambda}}(x_{k_1}+ \cdots +x_{k_{n-i}})\* L \v$. Now by Serre duality 
$h^0( \O_{D_{\lambda}}(x_{k_1}+ \cdots +x_{k_{n-i}})\* L \v= h^1(\O_{D_{\lambda}}(L+D)(-x_{k_1}- \cdots -x_{k_{n-i}})$. By the sequence
\begin{equation}   
0 \hpil L  \hpil \O_S(L+D) \* \I_Z  \hpil \O_{D_{\lambda}}(L+D)(-Z)  \hpil   0,
\end{equation}
with $Z$ the ideal defined by $x_{k_1}, \ldots ,x_{k_{n-i}}$, 
$h^1(\O_{D_{\lambda}}(L+D)(-x_{k_1}- \cdots -x_{k_{n-i}}))=h^1(\O_S(L+D) \* \I_Z)$. Since $h^1(L+D)=0$, this is equivalent to saying that $x_{k_1}, \ldots ,x_{k_{n-i}}$ pose independent conditions on $L+D$. 
but since $n-i \leq n$, the points pose independent conditions on $L$ by Lemma \ref{indepcond}(a), whence also on $L+D$, since $D$ is base point free. So
$h^0(\wedge ^{i+1} \Sigma)=0$ and by Riemann-Roch $h^1(\wedge ^{i+1} \Sigma)={n \choose i}(D^2+2+i)$.
Inserting for $n$, it follows that
\begin{equation}
  \label{eq:betti6}  h^1(\wedge ^{i+1} \E_{\lambda})= {c+\frac{1}{2}D^2 \choose i}(D^2+2+i)+{c+\frac{1}{2}D^2 \choose i+1}(\frac{1}{2}D^2+1+i).
\end{equation}
This improves (\ref{eq:betti1''}):
\begin{eqnarray}
  \label{eq:betti7}
 \dim \K^{\lambda} _{i,1}  & = & {c+\frac{1}{2}D^2 \choose
   i}(D^2+2+i)+{c+\frac{1}{2}D^2 \choose i+1}(\frac{1}{2}D^2+1+i) \\
 & & - {c+2+\frac{1}{2}D^2 \choose i+1}(\frac{1}{2}D^2+1)+
\nonumber  h^1(\wedge ^i \E_{\lambda} \*  L_{\lambda}). 
\end{eqnarray}
In particular, we see that
\begin{eqnarray}
\label{eq:betti8} 
\dim \K^{\lambda} _{i,1}  &  -  &  \dim \K^{\lambda} _{i-1,2} =  {c+\frac{1}{2}D^2 \choose i}(D^2+2+i) \\
\nonumber  & +  & {c+\frac{1}{2}D^2 \choose i+1}(\frac{1}{2}D^2+1+i) - {c+2+\frac{1}{2}D^2 \choose i+1}(\frac{1}{2}D^2+1) 
\end{eqnarray}
 is independent of $\lambda$.

Recall that the line bundle $L _{\lambda}$ on 
$D _{\lambda}$ is said to satisfy property $N _p$\index{property $N _p$} if the Betti-numbers satisfy the following:
\begin{equation} \label{eq:Np}
   \beta^{\lambda}_{0,j} = \left\{ \begin{array}{ll}
             1     & \mbox{ if $j=0$}, \\
             0     & \mbox{ if $j \not =0$ } 
          \end{array}
    \right .     
\hs \mbox{ and } \hs
\beta^{\lambda}_{i,j} \not = 0 \mbox{ if and only if } j=i+1, \mbox{ for } 0 < i \leq p.
\end{equation}

This means that $B ^{\lambda}$ has a resolution of the form
\begin{eqnarray*}
  \cdots  & \hpil & R_{\lambda}(-p-1) ^{\beta _{p,p+1}} \hpil  \cdots 
  \hpil  R_{\lambda}(-3) ^{\beta _{2,3}} \\
   & \hpil & R_{\lambda}(-2)
  ^{\beta _{1,2}} \hpil R_{\lambda} \hpil  B ^{\lambda} \hpil  0.
\end{eqnarray*}

In our case, we have

\begin{prop} \label{npprop}
Assume $c>0$. Then $L _{\lambda}$ satisfies property $N_{c-1}$ but not 
$N_c$.
\end{prop}

\begin{proof}
If ${D _{\lambda}}$ is smooth, then the second statement is immediate, since we have by \cite{kn1} and the conditions $(*)$ that $L _{\lambda}$ fails to be $(c+1)$-very ample, and the result follows from \cite[Thm. 2]{gl2}. By semicontinuity, $N_c$ fails for all $L _{\lambda}$.

The first statement is also immediate 
if ${D _{\lambda}}$ is smooth: Indeed,
 it follows from \cite[Thm. (4.a.1)]{gre}, since
 $\deg L _{D _{\lambda}} = 2g(D _{\lambda}) + c$ and $h^1 (L _{D _{\lambda}}) =0$.

We have to argue that the result still holds for the singular and reducible ${D _{\lambda}}$, in other words we have to show that
$\K^{\lambda} _{i,2}  =0$ for all $i \leq c-1$ and all $\lambda$.

By (\ref{eq:betti5}) we have to show that $h^1(\wedge ^{i+1} \E_{\lambda} \*  L_{\lambda})=0$ for all $i \leq c-1$.

Choose as above $n:=c+\frac{1}{2}D^2$ general points $x_1, \ldots, x_n$ of $D_{\lambda}$. We then get a sequence as (\ref{eq:betti5}), and tensoring this sequence with $L_{\lambda}$ yields
\begin{equation}
  \label{eq:betti9}
  0 \hpil \wedge ^i \Sigma \* \B \hpil  \wedge ^{i+1} \E_{L_{\lambda}} \* L_{\lambda} 
\hpil \wedge ^{i+1} \Sigma \* L_{\lambda} \hpil  0.
\end{equation}

The term on the right is a direct sum of ${n \choose i+1}$ line bundles of the form
$L_{\lambda}(-x_{k_1}- \cdots -x_{k_{i+1}})$. By Lemma \ref{indepcond}(a) it follows that
for all $i \geq 0$ we have 
$h^1(\wedge ^{i+1} \Sigma \* L_{\lambda})=0$.

The term on the left is a direct sum of ${n \choose i}$ line bundles of the form
$\O_{D_{\lambda}}(x_{k_1}+ \cdots +x_{k_{n-i}})$, whence of degrees 
$n-i \geq \frac{1}{2}D^2+1=p_a(D)$. By Lemma \ref{indepcond}(b) it follows that
$h^1(\wedge ^i \Sigma \* \B)=0$.

It follows that $h^1(\wedge ^{i+1} \E_{\lambda} \*  L_{\lambda})=0$ for all $i \leq c-1$.

An alternative proof of the fact that
$\K^{\lambda} _{i,2}  =0$ for all $i \leq c-1$ and all $\lambda$ goes as follows: Since $\varphi_L(D_{\lambda})$ is arithmetically normal and $D_{\lambda}$ is of pure dimension one, the Betti-numbers of $\varphi_L(D_{\lambda})$ are equal to the Betti-numbers of a general hyperplane section of it. This is a scheme $X$ of length $L.D=D^2+c+2$ in general position by Lemma \ref{gps}. We now argue as in the proof of \cite[Thm. 2.1]{gl2} to show that the scheme $X$ satisfies  $N_{c-1}$. Recall that $X$ either consists of 
distinct points  or at worst of a union of 
$L.D-2=D^2+c$ distinct points and a scheme of length two supported in one point, call it $Z$. The case of distinct points is exactly the statement in \cite[Thm. 2.1]{gl2}, so we have to show that the proof goes through in the other case. 
We leave it to the reader to verify that everything works as long as one writes the scheme $X$ as a disjoint union $X = X_1 \cup X_2$ as in the proof of \cite[Thm. 2.1]{gl2}, taking care that $Z \sub X_2$. This is possible, since $X_1$ should consist of $\frac{1}{2}D^2+c+1$ distinct points, which yields $\length X_2 =
\frac{1}{2}D^2+1 \geq 2$.
\hspace{0.07cm} $\square$ \end{proof}

From this proposition we therefore get
\begin{eqnarray}
   \label{eq:betti1a}  \K^{\lambda} _{i,2} & = & 0 \mbox{  for all } i \leq c-1. \\
   \label{eq:betti1b}  \K^{\lambda} _{c,2}& \neq & 0.
 \end{eqnarray}

Also, by the Theorem in \cite{glapp}, we have that for $D^2 > 0$:
\[  \beta^{\lambda}_{c,c+1} \not = 0 \mbox{  for all smooth irreducible } D_{\lambda}. \]
Indeed, $L _{D _{\lambda}} \iso F _{D _{\lambda}} + \omega _{D _{\lambda}} $, and $D^2 \leq 2c$ for $c >0$, and we calculate
\[ h^0 (F _{D _{\lambda}}) = h^0 (F) - \chi (F-D) = c+2 - \frac{1}{2}D^2 \geq 2, \]
\[ h^0 (\omega _{D _{\lambda}}) =  \frac{1}{2}D^2 +1 \geq 1 \]
and
\[ h^0 (F _{D _{\lambda}}) + h^0 (\omega _{D _{\lambda}}) -3 = c.\]

By semicontinuity it follows that
\begin{equation}
  \label{eq:betti1c}
 \K^{\lambda} _{c,1 } \neq  0 \mbox{  for all } \lambda. 
\end{equation}

Finally, recall that the line bundle $L _{\lambda}$ on 
$D _{\lambda}$ is said to satisfy property $M _q$\index{property $M _q$} if $\K^{\lambda} _{i,j }=0$ for all 
$i \geq h^0(L_{\lambda})-1-q=\frac{1}{2}D^2+c+1-q$ and $j \neq 2$.

We have

\begin{prop} \label{tjo}
  (a) If $c >0$, then $L _{\lambda}$ satisfies $M_1$.

  (b) If $D^2 \geq 4$ and $c \geq 3$, then $L _{\lambda}$ satisfies $M_2$. 
\end{prop}

\begin{proof}
The main ingredient in this proof is the proof of Green's $\K _{p,1}$ theorem \cite[(3.c.1)]{gre}.

Set $r:=h^0(L _{\lambda})=\frac{1}{2}D^2+c+1$.

To show (a), we argue as in the proof of statement (2) in \cite[(3.c.1)]{gre}, and assume that $\K^{\lambda} _{p,1 } \neq 0$, for $p=\frac{1}{2}D^2+c=r-1$. Taking a general hyperplane section $Z$ of $\varphi_L(D_{\lambda})$ we get that $\K^{\lambda} _{p,1 } \neq 0$ for $Z \sub \PP^{r-1}$. By Lemma \ref{gps} $Z$ is in general position, so if it consists of distinct points, then it follows from Green's {\it Strong Castelnuovo Lemma} \cite[(3.c.6)]{gre} that $Z$ lies on a rational normal curve,
whence the contradiction $D^2+c+2=L.D=\deg \varphi_L(D_{\lambda}) \leq r=\frac{1}{2}D^2+c+1$.

If $Z$ consists of $L.D-2$ distinct points and a scheme of length two with support in one point we have to show that Green's Strong Castelnuovo Lemma still can be used. The key point is where Green uses that any $r+2$ distinct points in general position in $\PP^{r-1}$ lie on a unique rational normal curve. This still holds true if we have $r+1$ distinct points with one additional tangent direction at one of them, when the whole scheme is in general position.

We leave it to the reader to verify that the Strong Castelnuovo Lemma holds in our case and that we can conclude as above that $Z$ lies on a rational normal curve, and get the same contradiction.

Now we prove (b). Once we have checked that the Strong Castelnuovo Lemma holds in our case, we can argue as in the proof of
(3) in \cite[(3.c.1)]{gre}, and find that either
$D^2+c+2=L.D=\deg \varphi_L(D_{\lambda}) \leq r+1 =\frac{1}{2}D^2+c+2$, which is not our case, or that 
$\varphi_L(D_{\lambda})$ lies on a surface of minimal degree, i.e. the Veronese surface in $\PP^5$, a ruled surface or a cone over a rational normal curve. 

In the first case we must have $r=5$, whence $c=2$ and $D^2=4$. 

In the two other cases, then if $\varphi_L(D_{\lambda})$ does not pass through the vertex of the cone
the ruling restricts to a Cartier divisor on $\varphi_L(D_{\lambda})$ and it cuts out a $g^1_2$ on $\varphi_L(D_{\lambda})$ which we can pull back by $\varphi_L$ to $S$. Then every element $Z$ in this $g^1_2$ on $S$ is a $0$-dimensional scheme of length $2$ failing to pose independent conditions on $|D|$.
Therefore $|D|$ must be hyperelliptic and by 
Proposition \ref{cliffgenD} we have $c=2$ and $D^2=4$.

We have left to treat the case where $\varphi_L(D_{\lambda})$ lies in a cone and passes through its vertex. Since $\varphi_L(D_{\lambda})$ cannot be a union of lines by (\ref{redpencil}) the ruling cuts out a $g^1_1$ on the component of $\varphi_L(D_{\lambda})$ obtained by removing the components which are lines of the ruling, if any. So this component  is an irreducible curve birational to $\PP^1$.
By (\ref{redpencil}) this curve is either smooth of genus $p_a(D)$ or $p_a(D)-1$ or has only one node or cusp and arithmetic genus $p_a(D)$ or $p_a(D)+1$. In all these cases we get that the curve has geometric genus $\geq p_a(D)-1 \geq 2$, since we assume $D^2 \geq 4$, a contradiction.
\hspace{0.07cm} $\square$ \end{proof}

The following lemma settles the remaining case $D^2=4$ and $c=2$, where in fact 
$L _{\lambda}$ does not satisfy $M_2$:

\begin{Lemma}
Let $(c,D^2)=(2,4)$. Then $\dim \K ^{\lambda}_{3,1}=3$ for all $\lambda$.
\end{Lemma}

\begin{proof}
We are in the case (Q) with $L \sim 2D$. 
By Proposition \ref{2-uple} either $\varphi_L$ is the $2$-uple embedding of $\varphi_D(S)$,
or there is  an elliptic pencil $|E|$ such that $E.D=2$. We will treat these two cases separately.

In the first case $\varphi_L(D _{\lambda})$ is the $2$-uple embedding of $\varphi_D(D _{\lambda})$,
for all $\lambda$. Now $\varphi_D$ maps $D _{\lambda}$ into $\PP^2$, so 
$\varphi_L(D _{\lambda})$ lies on the Veronese surface $V$ in $\PP^5$, i.e. the $2$-uple embedding of $\PP^2$.

We have $\Pic V \sim \ZZ l$, where $l^2=1$. The hyperplane class $H_V$ satisfies $H_V \sim 2l$, and since $\varphi_L(D _{\lambda})$ has degree $L.D=8$, we have $\varphi_L(D _{\lambda}) \sim 4l \sim 2H_V$. By \cite[(3.b.4)]{gre} we have
\[ \K ^{\lambda}_{3,1}= \K ^{\lambda}_{3,1}(V,H_V) \+ \K ^{\lambda}_{2,0}(V,H_V). \]
Both the latter are well-known, since  $V$ is a variety of minimal degree (see e.g. \cite[Lemma 5.2]{Sc}). In fact $\dim \K ^{\lambda}_{3,1}(V,H_V)=3$ and $\K ^{\lambda}_{2,0}(V,H_V)=0$.
Hence $\dim \K ^{\lambda}_{3,1}=3$, as asserted.

In the second case, any $D_{\lambda}$ has a $g^1_2$ given by $\O_{D_{\lambda}}(E)$.
Compare the two morphisms $f_E: D_{\lambda} \khpil \PP^1$ given by $|\O_{D_{\lambda}}(E)|$ and
$f_D: D_{\lambda} \khpil \PP^2$ given by $|\O_{D_{\lambda}}(D)|$. Since $h^1(E-D)=h^1(\O_S)=0$, these are the restrictions of $\varphi_E$ and $\varphi_D$ respectively. Since they both collapse every member of the $g^1_2$, we see that 
$f_D= g \circ f_E$, where $g: \PP^1 \khpil \PP^2$ is the $2$-uple embedding. It follows that
 $\O_{D_{\lambda}}(D) = \O_{D_{\lambda}}(D) ^{\* 2}$ and consequently $L_{D_{\lambda}} =4E_{D_{\lambda}}$.

The members of the $g^1_2$ sweep out a scrollar surface $S_0$ containing $\varphi_L(D _{\lambda})$. As before, we can compute its scroll type $(e_1, e_2)$ by first computing the ``dual scrollar invariants'' 
\begin{eqnarray*}
 d_i & = & h^0(L_{D_{\lambda}}-iE_{D_{\lambda}}) - h^0(L_{D_{\lambda}}-(i+1)E_{D_{\lambda}}) \\
     & = & h^0((4-i)E_{D_{\lambda}}) - h^0((3-i)E_{D_{\lambda}}). 
\end{eqnarray*}
We easily get $d_0=2$, $d_1=d_2=d_3=d_4=1$ and $d_{\geq 5}=0$. Recalling that
$e_i= \# \{j \hs | \hs d_j \geq i \}-1$ we get $(e_1, e_2)=(4,0)$, whence $S_0$ is a cone over a rational normal quartic. 

Note that since the $g^1_2$ is base point free $\varphi_L(D _{\lambda})$ does not intersect the vertex of $S_0$, so we can work with the desingularization $S_0$, which we by abuse of notation also denote by $S_0$.

We have $\mbox{Num} S_0 \iso \ZZ H_0 \+ \ZZ L_0$, where $H_0$ is the hyperplane class and $L_0$ is the class of the ruling, whence $H_0.L_0=1$, $H_0^2=4$ and $L_0^2=0$. Since $\varphi_L(D _{\lambda}).H_0= \deg \varphi_L(D _{\lambda})=8$ and 
$\varphi_L(D _{\lambda}).L_0=2$, we find $D_0 \sim 2H_0$. Moreover we have
$H^0(S_0, H_0-D_0)=H^1(S_0, qH_0-D_0)=H^1(S_0, qD_0)=0$ for all $q \geq 0$, so by \cite[(3.b.4)]{gre} we have
\[ \K ^{\lambda}_{3,1}= \K ^{\lambda}_{3,1}(S_0,H_0) \+ \K ^{\lambda}_{2,0}(S_0,H_0). \]
Again it is well-known that $\dim \K ^{\lambda}_{3,1}(S_0,H_0)=3$ and $\K ^{\lambda}_{2,0}(S_0,H_0)=0$ (see e.g. \cite[Lemma 5.2]{Sc}), so $\dim \K ^{\lambda}_{3,1}=3$, as asserted.
\hspace{0.07cm} $\square$ \end{proof}

Summing up, we have 
\begin{prop} \label{Bettiprop}
 Let $D$ be a free Clifford divisor on a polarized $K3$ surface
 $(S,L)$ of
 non-general Clifford index $c>0$ satisfying $D^2 >0$. Then the Betti-numbers of the 
$\varphi_L (D _{\lambda})$ satisfy:
\begin{itemize}
  \item[(a)] 
   $\beta^{\lambda}_{0,j} = \left\{ \begin{array}{ll}
             1     & \mbox{ if $j=0$}, \\
             0     & \mbox{ if $j \not =0$ } 
          \end{array}
    \right. $.
  \item[(b)] For $0 < i \leq c-1$, $\beta^{\lambda}_{i,j} \not = 0$ if and only if $j=i+1$,
  \item[(c)] $\beta^{\lambda}_{ij} =0$  for $i \geq c+1+ \frac{1}{2}D^2$.
  \item[(d)] $\beta^{\lambda}_{ij} =0$  for $j \geq i+3$.  
    \item[(e)] $\beta^{\lambda}_{i,i+1}-\beta^{\lambda}_{i-1,i+1} =
     {c+\frac{1}{2}D^2 \choose i}(D^2+2+i)+\\ {c+\frac{1}{2}D^2
     \choose i+1}(\frac{1}{2}D^2+1+i)- {c+2+\frac{1}{2}D^2 \choose
     i+1}(\frac{1}{2}D^2+1),$
     for $i >0$.
\item[(f)] $\beta^{\lambda}_{\frac{1}{2}D^2+c,\frac{1}{2}D^2+c+1} =0$.
  \item[(g)] $\beta^{\lambda}_{\frac{1}{2}D^2+c-1,\frac{1}{2}D^2+c } =0$ 
 for $D^2 \geq 4$ and $c \geq 3$. 
\item[(h)] $\beta^{\lambda}_{3,4} =3$ if $(c,D^2)=(2,4)$
\item[(i)] $\beta^{\lambda}_{c,c+1} \not = 0$ for $D^2 >0$.
  \item[(j)] $\beta^{\lambda}_{c,c+2} \not = 0$.
\end{itemize}
\end{prop}

So for $D^2 >0$, the $\varphi_L (D _{\lambda})$ all have a resolution of the form:
\begin{eqnarray*}
 0 \hpil R_{\lambda}(-\frac{1}{2}D^2-c-2) ^{\beta ^{\lambda}_{\frac{1}{2}D^2+c,\frac{1}{2}D^2+c+2}} \hpil \cdots  & \hpil  &  \\
 R_{\lambda}(-\frac{1}{2}D^2-c) ^{\beta ^{\lambda}_{\frac{1}{2}D^2+c-1,\frac{1}{2}D^2+c}} \+ 
R_{\lambda}(-\frac{1}{2}D^2-c-1) ^{\beta ^{\lambda}_{\frac{1}{2}D^2+c-1,\frac{1}{2}D^2+c+1}}  
& \hpil &    \\
\cdots  \hpil R_{\lambda}(-c-2) ^{\beta ^{\lambda}_{c+1,c+2}} \+  R_{\lambda}(-c-3) ^{\beta ^{\lambda}_{c+1,c+3}} 
&  \hpil & \\
R_{\lambda}(-c-1) ^{\beta ^{\lambda}_{c,c+1}} \+  R_{\lambda}(-c-2) ^{\beta ^{\lambda}_{c,c+2}} 
\hpil R_{\lambda}(-c) ^{\beta ^{\lambda}_{c-1,c}} & \hpil &   \\
 \cdots  \hpil R_{\lambda}(-3) ^{\beta ^{\lambda}_{2,3}} \hpil
      R_{\lambda}(-2) ^{\beta ^{\lambda}_{1,2}} \hpil R_{\lambda} \hpil B^{\lambda} \hpil 0.& & 
\end{eqnarray*}

It is easy to see that all the Betti-numbers for $D^2=2$ and $D^2=4$ are uniquely determined by the information above.
Combining with Example \ref{zeroint}, we get:

\begin{cor} \label{Bettis}
For $0 \leq D^2 \leq 4$ the Betti-numbers of the $\varphi _L(D_{\lambda})$ are the same for 
all $\lambda$ and uniquely given by the results above.
\end{cor}

We will now compute some concrete examples of projective models
of $K3$ surfaces contained in singular scrolls $\T$. We will use the
results above to obtain minimal resolutions of the
$\varphi_ L(D_{\lambda})$ in the projective spaces they span. We
will also give results (Proposition  \ref{upperres} and \ref{godeb})
showing how one can lift these resolutions to resolve $\O_{S''}$ and 
$\O_{S'}$ as $\O_{\T_0}$- and $\O_{\T}$-modules, respectively.

 \begin{exa} \label{simis}
{\rm As our first example we study the case (E0) with  $c=1$, $D^2=2$ 
and $g=6$. By Proposition \ref{Bettiprop} all the $\varphi_ L(D_{\lambda}) \sub \PP^3$
have minimal resolutions
\[  0 \hpil R(-4)^2 \hpil R(-2) \+ R(-3)^2 \hpil R \hpil B \hpil 0. \] 
This is the well-known resolution of a smooth curve of genus $2$ in $\PP^3$ (see e.g.
\cite{Sim}).}
\end{exa}

\begin{exa} \label{Dkvadrat2}
{\rm As another example we study the case when $D^2=2$ and $c \geq 2$, where
$\varphi_L(D_{\lambda}) \sub \PP^{c+2}$.

For $c=2$  a minimal resolution is of the following form:
\begin{eqnarray*}
 0 \hpil R(-5)^2  & \hpil & R(-3)^2  \+  R(-4)^{3} \\
& \hpil & R(-2)^4 \hpil R \hpil B^{\lambda}  \hpil  0,
\end{eqnarray*} 

For $c=3$  a minimal resolution is of the following form:
\begin{eqnarray*}
 0 \hpil R(-6)^2   \hpil  R(-5)^4  \+  R(-4)^{3}  \hpil  R(-3)^{12}\\
  \hpil  R(-2)^8 \hpil  \hpil B^{\lambda}  \hpil  0,
\end{eqnarray*}

For $c=4$  a minimal resolution is of the following form:
\begin{eqnarray*}
 0 \hpil R(-7)^2   \hpil  R(-6)^5  \+  R(-5)^{4}  \hpil  R(-4)^{25} \\
 \hpil  R(-3)^{30}  \hpil  R(-2)^{13} \hpil R \hpil B^{\lambda}  \hpil  0,
\end{eqnarray*}}
\end{exa}

\begin{exa} \label{Dkvadrat4}
{\rm As yet another example we study the case when $D^2=4$ and $c \geq 2$, where
$\varphi_L(D_{\lambda}) \sub \PP^{c+3}$.

For $c=2$  a minimal resolution is of the following form:
\begin{eqnarray*}
 0 \hpil R(-6)^3   \hpil  R(-5)^8 \+ R(-4)^3  \hpil \\
R(-4)^6  \+  R(-3)^{8} & \hpil & R(-2)^7 \hpil R \hpil B^{\lambda}  \hpil  0.
\end{eqnarray*}

For $c=3$  a minimal resolution is of the following form:
\begin{eqnarray*}
 0 \hpil R(-7)^3  \hpil  R(-6)^{10}  \hpil  R(-5)^6  \+  R(-4)^{15}\\
\hpil R(-3)^{25}  \hpil  R(-2)^{12} \hpil R \hpil B^{\lambda}  \hpil  0.
\end{eqnarray*}}
\end{exa}

\begin{rem} \label{projeul}
{\rm 
If we twist the resolution following Proposition \ref{Bettiprop}  with
$n$ and use the additivity of the Euler
characteristic, we obtain the following polynomial identity in the
variable $n$:
\begin{eqnarray*}
  (c+2+D^2)n-\frac{1}{2}D^2  = 
{n+c+1+\frac{1}{2}D^2 \choose c+1+\frac{1}{2}D^2}   +  \\
 \sum_{j=2}^{\frac{1}{2}D^2+c+2}  (-1)^j{n+c+\frac{1}{2}D^2+1-j \choose c+1+
\frac{1}{2}D^2}(\beta _{j-2,j} - \beta _{j-1,j}).
\end{eqnarray*}
It is easy to see that this identity alone determines the $\beta _{j-2,j} - 
\beta _{j-1,j}$ uniquely. Since $\beta _{j-2,j}=0$, for $j \leq c+1$, the 
$\beta _{j-1,j}$, for $j=2, \ldots,c+1$ are determined uniquely.
On the other hand this observation gives nothing that is not already
contained in the statement of Proposition \ref{Bettiprop}.}
\end{rem}

We now return to the general resolution, following Proposition \ref{Bettiprop}.
From  Corollary \ref{Bettis}, Example \ref{zeroint} for the case 
$D^2=0$ and \cite[Thm. (3.2)]{Sc} in general, we obtain the following:

\begin{prop} \label{upperres}
If $D^2=0$ and $c=1$, then the $\O_{\T _0}$-resolution $F_*$\index{$F_*$} 
of $\O_{S''}$ is
\[ 0 \hpil \O_{\T_0}(-3\H_0+(g-1)\F) \hpil \O_{\T_0} \hpil \O_{S''} \hpil 0.\]
           
If $D^2=0$ and $c \geq 2$, the resolution is of the following 
type:\index{$b_{i}^k$}
\begin{eqnarray*}
0 \hpil \O_{\T_0}(-(c+2)\H_0+(g-1)\F)  \hpil \+ _{k=1}^{\beta _{c-1}} 
\O_{\T_0}(-c\H_0+b_{c-1}^k\F)  \hpil  \\
\cdots  
 \hpil \+ _{k=1}^{\beta _1} \O_{\T_0}(-2\H_0 + b_1^k \F)  \hpil 
\O_{\T_0} \hpil \O_{S''} \hpil 0, 
\end{eqnarray*}
where ${\beta}_i = i {c+1 \choose i+1} - {c \choose i-1}$\index{${\beta}_i$}.

If $D^2 =2$ or $4$, or more generally if the Betti-numbers of all the
$\varphi_L(D_\lambda)$ are the same for all $\lambda$, then $\O_{S''}$ has a 
$\O_{\T_0}$-resolution $F_*$  of the following type:
\begin{eqnarray*}
& 0 & \hpil F _{\frac{1}{2} D^2+c}  \cdots \hpil F_{c+1} \hpil F _c \\
&    & \hpil \+ _{k=1}^{\beta _{c-1}} \O_{\T _0}(-c{\H_0} + b_{c-1}^k {\F}) 
\hpil  \cdots  \hpil \+ _{k=1}^{\beta _2} \O_{\T _0}(-3{\H_0} + b_2^k {\F}) \\
&    & \hpil \+ _{k=1}^{\beta _1} \O_{\T _0}(-2{\H_0} + b_1^k {\F})  \hpil 
\O_{\T _0} \hpil \O_{S''} \hpil 0. 
\end{eqnarray*}
Here  $\beta _i=\beta_{i,i+1}$\index{${\beta}_i$}, for $i=1, \ldots ,c$, and $F_c$ is an extension of 
the non-zero term 
\index{$b_{i,j}^k$}\[ \+ _{k=1}^{\beta _{c,c+2}} \O_{\T _0}(-(c+2){\H_0} + b_{c,c+2}^k {\F})\] 
by the non-zero term
\[ \+ _{k=1}^{\beta _{c,c+1}} \O_{\T _0}(-(c+1){\H_0} + b_{c,c+1}^k {\F}).\]

Moreover  $F_i$ is an extension of the term  
\[ \+ _{k=1}^{\beta _{i,i+2}} \O_{\T _0}(-(i+2){\H_0} + b_{i,i+2}^k {\F})\]
by the term
\[ \+ _{k=1}^{\beta _{i,i+1}} \O_{\T _0}(-(i+1){\H_0} + b_{i,i+1}^k {\F})\]
for $i=c+1, \ldots ,\frac{1}{2} D^2+c$.
\end{prop}

\begin{proof}
Since the Betti-numbers are the same for all $\lambda$ if $D^2=0$ by
Example \ref{zeroint}, the case $D^2=0$ is a direct
application of \cite[Thm. (3.2)]{Sc}. The case $D^2 \geq 2$ follows
from \cite[Thm. (3.2)]{Sc} and Corollary \ref{Bettis}.
\hspace{0.07cm} $\square$ \end{proof}

We recall the definition \label{seinbetti} $\beta _i=\beta _{i,i+1}$\index{$\beta _i$},
for $i=1,   \ldots   ,c-1$, and $d=c+2+\frac{1}{2}D^2=\dim\T$.

\begin{prop} \label{Euler}
The $b_{i,j}^k$ and $b_i^k$ in Proposition \ref{upperres}  
satisfy the following polynomial equation in $n$ (set $b_{i,i+2}^k=
\beta_{i,i+2}^k=0$ for all $i$ and $k$ if $D^2=0$, and set $b_{i,i+1}^k=
b_i^k$, for $i=1,   \ldots   ,c-1$ for all values of $D^2$):
\begin{eqnarray*}
{n+d-1 \choose d-1}(\frac{n(g+1)}{d}+1) - n^2(g+1+c+D^2)+\frac{1}{2}nD^2
-2  & = & \\
\sum _{i=1}^{c+\frac{1}{2}D^2}(-1)^{i+1}{n+d-i-2 \choose d-1}
(\frac{((n-i-1)(g+1)+d)\beta _{i,i+1}}{d}+\sum _{k=1}^{\beta _{i,i+1}} 
b_{i,i+1}^k) & + & \\
\sum _{i=c}^{c+\frac{1}{2}D^2}(-1)^{i+1}{n+d-i-3 \choose d-1}
(\frac{((n-i-2)(g+1)+d)\beta _{i,i+2}}{d} + 
\sum _{k=1}^{\beta _{i,i+2}} b_{i,i+2}^k).
\end{eqnarray*}
\end{prop}

\begin{proof}
Denote the term $i$ places to the left of $\O _{\T_0}$ in the resolution $F_*$
by $F_i$. The result follows, similarly as in the proof of \cite[Prop. 4.4(c)]
{Sc} from the identity 
\[ \chi (\O_{\T_0}(n\H_0)) - \chi (\O_{S''}(n\H_0)) = \sum_i(-1)^{i+1} 
\chi (F_i(n\H_0)). \]
To calculate  $\chi (\O_{S''}(n\H_0))$ one uses Riemann-Roch on $S''$ 
 and $\deg S'' = 2g+2c+2+2D^2$.

Moreover it is clear that for all large $n$, we have 
$\chi (F_i(n\H))=h^0(F_i(n\H))$, for all $i$, and  
$\chi (\O_{\T_0}(n\H_0))=h^0(\O_{\T_0}(n\H_0))$
since $\H_0$ is (very) ample on $\T_0$.
Then one uses (\ref{dyttned}) again.
\hspace{0.07cm} $\square$ \end{proof}

\begin{rem} \label{ikkehverb}
{\rm From the last result it is clear that the sums 
$\sum _{k=1}^{\beta _{i,j}}b_{i,j}^k$ are uniquely determined, but this does 
not necessarily apply to the 
$b_{i,j}^k$ individually. If $D^2 >0$, it is not even a priori clear that the 
$b_{i,j}^k$ are independent of the choice of pencil inside $|D|$ giving rise to $\T(c,D,\{ D _{\lambda} \})$. }
\end{rem}
\begin{cor} \label{maxbett}
  \begin{itemize}
\item [(a)] We have
\[\sum_{k=1}^{\beta_{1,2}}b_{1,2}^k=(\frac{1}{2}D^2+c-1)g+(1-c-D^2).\]  
\item [(b)] If $D^2 > 0$, 
 then 
\[ \sum _{k=1}^{\beta _{c+ \frac{1}{2}D^2,c+\frac{1}{2}D^2+2}}b_{c+ \frac{1}{2}
D^2,c+\frac{1}{2}D^2+2}^k= g(\frac{1}{2}D^2+1) +1. \]
  \end{itemize}
\end{cor}
\begin{proof}
We insert $n=2$ in Proposition \ref{Euler}. That gives part (a) directly.
Then we insert $n=0$ in Proposition \ref{Euler}. That gives
\[-1=g\beta _{c+ \frac{1}{2}D^2,c+\frac{1}{2}D^2+2} - \sum _{k=1}^
{\beta _{c+ \frac{1}{2}D^2,c+\frac{1}{2}D^2+2}}b_{c+ \frac{1}{2}D^2,c+\frac{1}{2}D^2+2}^k.\] 
This immediately gives the statement of part (b), since it follows from
Proposition \ref{Bettiprop} that 
\[\beta _{c+ \frac{1}{2}D^2,c+\frac{1}{2}D^2+2}=\frac{1}{2}D^2+1.\]  
\hspace{0.07cm} $\square$ \end{proof}

\begin{exa} \label{all2}
{\rm We return to the situation studied in Example \ref{Dkvadrat2}, with
$D^2=c=2$.

From that example and Corollary \ref{maxbett} we see 
that $\beta _{3,5}=2$ and 
\[ \sum_{k=1}^{\beta _{3,5}}b_{3,5}^k=2g-1.\]
We now  apply
Proposition \ref{Euler}:

Setting $n=1$, we get nothing, but setting $n=2$
we obtain 
\[\sum _{k=1}^{\beta _1}b_{1}^k= 2g+1-\beta _{1,2}=2g-3.\] 
Setting $n=3$ we obtain
\[\sum _{k=1}^{\beta _{2,3}}b_{2,3}^k= 2g-3-\beta _{2,3}=2g-5.\]
Continuing this way, we find the difference of the $b_{3,4}^k$ and the 
$b_{2,4}^k$ in terms of $\beta _1, \beta _{2,3}, \beta _{2,4}, \beta _{3,4}$, 
by setting $n=4$. This gives 
\[\sum _{k=1}^{\beta _{3,4}}b_{3,4}^k -\sum _{k=1}^
{\beta _{2,4}}b_{2,4}^k = (\beta _{2,3}-4)g +(\beta _{2,3}+\beta _{3,4}-
\beta _{2,4}-6)=-2g-7.\]}
\end{exa}

We will now ``push down'' results for $S''$ in $\T_0$ to results for 
$S'$ in $\T$.

\begin{defn} \label{push}
We define, for integers $a$ and $b$,
\index{$\O_{\T}(a\H+b\F)$}\[ \O_{\T}(a\H+b\F):=i_*\O_{\T_0}(a\H+b\F).\]
In particular, by the projection formula,
\[i_* \O_{\T_0}(a\H_0+b\F)=i_*\O_{\T_0}(a\H+(a+b)\F)=\O_{\T}(a) \otimes i_*((a+b)\F).\]
\end{defn}
We now return to the general situation. As a consequence of Proposition 
\ref{upperres} we have the following result:

\begin{prop} \label{godeb}
If $D^2=0$ and $c=1$, then the $\O_{\T}$-resolution $F_*$\index{$F_*$} 
of $\O_{S'}$ is
\[ 0 \hpil \O_{\T}(-3\H+(g-4)\F) \hpil \O_{\T} \hpil \O_{S'} \hpil
0.\]
In all other cases we assume $b_i^k \geq i $, for $i=1, \ldots ,c-1$
and all $k$ and $b_{i,j}^k \geq j-1$
for $j=i+1,i+2$, $i=c, \ldots ,\frac{1}{2} D^2+c$, and all $k$.

If $D^2=0$ and $c \geq 2$, then 
\begin{eqnarray*}
0 \hpil \O_{\T}(-(c+2)\H+(g-c-3) \F)  \hpil \+ _{k=1}^{\beta _{c-1}} 
\O_{\T}(-c\H + (b_{c-1}^k-c) \F)  \\ \hpil  \cdots  
 \hpil \+ _{k=1}^{\beta _1} \O_{\T}(-2\H + (b_1^k-2) \F)  \hpil 
\O_{\T} \hpil \O_{S'} \hpil 0, 
\end{eqnarray*}
is an $\O_{\T}$-resolution of $\O_{S'}$.

If $D^2 \geq 2$, 
and if there exists a resolution as described in
Proposition \ref{upperres}, in particular if
the Betti-numbers are the same for all 
$\{\varphi_L(D_{\lambda}) \}$, 
then 
$\O_{S'}$ has a $\O_{\T}$-resolution $F'_*$  of the following type:
\begin{eqnarray*} 
0  \hpil F' _{\frac{1}{2} D^2+c}  \hpil &  \cdots  & \hpil  F' _{c+1} \hpil F' _c \\
   \hpil \+ _{k=1}^{\beta _{c-1}} \O_{\T}(-c\H + (b_{c-1}^k-c)\F)  
   \hpil &  \cdots  & \hpil \+ _{k=1}^{\beta _2} \O_{\T}(-3\H +(b_2^k-3)\F) \\
   \hpil \+ _{k=1}^{\beta _1} \O_{\T}(-2\H + (b_1^k-2)\F)  
   \hpil & \O_{\T} & \hpil  \O_{S'} \hpil 0. 
\end{eqnarray*}
Here $F'_i=i_*(F_i)$, for all i. 
\end{prop}

\begin{proof}
See \cite[p.117]{Sc}.  The essential fact is that the map $i:\T_0 \iso  
\PP (\E) \khpil \T$ is a rational resolution of singularities, and that we 
therefore have $R^1i_*\O_{\T_0}=0$. Moreover $i_* \O_{S''}=\O_{S'}$, 
and $i_* \O_{\T_0}=\O_{\T}$. The condition on the $b_i$ and the 
$b_{i,j}$ gives that each term (except $\O_{S''}$) in the resolution of 
$\O_{S''}$ in Proposition  \ref{upperres} is an extension of terms of the 
form $\O_{\T_{0}}(a \H +b\F)$, with $b \geq -1$. As in 
\cite[(3.5)]{Sc} we then conclude that the
resolution therefore remains exact after pushing down.
\hspace{0.07cm} $\square$ \end{proof}

\begin{rem}
  {\rm By Proposition \ref{ideal} we already know that the ideal of
  $S'$ in $\T$ is the push-down by $i$ of the ideal of $S''$ in
  $\T_0$.

  If $b_2^k \geq 2$ for all $k$ when $c \geq 3$ or $D^2=0$ (resp. $b_{2,3}^k 
  \geq  2$ and $b_{2,4}^k \geq 3$ when $c \leq 2$ and $D^2>0$), then it
  automatically follows that $R^1 i_*F_2=R^1 i_*F_1=0$, so that we get
  an exact pushed-down right end
\[ i_*F_1 \hpil \O_{\T} \hpil \O_{S'} \hpil 0. \]
This means that the ideal of
  $S'$ in $\T$ is generated by the push-down by $i$ of the generators
  of the ideal of $S''$ in $\T_0$. }
\end{rem}

The next two results give the first examples of applications of the proposition.

\begin{cor} \label{c1zeroint}
Assume $D^2=0$. 

(a) If $c=1$ then $\O_{S'}$ has the  
following $\O_{\T}$-resolution:
\begin{equation}
  0 \hpil \O_{\T}(-3\H+(g-4)\F) \hpil \O_{\T} \hpil \O_{S'} \hpil 0.
\end{equation}

(b) If $c=2$, then $\O_{S'}$ has the following $\O_{\T}$-resolution:
\begin{eqnarray*}
  0 \hpil \O_{\T}(-4\H+(g-5)\F)  & \hpil &  \\
\O_{\T}(-2\H+a_1\F) \+ \O_{\T}(-2\H+a_2\F) & \hpil & \O_{\T} \hpil \O_{S'} 
\hpil 0,
\end{eqnarray*}
for two integers $a_1$ and $a_2$ such that $a_1 \geq a_2 \geq 0$ and $a_1+a_2=g-5$.
\end{cor}

\begin{proof}
Set $g_0=g+c+2$.

If $c=1$, then this is a part of Proposition \ref{godeb}. The essence
is as follows: By Proposition \ref{zeroint} a resolution of
$\O_{S''}$ as 
an $\O_{\T_0}$-module is
\begin{eqnarray*}
0 \hpil \O_{\T_0}(-3\H_0+(g_0-4)\F) \hpil \O_{\T_0} \hpil \O_{S''} 
\hpil 0.
\end{eqnarray*} 
Here $g_0-4 \geq 5$, whence (a) follows.

If $c=2$, we have $g_0 \geq 10$ and Proposition \ref{resolv}(a), gives 
that a resolution of $\O_{S''}$ as an 
$\O_{\T_0}$-module is:
\begin{eqnarray*}
0 \hpil \O_{\T_0}(-4\H_0+(g_0-5)\F)  &\hpil &  \\ 
\O_{\T_0}(-2\H_0+b_1\F) \+ \O_{\T_0}(-2\H_0+b_2\F) & \hpil & \O_{\T_0} 
\hpil \O_{S'} \hpil 0,
\end{eqnarray*}
for two integers $b_1$ and $b_2$ such that $b_1 \geq b_2 \geq 0$ and
$b_1+b_2=g_0-5$. From \cite[Thm. 5.1]{B} we have that $S''$ is 
singular along a curve if $b_1 > \frac{2(e_1+e_2+e_3)}{3}$, where
$(e_1,e_2,e_3)$ denotes the type of $\T_0$. This is
equivalent to  $b_2 < \frac{g_0-9+2e_4}{3}$. 
Hence $b_2 \geq 1$ and (b) follows.
\hspace{0.07cm} $\square$ \end{proof}

\begin{cor} \label{lowestexc}
  Let $c=1$, $D^2=2$ and $g=6$ as in Example \ref{simis}. 
Then $\O_{S''}$ has the following $\O_{\T_0}$-resolution:
\begin{eqnarray*} 
0 & \hpil & \O_{\T_0}(-4\H_0 + 6\F) \+\O_{\T_0}(-4\H_0 + 5\F)  \\
& \hpil & \O_{\T_0}(-2\H_0 + 4\F) \+ \O_{\T_0}(-3\H_0 + 4\F) \+
\O_{\T_0}(-3\H_0 + 3\F) \\
& \hpil &\O_{\T_0} \hpil \O_{S''} \hpil 0. 
\end{eqnarray*}

  In particular, $\O_{S'}$ has the following $\O_{\T}$-resolution:
\begin{eqnarray*} 
0 & \hpil & \O_{\T}(-4\H + 2\F) \+\O_{\T}(-4\H + \F)  \\
& \hpil & \O_{\T}(-2\H + 2\F) \+ \O_{\T}(-3\H + \F) \+
\O_{\T}(-3\H) \\
& \hpil &\O_{\T} \hpil \O_{S'} \hpil 0.
\end{eqnarray*}  
\end{cor}

\begin{proof}
In Example \ref{simis} the minimal resolutions of all the
$\varphi_L(D_{\lambda})$ are given. Corollary \ref{maxbett} gives $b^1_{1,2}=4$ and $b^1_{2,4}+
b^2_{2,4}=11$, while inserting $n=3$ in Proposition \ref{Euler} gives 
$b^1_{1,3}+b^2_{1,3}=7$. 
Then Proposition \ref{upperres}
gives the following resolution:
\begin{eqnarray} 
\label{elvis} 0  \hpil  \O_{\T_0}(-4\H_0 + b_{2,4}^1\F) \+\O_{\T_0}(-4\H_0 + (11-b_{2,4}^1)\F)   
& \hpil & F_1 \\
\nonumber  \hpil  \O_{\T_0} \hpil  \O_{S''} & \hpil & 0,
\end{eqnarray}
where $F_1$ is an extension
\begin{eqnarray} 
\label{presley} 0 & \hpil &  \O_{\T_0}(-2\H_0 + 4\F)   \hpil  F_1 \\
\nonumber  & \hpil &   \O_{\T_0}(-3\H_0 + b_{1,3}^1\F) \+
\O_{\T_0}(-3\H_0 + (7-b_{1,3}^1)\F) \hpil 0.
\end{eqnarray}

Without loss of generality we assume $b:=b_{1,3}^1 \geq 4$, and $a:=b_{2,4}^1
\geq 6$. The type of $\T_0$ is $(3,2,1,1)$.

Look at the composite morphism given by (\ref{elvis}) and
(\ref{presley})
\begin{eqnarray*} 
 \alpha &   :   & \O_{\T_0}(-4\H_0 + a \F) \+\O_{\T_0}(-4\H_0 + (11-a)\F) \\
        & & \hpil  \O_{\T_0}(-3\H_0 + b\F) \+ \O_{\T_0}(-3\H_0 + (7-b)\F). 
\end{eqnarray*}
Now $\alpha$ can be expressed by a matrix
\[  \left[ 
  \begin{array}{cc}
   \alpha_1 & \alpha_2 \\
   \alpha_3 & \alpha_4
    \end{array} \right],  \] 
with
\begin{eqnarray*}
  \alpha_1 & \in & H^0(\H_0 + (b-a)\F) \\
  \alpha_2 & \in & H^0(\H_0 + (a+b-11)\F) \\
  \alpha_3 & \in & H^0(\H_0 + (7-a-b)\F) \\
  \alpha_4 & \in & H^0(\H_0 + (a-b-4)\F), 
\end{eqnarray*}
whose determinant gives a section $g \in H^0 (2\H_0 - 4\F)$ whose zero
scheme contains $S''$.

If $(a,b) \not = (6,4)$, we have
\begin{eqnarray*}
H^0 (\H_0 + (7-a-b)\F) & = & \\
H^0 (\PP^1, \O_{\PP^1}(10-a-b) & \+ & \O_{\PP^1}(9-a-b) \+ \O_{\PP^1}(8-a-b)^2)=0, 
\end{eqnarray*}
whence $\alpha_3 =0$ and $g$ is  a product of two sections of $\H_0 +
(b-a)\F$ and $\H_0 + (a-b-4)\F$ respectively. But then $S''$ would have
degenerate fibers $S''_{\lambda}$, contradicting Proposition \ref{Bettiprop}(b).

So $(a,b) = (6,4)$ and we compute
\begin{eqnarray*}
  \Ext ^1 (\O_{\T_0}(-3\H_0 + 4\F) \+
\O_{\T_0}(-3\H_0 + 3\F), \O_{\T_0}(-2\H_0 + 4\F) & = & \\ 
   H^1 (\H_0) \+ H^1(\H_0+\F )  & = &  \\
   H^1 (\PP^1, \O_{\PP^1}(3) \+ \O_{\PP^1}(2) \+ \O_{\PP^1}(1)^2)  \+ & & \\
   H^1 (\PP^1, \O_{\PP^1}(4) \+ \O_{\PP^1}(3) \+ \O_{\PP^1}(2)^2) & = & 0,
\end{eqnarray*}
whence the sequence (\ref{presley}) splits and the first assertion
follows.

The second is then an immediate consequence of Proposition \ref{godeb}.
\hspace{0.07cm} $\square$ \end{proof}

Note that by this result, $S''$ is cut out in $\T_0$ by three
sections $q$, $c_1$ and $c_2$ of $\O_{\T_0}(2\H_0 -4\F)$, 
$\O_{\T_0}(3\H_0 -4\F)$ and  $\O_{\T_0}(3\H_0 - 3\F)$ respectively.

Now look at the three dimensional subvariety $V$ of $\T_0$ defined by
$q \in \O_{\T_0}(2\H_0 -4\F)$. Arguing as in the proof of Proposition
\ref{Chowring} we find that the class of $i (V)$ in the Chow group of
$\T$ is $2 \H_{\T} -2 \F_{\T}$, whence $i(V)$ has degree $4$ and
dimension $3$ in $\PP^{6}$. As in \cite[(7.12)]{S-D} we have that
$i(V)$ is a cone over the Veronese surface (whose vertex is the image
of $\Gamma$) and that this variety is the (reduced) intersection of all quadrics
containing $S'$.

A very useful result is the
following, involving so called  
``rolling factors'' coordinates\index{rolling factors coordinates} (see for example  \cite[p.59]{H},  
\cite[p.3]{St} or \cite{Re}):

\begin{Lemma} \label{roll}
The sections of $a\H-b\F$ on a smooth rational normal scroll of type 
$(e_1, \ldots ,e_d)$ can be identified
with weighted-homogeneous polynomials of the form 
\[P=\sum P_{i_1,   \ldots   ,i_d}(t,u)Z_1^{i_1} \ldots Z_d^{i_d},\]
where $i_1 + \cdots + i_d=a$, and $ P_{i_1, \ldots ,i_d}(t,u)$ is a
homogeneous polynomial of degree $-b+(i_1e_1+ \cdots +i_de_d)$. 

If we multiply $P$ by a homogeneous polynomial of degree $b$ in $t,u$, then
we get a defining equation of the zero scheme of the section, in term of
homogeneous coordinates of the projective space, within which the scroll is
embedded. Here $X_{k,j}=t^ju^{e_k-j}Z_k$, for $k=1, \ldots ,d$, and $j=0, \ldots ,e_k$,
are coordinates for this space. The equation is uniquely determined modulo
the homogeneous ideal of the scroll. 
\end{Lemma}

As a first application, we prove the analogue of Corollary \ref{c1zeroint}
for $c=3$.

\begin{cor} \label{c3}
  Let $D^2=0$ and $c=3$. Then a resolution of $\O_{S'}$ as an
 $\O_{\T}$-module is:
\begin{eqnarray*}
0 \hpil \O_{\T}(-5\H+(g-6)\F)  \hpil \+ _{i=1}^5 
\O_{\T}(-3\H+a_i\F)  & \hpil &  \\
 \+ _{i=1}^5 \O_{\T}(-2\H + b_i \F)  \hpil 
\O_{\T} \hpil \O_{S'} & \hpil & 0, 
\end{eqnarray*}
where all the $a_i$ and $b_i \geq -1$ and satisfy $\sum_{i=1}^5 b_i = 2g-12$
and $a_i=g-6-b_i$. 
\end{cor}

\begin{proof}
First recall that $g \geq 9$.

As a special case of Proposition \ref{upperres} we obtain that the resolution 
of $\O_{S''}$ as an $\O_{\T_0}$-module is:
\begin{eqnarray*}
0 \hpil \O_{\T_0}(-5\H_0+(g-1)\F)  \hpil \+ _{k=1}^5 
\O_{\T_0}(-3\H_0+b_2^k\F)  \hpil \\
 \+ _{k=1}^5 \O_{\T_0}(-2\H_0 + b_1^k \F)  \hpil 
\O_{\T_0} \hpil \O_{S''} \hpil 0. 
\end{eqnarray*}
From Corollary \ref{maxbett} we get $\sum_{i=1}^5 b_1^k = 2g-2$.
The self-duality of the resolution in this particular case gives 
$b_2^k=g-1-b_1^k$, if we for example order the $b_1^k$ in a non-increasing 
way, and the $b_2^k$ in a non-decreasing way. 

We will show that for all $g \geq 9$ all the $b_1^k \geq 1$ and all 
the $b_2^k \geq 2$, so that we can push down the resolution to one of 
$\O_{S'}$ as an $\O_{\T}$-module. 

Look at the map 
\[ \Phi: \+ _{k=1}^5\O_{\T_0}(-3\H_0+b_2^k\F)  \hpil 
 \+ _{k=1}^5 \O_{\T_0}(-2\H_0 + b_1^k \F).\]
Just like in the analysis of pentagonal curves in \cite{Sc}, it follows 
from \cite{B-E} that the map $\Phi$ is skew-symmetrical and that its 
Pfaffians generate the ideal of $S''$ in $\T_0$. See also \cite{W}.

Let the type of $\T_{0}$ be $(e_1, \ldots ,e_5)$, where $e_5=1$ and 
$\sum_{i=1}^4e_i=g$.

A key observation is the following: $b_1^1 \leq e_1+e_3$ and $b_1^1 \leq 2e_2$.
The first inequality holds, since otherwise we would have a quadratic relation
of the form $f(t,u,Z_1,Z_2)=0$ in each fiber. Hence the general fiber would be
reducible, a contradiction. The second inequality follows since its negation
implies that $Z_1$ is factor in one quadratic relation satisfied by the points 
of $S''$, a contradiction.  This gives $b_1^1 \leq \frac{2(e_1+e_2+e_3)}{3}
=\frac{2(g-e_4)}{3}$. Hence 
\[b_2^1=g-1-b_1^1 \geq g-1-\frac{2g}{3}+\frac{2e_4}{3}=\frac{g-3+2e_4}{3}
\geq 2,\] 
since $g \geq 9$. Hence $b_2^k \geq b_2^1 \geq 3$ for all $k$.

Another key observation is the following: $b_1^2 \leq e_1+e_4$ and 
$b_1^2 \leq e_2+e_3$. The first inequality holds, since otherwise the two-step
projection of the general fiber $D''$ of $S''$ into the $Z_1,Z_2,Z_3$-plane
from $P=(0,0,0,0,1)$ and $Q=(0,0,0,1,0)$ would be contained in 2 quadrics.
This in only possible if the projected image is a line, and in that case
the general fiber would be degenerate (contained in the $\PP^3$ spanned by
this line and $P$ and $Q$). This is impossible. The second inequality holds,
since otherwise there would be two independent relations of the form
\[Z_1f(t,u,Z_1, \ldots ,Z_5)+aZ_2^2=0\]
for each fiber. In that case we could eliminate the $Z_2^2$-term and obtain 
one relation with $Z_1$ as a factor, a contradiction.

These two inequalities for $b_1^2$ imply 
\[b_1^2 \leq \frac{e_1+e_2+e_3+e_4}{2}=\frac{g}{2}.\] 
Now we assume for contradiction that $b_1^5 \leq 0$. Then we get
\[b_2^5-b_1^k \geq b_2^5-b_1^2 = (g-1-b_1^5) - b_1^2 \geq g-1 -\frac{g}{2}
=\frac{g}{2}-1,\]
for $k=2,3,4$. In the matrix description of the map $\Phi$ there is one 
submaximal minor
with one column consisting of zero and sections of $\H_0-(b_2^5-b_1^k)\F$, for 
$k=2,3,4$. If all entries of this column have $Z_1$ as a factor, that would
lead to a contradiction, since the minor is the square of one of the
generators of the (Pfaffian) ideal of $S''$ on $\T_0$.
To avoid that $Z_1$ is a factor in each such entry, we must have $e_2 \geq 
\frac{g}{2}-1$. This gives $e_1+e_2 \geq g-2$, and $e_3+e_4 \leq g-(e_1+e_2) 
\leq g-(g-2)=2$. Hence $e_3=e_4=e_5=1$. But this implies that 
$D^2+h^1(R) \geq 3$, contradicting Proposition \ref{bound}.

Hence the assumption $b_1^5 \leq 0$ leads to a contradiction, and $b_1^k \geq
1,$ for all $k$. Hence the entire resolution can be pushed down to one of
$\O_{S'}$ as an $\O_{\T}$-module and the result follows. 
\hspace{0.07cm} $\square$ \end{proof}

\begin{rem} \label{f_*}
{\rm Assume that we are in the situation of Proposition
 \ref{upperres} (i.e. the Betti-numbers of all the
$\varphi_L(D_\lambda)$ are the same for all $\lambda$, for instance if $D^2 \leq 4$),
so that a finite set of sections of 
line bundles of type $a\H_0-b\F$ generate the ideal of the surface $S''$ on 
the smooth rational normal scroll $\T_{0}$ of type 
$(e_1, \ldots ,e_{d-r-1}, e_{d-r},   \ldots   ,e_d)$, where $e_{d-r}= \cdots =e_d=1$, $e_{d-r-1} \geq 2$, and
$V=\Sing \T \iso \PP^{r-1}$ for some $r \geq 0$. 
Let $W=i^{-1}(V)$. This is a subscroll of $\T_0$ of type $(1,   \ldots   ,1)$,
that is $\PP^r \times \PP^1$. 
The ideal generators can be classified into $3$ types:
\begin{itemize}
\item[(a)] Those that are sections of $a\H_0-b\F=a\H-(b-a)\F$, with $b > a$.
\item[(b)] Those that are sections of $a\H_0-b\F=a\H-(b-a)\F$, with $b = a$.
\item[(c)] Those that are sections of $a\H_0-b\F=a\H-(b-a)\F$, with $b < a$.
 \end{itemize}
For those of type (a) it is clear from Lemma \ref{roll} that their zero 
scheme contains $W$. Likewise one sees that those of type (b) can be written as
\[f(t,u,Z_1, \ldots ,Z_d)+ g(Z_{s+1}, \ldots ,Z_d),\]
where $Z_1$ or $Z_2$ or  \ldots  or $Z_s$ is a factor in every monomial of
$f(t,u,Z_1, \ldots ,Z_d)$, while $g$ is homogeneous of degree $a=b$.
Those of type (c) can be written as
\[f(t,u,Z_1, \ldots ,Z_d)+ h(t,u,Z_{s+1}, \ldots ,Z_d),\]
where $f$ is as described for type (b), while $h$ is bihomogeneous, of degree
$a-b>0$ in $t,u$ and degree $a$ in $Z_{s+1}, \ldots , Z_d$.

There is one fundamental difference between the sections of types (a) and (b)
on one hand, and those of type (c) on the other. Those of types (a) and (b)
are ``constant'' on the fibers of $i$, their zero scheme contains either the
whole fiber, or no point on the fiber, for each $\PP^1$, which is a fiber
of $i$. For the sections of type (c) this is only true if $h$ is the zero
polynomial, and then its zero scheme contains all of $W$. 

We therefore see (referring to the notation of Proposition \ref{upperres}) 
that if  $b_1^k \geq 2$, for $k=1,\ldots,\beta_1$ (we must use the formulation
$b_{1,j}^k \geq j$, for $j=2,3$ and $k=1,\ldots,\beta_{1,j}^k$ in the special 
case (E0) in Corollary \ref{lowestexc} with $c=1$, and $D^2=2$) in the resolution of $\O_{S''}$ as an 
$\O_{\T_0}$-module, then the ideal of $S''$
is generated by ``fiber constant'' equations, and if $Q$ is a point on $\T_0$ 
not on $S''$, then there is a fiber constant section of the type described, 
which does not contain $Q$ in its zero scheme. In short, fiber
constant equations cutting out $S''$ in $\T_0$ are also equations of
$S'$ in $\T$. }
\end{rem}

In Section \ref{Conc} we will classify the possible projective models
for $g \leq 10$,  and in particular the singular scrolls $\T$ appearing as
$\T(c,D,\{ D _{\lambda} \})$ in the various cases, and we will also show that 
projective models giving scrolls of all these types exist. 

In
Corollaries \ref{c1zeroint}, \ref{lowestexc} and \ref{c3} we showed
that in some particular cases we can push down the entire resolution
of $\O_{S''}$ to one of $\O_{S'}$.
In the rest of this section, through a series 
of additional examples, we take a closer look at the rest of the
singular scrolls appearing for $g \leq 10$, and using Lemma
\ref{roll} we will find restrictions on the $b_1^k$. 

Recall the $b_1^k$ and
$\beta_{i}$ described in Proposition \ref{upperres} (see also Remark \ref{ikkehverb}). To make the
notation simpler in the examples below, we give the following:

\begin{defn} \label{bene}
Let $b_k$ denote $b_1^k$, for $i=1, \cdots\ \beta_{1}$.
\end{defn}

In all the examples below we have $D^2 \leq 4$, so by Corollary \ref{Bettis} the Betti-numbers of the 
$\varphi_L(D_{\lambda})$  are independent
  of the $\lambda$.

\begin{exa} \label{merall2}
{\rm We return to the situation studied in Example \ref{Dkvadrat2} and 
Example \ref{all2}, with $c=2$, $D^2=2$ and $g \geq 7$. From Example 
\ref{all2} we see that the ideal of $S''$ in $\O_{\T_0}$ is generated
by four sections of the type $2\H_0-b_k \F$, where $\sum_{k=1}^4 b_k = 2g-3$.

The type of $\T_0$ is $(e_1+1,e_2+1,e_3+1,1,1)$,
where all $e_i \geq 0$ and $e_1+e_2+e_3=g-4$.
Let $Q$ be the subscroll of $\T_0$ formed by the two last directrices,
so $Q$ is the inverse image by $i$ of the line in $\PP^g$ spanned by the 
images by $\varphi_L$ of the basepoints of $D$.
We see that $Q$ is a quadric surface in $\PP^3$.
All the four sections of type $2\H_0-b_k\F$ must intersect $Q$ in, 
and therefore contain, the two lines
that form the inverse image by $i$ in $\T_0$ of the images 
by $\varphi_L$ of the basepoints of $D$ . But this is simply the two last 
directrices.
The intersection with $Q$ for one such section is obtained by using
Lemma \ref{roll} to express each of the sections, and set $Z_1=Z_2=Z_3=0$.
What remains must be a term of the type $P_{2-b_k}(t,u)Z_4Z_5$, where
$P_{2-b_k}$ is zero if $b_k \geq 3$, and a polynomial of degree $2-b_k$ 
otherwise.

We order the $b_k$ as $b_1 \geq b_2 \geq b_3 \geq b_4$. We see that if $b_4 
\leq 1$, in particular if 
$b_4 \leq 0$, then $b_3 \leq 2$, since
otherwise the total intersection of $Q$ with the four sections will consist
of $2-b_4$ lines transversal to the two directrices in addition to the two 
directrices.

If $g=7$, it is clear that
$\T_0$ has type $(2,2,2,1,1)$ or $(3,2,1,1,1)$. 
Then $Z_1$ is a factor in all sections of 
$2\H_0-b\F$ for $b \geq 5$ for both scroll types. Hence $b_1 \leq 4$.
If $b_1=3$, then the only possible combination is
$(b_1, \ldots ,b_4)=(3,3,3,2)$, since $\sum b_k=11$.
If $b_1=4$, then the only a priori possibilities are $(4,4,2,1)$ and
$(4,3,2,2)$. But $(4,4,2,1)$ is impossible for type $(2,2,2,1,1)$, 
since we then have two quadratic relations between $Z_1,Z_2,Z_3$ only. To see 
that this is impossible, let $D''$ be any smooth curve in $|f^* D-E|$, which can be identified 
with its image under $\varphi_H$.
The variables $Z_1, Z_2, Z_3$ restricted to $D''$ correspond to sections of 
$(H-E)_{D_0}$. Since this line bundle has degree $2g(D'')=4$, it is base 
point free and its sections map $D''$ into $\PP^2$ by a one-to-one or 
two-to-one map. This means that there is at most one quadratic relation between 
$Z_1,Z_2,Z_3$, whence $b_2 \leq 3$. So for type $(2,2,2,1,1)$
the only possibilities for $(b_1, \ldots ,b_4)$ are
\[(3,3,3,2) \hs \mbox{ and } \hs  (4,3,2,2).\]

For the type $(3,2,1,1,1)$, for each fiber of $\T_0$, the equations 
with $b_k \leq 2$, restricted to the subscroll $Z_1=Z_2=0$ with plane fibers, 
must cut out a subscheme of length $4$ (such that each 
subscheme of length $3$ spans a $\PP^2$) (these are the cases (E1) and (E2)). 
It takes 2 equations to do this.
Hence $b_3 \leq 2$. Moreover $b_2 \leq 3$. Assume $b_2 \geq 4$.
Then we would have two independent equations of type 
\[Z_1f(t,u,Z_1,Z_2,Z_3,Z_4,Z_5)+Z_2^2\]
for general $(t,u)$. From these equations we can eliminate the $Z_2^2$-term,
and derive one equation with $Z_1$ as a factor. This is a contradiction. 

Hence the only possibility for $(b_1, \ldots ,b_4)$ for the type 
$(3,2,1,1,1)$ is $(4,3,2,2)$.

If $g=8$, then $\T_0$ has type $(3,2,2,1,1)$. A similar argument as for $g=7$,
gives that the only possible combinations for $(b_1, \ldots ,b_4)$ are
\[(5,4,2,2), \hs (5,3,3,2) \hs \mbox{ and } \hs  (4,4,3,2).\]

For $g=9$ the type of $\T_0$ is a priori either $(3,3,2,1,1)$ or $(4,2,2,1,1)$.
In Section \ref{Conc} the type $(4,2,2,1,1)$ is ruled
out when $D$ is perfect. 
For the type $(3,3,2,1,1)$ we see that $b_1 \geq 6$ is 
impossible, since $2\H_0-6\F$ only has sections of the form $f(Z_1,Z_2)$. 
Moreover $b_4 \leq 2$, since we need to cut out the exceptional fibers.
Hence $b_1=5$, otherwise the sum of the $b_i$ would be at most $14$, and it is $15$. Any section of $2\H_0-5\F$ can be written in terms of 
$t,u,Z_1,Z_2,Z_3$ only. 
As for one case with $g=7$ we see that  we cannot have two quadratic relations 
between $Z_1,Z_2,Z_3$ for general fixed $(t,u)$, so $b_2 \leq 4$. We then see 
that the only possible combination is $(b_1,b_2,b_3,b_4)=(5,4,4,2)$.

For $g=10$ the type of $\T_0$ is a priori $(3,3,3,1,1)$, 
$(4,3,2,1,1)$ or $(5,2,2,1,1)$. In Section \ref{Conc} it is shown that only
the type $(4,3,2,1,1)$ occurs when $D$ is perfect.
In this case a more detailed analysis gives that the only
possible combinations for $(b_1,   \ldots   ,b_4)$ are
\[(6,5,4,2) \hs \mbox{ and } \hs (5,5,5,2).\]
 }
\end{exa}

\begin{exa} \label{non-prim}
{\rm Let us study the case $c=2$, $D^2=4$ and $g=9$. This gives a non-primitive
projective model, with $L \iso 2D$ (it is the case (Q) described in the text 
above Theorem \ref{mainsing}). The scroll $\T$ necessarily has type 
$(2,1,1,0,0,0)$ and $\T_{0}$ has type $(3,2,2,1,1,1)$.

We will now find the $b_k$. 

Lemma \ref{roll} gives $b_k \leq 4$, for all $k$, since
$Z_1$ is factor in every section of $2\H_0-b\F$, for $b \geq 5$. The complete
intersection $Z_1=Z_2=Z_3=0$ in $\T_0$ is a subscroll $N$ of type $(1,1,1)$ 
with a plane in each fiber. The $7$ equations cutting out $S''$ in $\T_0$ must
together cut out four points in each plane fiber, such that no three of 
these points are collinear, by Theorem \ref{mainsing}. It is clear that such a
configuration of points is contained in exactly two quadrics in each plane.
All sections of $2\H_0-b\F$ with $b \geq 3$ vanish on $N$, so we must have at least two of the $b_k$ less than $3$.  
Moreover, every section of $2\H_0-4\F$ can be written 
\[Z_1f(t,u,Z_1, \ldots ,Z_5)+aZ_2^2+bZ_2Z_3+cZ_3^2.\]
If $b_4 \geq 4$, there are four independent equations of this kind,
so we could eliminate the three last terms and obtain one relation with $Z_1$
as a factor. This is a contradiction, so $b_4 \leq 3$.
 
This leaves the unique possibility $(b_1, \ldots ,b_7)=(4,4,4,3,3,2,2)$.}   
\end{exa}

The following proposition describes the particular case in Proposition
\ref{2-uple}, that is we have $L \sim 2D$, $D^2=4$ and $c=2$ and
$D$ is hyperelliptic, which is also a particular case of the last
example. In this case there is a smooth curve $E$ (which is a perfect
Clifford divisor for $D$)
satisfying $E^2=0$ and $E.D=2$. Since $(D-E)^2 = 0$ and $(D-E).L =4$,
we have $D > E$, so $E$ does not satisfy the conditions (C6) and (C7). We will also see that $E$ is not always a perfect Clifford divisor for $L$.

\begin{prop} \label{non2-uple}
  Let $L$ and $D$ be as in the particular case of Proposition
  \ref{2-uple} (where $S' \sub \PP^9$ is not the $2$-uple embedding of
  $\varphi_D(S)$). 

  Then we are in one of the following three cases:
  \begin{itemize}
  \item[(i)] \hs $\R_{L,E} = \emptyset$ and $D \sim E + E'$, where $E'$ is
    a smooth elliptic curve such that $E.E'=2$. 
  \item[(ii)] \hs $\R_{L,E} = \{ \Gamma_1, \Gamma_2 \}$ and 
               $D \sim 2E + \Gamma_1 +\Gamma_2$. 
  \item[(iii)] \hs $\R_{L,E} = \{ \Gamma_0 \}$ and $D \sim 2E + \Delta_0$,
    where $\Delta_0$ has a configuration with respect to $E$ as in (E2). 
  \end{itemize}

  Let $\T = \T(2,E)$ be the scroll defined by $|E|$. 
  
In case (i), $\T$ is of type $(2,2,2,0)$ and $\O_{S'}$ has the following 
$\O_{\T}$-resolution:

\begin{eqnarray*}
  0 \hpil \O_{\T}(-4\H+4\F)   & \hpil & \O_{\T}(-2\H+4\F) \+ \O_{\T}(-2\H)  \\
                              & \hpil &  \O_{\T} \hpil \O_{S'} \hpil 0.
\end{eqnarray*}
In this case $\Sing \T \cap S'= \emptyset$, so $S' \iso S''$
where $S''$ sits in the smooth scroll $\T_0$ of type $(3,3,3,1)$. 

In the cases (ii) and (iii), $\T$ is of type $(4,2,0,0)$ and its singular locus 
is spanned by 
$< Z_{\lambda} >$ (using the same notation as in Theorem \ref{mainsing}). Moreover
$\O_{S'}$ has the following 
$\O_{\T}$-resolution:
\begin{equation}
  0 \hpil \O_{\T}(-4\H+4\F)  \hpil 
\O_{\T}(-2\H+4\F) \+ \O_{\T}(-2\H)  \hpil  \O_{\T} \hpil \O_{S'} 
\hpil 0.
\end{equation}
\end{prop}

\begin{proof}
 The three cases follow from Proposition \ref{singc=0}, by noting that
 we clearly have $\R_{L,E} = \R_{D,E}$.

 We have $h^0(L)=10$ and $h^0(L-E)=6$. We leave it to the reader to
 verify that in case (i) we have
\[ h^0(L-2E)=3 \hs \mbox{ and }  \hs h^0(L-3E)=0, \]
and that we in the cases (ii) and (iii) have
\[ h^0(L-2E)=3, \hs  h^0(L-3E)=2, \hs  h^0(L-4E)=1  \hs \mbox{ and }
\hs h^0(L-5E)=0. \]

 This yields the two scroll types $(2,2,2,0)$ and $(4,2,0,0)$
 respectively.

 In the cases (ii) and (iii), one can show as in the proof of Theorem
 \ref{mainsing} that $\Sing \T =  < Z_{\lambda} >$. 
 
The statement about the resolution in part (ii) and (iii) follows from Proposition 
\ref{upperres} and the upper (large) table in Section \ref{stypes} below. 
In case (i) the corresponding statement follows in part from these results.
Proposition \ref{upperres} and the table give that the resolution of 
$\O_{S''}$ in case (i) is
\[
  0 \hpil \O_{\T_0}(-4\H_0+8\F)  \hpil \O_{\T_0}(-2\H_0+4\F)^2  
\hpil \O_{\T_0} \hpil \O_{S''} 
\hpil 0,
\]
or
\begin{eqnarray*}
  0 \hpil \O_{\T_0}(-4\H_0+8\F)  & \hpil &  \O_{\T_0}(-2\H_0+6\F) \+ \O_{\T_0}(-2\H_0+2\F)  \\
                                 & \hpil &  \O_{\T_0} \hpil \O_{S''} \hpil 0.  
\end{eqnarray*}

On the other hand it is clear that there are no contractions across the
fibres in this case. Assume we have the upper of these two resolutions.
From Lemma \ref{roll} we then get that $S''$ is cut out in $\T_0$ by two 
equations of the form:
\[P_1(t,u)Z_1^2+P_2(t,u)Z_1Z_2+P_3(t,u)Z_1Z_3+P_4(t,u)Z_2^2+\]
\[P_5(t,u)Z_2Z_3+P_6(t,u)Z_3^2+c_1Z+1Z_4+c_2Z_2Z_4+c_3Z_3Z_4=0.\]
Here all the $P_i(t,u)$ are quadratic in $t,u$, and the $c_j$ are
constants. But both these quations contain the directrix line
$(Z_1,Z_2,Z_3,Z_4)=(0,0,0,1)$ of $\T_0$. This is precisely the inverse image 
of $\Sing \T$. Hence the inverse image of this line on $S$ is contracted,
a contradiction. Hence we are left with the lower of the two resolutions.
From Corollary \ref{c1zeroint} we then get the given resolution of $\O_{S'}$.
The last details in the description of case (i) follow from Remark 
\ref{Dnotperfect}. 
\hspace{0.07cm} $\square$ \end{proof}

\begin{rem}
{\rm   We see that in case (i) above, $E$ is not perfect, since $\Sing \T$ is a point, but there are no contractions across the fibers.

In the cases (ii) and (iii), $E$ is however perfect.

These two cases are therefore included in the table on p. \pageref{tablepage} (under scroll type $(4,2,0,0)$). However, also $D$ is a perfect Clifford divisor, so these cases can also be described as the case with scroll type $(2,1,1,0,0,0)$ in the same table.}
\end{rem}

\begin{exa} \label{c3D^22}

{\rm Let us study the case $c=3$, $D^2=2$ and $g=9$. In Section \ref{Conc} we
will show that projective models with such invariants occur, and that the
scroll type of $\T$ is either $(1,1,1,1,0,0)$ or $(2,1,1,0,0,0)$ when $D$ is perfect.
By Proposition \ref{upperres} we have $\beta_{1,3}=0$, and by
Corollary \ref{maxbett} we have $\beta _{1,2}=8$ and $\sum_{k=1}^8b_k=3g-4=23$. 

Assume first that the type is $(1,1,1,1,0,0)$, which implies 
that $\T_{0}$ has type $(2,2,2,2,1,1)$. 
Lemma \ref{roll} gives $b_k \leq 4$, for all $k$, since
$h^0 (2\H_0-b\F)=0$ for $b \geq 5$. The complete
intersection $Z_1=Z_2=Z_3=Z_4=0$ in $\T_0$ is a subscroll $Q$ of type $(1,1)$ 
with a line in each fiber. The $8$ equations cutting out $S''$ in $\T_0$ must
together cut out two points in each fiber of $Q$ (the inverse image in $S''$ of
$\Sing \T \cap S'$). For a general fiber, call these points $P_1$ and $P_2$.
Order the $b_k$ in a non-increasing way. To cut out the two points we must
have $b_8 \leq 2$, since all sections of $2\H_0-b\F$ vanish on $Q$ for 
$b \geq 3$. 

For general $(t,u)$, where the fiber $D''$ of $S''$
is smooth, $D''$ is a smooth curve of degree $7$ and genus $2$, which can be identified with a smooth curve in $|f^*D-E|$.
The complete linear system $|(H-E)_{D''}|$ is of degree $2g(D'')+1=5$
and in particular very ample. Now $Z_1, Z_2, Z_3, Z_4$ (restricted to $D''$) span 
$H^0((H-E)_{D''})$, which embeds $D''$ as a curve of degree $5$ and genus
$2$ in $\PP^3$.  As in Example \ref{simis} we conclude from \cite{Sim} that
this curve is contained in only one quadric surface.
On the other hand all sections of $2\H_0-4\F$ can be expressed in terms of 
$Z_1, Z_2, Z_3, Z_4$ only.
Hence no more than one of the $b_k$ can be $4$. This leaves us with only two 
possible cases:  
\[ (b_1, \ldots ,b_8)=(4,3,3,3,3,3,2,2) \hs \mbox{ or } \hs
(3,3,3,3,3,3,3,2). \] 
 
Assume now that the type is $(2,1,1,0,0,0)$, which implies 
that $\T_{0}$ has type $(3,2,2,1,1,1)$. 
Lemma \ref{roll} gives $b_k \leq 4$, for all $k$, since
all sections of $2\H_0-b\F$ have $Z_1$ as factor if  $b \geq 5$. The complete
intersection $Z_1=Z_2=Z_3=0$ in $\T_0$ is a subscroll $N$ of type $(1,1,1)$ 
with a plane in each fiber. The $8$ equations cutting out $S''$ in $\T_0$ must
together cut out three independent points in each fiber of $N$ (the inverse image in $S''$ of $\Sing \T \cap S'$).

Order the $b_k$ as above. To cut out the three points we must
have $b_6 \leq 2$, since all sections of $2\H_0-b\F$ vanish on $N$ for 
$b \geq 3$, and a net of three quadrics is needed to cut out three independent
points in a plane. 
On the other hand every section of $2\H_0-4\F$ can be written 
\[Z_1f(t,u,Z_1, \ldots ,Z_5)+aZ_2^2+bZ_2Z_3+cZ_3^2.\]
This gives $b_4 \leq 3$ as in Example \ref{non-prim}.
This leaves
\[ (b_1, \ldots ,b_8) = (4,4,4,3,2,2,2,2) \hs \mbox{ or } \hs
(4,4,3,3,3,2,2,2), \]    
as the only possibilities.}
\end{exa}

\begin{exa} \label{10c3D^22}
{\rm The case $c=3$, $D^2=2$ and $g=10$ is very similar to the analogous one 
for $g=9$, 
treated in Example \ref{c3D^22} and one can show in a similar way that }
\[(b_1, \ldots ,b_8)=(4,4,4,3,3,3,3,2).\]
{\rm We show in Section \ref{Conc} that the only
possible scroll type for $\T$ is $(2,1,1,1,0,0)$ when $D$ is perfect, corresponding to the type
$(3,2,2,2,1,1)$ for $\T_0$. }
\end{exa}

\begin{exa} \label{c3D^24}
{\rm In the case $c=3$, $D^2=4$ and $g=10$, it follows from
  Proposition \ref{Bettiprop} that the Betti-numbers of the 
$\varphi_L(D_{\lambda})$  are independent
  of the $\lambda$. Now we have a projective model of type 
  type (E0) (with $\beta _{1,2}=12$ and $\sum b_k=4g-6=34$) One can
  show that 
\[ (b_1, \ldots ,b_{12}) = (4,4,4,3,3,3,3,2,2,2,2,2). \]

In Section \ref{Conc} it is shown that the only scroll type occurring for $\T$
is $(2,1,1,0,0,0,0)$, which means $(3,2,2,1,1,1,1)$ for $\T_0$. }
\end{exa}

\begin{exa} \label{c3exa}
{\rm  Let $D^2=0$ and $c=3$, as in Corollary \ref{c3}. 

We will show in
Section \ref{Conc} that for $g=9$ the only smooth scroll occurring as
$\T=\T(3,D)$ is of type $(1,1,1,1,1)$, and the singular 
scrolls occurring are of types $(2,1,1,1,0)$,
$(2,2,1,0,0)$ and  $(3,1,1,0,0)$, corresponding to the smooth types 
$(3,2,2,2,1)$,
$(3,3,2,1,1)$ and $(4,2,2,1,1)$ for $\T_0$. We also show that all
these occur. 
By using similar techniques as in the previous examples, one can show
that $(b_1,b_2,b_3,b_4)=(2,1,1,1)$ for the type $(1,1,1,1,1)$,
$(b_1,b_2,b_3,b_4)=(4,3,3,3,3)$, $(4,4,3,3,2)$ or $(4,4,4,2,2)$ for
the type  $(3,2,2,2,1)$, and $(b_1,b_2,b_3,b_4)=(4,4,3,3,2)$ or 
$(4,4,4,2,2)$ for the types $(3,3,2,1,1)$ and $(4,2,2,1,1)$.

For $g=10$ we will show in Section \ref{Conc} that the only smooth scroll 
occurring as $\T=\T(3,D)$ is of type $(2,1,1,1,1)$, and the singular 
scrolls occurring are of types $(2,2,1,1,0)$, $(2,2,2,0,0)$, 
$(3,2,1,0,0)$, corresponding to the types $(3,3,2,2,1)$, $(3,3,3,1,1)$ and 
 $(4,3,2,1,1)$ for $\T_0$. We also show that all these occur. Again 
one can show
that $(b_1,b_2,b_3,b_4) =(2,2,2,1,1)$ for the type $(2,1,1,1,1)$,
$(b_1,b_2,b_3,b_4)=(4,4,4,3,3)$ or $(4,4,4,4,2)$
for the scroll type $(3,3,3,1,1)$, and 
$(b_1,b_2,b_3,b_4)=(5,5,4,2,2)$, $(5,5,3,3,2)$, $(5,4,4,3,2)$,
 $(4,4,4,3,3)$ or $(4,4,4,4,2)$ for the scroll types $(3,3,2,2,1)$ 
and $(4,3,2,1,1)$.
}  
\end{exa}

\section{More on projective models in smooth scrolls of $K3$ surfaces of low Clifford-indices} 
\label{lowind}
%\chaptermark{Models in smooth scrolls of low Clifford-indices}

In this section we will have a closer look at the situation described
in Section \ref{resol} for $c=1$, $2$ and $3$. We recall that $D$ is a free
Clifford divisor on a non-Clifford general polarized $K3$ 
surface $S$, and that $\T=\T(c,D)$ is smooth, which is equivalent to the 
conditions $D^2=0$ and $\R_{L,D}= \emptyset$ when $D$ is perfect. In any
case these two conditions are necessary to have $\T$ smooth, and the pencil 
$D _{\lambda}$ is uniquely determined.
The resolution of $\O_{S'}$ as an $\O_{\T}$-module was given in 
Proposition \ref{resolv}.

By Corollaries \ref{c1zeroint} and \ref{c3} such resolutions exist also if 
$\T$ is singular if  $D^2=0$. We will use this to take a closer look also at 
the situation for singular $\T(c,D)$ when $D^2=0$ and $c=1,2,3.$
We end the section with a statement valid for general $c$.

From the proof of Theorem \ref{exthm} it is  
clear that for each of the possible combinations of $c$ and $g$ there is an
$18$-dimensional family of isomorphism classes of polarized $K3$ surfaces with 
smooth scroll $\T(c,D,\{ D _{\lambda} \})$. Moreover it follows that there 
will be an $18+\dim(\Aut({\PP}^g)) = 17 + (g+1)^2$-dimensional family of such 
projective models of $K3$ surfaces. This is true, simply because there is 
only a finite number of linear automorphisms that leaves a smooth polarized 
$K3$ surface invariant. 

For each value of $c$ (and $g$) one can pose several questions about the
set (or subscheme of the Hilbert scheme) of projective models $S'$ of
$K3$ surfaces $S$ with
elliptic free Clifford divisor $D$ and such that $\T$ is smooth. 

All scrolls of the same type are projectively equivalent, and hence the 
configuration of projective models of $K3$ surfaces in one such scroll is a 
projectively equivalent copy of that in another. Some questions one can pose, 
are:
How many scrolls are there of a given type? How many projective models  $S'$ 
are there within each scroll? In how many scrolls of a given type is a given
$S'$ included? 

The answer to the first question is well-known, the remaining
ones we will study more closely.

We start with the following well-known result from \cite{H}:

\begin{prop} \label{harr}
The dimension of the set of scrolls of type $\bf{e}$ and dimension $d$
in ${\PP}^g$ is
\[ \dim (\Aut ({\PP}^g)) - \dim (\Aut({\PP}(\E))) = (g+1)^2 - 3 - d^2 - \delta_1, \] 
where $\delta_1 := \sum _{i,j} \max (0, e_i - e_j - 1)$\index{$\delta_1$}.
\end{prop}

If $D^2=0$, we recall that $d=c+2$ for the scroll $\T(c,D,\{ D _{\lambda} \})$.

\subsection{Projective models with $c=1$}
 \label{c=1}
We have $g \geq 5$. Let the 
projective model $S'$ and the smooth scroll $\T=\T(1,D)$ be 
given.  As is seen from Proposition 
\ref{resolv} the surface $S'$ corresponds to the divisor class $3\H-(g-4)\F$ 
on $\T$. 
Moreover, part (c) of the proposition can
be applied so we can obtain a resolution in $\PP^g$. This is even 
minimal, by the comment in \cite[Example 3.6]{Sc}.

Assume $\T$ has scroll type  $(e_1,e_2,e_3)$, with $e_1 \geq e_2 \geq e_3$

\begin{prop} \label{SiT}
The (projective) dimension of the set of sections of divisor type 
$3\H-(g-4)\F$ in $\T$ is equal to $29 + \delta _2$, where $\delta _2 := 
\sum \max (0,g-5- \sum _{i=1}^3a_ie_i)$\index{$\delta _2$}. Here the first summation is taken 
over those triples $(a_1,a_2,a_3)$
such that $a_i \geq 0$, for $i=1,2,3$, and $\sum _{i=1}^3 a_i =3$. 
If $S'$ is smooth, a general section is a smooth projective model of a 
$K3$ surface.
\end{prop}

\begin{proof}
We use the formula
$h^0(\PP (\E), a\H+b\F) = h^0({\PP}^1, \Sym ^a(\E ) \* \O_{{\PP}^1}(b))$,
with $a=3$ and $b= g-4$. This gives $30+ \delta _2$. Being a smooth 
model of a $K3$ surface is an open condition on
the set of sections of $3\H-(g-4)\F$, and since one section, the one giving 
$S'$, is smooth, a general section of the linear system is so, too.
\hspace{0.07cm} $\square$ \end{proof}

We also have:

\begin{prop} \label{onescroll}
Each projective model $S'$ of a $K3$ surface $S$ of Clifford index 1 in 
${\PP}^g$ for $g \geq 5$, with $\T(c,D,\{ D _{\lambda} \})$ smooth,
is contained in only one smooth rational normal scroll of dimension 3. 
\end{prop}

\begin{proof}
In \cite[Part (7.12)]{S-D} it is shown that the scroll $\T$ is the intersection
of all quadric hypersurfaces containing $S'$. Moreover, any other smooth 
scroll $\T_2$
containing $S'$ is an intersection of quadric hypersurfaces 
(\cite[Prop. 1.5(ii)]{S-D}), each of course 
containing $S'$. Hence $\T_2$ contains $\T$. Since the two scrolls have the 
same dimension and degree, they must be equal.
\hspace{0.07cm} $\square$ \end{proof}

From this we conclude:

\begin{cor} \label{totnuscroll}
Let a scroll type $(e_1,e_2,e_3)$ be given, and let  $\delta _1$ and 
$\delta _2$ be defined as above. Then $\delta _1 \geq \delta _2$, and
there is a set of dimension 
\[(g+1)^2 + 17+ \delta _2 - \delta _1= \dim (\Aut(\PP^g))+18 +\delta _2 - 
\delta _1,\] parametrizing projective models
of $K3$ surfaces in smooth scrolls of the given type. 
\end{cor}

\begin{proof}
From Proposition \ref{harr} we see that there is a 
$((g+1)^2 -12 - \delta _1)$-dimensional set of scrolls of the same type 
as $\T$ in  ${\PP}^g$.
We know that each $S'$ in each such scroll is contained in only one scroll.
In each scroll there is a $(29+\delta _2)$-dimensional set of projective models of $K3$ surfaces as
described. We have $\delta _1 \geq \delta _2$, since otherwise there would be 
too many models with Clifford index 1.
\hspace{0.07cm} $\square$ \end{proof} 

\begin{rem} \label{comparison}
{\rm For $c=1$ we have }
\index{$\delta _1$}\[
\delta _1=\max (0,e_1-e_2-1)+\max (0,e_1-e_3-1)+\max (0,e_2-e_3-1),
\]
{\rm and}
\index{$\delta _2$}\begin{eqnarray*}
\delta _2 & = & \max (0,e_1-e_2-3)+\max (0,e_1-e_3-3)+\max (0,e_2-e_3-3)+ \\
 & & \max (0,e_1+e_2-2e_3-3)+\max (0,e_1-2e_2+e_3-3). 
\end{eqnarray*}
\end{rem}

Moreover $\delta _1 =0$ if 
and only if the scroll type is maximally balanced, and $\delta_2=0$ if the 
scroll type is ``reasonably well balanced''. It is clear that $\delta _1=0$
implies $\delta _2=0$. We also see that if $5 \leq g \leq 8$, 
then $\delta _2=0$. Hence the cases $(e_1,e_2,e_3)=(3,1,1)$ or $(3,2,1)$ are
cases where the number  $(g+1)^2 + 17+ \delta _2 - \delta _1$ is strictly
less than  $(g+1)^2 + 17= \dim (\Aut ({\PP}^g))+18$, and it is clear 
that scrolls of these types cannot represent the general projective models 
with $c=1$ and fixed $g$ since by 
the construction as in Proposition \ref{exthm} we get an 18-dimensional family
of such models.

The inequality $\delta _1 \geq \delta _2$ does not follow directly from the
formulas in Remark \ref{comparison}, for example since $\delta _1 < 
\delta _2$, for scroll types $(g-4,1,1)$, when $g \geq 11$.
This enables us to conclude that these and other scroll types with
$\delta _1 <  \delta _2$ do not occur for the scrolls formed from
projective models of $K3$ surfaces as described. This statement will be
strengthened to apply for $g \geq 8$ below. On the other hand the 
mentioned type $(3,1,1)$ does 
occur for $g=7$. This can be seen by using Lemma \ref{roll}.

For $g=7$ and type $(3,1,1)$ one then gets a polynomial $P$ of the form:
\begin{eqnarray*}
P_{6,1}(t,u)Z_1^3 + P_{4,1}(t,u)Z_1^2Z_2 + P_{4,2}Z_1^2Z_3 + 
P_{2,1}(t,u)Z_1Z_2^2 + \\ 
P_{2,2}(t,u)Z_1Z_2Z_3 + 
P_{2,3}(t,u)Z_1Z_3^2 + c_1Z_2^3 + c_2Z_2^2Z_3 + c_3Z_2Z_3^2 + c_4Z_3^3,
\end{eqnarray*} 
where the $P_{i,j}$ are homogeneous of degree $i$, and the $c_k$
are constants. For any fixed $(t,u)$ and any fixed point in the $\PP^2$
thus obtained, we see that we can avoid that point lying on the zero
scheme of $P$ by choosing the $P_{i,j}$ and $c_k$ properly, so we conclude 
that the linear system $|3\H-(g-4)\F| = |3\H-3\F |$ is base point free, and 
hence its general section is smooth, by Bertini. 
Irreducibility also follows by a similar argument.

Using Lemma \ref{roll}, we see that for $g=8$
any section of  $3\H-(g-4)\F=3\H-4\F$ on a scroll of type $(3,2,1)$ can be 
identified with a $P$ of the form
\begin{eqnarray*}
P_{5,1}(t,u)Z_1^3 + P_{4,1}(t,u)Z_1^2Z_2 + P_{3,1}(t,u)Z_1^2Z_3 + 
P_{3,2}(t,u)Z_1Z_2^2 + \\
P_{2,1}(t,u)Z_1Z_2Z_3 + P_{1,1}(t,u)Z_1Z_3^2 + P_{2,2}(t,u)Z_2^3 + 
P_{1,2}(t,u)Z_2^2Z_3 + c_1Z_2Z_3^2.
\end{eqnarray*} 
So here there is no $Z_3^3$-term, and from that we conclude that any section
of $3\H-4\F$ on a scroll of type $(3,2,1)$ must have the directrix line, say 
$\Gamma$, corresponding to $e_3=1$ as a part of its zero scheme. The fact that 
$\delta _2 < \delta _1$ indicates that if we can form smooth projective models of $K3$ surfaces 
this way, the surface must have a Picard lattice of higher rank than two. 
We may check this. If $L$, $E$ and $\Gamma$ sit inside a lattice of rank 2, then
we can write $\Gamma=aL+bE$, for rational numbers $a,b$.
In addition we must have $\Gamma ^2=-2$ and $L.\Gamma = E.\Gamma =1$.
It is easy to check that this is impossible.

It is clear that the set of base points of the linear system $3\H-4\F$ is 
just the directrix line. This is true since for a fixed value of $(t,u)$ (each 
fixed fiber) and a point $Q$ outside $(0,0,1)$, we can just change the 
$P_{5,1}(t,u)Z_1^3$-term or the $P_{2,2}(t,u)Z_2^3$-term, if we want to 
avoid $Q$. 
If we choose $c_1 \neq 0$, then we obtain that the zero scheme of the 
corresponding section intersects all fibers in curves, smooth at $(0,0,1)$.
(Here we set $Z_3=1$ in order to write the equation of the curve in affine
coordinates around $(0,0,1)$. The existence of the non-zero linear term 
$c_1Z_2$ gives smoothness at this point.)
Hence the zero scheme of a general section is smooth on all of $\T$. We have
basically used the identities $3e_2 \geq g-4$ and $e_2+2e_3=g-4$, to conclude
as we do.  See Remark \ref{noZ1} for references to other authors who have 
already used this kind of reasoning.

Using Lemma \ref{roll} again we see that if $g \geq 8$, then any section of
$3\H -(g-4)\F$ on a scroll of type $(g-4,1,1)$ corresponds to a polynomial
$P$ with $Z_1$ as a factor, which means that its zero scheme must be 
reducible as a sum of sections $\H -(g-4)\F$ (a subscroll of type $(1,1)$) 
and $2\H$. Hence these scroll types cannot occur for $\T(c,D,\{ D _{\lambda} \})$.
In a similar way one can draw conclusions about sections on other scroll
types. The observation above also has an interesting consequence for the types
of singular scrolls $\T(c,D,\{ D _{\lambda} \})$ arising for the case $c=1$, $D^2=0$, $g \geq 5$.

\begin{cor} \label{no00}
The type of $\T(c,D,\{ D _{\lambda} \})$ is never $(g-2,0,0)$ for $c=1$ and $D^2=0$.
\end{cor}

\begin{proof}
Assume the type of $\T(c,D,\{ D _{\lambda} \})$ is $(g-2,0,0)$. Then the type of the associated
scroll $\T_0$ is $(g-1,1,1)=(g_0-4,1,1)$, and the divisor type of $S''$
in $\T_0$ is $3\H_0 -(g_0-4)\F$ (see Example \ref{zeroint} and the proof of 
Corollary \ref{c1zeroint}). But we just observed that this is impossible.
\hspace{0.07cm} $\square$ \end{proof}

In general we conclude in the same way:

\begin{prop}
If a type $(a,b,c)$ is impossible for a smooth scroll \linebreak $\T(c,D,\{ D _{\lambda} \})$ in $\PP^g$ 
with $c=1$ and $D^2=0$, then the type $(a-1,b-1,c-1)$ is impossible for any scroll
$\T(c,D,\{ D _{\lambda} \})$ in $\PP^{g-3}$ with $c=1$ and $D^2=0$.
\end{prop}

We will make a list including all possible scroll types for smooth \linebreak $\T(c,D,\{ D _{\lambda} \})$, 
for $g \leq 13$, with $c=1$ and $D^2=0$. By the previous lemma, this will give a list including all scroll types of $\T(c,D,\{ D _{\lambda} \})$, smooth or singular, for
$g \leq 10$ and $c=1$ and $D^2=0$. In the column with headline ``$\#$ mod.'' we give
the value of $18-\delta _1+\delta_ 2$.

The information in the list is essentially contained in \cite{Re} and 
\cite[p.8-10]{St}. We include it for completeness, and for the benefit of
the reader we also include, in Remark \ref{noZ1} below, a few words about how 
the information can be obtained.

\vspace{.5cm}
\begin{tabular}{|c|c|c||c|c|c||c|c|c|} \hline
$g$ & scroll type & $\#$ mod. & $g$ & scroll type &  $\#$ mod. & $g$ & scroll type & $\#$ mod.                                             \\ \hline
$5$  & $(1,1,1)$  &   $18$  & $9$ & $(3,2,2)$ & $18$ & 
     $12$ & $(5,3,2)$  &   $16$  \\  \hline
$6$  & $(2,1,1)$  &   $18$  & $10$ & $(5,2,1)$  &  $16$  
      & $12$ & $(4,4,2)$  &   $17$  \\  \hline
$7$  & $(3,1,1)$  &   $16$ & $10$ & $(4,3,1)$  &   $17$  
      & $12$ & $(4,3,3)$  &   $18$  \\  \hline
$7$  & $(2,2,1)$  &   $18$ & $10$ & $(4,2,2)$  &   $16$
      & $13$ & $(7,3,1)$  &   $18$  \\  \hline
$8$  & $(3,2,1)$  &   $17$ & $10$ & $(3,3,2)$  &   $18$ 
      & $13$ & $(6,3,2)$  &   $16$  \\  \hline
$8$  & $(2,2,2)$  &   $18$ & $11$ & $(5,3,1)$  &   $17$  
      & $13$ & $(5,4,2)$  &   $17$  \\  \hline
$9$  & $(4,2,1)$  &   $16$ & $11$ & $(4,3,2)$  &   $17$  
      & $13$ & $(5,3,3)$  &   $16$  \\  \hline
$9$  & $(3,3,1)$  &   $17$ & $11$ & $(3,3,3)$  &   $18$ 
      & $13$ & $(4,4,3)$  &   $18$  \\  \hline
 & & & $12$ & $(6,3,1)$  &   $17$ & & &   \\  \hline
\end{tabular}

\vspace{.5cm}

This gives the following possibilities for singular types:

\vspace{.5cm}

\begin{tabular}{|c|c|} \hline
$g$   & singular scroll types    \\ \hline
$5$   & $(2,1,0)$   \\ \hline
$6$   & $(3,1,0)$, $(2,2,0)$    \\ \hline
$7$   & $(4,1,0)$, $(3,2,0)$   \\ \hline
$8$   & $(4,2,0)$   \\ \hline
$9$   & $(5,2,0)$    \\ \hline
$10$  & $(6,2,0)$   \\ \hline
\end{tabular}
\vspace{.5cm}

The dimensions of the families on the singular scrolls of type $(e_1,e_2,e_3)$
in $\PP ^g$ are equal to those of type $(e_1+1,e_2+1,e_3+1)$  on the 
corresponding smooth scrolls in $\PP ^{g+3}$.

\begin{rem} \label{noZ1}
{\rm 
Among the smooth scroll types listed above, we may immediately conclude that
a general section of $3\H-(g-4)\F$ is smooth, and hence a smooth
projective model of a $K3$ surface, for
the types
\begin{eqnarray*}
(1,1,1), (2,1,1), (3,1,1), (2,2,1), (2,2,2), (3,2,2), \\
(4,2,2), (3,3,2), (3,3,3), (4,3,3), (5,3,3), (4,4,3).
\end{eqnarray*} 
These are the ones with $3e_3 \geq g-4$.
The last inequality implies that the complete linear system $3\H-(g-4)\F$ has
no base points, and hence a Bertini argument gives smoothness of the general
section. The remaining types have the third directrix (of degree $e_3$)
as base locus. We have seen above that for the type $(3,2,1)$ the zero
scheme of the general section of $3\H-(g-4)F$ is smooth, since 
$3e_2 \geq g-4$ and $e_2+2e_3=g-4$. The same identities hold for the types
$(4,3,2)$, $(4,4,2)$ and $(5,4,2)$ also, so the zero scheme of a general section
of $3\H-(g-4)\F$ is smooth for these types too.

A similar argument can be made for the types $(3,3,1)$, $(4,3,1)$, $(5,3,1)$,
$(6,3,1)$ and $(7,3,1)$. Here the identity $3e_3 < g-4$ gives that the third
directrix curve (a line) consists of base points for the linear system.
The identity $3e_2 \geq g-4$ gives that there are no other base points.
The identity $e_1+2e_3=g-4$ gives that in each fiber the curve that
arises as the intersection of that fiber and the zero scheme of a section of 
$3\H-(g-4)\F$ is smooth at $(0,0,1)$, provided we choose the section
with a non-zero $cZ_1Z_3^2$-term. The total zero scheme is also smooth then.

The remaining possible smooth scroll types for $g \leq 13$ on the list above
are different, in the sense that for a general section of $3\H-(g-4)\F$,
the zero scheme of the section is singular  at finitely many points.
It turns out that for these types, which are $(4,2,1)$, $(5,2,1)$, $(5,3,2)$
and $(6,3,2)$, the general zero schemes are singular at exactly one point each.

The reason is the following: Since $3e_3 < g-4$, the third directrix curve
consists of base points for the linear system. Since $3e_2 \geq g-4$ for these
types, there are no other base points. Since $e_2+2e_3 < g-4$, there is no
$Z_2Z_3^2$-term for any section. Since $e_1+2e_3 > g-4$, in fact $e_1+2e_3=
(g-4)+1$ for all these types, there is no $cZ_1Z_3^2$-term with $c$ a constant,
but there is an $L(t,u)Z_1Z_3^2$-term with $L(t,u)$ a linear expression in
$t$ and $u$. If the section is chosen general, $L(t,u)$ is not identically
equal to zero. For all fixed $(t,u)$ where $L(t,u) \neq 0$, the zero scheme
of the section of $3\H-(g-4)\F$ is then smooth. For the single zero of
$L(t,u)$, the zero scheme is however singular. 

A comment about the types not appearing on the list above:
The smooth scroll types $(4,4,1), (5,4,1), (6,4,1), (5,5,1)$ are 
eliminated the following way: A section of $3\H-(g-4)\F$ can have no term 
containing a $Z_3^2Z_i$-term for these scroll types. Hence all plane cubics 
in the pencil 
are singular where they meet the linear directrix. But this is a contradiction,
since the general element in the pencil $|D|$ is smooth. Here we have used
the identity $e_1+2e_3 < g-4$ for these types.
The other types are eliminated because $Z_1$ must be a factor in each relevant
section, a contradiction. These are the types with $3e_2 < g-4$.

The necessary and sufficient condition ``$e_1+2e_3 < g-4$ or $3e_2 < g-4$''
for eliminating scroll types is given in \cite{Re}, as quoted in 
\cite[Lemma 1.8.]{St}. In \cite[p.9]{St} one also describes on which scroll
types a general section of $3\H-(g-4)\F$ is singular in a finite number of
points.

In Section \ref{Conc} we will show that all the types listed above for
$g \leq 10$ actually occur.
}  
\end{rem}

\begin{rem} \label{nores}
{\rm 
One does not have to use the resolution from Proposition
\ref{resolv} to see that
a projective model $S'$ of a $K3$ surface of Clifford index one in a smooth scroll $\T$ as above must 
be of divisor type $3\H-(g-4)\F$ in $\T$.
Define the vector space $W=H^0(\I _{S'}(3))/ H^0(\I _{\T}(3))$. In a natural
way $W$ represents the space of cubic functions on $\T$ that vanish on
$S'$. 

Study the exact sequences:
\[ 0 \hpil H^0(\I _{S'}(3))  \hpil  H^0(\O _{\PP^g}(3)) \hpil  H^0(\O _{S'}
(3)) \hpil H^1(\I_ {S'}(3)) \hpil 0 \] and 
\[ 0 \hpil H^0(\I _{\T}(3))  \hpil  H^0(\O _{\PP^g}(3)) \hpil  H^0(\O _{\T}
(3))  \hpil H^1(\I_ {\T}(3)) \hpil 0. \]
One obtains  $\dim W= h^0(\I _{S'}(3)) - h^0(\I _{\T}(3))=g-3$,
since $h^0(\O _{\T}(3)) - h^0(\O _{S'}(3))=g-3$ and $h^1(\I _{S'}(3)) =
h^1(\I _{\T}(3))=0$ (see \cite[Prop 1.5(i)]{S-D}). Take $g-3$ arbitrary fibers $F$ of the 
ruling on $\T$, that is $g-3$ planes. For each plane it is one linear condition
on the elements in $W$ to contain it (since this is equivalent to contain an
point in the plane outside $S'$). Hence containing all the $g-3$ planes 
imposes $g-3$ conditions. These conditions must be independent, since 
otherwise there would be a cubic in  ${\PP}^g$, not containing $\T$, and 
containing the union of $S'$ and $g-3$ planes. This union has degree 
$(2g-2)+(g-3)=3g-5$. But by 
Bezout's theorem the cubic and $\T$ intersect in a surface of degree $3g-6$.
Hence, in particular, any choice of $g-4$ planes gives independent conditions,
and there is one, and only one, hypercubic (modulo the ideal of $\T$), which
contains $S'$ and $g-4$ planes in the pencil. By Bezout's theorem, it does not
contain more. Hence $S'$ in a natural way is a section of $3\H-(g-4)\F$.}
\end{rem}

\subsection{Projective models with $c=2$}
\label{c=2}
Let $\T=\T(2,D)$ with $D^2=0$. We have $g \geq 7$. 
Denote the type of $\T$ by $(e_1,e_2,e_3,e_4)$.
Proposition \ref{resolv}(a) (or Corollary \ref{c1zeroint} if $\T$ is not 
smooth) gives that $S'$ is a complete intersection in $\T$ of two divisors of 
type $2\H-b_1\F$ and $2\H-b_2\F$. 
By convention, we set $b_1 \geq b_2$. Part (d) of the same proposition gives 
the well-known fact that $b_1 + b_2=g-5$.
Such a situation has been thoroughly investigated in \cite{B}. 
As already mentioned in the proof of Corollary \ref{c3}, it follows
from \cite[Thm. 5.1]{B} that $S'$ is 
singular along a curve if $b_1 > \frac{2(e_1+e_2+e_3)}{3}$, or
equivalently  $b_2 < \frac{g-9+2e_4}{3}$. 
Hence $b_2 \geq 0$ for $g \geq 7$, for complete intersections with only
finitely many singularities. Since in particular $b_2 \geq -1$, it is clear 
that part (b) and (c) of Proposition \ref{resolv} can be used to give a resolution of $S'$ in ${\PP}^g$. This resolution is minimal because of the comment in 
\cite[Example 3.6]{Sc}. The fact that $b_2 \geq 0$ means that $S'$ can be 
viewed geometrically as a complete intersection of one hyperquadric containing 
$b_1$ three-planes in the pencil, and another containing $b_2$ three-planes
(throwing away the three-planes and taking the closure of what remains).

Let us study projective models in smooth scrolls for $c=2$ and  $D^2=0$ in general.
We see from Proposition \ref{harr} that the set parametrizing the 
scrolls having the 
same type as ${\T}$ has dimension $(g+1)^2-19 - \delta _1  = 
\dim (\Aut ({\PP^g}))-18- \delta _1$, where $\delta _1 := \sum _{i,j} 
\max (0, e_i - e_j - 1).$\index{$\delta_1$}
Therefore one expects the set of projective models of smooth $K3$ surfaces 
in each scroll to have dimension 36 if the scroll type is reasonably well 
balanced, to get a set of total dimension $\dim (\Aut({\PP^g}))+18$.
This ``expectation'' is based on the natural assumption that a set of
total dimension $\dim (\Aut ({\PP^g}))+18$ arises from $S'$ that sit inside 
maximally balanced scrolls. 
In this case there are two different sources of imbalance; that of 
the $e_i$, and that of the $b_k$. We will look more closely at this.

Set $\delta _2:=\max (0,b_1-b_2-1)$\index{$\delta_2$}. Assume first $b_1 > b_2$.  We calculate 
\[ \dim |\O_{\T}(2\H-b_1\F)| = 5g-6-10b_1+\delta _3, \]
where $\delta _3:=h^1(\Sym ^2 \E \* \O_{\PP^1} (-b_1))=0$\index{$\delta_3$} if and only if $e_4 \geq \frac{b_1-1}{2}$.

By the same sort of calculation we of 
course get 
\[ \dim |\O_{\T}(2\H-b_2\F)|= 5g-6-10b_2+\delta _4, \] 
where  $\delta_ 4 := h^1(\Sym ^2 \E \* \O_{\PP^1} (-b_2))=0$ \index{$\delta_4$}if and only if 
$e_4 \geq \frac{b_2-1}{2}$.

Now fix a zero scheme $Y$ of a section $s$ of $2\H-b_1\F$, and study the 
exact sequence
\[ 0 \hpil \O_{\T}((b_1-b_2)\F) \hpil  \O_{\T}(2\H-b_2\F) \hpil
 \O_Y(2\H-b_2\F) \hpil 0. \] 
This induces a sequence
\begin{eqnarray*}
0  & \hpil H^0(\O_{\T}((b_1-b_2)\F)) & \hpil  H^0(\O_{\T}(2\H-b_2\F))   \\ 
   & \hpil H^0(\O_Y(2\H-b_2\F))      & \hpil  H^1(\O_{\T}((b_1-b_2)\F)).
\end{eqnarray*}
Since $h^0(\O_{\T}((b_1-b_2)\F))= b_1-b_2+1$ and 
$h^1(\O_{\T}((b_1-b_2)\F))=0$, we obtain 
\begin{eqnarray*}
\dim |\O_{Y}(2\H-b_2\F)| & = & \dim |\O_{\T}(2\H-b_2\F)| -(b_1-b_2+1) \\
 & = & 5g-6-10b_2+\delta _4-(b_1-b_2+1) \\
 & = & 5g-8-10b_2+\delta _4-\delta _2.
\end{eqnarray*}
Summing up, we obtain
\begin{eqnarray*}
 & & \dim \O_{\T}(2\H-b_1\F) + \dim \O_{Y}(2\H-b_2\F) \\  
& = & 5g-6-10b_1 +\delta _3+ 5g-8-10b_2 + \delta _4-\delta _2 \\
& = & 10g-14-10(b_1+b_2) -\delta _2+ \delta _3+  \delta _4 \\ 
& = & 36-\delta_2+\delta _3+\delta _4.
\end{eqnarray*}
In particular, if $g$ is even, $b_1-b_2=1$ and $e_4 \geq \frac{b_1-1}{2}=
\frac{g-6}{4}$, we get 36.
We remark that $\frac{g-6}{4}$ is just $\frac{3}{4}$ less than the average
values of the $e_i$. In general there is thus a $(36-\delta _2+\delta _3+
\delta _4)$-dimensional set of complete intersections of type 
$(2\H-b_1\F,2\H-b_2\F)$, provided the section $s$ is uniquely determined. 
The latter fact follows, for example by the same kind of argument as in 
\cite[(6.2)]{Sc}, 
, where scrolls arising from tetragonal curves are treated (see also the proof of Proposition \ref{b1} below). 

We now assume $b_1=b_2(=b= \frac{g-5}{2})$, which can only occur if 
$g$ is odd. Then $h^0(2\H-b\F)=h^0({\PP^1},\Sym ^2(\E)\* \O_{\PP^1}(b))$.
This number is of the form $20+\delta _3=20+\delta _4$, where $\delta _3$
and $\delta _4$ are defined as in the case above. We see 
that  $\delta _3 = \delta _4=0$ if and only if $2e_4-b \geq -1$, or equivalently
$e_4 \geq \frac{g-7}{4}$. The average value of the $e_i$ is $\frac{g-3}{4}$,
which is just one more.
The set of complete intersections corresponds to an open set in the
Grassmannian $G(2,h^0(2\H-b\F))$, which has dimension $36+2\delta _3$.
Hence we get the expected number if $b_1=b_2$, and the $e_i$ are 
well balanced. Whether $b_1=b_2$ or not, we have now proved:
Let a scroll type $e_1,e_2,e_3,e_4)$ be given, and let $\delta _2,
\delta _3, \delta _4$ be defined as above in this section, and let
$\delta _1$ be defined as in Proposition  \ref{harr}

\begin{prop} \label{ineachscroll}
The set of complete intersections of type $(2\H-b_1\F,2\H-(g-5-b_1)\F)$
on a smooth rational normal scroll of dimension $4$ in 
$\PP^g$ of type $(e_1,e_2,e_3,e_4)$ is either empty or of dimension 
$36-\delta _2 +\delta _3+\delta _4$. 
\end{prop}

\begin{cor} \label{upfrom36}
The set of complete intersections of type $(2\H-b_1\F,2\H-(g-5-b_1)\F)$
with no or finitely many singularities on a smooth rational normal scroll of 
dimension $4$ in $\PP^g$ of 
 type $(e_1,e_2,e_3,e_4)$, is either empty or of 
dimension at least $36$ if $\delta _1 \geq 1$.
\end{cor}
\begin{proof}
Set $s=b_1-\frac{g-5}{2}$. If $s=0$, then $\delta_ 2=0$, and there is nothing
to prove. If $s \geq \frac{1}{2}$, then $\delta _2 =2s-1$. 
We split into $4$ subcases: $e_1+e_2+e_3+e_4= 4e+h$, where $h=0,1,2,3$.
If $h=0,1$, we have if $\delta _1 \geq 1$:
$2e_4 \leq  2(\lfloor \frac {g-3}{4} \rfloor -1) \leq \frac{g-7}{2}$, and
$e_3 + e_4 \leq \frac{g-3}{2} -1=\frac{g-5}{2}$.
This gives $b_1-2e_4-1 \geq \frac{g-5}{2}+s-\frac{g-7}{2}-1=s$, and
$b_1-(e_3+e_4)-1 \geq \frac{g-5}{2}+s-\frac{g-1}{2}-1=s-1$.
Hence 
\[\delta _3=h^1(\PP^1,\Sym ^2(\E)\* \O_{\PP^1}(b_1)) 
\geq s+(s-1)=2s-1=\delta _2,\]
and then the result follows from Proposition \ref{ineachscroll}.

If $h=3$, essentially the same method works (look at the three terms
$b_1-2e_3-1$, $b_1-e_3-e_4-1$ and $b_1-2e_4-1$).

If $h=2$, then essentially the same method works (look at the six terms
of the form $b_1-e_i-e_j-1$, for $i,j=2,3,4$), except in the case $s=1$,
and $(e_1,e_2,e_3,e_4)=(e+2,e,e,e)$. But in that case 
$b_1=\frac{g-1}{2}+1=2e+1$.
Using Lemma \ref{roll} we see that that $Z_1$ is a factor in every section of
$2\H-b_1\F$, so this case simply does not occur. See also the appendix of
\cite{B}.
\hspace{0.07cm} $\square$ \end{proof} 

If we add the assumption that there exists at least one smooth model $S'$ 
of a given intersection type, giving rise to a scroll of a given type, then
we can conclude that there is a set of dimension $(36-\delta _2+\delta _3+
\delta _4)$ parametrizing smooth projective models $S'$ in a scroll of the 
given type. We see that if the scroll type or $(b_1,b_2)$-type
is unbalanced, the dimension of the set of complete intersections (smooth or singular) can a priori go up, or it can go down from its ``expected'' value 36.
From the last corollary, however, we see that the dimension goes down only if 
the scroll type is maximally balanced, and the intersection type is not. 
 
If we just start with an arbitrary scroll type $(e_1, e_2, e_3, e_4)$, with
$e_4 \geq 1$, $\sum_{i-1}^{c+2}e_i=g-3$ and ``intersection type'' 
$b_1 \geq b_2 \geq 0$, with $b_1+b_2=g-5$, it is an intricate question to 
decide whether there are any that are smooth projective models of $K3$ surface,
or any that are not singular along a curve. This 
problem is studied in detail in \cite{B}, and we will study the cases for low
$g$ in Section \ref{stypes}.

The question whether a projective model $S'$ can be included in several 
scrolls of the same type simultaneously is not as simple to answer as in the 
case $c=1$.
In the case $c\geq 2$ the scroll $\T$ is no longer the intersection of the
hyperquadrics containing $S'$; in fact $S'$ itself is that intersection 
\cite{S-D}. The question is essentially how many divisor classes of
elliptic curves $E$ with $E.L=c+2$ there are on $S$.

On the other hand it is clear that a projective model cannot be contained
in a positive dimensional set of scrolls of the type in question. There is
a discrete family of divisor classes on $S$, so the dimension of the family
of smooth projective models of $K3$ surfaces on scrolls of the type in
question is now, as for $c=1$, equal to the sum of the dimension of the set 
of scrolls of a given type and the dimension of the set of smooth projective 
models of that type. This sum is 
\begin{eqnarray*}
(g+1)^2-19-\delta _1 + (36-\delta _2+\delta _3 +\delta _4) \\
= \dim (\Aut (\PP^g)) +18 - \delta _1 -\delta _2+ \delta _3 +\delta _4.
\end{eqnarray*}

To obtain the dimension of the set of projective equivalence classes,
since only a finite number of automorphisms of $\PP^g$ fixes a $K3$ surface,
we subtract the number $\dim (\Aut (\PP^g))$ and get
\[18 - \delta _1 -\delta _2+ \delta _3 +\delta _4.\]
By Theorem \ref{exthm} this number is equal to $18$ for at least one scroll 
type, where there exists smooth complete intersections of that type, and where
the fibers of the complete intersection represent a free Clifford divisor
on $S$ with $c=2$. It is clear that if we choose the most balanced scroll
type for a fixed $g$, then $\delta _1=0$. If in addition we choose the
most balanced $(b_1,b_2)$-type, then $\delta _2=\delta _3=\delta _4=0$.
 Using Lemma \ref{roll} it is also
easy to prove that for all $g$ the most balanced scroll and intersection type
(the unique combination for fixed $g$ with $\delta _1= \delta _2=0)$ then
a general complete intersection will be a smooth projective model of a $K3$ surface.

\subsection{An interpretation of $b_1$ and $b_2$} 
\label{theb}

We will briefly study the case $c=2$ in an
analogous manner as the case $c=1$ was studied in Remark \ref{nores} when we 
showed that a projective model of a $K3$ surface of 
Clifford index one with smooth asssociated scroll $\T$ must be of divisor type $3\H+(g-4)\F$ in ${\T}$ without 
using the resolution from Proposition \ref{resolv}. 

Define the vector space $W=H^0(\I _{S'}(2))/ H^0(\I _{\T}(2))$. In a natural
way $W$ represents the space of quadric functions on $\T$ that vanish on
$S'$. 

As in Remark \ref{nores} one obtains $\dim W= h^0(\I _{S'}(2)) - h^0(\I _{\T}(2))=g-3$,
since $h^0(\O _{\T}(2)) - h^0(\O _{S'}(2))=g-3$ and $h^1(\I _{S'}(2)) =
h^1(\I _{\T}(2))=0$ (see \cite[Prop. 1.5(i) and Theorem 6.1(ii)]{S-D}. 
Assume $g$ is odd. Take $b=\frac{g-5}{2}$ arbitrary fibers $F$ of the 
ruling on $\T$, that is three-planes. For each three-plane we have two 
independent linear conditions on the quadric hypersurfaces to contain
it. (First, take one point in the three-plane, not on $S'$. There is only one
quadric surface in the threespace containing this point and the intersection 
with $S'$. Then take another point outside this quadric surface. To contain 
these two points and $S'$ is equivalent to containing the threespace and $S'$.)
So, one naively expects there to be $2b=g-5$ conditions to contain all 
the $b$ three-planes. Hence there should be a pencil, and only a pencil, of 
elements of $W$ doing so. Intersecting the elements of the pencil, one would 
expect to get the projective model of the $K3$ surface. If it really were so 
simple, however, all intersection types $(b_1,b_2)$ would be completely balanced. This is not always true, and one reason is that two different elements of 
$W$ may intersect $\T$ in a common threedimensional component dominating 
${\PP^1}$ in the fibration on $\T$.

We are therefore not able to imitate the reasoning of Remark \ref{nores}, and 
thereby establish the fact that the ideal of $S'$ in $\T$ is generated as it 
is, without using Proposition \ref{resolv}.
On the other hand we may use the knowledge that 
we have from Proposition \ref{resolv}, that $S'$ is indeed of intersection type
$(2\H-b_1\F,2\H-b_2\F)$ in its scroll. Make no assumption on the parity of $g$.

\begin{prop} \label {b1}
The invariant $b_1$ is equal to the largest number $k$, such that there exists
a non-zero element $Q$ of $W$ (a hyperquadric in ${\PP^g}$ containing $S'$, 
but not ${\T}$) containing $k$ three-planes in the pencil. 

Moreover, 
$b_2$ is the largest number $m$, such that there exists a non-zero element 
of $W$ containing $m$ three-planes in the pencil, and intersecting $\T$ in a 
different 3-dimensional dominant component than $Q$ does. 
\end{prop}

\begin{proof}
Define two elements in $W$ to be congruent if they have the same 
dominating three-dimensional component (but possibly differ in which 
three-planes they 
contain). Define the index of an element of $W$ as the number of three-planes 
it contains (if necessary, counted with multiplicity). It is clear that if 
two elements of $W$ are congruent, then they have the same index (they 
correspond to a well defined divisor class $2\H-i\F$, where $i$ is the index). 
It follows from a Bezout argument that if the sum of the indices of 
two elements is larger than $2b=g-5$, then they must be congruent.
This shows both assertions. It also shows that the element in $2\H-b_1\F$ which
any element in $W$ with index larger than $b$ gives rise to, is the same.
(This element is nothing but the congruence class of the element in $W$.)
\hspace{0.07cm} $\square$ \end{proof}

\subsection{Possible scroll types for $c=2$} \label{stypes}
\label{type2}

Almost all the information in this subsection can also be found in \cite{B}
and \cite{St}, taken together, but we include it for completeness, and present
it in our own way, for the sake of the reader.

Also in the case $c=2$ it is possible to use Lemma \ref{roll} to obtain
useful conclusions for many concrete scroll types. First we will look at 
possible smooth scroll types for $g \leq 10$. 

For $g=7$ and $g=8$, 
the only possible smooth scroll types are $(1,1,1,1)$ and $(2,1,1,1)$,
respectively. For $g=7$, Lemma \ref{roll} immediately gives $b_1 \leq 2$.
Here $b_1+b_2=2$, so $b_2 \geq 0$. For $g=8$ we see that $Z_1$ is a factor
in every section of $2\H-b_1\F$, for $b_1 \geq 3$. Since a reducible section
which is the sum a section $\H-2\F$ and a section $\H-(b-2)\F$  would 
intersect another section of type
$2\H-b_2\F$ in something reducible, alternatively  since $S'$ is non-degenerate, we 
must have $b_1=2$ and $b_2=1$.

For $g=9$ we have  $b_1+b_2=4$ and two smooth scroll types $(2,2,1,1)$ and 
$(3,1,1,1)$. For the latter type Lemma \ref{roll} gives $b_1=b_2=2$,
since any section of $2\H-b\F$, with $b \geq 3$ must be a product of a section
$\H-b\F$ and a section of  $\H$. For the type $(2,2,1,1)$ any section of $2\H-3\F$ has a 
zero scheme containing the subscroll generated by the 
two linear directrices (a quadric surface $Q$).
If $b_1=3$, then $b_2=1$, and the section $2\H-b_2\F$ intersects $Q$ as a
curve of type $(1,2)$ on $Q$. This is a rational twisted cubic $\Gamma$.
This corresponds to the fact that $\delta _1=\delta _3= \delta _4=0$ and 
$\delta _2=1$, so that the set of complete intersections inside $\T$ has 
dimension at most $36- \delta _2+\delta _3+\delta _4=35$, and taking the union
over all $\T$ of the same type we get dimension at most $\dim (\Aut (\PP^g))+
17$.

Any section of $2\H-4\F$ has a zero scheme, which restricts to two lines in
each fiber of $\T$. This is impossible if this scheme shall contain
a (necessarily non-degenerate) model $S'$. We also have 
$h^0(2\H-b\F)=0$ for $b \geq 5$. Hence $b_1$ is $2$ or $3$.

For $g=10$, we have $b_1+b_2=5$ and a priori three possible scroll types
$(4,1,1,1)$, $(3,2,1,1)$ and $(2,2,2,1)$. But the first cannot occur, since
any section of $2\H-b_1\F$ must have total weight $-b_1 \leq -3$, and then
$Z_1$ must be a factor, using Lemma \ref{roll}.
This is impossible.
Likewise, if $\T$ has type $(3,2,1,1)$ and $b_1 \geq 5$, we conclude that
$Z_1$ must be a factor, again impossible. The cases $b_1=4$ and
$b_1=3$ are 
however possible.

If $\T$ has type $(2,2,2,1)$, then $b_1 \leq 4$, since $h^0(2\H-b\F)=0$ if 
$b \geq 5$. If $b_1=4$, then the zero scheme of any section of 
$2\H-b_1\F=2\H-4\F$
contains the linear directrix of $\T$ twice (its equation is a homogeneous 
quadric in $Z_1, Z_2, Z_3$ involving neither $Z_4, t$ nor $u$).
If we intersect with a section of $2\H-b_2\F=2\H-\F$, and interpret it as the
intersection with a quadric containing a fiber, and throw away the fiber,
the residual intersection with $\T$ must contain one point of its linear
directrix. This must then be a singular point of $S'$.
Hence only $b_1=3$, $b_2=2$ gives a smooth $S'$ for this scroll type.
For these invariants $\delta _i=0$, for $i=1,2,3,4$, so we get a family of
total dimension $\dim (\Aut (\PP^g))+18$.

So far we have studied smooth scroll types for $7 \leq g \leq 10$. At this 
point we could either proceed with smooth scroll types for $g \geq 11$, or look
at singular scroll types for $g \geq 7$. These topics are closely related.
Assume we have a singular scroll $\T(c,D,\{ D _{\lambda} \})$ for $g \geq 7$, $D^2=0$, $c=2$.
Then the associated smooth scroll $\T_0$ is contained in $\PP^{g+4}$, with
a resolution as in Proposition \ref{upperres}:
\begin{eqnarray*}
0 \hpil \O_{\T_0}(-4\H_0+(g-1)\F)  & \hpil & \+ _{k=1}^2 \O_{\T_0}(-2\H_0 + b_k \F)  \\
&  \hpil & \O_{\T_0} \hpil \O_{S''} \hpil 0.
\end{eqnarray*}
So, at this point we can use the method of rolling factors to check what
scroll types in $\PP^{g+4}$ that may contain a surface like $S''$.
Scroll types $(e_1+1, \ldots ,e_4+1)$ for $\T_0$ correspond to types 
$(e_1, \ldots ,e_4)$ for $\T$.

Let us study complete intersection surfaces in smooth scroll types for $g=11$ 
with this dual viewpoint.
Now $b_1+b_2=6$ and $\deg \T=8$ and there are a priori
$5$ different possible scroll types: 
\[(5,1,1,1), \hs (4,2,1,1), \hs (3,2,2,1), \hs (3,3,1,1) \hs \mbox{ and } \hs (2,2,2,2).\]
The type $(5,1,1,1)$ cannot occur, for the same
reason that $(4,1,1,1)$ cannot occur for $g=10$.

For $(4,2,1,1)$ and $(3,3,1,1)$ we can conclude
that $Z_1$ or $Z_2$ is a factor in every term of every section of $2\H-3\F$. 
This gives that the subscroll formed by the
two linear directrices is contained in $S'$ (or $S''$). This is clearly
impossible. Hence $b_1 \geq 4$ for these types.

If $b_1 \geq 5$, then $Z_1$ is a factor in every section of $2\H-b_1\F$,
for each of the types $(4,2,1,1), (3,2,2,1)$ and $(2,2,2,2)$, which gives a 
contradiction.   
For the type $(3,3,1,1)$ we argue as follows:
If $b_1 \geq 5$, then no term of the form $Z_iZ_3$ or $Z_iZ_4$ can occur as
factor in a monomial of a section of $2\H-b_1\F$, so for each fixed value
of $(t,u)$ we get a quadric in $Z_1, Z_2$ only. This defines a union of two
planes in each $\PP^3$ which is a fiber of $\T$ (or $\T_0$). Hence each 
fiber of $S'$ (or $S''$) is degenerate, a contradiction. So $b_1=4$.

We make the same kind of considerations for all $g \leq 14$.
We end up with the following a priori possible combinations of smooth scroll 
type and $b_1$, for $7 \leq g \leq 14$ (of course $b_2=g-5$). For each scroll 
type and intersection type $(b_1,b_2)$ we
indicate whether the general zero scheme of a complete intersection of
type $(2\H-b_1\F,2\H-b_2\F)$ is smooth or singular. See also 
Remark \ref{posstype}.  In the column with headline `` $\#$ mod.'' we give
the value of $18-\delta_ 1-\delta _2+\delta _3 +\delta _4$.
This table contains information that can also be found in \cite{B} and 
\cite{St}, taken together. In \cite{B} all possible scroll and intersection 
types for smooth scrolls for all $g \geq 7$ are listed, and in \cite{St} the
information on the moduli of the corresponding families are given. We include 
the list for $g \leq 14$, since it will be useful in the study of projective 
models on singular scrolls, and for the lattice-theoretical considerations in 
Section \ref{Conc}.

\vspace{.5cm}
\begin{tabular}{|c|c|c|c|c||c|c|c|c|c|} \hline
$g$ & scroll type & $b_1$ & comp. int. & $\#$ mod. & 
                       $g$ & scroll type & $b_1$ & comp. int. & $\#$ mod. \\ \hline
$7$  & $(1,1,1,1)$ & $1$ & Smooth  & $18$ & 
$12$ & $(4,3,1,1)$ & $5$ & Smooth  & $16$ \\ \hline
$7$  & $(1,1,1,1)$ & $2$ & Smooth  & $17$ & 
$12$ & $(4,2,2,1)$ & $4$ & Singular& $15$ \\ \hline
$8$  & $(2,1,1,1)$ & $2$ & Smooth  & $18$ & 
$13$ & $(3,3,2,2)$ & $4$ & Smooth  & $18$ \\ \hline
$9$  & $(2,2,1,1)$ & $2$ & Smooth  & $18$ & 
$13$ & $(3,3,2,2)$ & $5$ & Smooth  & $17$ \\ \hline
$9$  & $(2,2,1,1)$ & $3$ & Smooth  & $17$ & 
$13$ & $(3,3,3,1)$ & $4$ & Smooth  & $17$ \\ \hline
$9$  & $(3,1,1,1)$ & $2$ & Smooth  & $15$ & 
$13$ & $(3,3,3,1)$ & $6$ & Smooth  & $18$ \\ \hline
$10$ & $(2,2,2,1)$ & $3$ & Smooth  & $18$ & 
$13$ & $(4,2,2,2)$ & $4$ & Smooth  & $15$ \\ \hline
$10$ & $(2,2,2,1)$ & $4$ & Singular& $17$ & 
$13$ & $(4,3,2,1)$ & $4$ & Singular& $16$ \\ \hline
$10$ & $(3,2,1,1)$ & $3$ & Smooth  & $16$ & 
$13$ & $(4,3,2,1)$ & $5$ & Smooth  & $16$ \\ \hline
$10$ & $(3,2,1,1)$ & $4$ & Singular& $17$ & 
$13$ & $(4,3,2,1)$ & $6$ & Smooth  & $18$ \\ \hline
$11$ & $(2,2,2,2)$ & $3$ & Smooth  & $18$ & 
$13$ & $(5,3,1,1)$ & $6$ & Smooth  & $17$ \\ \hline
$11$ & $(2,2,2,2)$ & $4$ & Smooth  & $17$ & 
$14$ & $(3,3,3,2)$ & $5$ & Smooth  & $18$ \\ \hline
$11$ & $(3,2,2,1)$ & $3$ & Smooth  & $17$ & 
$14$ & $(3,3,3,2)$ & $6$ & Singular& $17$  \\ \hline
$11$ & $(3,2,2,1)$ & $4$ & Smooth  & $17$ & 
$14$ & $(4,3,2,2)$ & $5$ & Smooth  & $16$ \\ \hline
$11$ & $(4,2,1,1)$ & $4$ & Singular& $15$ & 
$14$ & $(4,3,2,2)$ & $6$ & Singular& $17$  \\ \hline
$11$ & $(3,3,1,1)$ & $4$ & Smooth  & $16$ & 
$14$ & $(4,3,3,1)$ & $5$ & Smooth  & $17$  \\ \hline
$12$ & $(3,2,2,2)$ & $4$ & Smooth  & $18$ & 
$14$ & $(4,4,2,1)$ & $5$ & Singular& $16$  \\ \hline
$12$ & $(3,3,2,1)$ & $4$ & Smooth  & $17$ & 
$14$ & $(5,3,2,1)$ & $5$ & Singular& $15$  \\ \hline
$12$ & $(3,3,2,1)$ & $5$ & Smooth  & $17$ & 
$14$ & $(5,3,2,1)$ & $6$ & Singular& $16$  \\ \hline
\end{tabular}

\vspace{.5cm}

For perfect Clifford divisors $D$ with $D^2=0$ and singular scrolls 
$\T=\T(2,D)$ we get 
the following list of a priori possible cases in $\PP^g$, for 
$7 \leq g \leq 10$ (subtracting $2$ from all values of $b_i$ for $S''$ in 
$\T_{0}$ in $\PP^{g+4}$, for $i=1,2$):

\vspace{.5cm}

\begin{tabular}{|c|c|c|c||c|c|c|c|} \hline \label{table2}
$g$ & sing. scroll type & $b_1$ & $(S'')_{virt}$ & 
                  $g$ & sing. scroll type & $b_1$ & $(S'')_{virt}$  \\ \hline
$7$ & $(2,1,1,0)$ & $1,2$ & Smooth  & 
$9$ & $(4,2,0,0)$ & $4$   & Smooth   \\ \hline
$7$ & $(2,2,0,0)$ & $2$   & Smooth  & 
$9$ & $(3,2,1,0)$ & $2,3$ & Smooth only for $b_1=3$  \\ \hline
$7$ & $(3,1,0,0)$ & $2$   & Singular& 
$9$ & $(2,2,2,0)$ & $2$   & Smooth   \\ \hline
$8$ & $(3,2,0,0)$ & $3$   & Smooth& 
$10$& $(4,2,1,0)$ & $3,4$ & Singular \\ \hline
$8$ & $(3,1,1,0)$ & $2$   & Singular& 
$10$& $(3,3,1,0)$ & $3$   & Singular \\ \hline
$8$ & $(2,2,1,0)$ & $2$   & Smooth  & 
$10$& $(3,2,2,0)$ & $3$   & Smooth \\ \hline
\end{tabular}

\vspace{.5cm}

\begin{rem} \label{posstype}
{\rm For each smooth scroll $\T$ and intersection type where the general 
element $S'$ is smooth (on the upper list of types for $7 \leq g \leq 14$) 
it is clear that we have a smooth projective model of a $K3$ surface. For the 
remaining cases (on that list) it is natural to interpret them as 
projective models $S'$ of $K3$ surfaces by non-ample linear systems.
The types on the upper list are the only ones for $g \leq 14$ where a general 
complete 
intersection $(2\H-b_1\F,2\H-b_2\F)$, with $b_1+b_2=g-5$, is either smooth, 
or singular in a finite number of points. 
In Section \ref{Conc}, moreover, we describe all projective 
models for low $g$, including those with $c=2$. All scroll types listed above 
(smooth as in the upper list or singular as in the lower list) for $g \leq 10$
reappear in the description in Section \ref{Conc}.

To decide which intersections that are in general smooth, which intersections 
that are in general singular in finitely many points, and which intersections
that are in general singular along a curve (or even reducible) one uses 
Lemma \ref{roll}, similarly as in \cite{B} and 
\cite{St}. In particular we have checked with the Appendix in \cite{B}, which 
gives a list of smooth complete 
intersection $K3$ surfaces in $4$-dimensional smooth rational normal scrolls 
and also a list of relevant intersections with only finitely many 
singularities.
 
On the lower list, concerning singular scrolls for $\leq 10$, we have listed all scroll and intersection types which might a priori appear as ``images'' by the map $i$ of scrolls $\T_0$ and surfaces $S''$ on the upper list, provided
that the Clifford divisor $D$ is perfect.

In the columns with heading $(S'')_{virt}$ we have indicated whether a
general complete intersection of type in question on $\T_0$, which contains
the exceptional divisor of the map $i$ from $\T_0$ to $\T$, is smooth or
singular. We call such a complete intersection $(S'')_{virt}$, since we do not
a priori know that it is an $S''$. If $D$ is perfect and the scroll $\T$ is singular, each occurring
projective model $S'$ on $\T$ is  of course also singular, but $S''$ smooth 
means that all singularities of $S'$ are due to contractions across the 
fibers; there are no contractions in the individual fibers.

We will illustrate that the issues whether a smooth scroll $\T(2,D)$ and
an associated intersection type for a model $S'$ appears on the upper list, 
is different from the issue whether the scroll and intersection type appears
for a $\T_0$ and a $S''$. If $D$ is perfect, it is a priori possible that all complete 
intersections, or a general one, represents a model $S'$, but not an $S''$.
As an example, look at the type $(3,3,3,1)$, with $(b_1,b_2)=(6,2)$.
Then $S'$ consists of the common zeroes of two sections of the form
\[
 c_1Z_1^2 + c_2Z_1Z_2 + c_3Z_1Z_3+ c_4Z_2^2+ c_5Z_2Z_3+c_6Z_3^2 
\]
and
\[
f(t,u,Z_1,Z_2,Z_3,Z_4) + c_7Z_4^2,  
\]
where $f(t,u,Z_1,Z_2,Z_3,Z_4)$ is contained in the ideal generated by
$Z_1,Z_2,Z_3$. The  general such intersection is smooth, and does not 
intersect the last directrix ($Z_1=Z_2=Z_3=0$) at all. But in order to be a
surface of the form $S''$, associated to a perfect Clifford divisor $D$, the 
intersection must contain the last directrix.
This forces $c_7$ to be zero. In that case the intersection is no longer
smooth, in fact it contains the directrix in its singular locus, and hence
it cannot be an $S''$. Hence scroll type $(3,3,3,1)$ with $b_1=6$ appears on
the upper list, but the corresponding ``pushed down'' type $(2,2,2,0)$ does
not appear in combination with (the revised) $b_1=4$.

A similar, but slightly different case, is the scroll type $(3,2,2,1)$ and
intersection type $(4,2)$.
Then a surface of the form $S'$ would consist of the common zeroes of two 
sections of the form:
\begin{eqnarray*}
P_{2,1}(t,u)Z_1^2 + P_{1,1}(t,u)Z_1Z_2 + P_{1,2}(t,u)Z_1Z_3+  \\
     c_1Z_2^2+c_2Z_2Z_3 + c_3Z_3^2+c_4Z_1Z_4 
\end{eqnarray*}
and
\begin{eqnarray*}
f(t,u,Z_1,Z_2,Z_3) + P_{2,2}(t,u)Z_1Z_4 + \\
P_{1,3}(t,u)Z_2Z_4 + P_{1,4}(t,u)Z_3Z_4+c_5Z_4^2.
\end{eqnarray*}
If this is an $S''$ for a perfect Clifford divisor $D$, then it contains the 
last directrix, which means $c_5=0$.
Even if $c_5=0$, the intersection will in general be smooth if  
$c_4 \neq 0$, and $P_{1,3}$ and $P_{1,4}$ have no common roots.
In this example only a subfamily of positive
codimension of the $(\dim \Aut (\PP^g) + 18- \delta_ 1-\delta_ 2 +\delta _3 +
\delta _4)$-dimensional family  of all complete intersections of that type are
of the form $S''$. }
\end{rem}

\begin{rem}   \label{Dnotperfect}
{\rm If we only assume that $D$ is free (and not perfect), we get the 
following additional 
a priori possible cases:

\vspace{.5cm}
\begin{tabular}{|c|c|c|c|} \hline \label{table3}
$g$ & sing. scroll type & $b_1$ & $(S'')_{virt}$  \\ \hline
$8$ & $(2,2,1,0)$ & $3$   & Smooth  \\ \hline
$9$ & $(2,2,2,0)$ & $4$   & Smooth   \\ \hline
$9$ & $(3,2,1,0)$ & $4$   & Singular \\ \hline
\end{tabular}

\vspace{.5cm}
As proven in Remark \ref{posstype} above, if $\T(2,D)$ has type 
$(2,2,2,0)$ with $b_1=4$, then $S''$ cannot contain the inverse image by 
$i: \T_0 \khpil \T$ of the point singular locus of $\T$. Therefore $S'$ 
cannot contain
the point singular locus of $\T(2,D)$, and $S'' \iso S'$, and $D$ is not perfect.
This completes the proof of Proposition \ref{non2-uple}.
A similar conclusion can be drawn about the two other cases in the last
table, if they occur.}
\end{rem}

\subsection{Projective models with $c=3$}
\label{c=3}
         
Assume $\T=\T(3,D)$ for a free Clifford divisor $D$ with $D^2=0$.
If $\T$ is smooth, we get from Proposition \ref{resolv}(a) that $\O_{S'}$ 
has a resolution (as an $\O_{\T}$-module) of the following form:
\[0 \khpil \O_{\T}(-5\H+(g-6)\F)  \khpil \+ _{k=1}^5 
\O_{\T}(-3\H+b_k\F)  \khpil \]
\[ \+ _{k=1}^5 \O_{\T}(-2\H + a_k\F)  \khpil 
\O_{\T} \khpil \O_{S'} \khpil 0. \]
From Corollary \ref{c3} we conclude that we have such a resolution even if
$\T$ is non-smooth.
We see from \cite{Sc} that we are in a situation very similar to that
of a pentagonal canonical curve, which is natural, since a general hyperplane 
section of $S'$ is such a curve. We do not intend to say as much about this 
situation as about the cases $c=1$ and $2$. Study the skew-symmetrical map $\Phi$ in the resolution
above, already introduced in Corollary \ref{c3}:
\[ \Phi:  \+ _{k=1}^5\O_{\T}(-3\H+b_k\F)  \khpil  \+ _{k=1}^5 \O_{\T}
(-2\H + a_k\F) .\]
Recall that the 
Pfaffians of this map generate the ideal of $S'$ in $\T$. Clearly $\T$ is a 
rational normal scroll of degree $g-4$ in ${\PP^g}$. Let its type be 
$\bf{e}$$=(e_1, \ldots ,e_5)$. 

From Proposition \ref{harr} the dimension of the set of scrolls of type 
$\bf{e}$ in ${\PP}^g$ is equal to $(g+1)^2 -28 - \delta_1= \dim (\Aut ({\PP^g})
-27-\delta _1$, where $\delta_1 := \sum _{i,j} \max (0, e_i - e_j - 1)$\index{$\delta_1$}.
To obtain the number  $18+\dim (\Aut ({\PP^g})$ for the dimension of the set
of projective models of $K3$ surfaces in scrolls of some type, one expects a 
$45$-dimensional set of such models in a given scroll, provided the scroll
type is reasonably well balanced. We will look into this issue, but we 
will not give a rigorous proof that we can find such a $45$-dimensional set.

A given projective model $S'$ is characterized by the ten above-diagonal 
entries of a five-by-five matrix description of the map $\Phi$. 
These entries are sections of: 
\begin{eqnarray*}                
\H-(b_2-a_1)\F, \H-(b_3-a_1)\F, \H-(b_4-a_1)\F, \H-(b_5-a_1)\F, \\
\H-(b_3-a_2)\F, \H-(b_4-a_2)\F, \H-(b_5-a_2)\F, \\
\H-(b_4-a_3)\F, \H-(b_5-a_3)\F, \H-(b_5-a_4)\F. 
\end{eqnarray*}
We have $h^0({\T}, \H-(b_i-a_j)\F)=g+1-5(b_i-a_j)+\delta _{2,i,j}$,
where $\delta _{2,i,j}:= h^1(\PP^1, \E \* \O_{\PP^1} (a_j-b_i))$\index{$\delta _{2,i,j}$} and is zero if and only
if $e_5-(b_i-a_j) \geq -1$.
In all, there set of choices of the ten linear terms has dimension 
\[10(g+1)-\sum _{i>j}(5(b_i-a_j)+\delta _{2,i,j}).\]
Moreover, we have 
\begin{eqnarray*}
\sum _{i>j}(b_i-a_j) & = & b_2+2b_3+3b_4+4b_5 -4a_1-3a_2-2a_3-a_4 \\
      & = & b_2+2b_3+3b_4+4b_5 -4(g-6-b_1)-3(g-6-b_2) \\ 
      &   &         -2(g-6-b_3)-(g-6-b_4)\\
      & = & 4(b_1+ \cdots +b_5)-10(g-6) \\
      & = & 2g-12,
\end{eqnarray*}
where we have used the self-duality of the resolution (Proposition 
\ref{resolv}(b)) which gives $a_i=g-6-b_i$, for $i=1,   \ldots   ,5$) and Proposition 
\ref{resolv}(d) (which gives $b_1+ \cdots +b_5=3g-18$).
Inserting this in 
the expression above we obtain the number
\[10(g+1)-5(2g-12)+\delta _2=70+ \delta _2,\]
where $\delta _2:=\sum_{i>j}\delta _{2,i,j}$\index{$\delta _2$} for the dimension of the set of 
choices of entries in the matrix determining the map $\Phi$.

We see from  $a_i=g-6-b_i$, for $i=1,   \ldots   ,5$, that $\sum_{i=1}^5a_i=2g-12$,
so the average value of the $b_i-a_j$ is $\frac{g-6}{5}$.
The average value of the $e_i$ is $\frac{g-4}{5}$, so if both the
$e_i$ and the $b_j$ (and therefore the $a_j$) are maximally balanced, we will
in fact have $e_5-(b_i-a_j) \geq -1$, so $\delta _{2,i,j}=0$, for each
$i>j$. 

To obtain the desired value 45 in the maximally balanced situation, one needs
to argue that it is correct to subtract 25, in the sence that there is 
typically a 25-dimensional family of matrix decriptions giving rise to 
each projective model of a $K3$ surface as described. We do not know how to do this
in a rigorous way, but the problem is related to the one mentioned in 
\cite[p. 457]{B-E}, where one treats matrix descriptions of maps between two 
free modules of rank 5 over a ring (see also \cite{Bea}). Translating the 
discussion in 
\cite{B-E} into our situation, the issue is: Do two matrices $A'$ and
$A$ have 
the same Pfaffian ideal if and only if there is a matrix $B$, such that
$A'=BAB^t$? In an extremely simple case, take
$g=11$ and $a_i=2$ for all $i$ (and consequently $b_j=3$ for all $j$), so that 
all entries in the matrix representation $A$ of $\Phi$ are sections of the 
same line bundle on $\T$ (in this case 
$\H-\F$). One can imagine the set of five-by-five matrices acting on the
matrix $A$ representing $\Phi$ as $A \khpil BAB^t$, for all $B$ in $GL(5)$. In a situation where the $a_i$ are less balanced, 
one can imagine an analogous matrix $B$ with entries in suitably manufactured 
line bundles, so that the ``shape'' of $A$ is preserved under a similar action.
By this we mean that if entry $A_{i,j}$ of $A$ is a section of a line bundle $L_{i,j}$, then entry $A'_{i,j}$ of $BAB^t$ is also a section of $L_{i,j}$.
One must then count the sections in the entries of $B$, control the
stabilizers of the action, and show that all $A$ with the same Pfaffian ideal 
are in the same orbit by the action.

A natural candidate for such a matrix $B$ is one where the entry $B_{ij}$ is
chosen as a general section of $\O_{\T}((a_j-a_i)\F)=\O_{\T}((b_i-b_j)\F)$, 
for all $(i,j)$. Since $h^0(\O_{\T}((a_j-a_i)\F))+h^0(\O_{\T}((a_i-a_j
)\F)=
2+\max (0, |a_i-a_j|-1)$, we see that $\sum_{i,j}h^0(\O_{\T}((a_j-a_i)\F)=25$
if and only if the $a_i$ are chosen in a maximally balanced way.
Set $\delta _3=\sum_{i>j}\max(0,|a_i-a_j|-1)$. Then the dimension of the set of
choices of matrix $B$ as described is $25+\delta _3$ (we see that $\det B$ 
is a constant, and we look at the closed subset of those $B$
with non-zero determinant). One checks that
$BAB^t$ is antisymmetric, and has entries that are sections in the same
line bundles as the corresponding ones for $A$. This leads to the 
following:

\begin{conj} \label{c3dim}
Let $\T$ be a fixed rational normal scroll of maximally balanced type
and dimension $5$ in $\PP^g$, for $g \geq 9$. Let $\M(\T,c)$\index{$\M(\T,c)$} be the set of projective models of 
$K3$ surfaces $S$ of Clifford index $3$, with a perfect, elliptic Clifford 
divisor  $D$, such that $\T=\T(c,D)$. Then $\dim M=45.$
For an arbitrary scroll type (not necessarily smooth), and given 
combination $(a_1, \ldots ,a_5)$ the corresponding set $\M(\T,c)$ is empty, or it has dimension $45+ \delta _2 -\delta _3$.
We have $\delta_2 \geq \delta_3$ if $\delta _1 \geq 1$.
\end{conj}
              
\begin{rem} \label{c3rem1}
{\rm 
The first statement of the conjecture will be proved in Proposition \ref{genc} 
below. For the second statement, see the discussion above.
The last statement ( $\delta_2 \geq \delta_3$ if $\delta _1 \geq 1$) of the 
conjecture does not follow directly from purely
numerical considerations. As an example, take the case $g=11$, scroll type
$(3,1,1,1,1)$ and $(a_1, \ldots ,a_5)=(1,2,2,2,3)$, which gives 
$(b_1, \ldots ,b_5)=(4,3,3,3,2)$.
Here $\delta _1=4$, and $\delta _3=1$. For all $(i,j)$ with $i>j$, we have
$e_5-(b_i-a_j) \geq -1$, so $\delta _2=0$.

On the other hand the entries outside the diagonal in the first row of a 
matrix 
description of $\Phi$ are sections of 
$\H-(b_2-a_1)\F, \H-(b_3-a_1)\F, \H-(b_4-a_1)\F$ and $\H-(b_5-a_1)\F$, which
here are  $\H-2\F, \H-2\F, \H-2\F$ and  $\H-\F$.
Taking the submaximal minor where we disregard the term $\H-\F$, we see that
$Z_1$ is a factor, since $Z_1$ is a factor in every section of $\H-2\F$.
This is a contradiction, and hence the case does not occur. 
}
\end{rem}

\begin{rem} \label{c3rem2}
{\rm We have now seen (as a special case) that one way to prove the (well 
known) formula $\dim (\Aut(\PP^g))+18$ for the dimension of the set of 
projective models on rational normal scrolls of maximally balanced types in $\PP^g$ (with elliptic Clifford divisor $D$), at least in each of the cases $c=1,2,3$, is to first compute 
the dimension of the set of scrolls, and then add the dimension of the set of 
projective models in each scroll. Using the same method, one deduces the 
well-known fact that the set of $k$-gonal curves in $\PP^g$ on
rational scrollar surfaces of maximally balanced types is empty or has 
dimension 
\[\dim (\Aut(\PP^{g-1}))+2g+2k-5\] 
in each of the cases $k=3,4,5$. But the set is not empty, as is shown for 
example in \cite{bal}, where one shows that for all $k$, the general canonical 
$k$-gonal curve has maximally balanced scroll type (for its gonality scroll).
For canonical curves, the scroll type is determined by the dual scrollar 
invariants $h^0(K-rD)$, in other words by $h^0(rD)$, for $r=1,2,\ldots$ for the 
gonality divisor $D$. One sees that for $k=3,4,5$ one can find the dimension of
the sets of $k$-gonal curves with fixed scrollar invariants (if non-empty) in
$\PP^{g-1}$, corresponding to sets of curves with prescribed values of 
$h^0(rD)$, for $r=1,2,\ldots$, by using similar methods as in the subsections 
above.}
\end{rem}

\subsection{Higher values of $c$} 
\label{highc}

From Proposition \ref{harr} we see that the dimension of 
the set of scrolls of a given type in $\PP^g$ is
$(g+1)^2-3-(c+2)^2 - \delta _1$, where $\delta _1$ is a non-negative
number, which is zero if and only if the scroll type is maximally balanced.
We recall the exact value: 
\[\delta_1 = \sum _{i,j} \max (0, e_i - e_j - 1).\]
Since we know that for all $c$ in the range in question there exists a set of 
dimension
$\dim (\Aut(\PP^g))+18=(g+1)^2+17$ parametrizing projective model of $K3$ 
surfaces in $\PP^g$ with Clifford-index $c$ fibered by elliptic curves on a scroll of some type, we know that for this type, the set of projective models of 
$K3$ surfaces of Clifford 
index $c$, with smooth associated scrolls $\T$, has dimension at least
\[(g+1)^2+17 - ((g+1)^2-3-(c+2)^2)= (c+2)^2+20.\] 
A scroll type with $\delta _1=0$ is then a natural candidate.
We have:

\begin{prop} \label{genc}
Let $g \geq 5$ and $1 \leq c <  \lfloor \frac{g-1}{2} \rfloor$. Let $\T$ be a
fixed rational normal scroll of maximally balanced type of dimension $c+2$ in 
$\PP^g$. Let $\M(\T,c)$ be the set of projective models of $K3$ 
surfaces $S$ of Clifford index $c$, with a perfect, elliptic Clifford divisor 
$D$, such that $\T=\T(c,D)$. Then 
\[\dim \M(\T,c)=(c+2)^2+20.\]
For (not necessarily smooth) scrolls $\T$ with types with $\delta _1 >0$, 
the corresponding set $\M(\T,c)$ is empty, or 
\[\dim \M(\T,c) \leq (c+2)^2+20 + \delta _1.\] 
\end{prop}

\begin {proof}
Let $S$ be a $K3$ surface with Picard group as in Lemma \ref{help}, that is such that $\Pic S \iso \ZZ L + \ZZ D$, with $L^2=2g-2$, $D^2=0$ and $LD=c+2$.
Let us study the scroll $\T(c,D)$. By Proposition \ref{exprop}, we have that $D$ is a free Clifford divisor and the ``dual scrollar invariants'' $d_r$ 
(see Section \ref{sins}) have the form: 
\[d_r=h^0(L-rD)-h^0(L-(r+1)D).\]

Assume that $S$ contains a smooth rational curve $\Gamma$. Then 
$\Gamma=aL+bD$, for integers $a$ and $b$. This gives $a^2(2g-2)+2ab(c+2)=-2$, 
which gives
$a(a(g-1)+b(c+2))=-1$. This, together with $D.\Gamma \geq 0$   gives $a=1$
and $b=\frac{-g}{c+2}$. Hence $S$ contains a rational curve $\Gamma$ if and only if $(c+2) | g$, in which case $\Gamma \sim L-nD$, for $n:= \frac{g}{c+2}$.

We will show that the scroll 
$\T(c,D)$ will be of maximally balanced type.
From the way the scrollar invariants $e_1, \ldots ,e_{c+2}$ are formed from the 
dual scrollar invariants $d_1, d_2, \ldots $ we see that the scroll 
type is maximally balanced if and only if 
\[h^0(L-rD)=(g+1)-r(c+2),\] 
for all $r \geq 0$, such that $L-rD$ is effective.
By Riemann-Roch we see that this happens if and only if $h^1(L-rD)=0$ for these $r$.

Set $B_r:=L-rD$. 

Assume first that $B_r$ is not nef. Then $|B_r|$ has a fixed component $\Sigma$ supported on a union of smooth rational curves. But we have just seen that the only such curve is of the form $\Gamma \sim L-nD$, with $n:= \frac{g}{c+2} \in \ZZ$. So we can write $\Sigma = m\Gamma$, for an integer $m \geq 1$, and denoting the (possibly zero) moving part of $|B_r|$ by $B_r^0$, we have
\[ B_r \sim B_r^0 + m\Gamma. \]
Furthermore, by our assumptions that $B_r$ is not nef, we have $B_r.\Gamma <0$.

We have 
\begin{equation}
  \label{eq:br}
 B_r^0 \sim B_r -m\Gamma \sim L-rD-m\Gamma \sim (1-m)\Gamma + (n-r)D.
\end{equation}

Since $D$ is nef, we have $\Gamma.D \geq 0$ and $B_r^0.D=(1-m)\Gamma. D \geq 0$, whence $m=1$.  

By (\ref{eq:br}) this implies that $B_r^0 \sim (n-r)D = (\frac{g}{c+2}-r)D$, whence
\[ 
\Gamma.B_r^0 = (\frac{g}{c+2}-r)\Gamma.D = (\frac{g}{c+2}-r)L.D = (\frac{g}{c+2}-r)(c+2)=g-r(c+2) \geq 0, \]
and since $\Gamma.B_r = \Gamma.B_r^0 -2 <0$, we must have 
\[ g-r(c+2) = 0 \hs \mbox{or} \hs 1. \]
In the first case, we get $r=\frac{g}{c+2}=n$, whence $B_r^0 = 0$ and $B_r=\Gamma$. In the second case we get the contradiction
\[ n = \frac{g}{c+2} = r + \frac{1}{c+2}. \]

So if $B_r$ is not nef, then $B_r=\Gamma$ and $h^1(B_r)=h^1(\Gamma)=0$.

Now assume $B_r$ is nef.

By Proposition \ref{sd1} and Lemma \ref{sd2}, we have that $h^1(B_r) >0$ if and only if $B_r \sim mE$ for an integer $m \geq 2$ and a smooth elliptic curve $E$. By $B_r ^2=(L-rD)^2=2g-2-2r(c+2)=0$, we get 
\[ r=\frac{g-1}{c+2}. \]
Furthermore $L^2=2g-2=2rmD.E = \frac{2m(g-1)}{c+2} D.E >0$, whence $D.E >0$ and
 $c+2=mD.E$. But this gives
\[ 0 < | \disc (D,E)| = (D.E)^2= \frac{(c+2)^2}{m^2} < (c+2)^2= | \disc (L,D)|, \]
a contradiction, since $L$ and $D$ generate $\Pic S$.

This shows that $h^1(L-rD)=0$ for all $r$ such that
$L-rD \geq 0$.

From Lemma \ref{help} we then have an 
abstract $18$-dimensional family of $K3$ surfaces. From the argument above we
know that these $K3$ surfaces give rise to projective models with balanced
$(c+2)$-dimensional  scrolls, i.e. an  $(\Aut (\PP ^{g})+18)$-dimensional set of
projective models of such surfaces. Hence the first part of the statement of 
the proposition follows, since there is an 
$(\Aut (\PP ^{g})-2 -(c+2)^2)$-dimensional family of
$(c+2)$-dimensional rational normal scrolls
of maximally balanced type in 
$\PP ^{g}$, and all projective models are contained in finitely many such 
scrolls, and all scrolls of the same type are projectively equivalent. We see 
that we can construct a concrete family of dimension $(c+2)^2+20$ in each 
scroll of maximally balanced type, by using the surfaces from 
Lemma \ref{help}.

Assume the scroll $\T$ is not maximally balanced, that is $\delta _1 > 0$. Then
the statement
\[\dim \M(\T,c) \leq (c+2)^2+20 + \delta _1\]
follows from the fact that there is no abstract $19$-dimensional family of 
$K3$ surfaces in $\PP^g$ with Clifford index $c$, and perfect elliptic
Clifford divisor.
Assume $\dim \M(\T,c) \geq (c+2)^2+21+ \delta _1$. Then, by taking the union over
the $(\Aut (\PP ^{g})-2 -(c+2)^2-\delta _1)$-dimensional family of rational 
normal scrolls in $\PP^g$ of the same type as $\T$, we obtain an  
$\Aut (\PP ^{g})+19$-dimensional set of projective models of $K3$ surfaces
in question. Here we use again that all projective models are contained in 
finitely many such scrolls, and all scrolls of the same type are projectively 
equivalent.
\hspace{0.07cm} $\square$ \end{proof} 

\begin{rem} \label{exctypes}
{\rm We also conjecture that $\M(\T,c)$ (defined as above) is empty or:
\[(c+2)^2+20 \leq \dim \M(\T,c)\] 
even if the scroll type of $\T$ is not maximally balanced. 
This conjecture is inspired by Proposition \ref{SiT}, Corollary \ref{upfrom36}
and Remark \ref{c3rem1}. (In many examples for $c=1,2$ with non-zero 
$\delta _1$ a strict inequality is impossible.)

Set $M'=$ (the largest component of) $\Hilb_ {\T}^{(g-1)x^2+2}$. Then it is 
clear that 
\[(c+2)^2+20 \leq \dim M'.\] 
This is true because we can define the relative Hilbert scheme
 \[\M _H'= \Hilb _{\T _H}^{(g-1)x^2+2},\]
where $H$ is the (parameter) Hilbert scheme of rational curves of degree 
$g-c-1$ in $G(c+1,g)$, that is: The parameter space of rational normal 
$(c+2)$-dimenional scrolls in $\PP^g$. Here $\T _H$ is the "universal scroll",
such that the fibre $\T _{[t]}$ is $\T$ if $[t]$ is the parameter point in $H$ 
corresponding to $\T$. It is well known, and follows from for example
\cite{Str}, \cite{RRW}, and \cite[p. 62]{H}, that $H$ is irreducible, and 
that the maximally balanced scrolls correspond to an open dense stratum of $H$.
Since the fibre $M_{[t]}'$ of $\M _H'$ has dimension at least 
$\dim \M(\T,c) =(c+2)^2+20$ for all $[t]$ corresponding to scrolls of 
maximally balanced type, we have $\dim M_{[t]}' \geq =(c+2)^2+20$ for 
the $[t]$ corresponding to scrolls of less balanced types. 

In order to prove the conjecture, we have to pass from $M'$ to $\M(\T,c)$. It is not
entirely clear to us how to do this. If the conjecture is true, we get 
\[(c+2)^2+20 \leq \dim \M(\T,c) \leq c+2)^2+20 + \delta _1.\] 
Moreover the cases $c=1$, scroll type $(7,3,1)$, and $c=2$, scroll types 
$(3,3,3,1)$ with $b_1=6$, and $(4,3,2,1)$ with $b_1=6$, reveal that a strict inequality
$\dim \M(\T,c) < (c+2)^2+20 + \delta _1$ is not always correct, even if 
$\delta_ 1 > 0$. In Section \ref{Conc} one
sees that these cases indeed occur with the fiber $D$ a perfect Clifford 
divisor. In these cases both the most balanced scroll/intersection type and
the mentioned non-balanced types give families of dimension 
$\dim (\Aut(\PP^g))+18$.
}
\end{rem}

\section{$BN$ general and Clifford general $K3$ surfaces}
\label{bncl}

It is shown in \cite{Mu1} that a projective model of a general $K3$ surface in
$\PP ^g$, for $g=6$, $7$, $8$, $9$ and $10$, is a complete
intersection in a homogeneous spaces described below. 

We recall the following definition of Mukai:

\begin{defn}[Mukai \cite{Mu2}]  \label{BNgendef}
 A polarized $K3$ surface $(S,L)$ of genus $g$ is said to be {\rm Brill-Noether ($BN$) general}\index{Brill-Noether (BN) general surface} if the inequality $h^0 (M) h^0(N) < h^0(L)=g+1$ holds for any pair $(M,N)$ of non-trivial line bundles such that $M \* N \iso L$.
\end{defn}

\begin{rem} \label{BNgenrem}
{\rm  One easily sees that this is for instance satisfied if any smooth curve $C \in |L|$ is Brill-Noether general\index{Brill-Noether general curve}, i.e. carries no line bundle $\A$ for which $\rho(A) := g-h^0(\A)h^1(\A) <0$. This is because any nontrivial decomposition $L \sim M + N$ with $h^0 (M) h^0(N) \geq g+1$ yields $h^0(M_C)h^1(M_C) \geq h^0(M)h^1(N) >g$. It is an open question whether the converse is true.} 
\end{rem}

Clearly the polarized $K3$ surfaces which are $BN$ general form a $19$-dimensional Zariski open subset in the moduli space of polarized $K3$ surfaces of a fixed genus $g$.

The following theorem is due to Mukai. We use the following convention: For a vector space $V^i$ of dimension
$i$, we write $G(r,V^i)$ (resp. $G(V^i,r)$) for the Grassmann variety of $r$-dimensional subspaces (resp. quotient spaces) of $V$.

The variety $\Sigma ^{10} _{12} \sub \PP ^{15}$ is a $10$-dimensional spinor variety of degree $12$. Let $V^{10}$ be a $10$-dimensional vector space with a nondegenerate second symmetric tensor $\lambda$. Then $\Sigma ^{10} _{12}$ is one of the two components of the subset of $G(V^{10}, 5)$ consisting of $5$-dimensional totally isotropic quotient spaces \footnote{A quotient $f: V \khpil V'$ is totally isotropic with respect to $\lambda$ if $(f \* f)(\lambda)$ is zero on $V' \* V'$.}.

The variety $\Sigma ^6 _{16} \sub \PP ^{13}$ is the Grassmann variety of $3$-dimensional totally isotropic quotient spaces of a $6$-dimensional vector space $V^6$ with a nondegenerate second skew-symmetric tensor $\sigma$. It has dimension $6$ and degree $16$.

Also, $\Sigma ^5 _{18}= G/P \sub \PP ^{13}$, where $G$ is the automorphism group of the Cayley algebra over $\CC$ and $P$ is a maximal parabolic subgroup. The variety has dimension $5$ and degree $18$. 

Finally, in the case $g=12$, let $V^7$ be a $7$-dimensional
vector space 
and $N \sub \wedge ^2 V \v $ a $3$-dimensional vector space of
skew-symmetric
bilinear forms, with basis $\{ m_1, m_2, m_3 \}$. We denote by 
$\Grass (3, V^7, m_i)$ the subset of $\Grass (3, V^7)$ consisting of 
$3$-dimensional subspaces $u$ of $V$ such that  the restriction of
$m_i$ to $U \x U$ is zero. Then  
$\Sigma ^3_{12} = \Grass (3, V^7, N) := \cap \Grass (3, V^7, m_i)$. It
has dimension $3$ and degree $12$.

\begin{thm}[Mukai \cite{Mu2}] \label{BNgenthm}
 The projective models of $BN$ general polarized $K3$
  surfaces of small genus are as follows:

\vspace{.4cm}
\begin{tabular}{|c|l|} \hline
genus & projective model of BN general polarized $K3$ surface \\ \hline
$2$  &  \hspace{.5cm} $S_2 \hpil \PP ^2$ double covering with branch sextic  \\  \hline
$3$  &  \hspace{.5cm} $(4) \sub \PP ^3$                \\    \hline
$4$  &  \hspace{.5cm} $(2,3) \sub \PP ^4$              \\    \hline
$5$  &  \hspace{.5cm} $(2,2,2) \sub \PP ^5$            \\    \hline
$6$  &  \hspace{.5cm} $(1,1,1,2) \cap G(2, V^5) \sub \PP ^6$      \\  \hline
$7$  &  \hspace{.5cm} $(1 ^8)\cap \Sigma ^{10} _{12} \sub \PP ^7$ \\  \hline
$8$  &  \hspace{.5cm} $(1^6) \cap G(V^6, 2) \sub \PP ^8$          \\  \hline
$9$  &  \hspace{.5cm} $(1^4) \cap \Sigma ^6 _{16} \sub \PP ^9$    \\  \hline
$10$  & \hspace{.5cm} $(1^3) \cap \Sigma ^5 _{18} \sub \PP ^{10}$ \\  \hline
$12$  & \hspace{.5cm} $S_{12} =  (1)     \sub \Sigma ^3_{12}$   \\    \hline
\end{tabular}
\end{thm} 
\vspace{.4cm}

In this section we will compare the notion of $BN$ generality with
our notion of Clifford generality as given in Section \ref{cliff} . We only treat the cases $g \leq 10$. 

It is an easy computation to check that a $BN$ general $K3$ surface is also Clifford general:

\begin{prop} \label{C=BN}
  Let $(S,L)$ be a polarized $K3$ surface of genus $g$. If
  $(S,L)$ is $BN$ general, then it is Clifford general.
\end{prop}

\begin{proof}
Assume that $(S,L)$ is not Clifford general, and let 
$c = \Cliff L < \lfloor \frac{g-1}{2} \rfloor$ and $D$ any Clifford divisor with $F:=L-D$.
Using (C1) and (C3) together with Riemann-Roch one easily computes
\begin{eqnarray}
  \label{eq:bn1} h^0(D)+ h^0(F) & = & \frac{1}{2}D^2+2+ \frac{1}{2}F^2+2 \\
  \nonumber              & = & \frac{1}{2}L^2+2-D.F+2=g+1-c \geq \frac{g+5}{2}.
\end{eqnarray}
Since $h^0(F) \geq h^0(D) \geq 2$ and for fixed $d \geq 2$ the function $f_d(x)=x(d-x)$ obtains
its maximal value in $[2,d]$ at $x=2$, we get
\[h^0(D)h^0(F) \geq 2(h^0(D)+h^0(F)-2) \geq 2(\frac{g+5}{2}-2) =g+1=h^0(L). \]
Hence $(S,L)$ is not $BN$ general.\hspace{0.07cm} $\square$ 
\end{proof}

For low genera we have:

\begin{prop} \label{bn=cl}
  Let $(S,L)$ be a polarized $K3$ surface of genus $g=2,3 \ldots ,7$ or $9$. Then 
  $(S,L)$ is $BN$ general if and only if it is Clifford general.

  If $g=8$ resp. $10$, then $(S,L)$ is Clifford general but not $BN$ 
   general if
  and only if there is an effective divisor $D$ satisfying $D^2=2$ and 
  $D.L =7$ resp. $8$, and there are no divisors satisfying the conditions 
  $(*)$ for $c < 3$ resp. $4$.
\end{prop}

\begin{proof}
  We must investigate the condition that there exists an effective decomposition $L \sim D+F$ such that $h^0(F)h^0(D) \geq g+1$, but 
$\Cliff \O_C(D) \geq \lfloor \frac{g-1}{2} \rfloor$ for any smooth curve $C \in |L|$.

By Riemann-Roch, we have 
$\Cliff \O_C(D) = g+1-h^0(\O_C(D))-h^1(\O_C(D)) \linebreak \leq g+1-h^0(D)-h^0(F)$, so we easily see that we must be in one of the two cases above.
\hspace{0.07cm} $\square$ 
\end{proof}

Since the divisor $D$ in the proposition satisfies $\Cliff L = \Cliff \O_C(D)$, we have $D \in \A ^0(L)$, so we get the following from Propositions \ref{mu1} and \ref{mu2}:

\begin{cor} \label{bn=cl:cor}
 Any divisor $D$ as in Proposition \ref{bn=cl} must satisfy $h^1(D)=h^1(L-D)=0$, and among all such divisors we can find one satisfying the conditions (C1)-(C5).
\end{cor}

By arguing as in the proof of Proposition \ref{exprop}, with the lattice
$\ZZ L \+ \ZZ D$, with
\[  \left[ 
  \begin{array}{cc}
  L^2  &  L.D   \\ 
  D.L  &  D^2 
    \end{array} \right]  = 
    \left[
  \begin{array}{cc}
  2(g-1)  &  c+4   \\ 
  c+4   &  2      
    \end{array} \right]     \]
for $g=8$ and $10$ and $c=\lfloor \frac{g-1}{2} \rfloor =3$ and $4$ respectively, we see that there exists an $18$-dimensional family of isomorphism classes of polarized $K3$ surfaces that are Clifford general but not $BN$ general for both $g=8$ and $g=10$.

We will in the next two sections investigate these two cases. A choice of a subpencil $\{ D _{\lambda} \}$ of $|D|$ gives as before a rational normal scroll $\T$ within which $\varphi_L(S)= S'$ is contained. Unfortunately, as we will see, we no longer have such a nice result about $V = \Sing \T$ as Theorem \ref{mainsing}, since the Clifford index $c$ is now the general one. We will however be able to describe these particular cases in a similar manner, too.

Let us first consider the case $g=8$, where $c=3$. We have $D^2=2$, $D.L=7$, $h^0(L)=9$ and $h^0(L-D)=3$. Since $(L-2D)^2=-6$ and $(L-2D).L=0$, we have $h^0(L-2D)=0$ or $1$ and $h^0(L-3D)=0$.

Recall that the type 
$(e_1,  \ldots  ,e_d)$ of the scroll $\T$, with $d= d_0= \dim \T$, is given by 
\begin{equation} 
  e_i = \# \{ j \hs | \hs d_j \geq i \}-1 ,
\end{equation}
where 
\[ d_i := h^0(L-iD)-h^0(L-(i+1)D).\]
We have $d _{\geq 3} =0$ and $(d_0,d_1,d_2) = (6, 3-h^0(L-2D), h^0(L-2D))$ and the two possible scroll types

\[ (e_1, \ldots ,e_6) = \left\{ \begin{array}{ll}
             (1,1,1,0,0,0) & \hs \mbox{if} \hs h^0(L-2D)=0  \\
             (2,1,0,0,0,0) & \hs \mbox{if} \hs h^0(L-2D)=1
          \end{array}
    \right . \]

We first study the case $h^0(L-2D)=0$, that is $h^1(L-2D)=1$. We have $V= \Sing \T= \PP^2$. Here we already see that Theorem \ref{mainsing} will not apply, since it is clear by the examples given by the lattice above that there are such cases with no contractions across the fibers. Denote the two base points of the pencil $\{ D _{\lambda} \}$ by $p_1$ and $p_2$ and their images under $\varphi_L$ by $x_1$ and $x_2$. We have the following result:

\begin{Lemma} \label{intersectBN}
  Either 

(i) $\R_{L,D} = \emptyset$, or 

(ii) $\R_{L,D} =\{ \Gamma \}$ and $V$ intersects $S'$ in $x_1$, $x_2$ and $y := \varphi_L(\Gamma)$, and $V = <x_1, x_2, y >$.
\end{Lemma}

\begin{proof}
  We first show that $\R_{L,D}$ is either empty or contains at most one curve.

  Choose any smooth $D_0 \in |D|$. Set $F :=L-D$ as usual. Since $\deg F_{D_0}= c+2=5=2p_a(D)+1$, one has that $F_{D_0}$ is very ample, and by arguing as in Lemma \ref{delta0}, we get that $D.\Delta=0$ or $1$. This shows the assertion.

  By arguing as in the proof of Theorem \ref{mainsing}, we get that $V$ intersects $S'$ in at most three points (two of which must of course be $x_1$ and $x_2$) and that these three points 
are independent.
\hspace{0.07cm} $\square$ \end{proof}

By this lemma, there are only two cases occurring for $h^0(L-2D)=0$, which we denote by (CG1)\index{(CG1)} and (CG2)\index{(CG1)}, since they are Clifford general:
\begin{itemize}
  \item[(CG1)] \hs \hs $\R_{L,D} = \emptyset$, 

  \item[(CG2)] \hs \hs $\R_{L,D} = \{ \Gamma \}$. 
\end{itemize}

If $h^0(L-2D)=1$, then since $F^2=2$ we have
\[ L \sim 2D + \Delta, \]
where $\Delta$ is the base divisor of $|F|$ and satisfies $\Delta ^2=-6$, 
$\Delta.L=0$ and $\Delta.D=3$. By arguing as in the proof of Proposition 
\ref{E0-E2}, we find that $L$ is as in one of the five following cases 
(where all the $\Gamma$ and $\Gamma_i$ are smooth rational curves): 
\begin{itemize}
  \item[(CG3)]  \hs \hs \index{(CG3)}$L \sim 2D + \Gamma_1 + \Gamma_2 + \Gamma_3$, with the 
                following configuration:
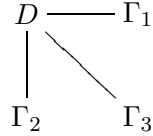
\[ \xymatrix{
{D} \ar@{-}[d] \ar@{-}[dr] \ar@{-}[r] & {\Gamma_1} \\
{\Gamma_2} & {\Gamma_3} 
} \]
 and $\R_{L,D} = \{ \Gamma_1, \Gamma_2, \Gamma_3 \}$, 
  \item[(CG4)] \hs \hs \index{(CG4)}$L \sim 2D + \Gamma + 2\Gamma_0 + 2\Gamma_1+ \cdots 2\Gamma_N + 
               \Gamma_{N+1} + \Gamma_{N+2}$, with the following configuration:
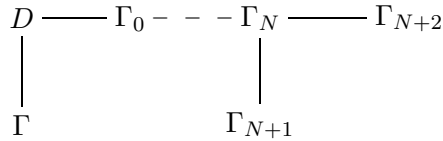
\[ \xymatrix{
{D} \ar@{-}[r] \ar@{-}[d] & {\Gamma_0} \ar@{--}[r] & 
          {\Gamma_N} \ar@{-}[d] \ar@{-}[r] & {\Gamma_{N+2}}   \\
     {\Gamma} & & {\Gamma_{N+1}} 
} \]
 and $\R_{L,D} = \{ \Gamma, \Gamma_0 \}$, 
  \item[(CG5)]  \hs \hs \index{(CG5)}$L \sim 2D + 3\Gamma_1 + 2\Gamma_2 + 2\Gamma_3+ 
                \Gamma_4 + \Gamma_5$, with the 
                following configuration:
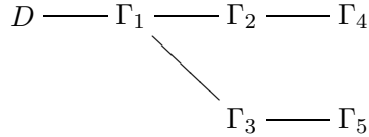
\[ \xymatrix{
{D} \ar@{-}[r] & {\Gamma_1} \ar@{-}[dr] \ar@{-}[r] & {\Gamma_2} \ar@{-}[r] 
         & {\Gamma_4}  \\
  & & {\Gamma_3} \ar@{-}[r] & {\Gamma_5}
} \]
 and $\R_{L,D} = \{ \Gamma_1 \}$.
  \item[(CG6)]  \hs \hs \index{(CG6)}$L \sim 2D +3\Gamma_0 + 4\Gamma_1 +
    2\Gamma_2 + 3\Gamma_3+ 2\Gamma_4 + \Gamma_5$, with the 
                following configuration:
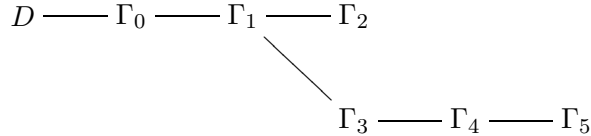
\[ \xymatrix{
{D} \ar@{-}[r] & {\Gamma_0} \ar@{-}[r] & {\Gamma_1} \ar@{-}[dr] \ar@{-}[r]
& {\Gamma_2}  \\
  & & & {\Gamma_3} \ar@{-}[r] & {\Gamma_4} \ar@{-}[r] &{\Gamma_5}
} \]
 and $\R_{L,D} = \{ \Gamma_0 \}$.
  \item[(CG7)]  \hs \hs \index{(CG7)}$L \sim 2D +3\Gamma_0 + 4\Gamma_1 +
    5\Gamma_2 + 6\Gamma_3+ 4\Gamma_4 + 2\Gamma_5 + 3\Gamma_6$, with the 
                following configuration:
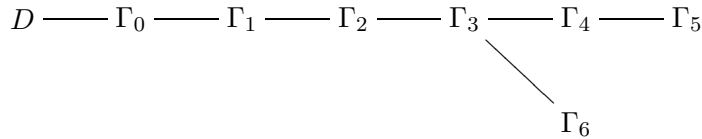
\[ \xymatrix{
{D} \ar@{-}[r] & {\Gamma_0} \ar@{-}[r] & {\Gamma_1}  \ar@{-}[r]
& {\Gamma_2} \ar@{-}[r] & {\Gamma_3} \ar@{-}[dr] \ar@{-}[r]&
{\Gamma_4} \ar@{-}[r] & {\Gamma_5}   \\
  & & & & &{\Gamma_6} 
} \]
 and $\R_{L,D} = \{ \Gamma_0 \}$.

\end{itemize}
Defining $\Z_{\lambda}$ as in (\ref{eq:att})-(\ref{eq:dtt}), we see that 
$\length \Z_{\lambda}=5$ and by arguing 
as in the proof of Theorem \ref{mainsing} in these 
five cases we get that for any $D \in \D$:
\[ V = < \Z_{\lambda} > = \PP^3, \]
any subscheme of length $4$ spans a $\PP^3$ and $V \cap S'$ has support only 
on this scheme.

For the cases (CG1)-(CG7) we can now argue as in Section \ref{singscrolls}.
In particular, we get a commutative diagram as on page 
\pageref{commutativepage}, and Proposition \ref{Bettiprop}, Corollary \ref{Bettis}, 
Propositions \ref{upperres} and \ref{Euler} and Corollary \ref{maxbett} 
still apply. All the $\varphi_L(D _{\lambda})$ have the same Betti-numbers and their resolutions are given in Example \ref{Dkvadrat2}.

In the cases (CG1) and (CG2) the type of $\T_0$ is $(2,2,2,1,1,1)$.
We leave it to the reader to use Lemma \ref{roll} to show that
the only possible combinations of the $b_i$'s (defined in Definition \ref{bene}) are
\[(b_1, \ldots ,b_8)=(4,3,3,3,2,2,2,1), \hs (4,3,3,2,2,2,2,2), \hs (3,3,3,3,2,2,2,2).\]

In the cases (CG3)-(CG7) the type of $\T_0$ is $(3,2,1,1,1,1)$.
We again leave it to the reader to show that $(b_1, \ldots ,b_8)=(4,3,3,2,2,2,2,2)$ is the only possibility. 

We conclude this section by showing that all the cases (CG1)-(CG7) actually exist, by arguing with the help of Propositions 
\ref{morrison} and \ref{plr}.

The case (CG1) can be realized by the lattice just below Corollary 
\ref{bn=cl:cor} and therefore has number of moduli $18$.

We now show that the case (CG2) can be realized by the lattice 
$\ZZ D \+ \ZZ F \+ \ZZ \Gamma$, with intersection matrix:
\[  \left[ 
  \begin{array}{ccc}
  D^2       &  D.F       & D.\Gamma  \\ 
  F.D       &  F^2       & F.\Gamma  \\
  \Gamma.D  &\Gamma.F    & \Gamma ^2
\end{array} \right]  = 
    \left[
  \begin{array}{ccc}
  2 & 5  & 1    \\ 
  5 & 2  & -1   \\
  1 & -1 & -2
\end{array} \right]     \] 
 
One easily checks that this matrix has signature $(1,2)$, so by Proposition 
\ref{morrison} there is an algebraic $K3$ surface with this lattice as its 
Picard lattice.

Set $L:=D+F$. By Proposition \ref{plr} we can assume that $L$ is nef, whence 
by Riemann-Roch $D$, $F >0$.

We first show that $L$ is base point free and that $\Cliff L=3$. Since 
$D.L-D^2-2=3$, we only need to show that there is no effective divisor $B$ on 
$S$ satisfying either
\begin{eqnarray*}
  B^2=0, & B.L=1,2,3,4, \hs \mbox{or} \\
  B^2=2, & B.L=6.
\end{eqnarray*}

Setting $B \sim xD+yF+z\Gamma$, one finds
\[ B.L =7(x+y), \]
which is not equal to any of the values above. Furthermore, $D$ forces 
$(S,L)$ to be non-$BN$ general. Since one easily sees that we cannot be in 
any of the cases (CG1), (CG3)-(CG7), we must be in case (CG2).

We can argue in the same way for the cases (CG3)-(CG7), with the
obvious lattices. The number of moduli of these cases are $16$, $15$,
$14$, $13$ and $12$, respectively. We leave these cases to the reader.

The case $g=10$ is very similar to the previous case. We have $c=4$, $D.L=8$, $h^0(L)=11$ and 
$h^0(L-D)=4$. Since $(L-2D)^2=-6$ and $(L-2D).L=2$, we have $h^0(L-2D)=0$ or 
$1$ and $h^0(L-3D)=0$. This gives as before $d _{\geq 3} =0$ and 
$(d_0,d_1,d_2) = (7, 4-h^0(L-2D), h^0(L-2D))$ and the two possible scroll types
\[ (e_1, \ldots ,e_6) = \left\{ \begin{array}{ll}
             (1,1,1,1,0,0,0) & \hs \mbox{if} \hs h^0(L-2D)=0  \\
             (2,1,1,0,0,0,0) & \hs \mbox{if} \hs h^0(L-2D)=1
          \end{array}
    \right . \]
We now get exactly analogous cases (CG1)'\index{(CG1)'} and 
(CG2)'\index{(CG2)'} as for $g=8$, corresponding to the scroll type $(1,1,1,1,0,0,0)$. If $h^0(L-2D)=1$, write as usual $F:=L-D$ and denote by $\Delta$ the base divisor of $|F|$, so that we have 
\[ L \sim 2D + A + \Delta. \]
for some $A >0$ satisfying $A.L =(L-2D).L=2$ and $A.\Delta=0$. We can now 
show that $2=h^1(R)=\Delta.D$, so that $A^2=-2$ and $A.D=2$. 
By arguing as in the proof of Proposition 
\ref{E0-E2} again, we find that $L$ is as in one of the two following cases 
(where all the $\Gamma_i$ are smooth rational curves such that $\Gamma_i.A=0$): 
\begin{itemize}
  \item[(CG3)'] \hs \hs \hs  \index{(CG3)'}$L \sim 2D + A + \Gamma_1 + \Gamma_2$, with 
                $\Gamma_1.D = \Gamma_2.D=1$, $\Gamma_1.\Gamma_2=0$ and 
                $\R_{L,D} = \{ \Gamma_1, \Gamma_2\}$, 
  \item[(CG4)'] \hs \hs \hs \index{(CG4)'}$L \sim 2D + A + 2\Gamma_0 + 2\Gamma_1+ \cdots +2\Gamma_N + 
               \Gamma_{N+1} + \Gamma_{N+2}$, with all the $\Gamma_i$ having a 
                configuration as in (E2), $\Gamma_i.A=0$,
                $\R_{L,D} = \{ \Gamma_0 \}$, 
\end{itemize}
Defining $\Z_{\lambda}$ as in (\ref{eq:att})-(\ref{eq:dtt}), we see that 
$\length \Z_{\lambda}=4$. By arguing 
as in the proof of Theorem \ref{mainsing} in these 
two cases we get that for any $D \in \D$:
\[ V = < \Z_{\lambda} > = \PP^3 \]
and $V \cap S'$ has support only on this scheme.

As above, we can argue as in Section \ref{singscrolls}, and find that for the cases 
(CG1')-(CG4') the $\varphi_L(D _{\lambda})$ have the same Betti-numbers and their resolutions are given in Example \ref{Dkvadrat2}.

For the cases (CG1)' and (CG2)' the type of $\T_0$ is $(2,2,2,2,1,1,1)$.
Again one can use Lemma \ref{roll} to show that
the only possible combinations of the $b_i$'s (defined in Definition \ref{bene}) are 
\[(4,3,3,3,3,3,3,3,2,2,2,2,2)\hs \mbox{ and } \hs (3,3,3,3,3,3,3,3,3,2,2,2,2)\]
for the case (CG1)', and 
\[
(4,3,3,3,3,3,3,3,2,2,2,2,2), \hs (3,3,3,3,3,3,3,3,3,2,2,2,2) \hs \mbox{
  and } \]
\[ (4,3,3,3,3,3,3,3,3,2,2,2,1)  
\]
for the case (CG2)'.

The type of $\T_0$ for the cases (CG3)' and (CG4)' is 
$(3,2,2,1,1,1,1)$. The only possible cases for the $b_i$'s are found to be

\[ (4,4,4,3,3,3,2,2,2,2,2,2,2) \hs \mbox{ and } \hs (4,4,3,3,3,3,3,2,2,2,2,2,2),\]

In the same way as for the cases (CG1)-(CG7), we can show the existence of each of the types
(CG1)'-(CG4)' for $g=10$ and show that their number of moduli $18$, $17$, $16$ and $15$ respectively.

These results will all be summarized in the next section, together with all non-Clifford general projective models for $g \leq 10$.

\section{Projective models of $K3$ surfaces of low genus}
\label{Conc}

In this section  we will use the results obtained in the previous
ones to classify all
projective models of non-BN-general $K3$ surfaces of genus at most $10$. 
Together with Mukai's description of the general models we are then 
able to give a complete classification and characterization for these
genera.
The central part of the section is Section \ref{tables} where we 
give tables summing up the essential information concerning  the
various projective models appearing of non-$BN$ general $K3$ surfaces for 
$5 \leq g \leq 10$.

An important intermediate step is performed in Section  \ref{perfe} 
where we  describe the possible perfect Clifford
divisors for $c=1$,
$2$ and $3$ and also in some more detail the cases where $h^1(L-2D) >0$, 
since this last number determines the singular locus of the scroll $\T$ by (\ref{eq:sing}). The description is valid for all genera, not
only the small ones, but for  $g \geq 11$ cases with $c \geq 4$
appear, even for $(S,L)$ non-Clifford general. For $g \leq 10$ we always have $c \leq 3$ for
the non-Clifford general models.

The reason why we concentrate on perfect Clifford divisors is 
to make the classification in Section \ref{tables} simpler. If we did
not restrict to perfect Clifford divisors, we would get more
projective models, but the extra projective models would also have
been possible to describe with a perfect Clifford divisor, whence they
would belong to our list.

Section \ref{decomp} is purely technical, and devoted to a new
decomposition of the divisor $R=L-2D$
for each free Clifford divisor $D$. The new decomposition, with the
added property  (\ref{eq:noe4}) is necessary to make the description 
in Section \ref{perfe} work.

In Section  \ref{scroty} we show how one can calculate the scroll
types of the relevant ambient scrolls appearing in the various cases.

The exposition  in Section  \ref{tables} contains detailed information
about the Picard lattice of $S$, and the singularity type of
$S'=\varphi_L(S)$ in many subcases. In Section   \ref{sce}
we show how this information can be obtained in some typical cases,
and leave the arguments in the remaining ones to the reader.

\subsection {A new decomposition of $R$}
\label{decomp}

Assume that $D$ is a free Clifford divisor. We
recall from  Section \ref{tedious} that $R=L-2D$, and that $\Delta=0$
if $H^0(R)=0$ by Lemma \ref{delta0}. If $R>0$, we have 
$L=2D+A+\Delta$, where $D+A$ is the moving part of $F:=L-D$, and
$\Delta$ is the base divisor of $F$. So $R \sim A+\Delta$ is an effective
decomposition of $R$. Recall from Lemma \ref{deltaA} that $\Delta.A=0$,
except for the cases (E3) and (E4). To make the classification
simpler,
we would like to find a new effective decomposition of $R$, say 
$R\sim  A'+\Delta'$, with a stronger property than the one in  
Lemma \ref{deltaA},
namely that $\Delta''.A'=0$ for every effective $\Delta'' \le \Delta'$.
At the same time we would like $A'$ and $\Delta'$ to enjoy the same
intersection properties and cohomological properties as $A$ and
$\Delta$, so that the results in Section \ref{tedious} are still valid.
(We are grateful to Gert M. Hana for pointing out the need for such a 
new decomposition)

\begin{prop} \label{newpair}
Let $(S,L)$ be a polarized $K3$ surface of non-general Clifford index,
with free Clifford index $D$ not as in (E3) or (E4), and such that
$R:=L-2D>0$. Let $A$ and $\Delta$ be defined as above. Then there
exists an effective decomposition $R=A'+\Delta'$\index{$A'$}\index{$\Delta'$} such that the
following properties hold:

\begin{eqnarray}  
\label{eq:noe1} \Delta' \leq \Delta  \hs \mbox{ and } \hs 
  A' \geq A \\
\label{eq:noe2}
\left[ \begin{array}{ccc}
       D^{2}        & D.A
       & D.\Delta                    \\
       D.A          & A^2
       & A.\Delta                    \\
       D.\Delta        & A.\Delta
       & \Delta^2                    \\
       \end{array}         \right]
  =
  \left[ \begin{array}{ccc}
       D^{2}        & D.A'
       & D.\Delta'                   \\
       D.A'          & A'^2
       & A'.\Delta'                  \\
       D.\Delta'        & A'.\Delta'
       & \Delta'^2                   \\
       \end{array}         \right]   \\
\label{eq:noe3} h^i(A')=h^i(A) \hs \mbox{and} \hs
  h^i(\Delta')=h^i(\Delta) \hs \mbox{for } i=0,1,2. \\
\label{eq:noe4} \Delta''.A'=0 \hs \mbox{for every effective} \hs
  \Delta'' \leq \Delta'. 
\end{eqnarray}
\end{prop}

\begin{rem} \label{wb}
{\rm Note that $R \sim A+\Delta$ always satisfies
(\ref{eq:noe1})-(\ref{eq:noe3}), 
so that property
(\ref{eq:noe4}) is the reason why we want to find a new decomposition. Moreover
note that (\ref{eq:noe1})-(\ref{eq:noe4}) ensure that all the important
results in Sections \ref{sing} and \ref{tedious} for $A$ and $\Delta$
are still valid for 
$A'$ and $\Delta'$.
To be more precise, Proposition \ref{Gset}, Remark \ref{perf1}, 
Proposition \ref{h1delta},
Lemmas \ref{delta0}, \ref{deltaA}, \ref{lemA}, \ref{Adec} and
Proposition \ref{h1A} are valid with $A$ and $\Delta$ replaced by $A'$
and $\Delta'$.}
\end{rem}

We will give an algorithmic proof of Proposition \ref{newpair}. First
we will state and prove the following result.

\begin{Lemma} \label{mer}
Assume we are neither in case (E3) nor (E4), and that we have an
effective decomposition $R \sim A_i+\Delta_i$ such that (\ref{eq:noe1})-(\ref{eq:noe3}) hold.
If there exists a smooth rational curve $\Gamma \leq \Delta_i$ such
that
$\Gamma.A_i >0$, then $\Gamma.A_i=1$, and $\Gamma.D=0$.
\end{Lemma}

\begin{proof}
Remember that $F_0 \sim D+A$ is the moving part of $F$.
We write $F_i := D+A_i$.
Then $F_0 \leq F_i \leq F$.
Hence we have $h^0(F_i)=h^0(F_0)$.
Since $(A_i,\Delta_i)$ satisfies (\ref{eq:noe2}), we have $F_0^2=F_i^2$.
Riemann-Roch then gives $h^1(F_i)=h^1(F_0)=0$ by Lemma \ref{F0smooth}.
Here we use that we are not in any of the cases (E3) or (E4).

Let now $\Gamma \leq \Delta_i$ be a smooth rational curve such that
$\Gamma.A_i >0$.

Using Riemann-Roch yet another time gives
\[
h^0(F_i+\Gamma)-h^0(F_i)=F_i.\Gamma -1 + h^1(F_i+\Gamma)=0.
\]
Hence $F_i.\Gamma \leq 1$.

Since $D.\Gamma \geq 0$ we get $\Gamma.A_i \leq 1$.
So if $\Gamma.A_i > 0$,
then $\Gamma.A_i=1$ and $\Gamma.D=0$.
\hspace{0.07cm} $\square$ \end{proof}

\renewcommand{\proofname}{Proof of Proposition  \ref{newpair} }
\begin{proof}
Write $\Delta_0 := \Delta$ and $A_0 := A$.
Given an effective decomposition $R \sim A_i+\Delta_i$ satisfying
(\ref{eq:noe1})-(\ref{eq:noe3}),
 assume that there exists
a smooth rational curve $\Gamma \leq \Delta_i$ such
that $\Gamma.A_i > 0$.
Write $A_{i+1} := A_i + \Gamma$ and $\Delta_{i+1} := \Delta_i -
  \Gamma$.
Then $R \sim A_{i+1}+\Delta_{i+1}$ satisfies
(\ref{eq:noe1})-(\ref{eq:noe2}) by 
the previous lemma.
Clearly $h^2(A_{i+1})=h^2(\Delta_{i+1})=0$, and since
$A_{i+1}^2=A_i^2$ and $\Delta_{i+1}^2=\Delta_i^2$, it suffices to show
that $h^0(A_{i+1})=h^0(A_i)$ and $h^0(\Delta_{i+1})=h^0(\Delta_i)$ to
show that $R \sim A_{i+1}+\Delta_{i+1}$ satisfies
(\ref{eq:noe3}). It is obvious that $h^0(\Delta_{i+1})=h^0(\Delta_i)=1$
since $\Delta_{i+1} \leq \Delta_i$. Furthermore
$h^0(A_{i+1})=h^0(A_i)$ since $\Gamma$ is fixed in $A_i+\Gamma$, as 
$\Gamma.(A_i+\Gamma)=-1$. 
Hence $R \sim A_{i+1}+\Delta_{i+1}$ satisfies (\ref{eq:noe3}).

We repeat this process if necessary,
and it is obvious  that the procedure will stop after finitely many
steps,
say for $i=n \geq 0$, since $\Delta_0 > \Delta_1 > \ldots > \Delta_n$.
For the effective decomposition $R \sim A_n+\Delta_n$ there exists no
smooth rational curve $\Gamma \leq \Delta_n$ such that $\Gamma.A_n>0$,
whence the decomposition satisfies (\ref{eq:noe4}) as well.
\hspace{0.07cm} $\square$ \end{proof}

\renewcommand{\proofname}{Proof}
\begin{Lemma}  \label{l:gamma.A=0}
Assume we are neither in case (E3) nor (E4), and that
for every $\Gamma \in \mathcal{R}_{L,D}$ we have $\Gamma.A=0$.
Then $R \sim A+\Delta$ satisfies (\ref{eq:noe1})-(\ref{eq:noe4}).
\end{Lemma}
\begin{proof}
If an effective divisor $B \leq \Delta$ satisfies $A.B \neq 0$, then
some smooth rational curve $\Gamma \leq \Delta$ (possibly equal to
$B$), must satisfy $A.\Gamma < 0$. But $(D+A).\Gamma=0$ or $1$.
Hence $\Gamma \in \mathcal{R}_{L,D}$.
But then $A.\Gamma=0$ by the assumptions, a contradiction.
\hspace{0.07cm} $\square$ \end{proof}

\subsection{Perfect Clifford divisors for low $c$}
\label{perfe}

From now on $(A',\Delta')$ will be a pair of divisors 
satisfying (\ref{eq:noe1})-(\ref{eq:noe4}).

Furthermore, in the list below  we have: 
\begin{itemize}
\item $\Gamma$ is a smooth rational curve such that $\Gamma.D=1$ and 
        $\Gamma.A'=0$.
\item $\Gamma_1$ and $\Gamma_2$ are 
        smooth rational curves such that 
        $\Gamma_1.D=\Gamma_2.D=1$ and 
        $\Gamma_1.A'=\Gamma_2.A'=\Gamma_1.\Gamma_2=0$.
\item \index{$\Delta_0$} $\Delta_0 := 2{\Gamma_0} + 2{\Gamma_1} +  \cdots  + 2{\Gamma_N} + {\Gamma_{N+1}} + {\Gamma_{N+2}}$, for $N \geq 0$, 
with a configuration with respect to $D$ as in (E2) and such that 
$A'.{\Gamma_i} =0$ for $i=0, \ldots ,N+2$.
\end{itemize}

Also we denote the different cases by $\{c,D^2\}$\index{$\{c,D^2\}$}.

Here is the list of all possible perfect Clifford divisors for $c=1$,
$2$ and $3$, and the cases where $h^1(R) >0$:

 \[ \mathbf{c=1, L^2 \geq 8} \]
\begin{itemize}
  \item[$\{1,0\}$] \hspace{.5cm}  $D^2=0$, $D.L=3$, $\dim \T=3$.
  \item[$\{1,2\}$] \hspace{.5cm} $D^2=2$, $L^2= 10$, $L \sim 2D+\Gamma$ as in (E0), 
                   $\dim \T=4$.
\end{itemize}

Moreover,  $h^1(R) \not =0$ if and only if $L$ is as in the following case:
\begin{itemize}
  \item[$\{1,0\}^a$] \hspace{.5cm} $L \sim 2D+A'+\Gamma$, ${A'}^2 \geq -2$, $D.A'=2$,  
                 $L^2=A'^2+10 \leq 18$ with equality if and only if $L
                 \sim 6D+3\Gamma$, $h^1(R)=1$, $\R_{L,D}=\{ \Gamma \}$. 
\end{itemize}

 \[ \mathbf{c=2, L^2 \geq 12} \]
\begin{itemize}
  \item[$\{2,0\}$] \hspace{.5cm} $D^2=0$, $D.L=4$, $\dim \T=4$.
  \item[$\{2,2\}$] \hspace{.5cm} $D^2=2$, $D.L=6$, $L^2 \leq 18$ with equality if and only if 
               $L \sim 3D$, $\dim \T=5$.
  \item[$\{2,4\}$] \hspace{.5cm} $D^2=4$, $L^2 =16$, $L \sim 2D$ as in (Q), $\dim \T=6$.
\end{itemize}

Moreover,  $h^1(R) \not =0$ if and only if $L$ is as in one of the following cases:
\begin{itemize}
  \item[$\{2,0\}^a$] \hspace{.5cm} $L \sim 2D+A'+\Gamma$, ${A'}^2 \geq -2$, $D.A'=3$,  
                 $L^2={A'}^2+14 \leq 32$ with equality if and only if 
              $L \sim 8D + 4 \Gamma$, $h^1(R)=1$, $\R_{L,D}=\{ \Gamma \}$.
  \item[$\{2,0\}^b$] \hspace{.5cm} $L \sim 2D+A'+\Gamma_1+\Gamma_2$, ${A'}^2 \geq 0 $, $D.A'=2$, 
                 $L^2= {A'}^2 + 12 \leq 16$ with equality if and only if
                 $L \sim 4D +2\Gamma_1 + 2\Gamma_2$, $h^1(R)=2$, $\R_{L,D}=\{ \Gamma_1
                 ,\Gamma_2  \}$ \label{footpage}
                 \footnote{If $L^2=14$, then the moving part
                   of $A'$ is a perfect Clifford divisor of type
                   $\{2,2\}$ containing $D$, and if $L^2=16$, then $A'$
                   is a perfect Clifford divisor of type
                   $\{2,4\}$ containing $D$.}.
  \item[$\{2,0\}^c$] \hspace{.5cm} $L \sim 2D+A'+ \Delta_0$, ${A'}^2 \geq 0 $, $D.A'=2$, 
                 $L^2={A'}^2 + 12 \leq 16$ with equality if and only if
                 $L \sim 4D +2\Delta_0$, $h^1(R)=2$, 
                 $\R_{L,D}=\{ \Gamma_0 \}$
                 \footnote{Same comment as above.}.
  \item[$\{2,2\}^a$] \hspace{.5cm} $L \sim 2D+\Gamma_1+\Gamma_2$ as in (E1), $L^2=12$, $h^1(R)=1$,
                 $\R_{L,D}=\{ \Gamma_1 ,\Gamma_2 \}$.
\item[$\{2,2\}^b$] \hspace{.5cm} $L \sim 2D + \Delta_0$ as in (E2), $L^2=12$, $h^1(R)=1$, $\R_{L,D}=\{ \Gamma_0 \}$.
\end{itemize}

 \[ \mathbf{c=3, L^2 \geq 16} \]
\begin{itemize}
  \item[$\{3,0\}$] \hspace{.5cm} $D^2=0$, $D.L=5$, $\dim \T=5$.
  \item[$\{3,2\}$] \hspace{.5cm} $D^2=2$, $D.L=7$, $L^2 \leq 22$, $\dim \T=6$.
  \item[$\{3,4\}$] \hspace{.5cm} $D^2=4$, $L^2 =18$, $L \sim 2D+\Gamma$ as in (E0), $\dim \T=7$.
\end{itemize}

Moreover,  $h^1(R) \not =0$ if and only if $L$ is as in one of the following cases:
\begin{itemize}
  \item[$\{3,0\}^a$] \hspace{.5cm} $L \sim 2D+A'+\Gamma$, ${A'}^2 \geq -2$, $D.A'=4$,  
                 $L^2={A'} ^2+18 \leq 50$ with equality if and only if $L
                 \sim 10D + 5\Gamma$, $h^1(R)=1$, $\R_{L,D}=\{ \Gamma \}$.
  \item[$\{3,0\}^b$] \hspace{.5cm}  $L \sim 2D+A'+\Gamma_1+\Gamma_2$, ${A'}^2 \geq 0$, $ D.A'=3$, 
                 $L^2={A'}^2+16 \leq 24$, $h^1(R)=2$, $\R_{L,D}=\{ \Gamma_1 ,\Gamma_2 \}$.
  \item[$\{3,0\}^c$] \hspace{.5cm} $L \sim 2D+A'+ \Delta'$, ${A'}^2 \geq 0 $, $D.A'=3$, 
                 $L^2={A'}^2+16 \leq 24$, $h^1(R)=2$, $\R_{L,D}=\{ \Gamma_0 \}$.
  \item[$\{3,2\}^a$] \hspace{.5cm} $L \sim 2D+{A'}'+\Gamma$, ${A'}^2 =-2$, $D.A'=2$,  
                 $L^2=16$, $h^1(R)=1$, $\R_{L,D}=\{ \Gamma \}$.
\end{itemize}

This list is obtained by using the relations $(*)$ and (\ref{eq:index}) in 
Section \ref{cliff} together with Propositions \ref{Gset}, \ref{h1delta} and
\ref{bound}. We now show how it works for $c=3$. 

The three cases $\{3,0\}$, $\{3,2\}$ and $\{3,4\}$ follow directly from the relations 
$(*)$. If $D^2>0$, then by (\ref{eq:index}) we must have $L^2 \leq 24$. Assume
$L^2=24$ and consider the divisor $E:=L-3D$. This satisfies $E^2=0$ and 
$E.L = 3$, thus inducing a Clifford index $1$ on $L$, a contradiction. So 
$L^2 \leq 22$.

Now assume we are in case  $\{3,0\}$ and $h^1(R) >0$. By Propositions 
\ref{h1delta}, \ref{bound} and \ref{newpair}  we have 
$1 \leq D.\Delta'= D.\Delta \leq 2$. Since 
\[ 5=D.L = A'.D + \Delta'.D, \]
we have the two possibilities:
\begin{itemize}
  \item[(a)] $\Delta'.D=1$ and $D.A'=4$,
  \item[(b)] $\Delta'.D=2$ and $D.A'=3$.
\end{itemize}

In case (a), there has to exist a smooth rational curve $\Gamma$ in
the support of $\Delta'$ such that $\Gamma.D=1$ and $\Gamma.A'=0$ (the
last equality follows from (\ref{eq:noe4}) of Proposition \ref{newpair}). Write 
\[ L \sim 2D + A' + \Gamma + \Delta''. \]
Clearly $D.\Delta''=A'.\Delta''=0$, and by $0=\Gamma.L = 2-2 +\Gamma.\Delta''$, we also get $\Gamma.\Delta''=0$, whence
\[ (2D + A' + \Gamma).\Delta'' =0, \]
and we must have $\Delta'' =0$ since $L$ is numerically $2$-connected. This 
establishes case $\{3,0\}^a$. From the Hodge index theorem on $L$ and
$A'$ it follows that $L^2 \leq 50$ with equality if and only if $4L \sim 5A'$.

In case (b), there either exist two (and only two) disjoint smooth rational 
curves $\Gamma_1$ and $\Gamma_2$ in the support of $\Delta'$ such that 
$\Gamma_1.D=\Gamma_2.D=1$ and $\Gamma_1.A'=\Gamma_2.A'=0$, or there exists
one and only one smooth rational curve 
$\Gamma_0$ in the support of $\Delta'$ (necessarily with multiplicity $2$) such that $\Gamma_0.D=1$ and 
$\Gamma_0.A'=0$. Arguing as in the proof of Proposition \ref{E0-E2},
these two cases give the cases $\{3,0\}^b$ and $\{3,0\}^c$
respectively. Again it follows from the Hodge index theorem on $L$ and
$A'$ that $L^2 \leq 24$.

Assume we are in case $\{3,2\}$ and $h^1(R) >0$. By Propositions 
\ref{h1delta}, \ref{bound}, and \ref{newpair} we have 
$L^2 = 16$, $D.\Delta' =1$ and $\Delta'^2=-2$. There has to exist a
smooth rational curve $\Gamma$ in the support of $\Delta'$  such that 
$\Gamma.D=1$ and 
$\Gamma.A'=0$. Arguing as above, we easily find that $L \sim 2D+A'+\Gamma$. 
Since $7=D.L = 2D^2 + A'.D + \Delta'.D$, we have $A'.D=2$, and since $16= L^2= 18 +{A'}^2$, we must have ${A'}^2=-2$. This is case $\{3,2\}^a$. 

We leave the easier cases $c=1$ and $2$ to the reader, but make a
comment on the cases $\{2,0\}^b$ and $\{2,0\}^c$.

From  the Hodge index theorem on $L$ and
$A'$ we get that $L^2 \leq 16$ with equality if and only if $L \sim
2A'$. If ${A'}^2=2$ or $4$, one calculates
\[ A'.L - A^2-2 = 2, \]
$(A'-D)^2 \geq -2$ and $(A'-D).D=2$, whence by Riemann-Roch $A' \geq D$,
so $D$ does not satisfy the condition (C6). However, since $A'$
computes the Clifford index
of $L$, we have $h^1(A)=h^1(A')=0$ by Proposition \ref{newpair}, whence $D$ is
perfect by Lemma \ref{pseudo-perfect}. If $L \sim 2A'$, one easily
sees that $A'$ is base point free, whence perfect.

These cases are particularly interesting, since $S'$ is contained in
two scrolls of different types. 

Note that for $g \leq 10$ (equivalently $L^2 \leq 18$) a polarized $K3$ 
surface of non-general Clifford index must have $c \leq 3$, so the above cases are sufficient to consider these surfaces. We know that the general $K3$ surface has general Clifford index. The following proposition considers the 
dimension of the families in the list above.

\begin{prop} \label{dimfam}
 The number of moduli of polarized $K3$ surfaces of genus $g$, with $5 \leq g  \leq 10$,
 and non-general Clifford index $c >0$ of each of the types $\{1,0\}$, 
 $\{1,2\}$, $\{2,0\}$, $\{2,2\}$ with $g \leq 9$, $\{3,0\}$, $\{3,2\}$ and 
 $\{3,4\}$ is $18$, and of each 
 of the types $\{2,2\}$ with $g=10$ and $\{2,4\}$ is $19$. 

 Furthermore the general projective model of each of these types satisfies 
 $h^1(L-2D)=0$, and the general projective model of each of these types except
 for the types $\{1,2\}$ and $\{3,4\}$ is smooth.

 The number of moduli of each of these types with $h^1 (L-2D) >0$ is 
 $\leq 17$, except for the type $\{1,0\}^a$ for $g=10$, whose number
 is $18$.
\end{prop}

\begin{proof}
 In the cases $\{2,2\}$ with $g=10$ and $\{2,4\}$ we have $L \sim 3D$ and $L \sim 2D$
 respectively, so it is clear that those cases can be realized with a Picard group of rank $1$ and hence live in $19$-dimensional families.

 In the other cases, one easily sees that $L$ and $D$ are linearly independent, and we will show that these cases can all be realized with a Picard group of rank $2$. Arguing as in the proof of Proposition \ref{exprop} we easily see that there is a $K3$ surface $S$ such that $\Pic S \iso \ZZ L \+ \ZZ D$ such that $L^2$, $L.D$ and $D^2$ have the values corresponding to the different cases in question and such that $D$ is a perfect Clifford divisor for $L$. This has already been done for $D^2=0$ in the proof of Proposition \ref{exprop}, and a case by case study establishes the proof in the other cases.

Recall now that $h^1(L-2D) >0$ if and only if there exists a smooth rational curve $\Gamma$ such that $\Gamma.L=0$ and $\Gamma.D=1$, and (since $c >0$) $\varphi_L (S)$ is singular if and only if there exists a smooth rational curve $\Gamma$ such that $\Gamma.L=0$ and $\Gamma.D=0$ or $1$.

Assuming that the rank of the Picard group is two, we can write 
 $\Gamma = aL + bD$, for $a,b \in \QQ$. The conditions $\Gamma ^2=-2$, 
 $\Gamma .L=0$ and $\Gamma.D=0$ or $1$ give the equations:
\begin{eqnarray*}
  a^2(g-1)+ab(c+2+D^2)+\frac{b^2D^2}{2} & = & -1, \\
  a(2g-2)+b(c+2+D^2) & = & 0 \mbox{ and } \\
  a(c+2+D^2) + bD^2 & = & 0 \mbox{ or } 1.
\end{eqnarray*}
  A case by case check reveals that we have a solution only when $\Gamma.D=1$
  and then in the following cases:
\begin{itemize}
  \item[(a)] $\{1,2\}$, with $\Gamma \sim L-2D$,
  \item[(b)] $\{3,4\}$, with $\Gamma \sim L-2D$,
  \item[(c)] $\{1,0\}$ for $g=10$, with $L \sim 6D + 3\Gamma$.  
\end{itemize}
One can easily show that case (c) can be realized with a lattice of the form 
$\ZZ D \+ \ZZ \Gamma$, with $D^2=0$, $D.\Gamma=1$ and $\Gamma^2=-2$.

This concludes the proof of the Proposition.
\hspace{0.07cm} $\square$ \end{proof}

\subsection{The possible scroll types}
 \label{scroty}

We now would like to study which scroll types are possible for each value of 
$(g,c,D^2)$ with $5 \leq g \leq 10$ and $1 \leq c \leq 3$. Recall that the type 
$(e_1,  \ldots  ,e_d)$ of the scroll $\T$, with $d= \dim \T$, is given by 
\begin{equation} 
  e_i = \# \{ j \hs | \hs d_j \geq i \}-1 ,
\end{equation}
where 
\begin{eqnarray*}
  d =  d_0 & := & h^0(L)-h^0(L-D) = c+2+\frac{1}{2}D^2,    \\
       d_1 & := & h^0(L-D)-h^0(L-2D)= d_0 -r,              \\
       \vdots   &    &                                          \\ 
       d_i & := & h^0(L-iD)-h^0(L-(i+1)D), \\
       \vdots   &    &                                          \\ 
\end{eqnarray*}
with
\begin{equation} 
  r = \left\{ \begin{array}{ll}
             D^2 + h^1 (L-2D) & \mbox{ if $ L \not \sim 2D$ (equiv. 
                                       $D^2 \not = c+2$)}, \\
             D^2-1            & \mbox{ if $ L \sim 2D$ (equiv. 
                                       $D^2 = c+2$)}
          \end{array}
    \right .     
\end{equation}

In the cases $\{1,2\}$ and $\{3,4\}$, which are both of type (E0), and the 
case $\{2,4\}$, which is of type (Q), we have $h^0(L-2D)=1$ and $h^0(L-iD)=0$ 
for all $i \geq 3$, so the scroll types are immediately given.

In the case $\{2,2\}$ with $g=10$, we have $L \sim 3D$, so $h^0(L-2D)=h^0(D)=2$, $h^0(L-3D)=1$ and $h^0(L-iD)=0$ for all $i \geq 3$.

We will now consider one by one the remaining cases and gather the result in 
the tables in Section \ref{tables} below.

If $c=1$ or $2$ and $D^2=0$ the possible scroll types are given in
Section \ref{lowind}. We now briefly review these cases.

Let us first consider the case $c=1$ and $D^2=0$ (case $\{1,0\}$).

For $g=5$ the two possible scroll types are $(1,1,1)$ and $(2,1,0)$. One 
easily sees that the first case corresponds to
$(d_0,d_1,d_2)=(3,3,0)$, whence $h^0 (L-2D)= h^1 (L-2D)=0$ and the
second corresponds to $(d_0,d_1,d_2)=(3,2,1)$, whence $h^1 (L-2D)=1$
and we are in case $\{1,0\}^a$. 

For $g=6$ we have three possible scroll types: $(2,1,1)$, $(2,2,0)$ and 
$(3,1,0)$. Comparing with the possible values of the $d_i$, one finds that the first case corresponds to $h^0 (L-2D)=1$ (and $h^1 (L-2D)=0$). Moreover, the two last cases corresponds to the case  
$\{1,0\}^a$ with $A' \not > D$ and $A' >D$ respectively. 

For $g=7$ there are four possible scroll types: $(2,2,1)$, $(3,1,1)$,
$(3,2,0)$ and $(4,1,0)$. We see that the two first cases correspond to 
$h^1 (L-2D)=0$, with $h^0 (L-3D)=0$ and $1$ respectively.  The two last cases have $h^1 (L-2D)=1$ and therefore correspond to $\{1,0\}^a$ with $A' \not >2 D$ and $A' > 2D$ respectively. 

We now leave the cases $g=8$, $9$ and $10$ to the reader.

If $c=2$ and $D^2=0$ (case $\{2,0\}$) then $12 \leq L^2 \leq 18$.

We leave the easiest case $g=7$ to the reader.

If $g=8$ we have seen that the four possible scroll types are $(2,1,1,1)$, $(2,2,1,0)$, $(3,1,1,0)$ and $(3,2,0,0)$. The scroll $(2,1,1,1)$
corresponds to $h^1(R)=0$,
whereas the scrolls $(2,2,1,0)$ and $(3,1,1,0)$ correspond to the case
$\{1,0\}^a$ with $A' \not > D$ and $A' > D$ respectively.  The type $(3,2,0,0)$ corresponds to a polarized surface that also has a different perfect Clifford divisor, and is hence contained in another scroll as well, by the footnote on page \pageref{footpage}.

If $g=9$ we have seen that the five possible scroll types are
$(2,2,1,1)$, $(3,1,1,1)$, $(2,2,2,0)$, $(3,2,1,0)$ and $(4,2,0,0)$.
The types $(2,2,1,1)$ and \linebreak $(3,1,1,1)$ correspond to $h^1(R)=0$ with
$h^0 (L-3D)=0$ and $1$ respectively. (One easily sees that the scroll type
$(3,1,1,1)$ can be realized by a $K3$ surface $S$ with Picard group
$\Pic S \iso \ZZ D \+ \ZZ \Gamma$, for a smooth rational curve
$\Gamma$ satisfying $\Gamma.D=2$, and with $L \sim 3D + \Gamma$.
Therefore, it has number of moduli  $18$.) The scroll types
$(2,2,2,0)$ and  $(3,2,1,0)$ correspond to the case $\{1,0\}^a$ with
$A' \not > D$ and $A' > D$ respectively. The type $(4,2,0,0)$ corresponds to a polarized surface that also has a different perfect Clifford divisor, and is hence contained in another scroll as well, by the footnote on page \pageref{footpage}.

If $g=10$ there are again five possible scroll types: $(2,2,2,1)$, $(3,2,1,1)$, $(3,2,2,0)$, $(3,3,1,0)$ and $(4,2,1,0)$. Again the two first correspond to 
$h^1(R)=0$ with $h^0 (L-3D)=0$ and $1$ respectively. 
The three last cases correspond to the case $\{1,0\}^a$ with $h^0 (A'-D)=1$ and $2$ respectively, but $A' \not > 2D$ for the two first cases, and 
$A' > 2D$ for the last case.

If $c=2$ and $D^2=2$ (case $\{2,2\}$), then $12 \leq L^2 \leq 16$ (the case $L^2=18$ being already treated). We have 
\[ (L-3D).L = L^2-18 < 0. \]
By the nefness of $L$ we must have $h^0(L-3D)=0$. Since $(L-2D)^2 =L^2-16$, we get by Riemann-Roch $h^0(L-2D) = \frac{1}{2}L^2 -6 +h^1(R)$. This gives $d_{\geq 3} =0$ 
and the two possibilities $(d_0, d_1, d_2) = (5, 3,\frac{1}{2}L^2 -6)$ or 
$(5,2,1)$, the latter occurring if and only if $L^2=12$ and $L$ is of type (E1)
or (E2) (the special cases $\{2,2\}^a$ and $\{2,2\}^b$). The corresponding scroll types in 
the first situation are then $(1,1,1,0,0)$ for $g=7$, $(2,1,1,0,0)$ for $g=8$ 
and $(2,2,1,0,0)$ for $g=9$. For $g=7$ and $L$ of type (E1)
or (E2) the scroll type is $(2,1,0,0,0)$.

If $c=3$ and $D^2=0$ (case $\{3,0\}$) then $L^2=16$ or $18$. We have 
\[ (L-3D).L = L^2-15  \leq 3 \]
and
\[ (L-4D).L = L^2-20 < 0. \]  
This gives immediately 
$h^0(L-iD)=0$ for all $i \geq 3$, whence $d_{\geq 4}=0$. Also, since $c=3$, 
we must have $h^0(L-3D) \leq 1$. We also have by Riemann-Roch $h^0(L-2D)= 
\frac{1}{2}L^2 - 8 + h^1(R)$. 

Let us first consider the case $g=9$. Then we have 
$(d_0, d_1, d_2, d_3) = (5, 5-h^1(R), h^1(R)-h^0(L-3D), h^0(L-3D))$.
If $h^1(R)=0$, then $h^0(L-2D)=h^0(L-3D)=0$ and 
$(d_0, d_1, d_2, d_3)=(5,5,0,0)$. The corresponding scroll type is
$(1,1,1,1,1)$. The cases with $h^0(R)= h^1(R) >0$ are $\{3,0\}^a$,  $\{3,0\}^b$ and $\{3,0\}^c$. In the first we have $h^0(R)= h^1(R)=1$, whence $h^0(L-3D)=0$ and 
$(d_0, d_1, d_2, d_3)=(5,4,1,0)$. The corresponding scroll type is 
$(2,1,1,1,0)$. In the cases $\{3,0\}^b$ and $\{3,0\}^c$, we have 
$h^0(R)= h^1(R)=2$. If $h^0(L-3D)=0$ (eqv. $A' \not >D$), we get 
$(d_0, d_1, d_2, d_3)=(5,3,2,0)$ and the scroll type is $(2,2,1,0,0)$. If 
$h^0(L-3D)=1$ (eqv. $A' >D$), we get 
$(d_0, d_1, d_2, d_3)=(5,3,1,1)$ and the scroll type is $(3,1,1,0,0)$.

If $g=10$, we have 
$(d_0, d_1, d_2, d_3) = (5, 5-h^1(R), 1+ h^1(R)-h^0(L-3D),
h^0(L-3D))$. If $h^1(R)=0$, then $h^0(L-2D)=1$ and $h^0(L-3D)=0$ and 
$(d_0, d_1, d_2, d_3)=(5,5,1,0)$. The corresponding scroll type is
$(2,1,1,1,1)$. The cases with $h^1(R) >0$ are $\{3,0\}^a$,  $\{3,0\}^b$ and $\{3,0\}^c$ as in the case $g=9$. Arguing as in that case, we get $(d_0, d_1, d_2, d_3)=(5,4,2,0)$ 
and scroll type $(2,2,1,1,0)$ in the case $\{3,0\}^a$ (where $h^1(R)=1$), and we get $(d_0, d_1, d_2, d_3)=(5,3,3,0)$ and scroll type $(2,2,2,0,0)$ if $h^0(L-3D)=0$ (eqv. $A' \not >D$), and $(d_0, d_1, d_2, d_3)=(5,3,2,1)$ and scroll type 
$(3,2,1,0,0)$ if $h^0(L-3D)=1$ (eqv. $A' >D$) in the two latter cases 
(where $h^1(R)=2$). 

If $c=3$ and $D^2=2$ (case $\{3,2\}$), then $L^2=16$ or $18$. We have 
\[ (L-3D).L = L^2-21  <0, \]
whence $h^0(L-iD)=0$ for all $i \geq 3$, whence $d_{\geq 3}=0$. By 
Riemann-Roch, $h^0(L-2D)= \frac{1}{2}L^2 - 8 + h^1(R)$ and we have 
$(d_0, d_1, d_2)=(6,4-h^1(R),\frac{1}{2}L^2 - 8 + h^1(R),0)$. 

If $g=9$ and $h^1(R)=0$, then $(d_0, d_1, d_2)=(6,4,0)$ and the 
corresponding scroll type is $(1,1,1,1,0,0)$. The case with 
$h^1(R) >0$ is given by $\{3,2\}^a$. In this case we have 
$(d_0, d_1, d_2)=(6,3,1)$ and the corresponding scroll type is $(2,1,1,0,0,0)$.

If $g=10$, then we automatically have $h^1(R)=0$, whence 
$(d_0, d_1, d_2)=(6,4,1)$ and the 
corresponding scroll type is $(2,1,1,1,0,0)$.

We will summarize these results below.

Furthermore, we can prove, by arguing with lattices that all the cases mentioned above exist, and calculate the number of their moduli. In many cases, we can also explicitly find an expression for $L$ in terms of $D$ and some smooth rational curves on the surface. Also, by studying the Picard lattices, we can find the curves that are contracted by $L$, and hence find the 
singularities of the generic surfaces in question.

All these informations are also summarized below, in section \ref{tables}.

\subsection{Some concrete examples}
 \label{sce}

In this section, we focus on some concrete examples, to give the reader an idea of the proofs. We then leave all the other cases to the reader, and conclude the section by giving the list of all projective models of genus $\leq 10$ in section \ref{tables}.

\begin{exa}
{\rm
  We start with an easy case: $g=6$, $c=1$, $D^2=0$ and the scroll type 
$(3,1,0)$. This occurs if $L$ is of type $\{ 1,0 \} ^a$ with $A' > D$ (and also $\R_{L,D}= \{ \Gamma \}$). By 
arguing as in the proof of Proposition \ref{E0-E2} we find that 
$L \sim 3D + 2\Gamma + \Gamma_0 + \Gamma_1$, where $\Gamma$, $\Gamma_0$ and
$\Gamma_1$ are smooth rational curves, with the following configuration:
\[ \xymatrix{
{D} \ar@{-}[d] \ar@{-}[r] & {\Gamma} \ar@{-}[r]  &{\Gamma_1} \\
{\Gamma_0}. & &
} \]

By Propositions \ref{morrison} and \ref{plr} there is an algebraic $K3$ surface $S$ with Picard group
$\Pic S = \ZZ D \+ \ZZ \Gamma \+ \ZZ \Gamma_0 \+ \ZZ \Gamma_1$ and 
 intersection matrix
\[  \left[ 
  \begin{array}{cccc}
   D^2      & D.\Gamma        & D.\Gamma_0        & D.\Gamma_1     \\ 
 \Gamma.D   & \Gamma^2        & \Gamma.\Gamma_0   & \Gamma.\Gamma_1    \\    
 \Gamma_0.D & \Gamma_0.\Gamma & \Gamma_0  ^2      & \Gamma_0.\Gamma_1  \\   
 \Gamma_1.D & \Gamma_1.\Gamma & \Gamma_1.\Gamma_0 & \Gamma_1 ^2     
    \end{array} \right]  = 
    \left[
  \begin{array}{cccc}
  0  & 1  & 1  & 0     \\ 
  1  & -2 & 0  & 1      \\  
  1  & 0  & -2 & 0     \\  
  0  & 1  & 0  & -2             
    \end{array} \right],     \] 
and such that $L:= 3D + 2\Gamma + \Gamma_0 + \Gamma_1$ is nef (whence by Riemann-Roch $D$ and $\Gamma_0 >0$).

We have $D.L-D^2-2=1$. To show that $L$ is base point free and of Clifford index $1$, it suffices to show that there is no effective divisor $E$ such that $E^2=0$ and $E.L=1$ or $2$.

Set $E \sim xD + y\Gamma + z\Gamma_0 + w\Gamma_1$. Since we can assume $E \in \A^0(L)$, and $E$ base point free, we easily see that 
\[ E.\Gamma_0 =x-2z =0 \hs \mbox{or} \hs 1, \]
whence
\[ E.L =3x +z = 7z \hs \mbox{or} \hs 7z+3, \]
which can never be equal to $1$ or $2$.

By Riemann-Roch either $\Gamma >0$ or $-\Gamma >0$. If the latter is the case, write $\Gamma = - \gamma$, and we then have $D=D_0+ \gamma$ with $D_0 >0$, 
since $D.\gamma=-1$. Therefore, we can write
\[ L \sim 3(D_0+ \gamma) - 2\gamma + \Gamma_0 + \Gamma_1 =  3D_0 + 2\gamma 
+ \Gamma_0 + \Gamma_1.\]

We can use the same argument if $-\Gamma_1 >0$, so possibly after a change of 
basis, we can assume $D$, $\Gamma$, $\Gamma_0$ and $\Gamma_1 >0$. It is then easy to check that $D$ is nef, whence a perfect Clifford divisor.
}
\end{exa}
%\vspace{.2cm}

\begin{exa}
{\rm
Let us consider the case $g=9$, $c=2$, $D^2=0$ and the scroll type 
$(3,2,1,0)$. This occurs if $L$ is of type $\{ 2,0 \} ^a$ with ${A'} > D$ (and also $\R_{L,D}= \{ \Gamma \}$). An analysis as in the proof of Proposition \ref{E0-E2} shows that $L$ is one of the following three types:
\begin{itemize}
 \item[(a)] $L \sim 3D + 2\Gamma + \Gamma_1 + \Gamma_2$, with 
  the following configuration:
\[ \xymatrix{
{D} \ar@{-}[d] \ar@{-}[dr] \ar@{-}[r] & {\Gamma}  \ar@{-}[d] \\
{\Gamma_2}   & {\Gamma_1} 
} \]
\item[(b)] $L \sim 3D + 2\Gamma + \Gamma_1 + \Gamma_2$, with 
  the following configuration:
\[ \xymatrix{
{D} \ar@{=}[d] \ar@{-}[r] & {\Gamma}  \ar@{-}[r] & {\Gamma_1} \\
{\Gamma_2}   & & 
} \] 
\item[(c)] $L \sim 3D + 2\Gamma + \Gamma_1 + \Gamma_2 +  \cdots  + \Gamma_{n+3}$, 
for $n \geq 0$ (in general $n=0$) with 
  the following configuration:
\[ \xymatrix{
          & {D} \ar@{-}[d] \ar@{-}[dr] \ar@{-}[r] & {\Gamma_2} \ar@{--}[d] \\
{\Gamma_1} \ar@{-}[r] & {\Gamma}  & {\Gamma_{n+3}.}    
} \]
\end{itemize}
(Actually case (b) can be looked at as a special case of case (c), 
with ``$n=-1$''.)

  One can easily show that both cases (a) and (b) do not occur with a Picard 
group of rank $<4$, and case (c) does not occur with a Picard group of rank 
$<5$.
We now show that both case (a) and (b) occur with a Picard group of rank $4$.

We first consider case (a).

By Propositions \ref{morrison} and \ref{plr} there is an algebraic $K3$ surface $S$ with Picard group
$\Pic S = \ZZ D \+ \ZZ \Gamma \+ \ZZ \Gamma_1 \+ \ZZ \Gamma_2$ and 
 intersection matrix corresponding to the configuration above, and such that 
$L:= 3D + 2\Gamma + \Gamma_1 + \Gamma_2$ is nef (whence $D$, $\Gamma_1$, 
$\Gamma_2 >0$ by Riemann-Roch).

We calculate $D.L -D^2-2=2$. To show that $L$ is base point free and that 
$\Cliff L=2$ with $D$ as a perfect Clifford divisor, it will suffice to show 
that there are no divisor $B$ on $S$ satisfying $B^2=0$, $B.L=1,2,3$ or 
$B^2=2$, $B.L=6$, and that $D$ is nef.

Write $B \sim xD + y\Gamma + z\Gamma_1 + w\Gamma_2$. Since we can assume 
$B \in \A^0(L)$, and $B$ base point free, we easily see that 
\[ B.\Gamma_2 =x-2w =0 \hs \mbox{or} \hs 1, \]
and
\[ B.\Gamma_1 =x+y-2z =0, 1 \hs \mbox{or} \hs 2. \]
By the Hodge index theorem one also finds
\[ B.D =y+z+w = \left\{ \begin{array}{ll}
             1 & \mbox{ if $B^2=0$,} \\
             2 & \mbox{ if $B^2=2$.}
          \end{array}
    \right . \]
Also, we have
\[B.L=4x+3z+w = \left\{ \begin{array}{ll}
             1, 2, 3 & \mbox{ if $B^2=0$,} \\
             6 & \mbox{ if $B^2=2$.}
          \end{array}
    \right . \]
One checks by inspection that these four equations have no integer solutions.

By Riemann-Roch, either $\Gamma>0$ or $-\Gamma >0$. As in the previous example, possibly after a change of basis one can assume that $\Gamma>0$ and that $D$ is nef, whence perfect.

We now consider case (b).

Again by Propositions \ref{morrison} and \ref{plr} there is an algebraic $K3$ 
surface $S$ with Picard group
$\Pic S = \ZZ D \+ \ZZ \Gamma \+ \ZZ \Gamma_1 \+ \ZZ \Gamma_2$ and 
 intersection matrix corresponding to the configuration for (b) above, and such 
that 
$L:= 3D + 2\Gamma + \Gamma_1 + \Gamma_2$ is nef (whence $D$ and $\Gamma_2 >0$ 
by Riemann-Roch).

We calculate $D.L -D^2-2=2$. To show that $L$ is base point free and that 
$\Cliff L=2$ with $D$ as a perfect Clifford divisor, it will again suffice to 
show 
that there are no divisor $B$ on $S$ satisfying $B^2=0$, $B.L=1,2,3$ or 
$B^2=2$, $B.L=6$, and that $D$ is nef.

Write $B \sim xD + y\Gamma + z\Gamma_1 + w\Gamma_2$ as before. Again by the Hodge index theorem and since we can assume 
$B \in \A^0(L)$, and $B$ base point free, we get 
\[ B.D =y+2w = \left\{ \begin{array}{ll}
             1 & \mbox{ if $B^2=0$,} \\
             2 & \mbox{ if $B^2=2$,}
          \end{array}
    \right . \]
\[ B.\Gamma_2 =2(x-w) =0 \hs \mbox{or} \hs 2, \]
and
\[ B.\Gamma_1 =y-2z = -1,0, \hs \mbox{or} \hs 1 \]
(since we do not know whether it is $\Gamma_1$ or $-\Gamma_1$ which is 
effective). Combining these equations with
\[B.L= 2(2x+y+2w), \]
we find no integer solutions. Again, possibly after a change of basis, we get 
that $D$ is perfect and that all $D$, $\Gamma$, $\Gamma_1$ and $\Gamma_2>0$.

We can also check which curves are contracted by $L$. 

In case (a), the only 
contracted curve is in general $\Gamma$, so all surfaces in that family has an
$A_1$ singularity, and the general surface has only such a singularity. 
Furthermore $S''$ is then in general smooth.

In case (b), the only 
contracted curves are in general $\Gamma$ and $\Gamma_1$ , so all surfaces 
in that family has an
$A_2$ singularity, and the general surface has only such a singularity. 
Furthermore $S''$ is then necessarily singular.

By comparing with the table on page \pageref{table2}, we then find that case (a) has 
$b_1=3$ and case (b) has $b_1=2$.
}
\end{exa}

\begin{exa}
{\rm
As an easy example we consider the case  $g=10$, $c=1$, $D^2=0$ and the 
scroll type $(5,2,1)$. This occurs if $\R_{L,D}=\emptyset$ and $h^0(L-5D)=1$.
An analysis as in the proof of Proposition \ref{E0-E2} shows that 
$L \sim 5D + 3\Gamma_1 + 2\Gamma_2 + \Gamma_3$, with 
the following configuration:
\[ \xymatrix{
{D} \ar@{-}[r] & {\Gamma_1} \ar@{-}[r] & {\Gamma_2} \ar@{-}[r] & {\Gamma_3.} 
} \]
One can easily show that this cannot be achieved with a Picard group of rank 
$<4$.

By Propositions \ref{morrison} and \ref{plr} again there is an algebraic $K3$ 
surface $S$ with Picard group
$\Pic S = \ZZ D \+ \ZZ \Gamma_1 \+ \ZZ \Gamma_2 \+ \ZZ \Gamma_3$ and 
 intersection marix corresponding to the configuration above, and such that 
$L:= 5D + 3\Gamma_1 + 2\Gamma_2 + \Gamma_3$ is nef (whence $D$ and 
$\Gamma_1 >0$ by Riemann-Roch).

We calculate $D.L -D^2-2=1$. To show that $L$ is base point free and that 
$\Cliff L=1$ with $D$ as a perfect Clifford divisor, it will suffice to show 
that there is no divisor $E$ on $S$ satisfying $E^2=0$, $E.L=1,2$ and that 
$D$ is nef.

By the Hodge index theorem $36 E.D \leq (E+D)^2 L^2 \leq ((E+D).L)^2 \leq 25$, whence $E.D=0$. Writing $E \sim xD + y\Gamma_1 + z\Gamma_2 + w\Gamma_3$, we get
\[E.D =y=0,\]
whence
\[E.L =3x+y=3x \not = 1 \hs \mbox{or} \hs 2. \]

Possibly after a change of basis, we get 
that $D$ is perfect and that all $D$, $\Gamma_1$, $\Gamma_2$ and $\Gamma_3>0$.

One finds that the only contracted curves with this Picard group are 
$\Gamma_2$ and $\Gamma_3$, so the general surface in this family has an $A_2$ singularity.
}
\end{exa}

\begin{exa}
{\rm
  We give a more involved example: $g=10$, $c=2$, $D^2=0$ and the 
  scroll type $(3,2,1,1)$. This occurs if $\R_{L,D}=\emptyset$ and 
  $h^0(L-3D)=1$. By the table on page 62, we must have $b_1=3$ or $4$, and we will now show that both these cases exist (with the number of moduli $17$ and $16$ respectively).

  One easily sees that there is no way to achieve this situation with a 
  Picard group of rank $<3$. We will now show that it is possible with a 
  Picard group of rank $3$.

By Propositions \ref{morrison} and \ref{plr} there is an algebraic $K3$ 
surface $S$ with Picard group
$\Pic S = \ZZ D \+ \ZZ \Gamma_1 \+ \ZZ \Gamma_2$ and 
 intersection matrix
\[  \left[ 
  \begin{array}{ccc}
   D^2        &   D.\Gamma_1      & D.\Gamma_2     \\ 
 \Gamma_1.D   &   \Gamma_1^2      & \Gamma_1.\Gamma_2     \\    
 \Gamma_2.D   & \Gamma_2.\Gamma_1 & \Gamma_2 ^2         
    \end{array} \right]  = 
    \left[
  \begin{array}{ccc}
  0  & 2  & 0       \\ 
  2  & -2 & 1        \\  
  0  & 1  & -2       
    \end{array} \right],     \] 
and such that $L:= 3D + 2\Gamma_1 + \Gamma_2$ is nef (whence by 
Riemann-Roch $D$ and $\Gamma_1 >0$).

We have $D.L-D^2-2=2$. To show that $L$ is base point free and of Clifford index $2$ with $D$ as a perfect Clifford divisor, it suffices to show that there is no effective divisor $B$ such that $B^2=0$, $B.L=1, 2, 3$, or $B^2=2$, 
$B.L=6$.

By the Hodge index theorem one has
\[ 18(B^2+2B.D) = L^2(B+D)^2 \leq ((B+D).L)^2 = (B.L+4)^2, \]
which gives $B.D \leq 1$.

Writing $B \sim xD + y\Gamma_1 + z\Gamma_2$, we have $B.D =2y$, whence
$y=0$.

Since either $\Gamma_2>0$ or $-\Gamma_2>0$ and we can assume $B \in \A^0(L)$, we must have $B.\Gamma_2=y-2z =-2z = -1,0,1$. We therefore get $z=0$.

So $B$ is a multiple of $D$, a contradiction.

Possibly after a change of basis, we get 
that $D$ is perfect and that also $\Gamma_2>0$.

One finds that the only contracted curve with this Picard group is 
$\Gamma_2$, so that all surfaces in this family have at least an 
$A_1$ singularity, and the general such surface has such a singularity. By comparing with the table on page \pageref{table2}, we see that we must have $b_1=4$.

But there is also another family of surfaces. Again we find that there is an 
algebraic $K3$ surface $S$ with Picard group
$\Pic S = \ZZ D \+ \ZZ \Gamma_1 \+ \ZZ \Gamma_2 \+ \ZZ \Gamma_3$ and 
 intersection matrix
\[  \left[ 
  \begin{array}{cccc}
   D^2      & D.\Gamma_1        & D.\Gamma_2        & D.\Gamma_3  \\ 
 \Gamma_1.D & \Gamma_1^2        & \Gamma_1.\Gamma_2 & \Gamma_1.\Gamma_3  \\    
 \Gamma_2.D & \Gamma_2.\Gamma_1 & \Gamma_2 ^2       & \Gamma_2.\Gamma_3    \\ 
 \Gamma_3.D & \Gamma_3.\Gamma_1 & \Gamma_3.\Gamma_2 & \Gamma_3 ^2       
  \end{array} \right]  = 
    \left[
  \begin{array}{cccc}
  0  & 2  &  1  & 1     \\ 
  2  & -2 &  0  & 0      \\  
  1  & 0  & -2  & 0   \\ 
  1  & 0  &  0  &  -2 
    \end{array} \right],     \] 
and such that $L:= 3D + \Gamma_1 + \Gamma_2 + \Gamma_3$ is nef.

One can show that $\Cliff L=2$ with $D$ as a perfect Clifford divisor (again 
after possibly changing the basis). Furthermore, one finds that with this 
lattice, there are no contracted curves, whence $S'$ is smooth. By comparing 
with the table on page \pageref{table2}, we see that we must have $b_1=3$.
}
\end{exa}

\subsection{The list of projective models of low genus}
\label{tables}

We will now summarize essential information about birational projective
models $S'$ of $K3$ surfaces of genera $5 \leq g \leq 10$.
In some cases we are able to give a resolution of $S'$ in its scroll
$\T$. When we are not able to do this, we give the vector bundle a
section of which cuts out $S''$ in $\T_0 \iso \PP (\E)$ (which is the
dual of the vector bundle $F_1$ in the resolution
\[ \cdots \hpil F_2 \hpil F_1 \hpil \O_{\T_0} \hpil \O_{S''} \hpil 0.)\]
This vector bundle is a direct sum of line bundles, which we write as
a linear combination of the line bundles $\H$ and $\F$ on $\PP(\E)$,
where $\H= i^*\O_{\PP^g}(1)$ and $\F = \pi^*\O_{\PP^1}(1)$, with
\[ 
\xymatrix{
   \PP (\E) \ar[r]^{i} \ar[d]^{\pi} & \T   \sub \PP ^{g}   
  \\ 
   \PP ^1.  & 
} \]
Also note that we in all cases have $\I_{S'/\T} = i_*\I _{S''/\T_0}$ by
   Proposition \ref{ideal}, and that in most cases, by Remark
   \ref{f_*}, the sections of $F_1 \v$ are constant on the fibers of
   $i$, whence they also give ``equations'' cutting
   out $S'$ in $\T$ set-theoretically. 

The singularity type listed in the
rightmost column of the tables below indicates that for ``almost all''
$K3$ surfaces in question its projective model $S'$ has singularities
exactly as indicated, and that none have milder singularities. By ``almost all'' we here mean that the moduli of
the exceptional set of $K3$ surfaces in question with different singularity 
type(s) is strictly smaller than the number of moduli listed in the middle 
column.
These exceptional $K3$ surfaces will have ``worse'' singularities
than the one(s) listed in the rightmost column.  

In the tables below, $c$ is as usual the Clifford index of $L$, $D$ is a perfect Clifford divisor and $A$ is as defined in (\ref{eq:F}). To find the tables we use $A'$ and $\Delta'$ as above, but since $A$ and $A'$ (resp. $\Delta$ and $\Delta'$) enjoy the same intersection and cohomology properties, we can then reintroduce $A$ (resp. $\Delta$). In particular, the tables below are still valid if one exchanges $A$ with $A'$.

 \[ \mathbf{g=5} \]
The general projective model is a complete intersection of three hyperquadrics. The others are as follows:
\vspace{.4cm}

\begin{tabular}{|c|c|c|c|c|c|} \hline
$c$ & $D^2$ & scroll type & $\#$ mod. & type of $L$ & sing.  \\ \hline
$1$ & $0$   & $(1,1,1)$   & $18$   & $h^0(L-2D)=0$ & sm. \\  \hline
$1$ & $0$   & $(2,1,0)$   & $17$   & $\{1,0\}^a$, $A^2=-2$ & $A_1$                                  \\  \hline
\end{tabular}

\vspace{.4cm}
In these cases $\O_{S'}$ has the following $\O_{\T}$-resolution:
\[
  0 \hpil \O_{\T}(-3\H+\F) \hpil \O_{\T} \hpil \O_{S'} \hpil 0.
\]

 \[ \mathbf{g=6} \]
The general projective model is a hyperquadric section of a Fano $3$-fold of 
index $2$ and degree $5$. 
The others are as follows:
\vspace{.4cm}

\begin{tabular}{|c|c|c|c|c|c|} \hline
$c$ & $D^2$ & scroll type & $\#$ mod. & type of $L$ & sing.  \\ \hline
$1$ & $0$   & $(2,1,1)$   & $18$ & $h^0(L-2D)=0$ & sm. \\  \hline
$1$ & $0$   & $(2,2,0)$   & $17$ & $\{1,0\}^a$, $A^2=0$, 
                                    $A \not >D$ & $A_1$  \\  \hline
$1$ & $0$   & $(3,1,0)$   & $16$ & $\{1,0\}^a$, $A^2=0$, 
                                    $A >D$ $^{(i)}$ & $A_2$ \\  \hline
$1$ & $2$   & $(2,1,0,0)$   & $18$ &  (E0) & $A_1$  \\  \hline
\end{tabular}

\vspace{.4cm}
In the three first cases  $\O_{S'}$ has the following $\O_{\T}$-resolution:
\[
  0 \hpil \O_{\T}(-3\H+2\F) \hpil \O_{\T} \hpil \O_{S'} \hpil 0.
\] 

In the last case, $S'$ has a resolution:
\begin{eqnarray*} 
0 & \hpil & \O_{\T}(-4\H + 2\F) \+\O_{\T}(-4\H + \F)  \\
& \hpil & \O_{\T}(-2\H + 2\F) \+ \O_{\T}(-3\H + \F) \+
\O_{\T}(-3\H) \\
& \hpil &\O_{\T} \hpil \O_{S'} \hpil 0
\end{eqnarray*}

\vspace{.4cm}
{\bf Comments on the types of $L$:}
\begin{itemize}
 \item[(i)] $L \sim 3D + 2\Gamma + \Gamma_0 + \Gamma_1$, with the following 
configuration:

\[ \xymatrix{
{D} \ar@{-}[d] \ar@{-}[r] & {\Gamma} \ar@{-}[r]  &{\Gamma_1} \\
{\Gamma_0} & &
} \]
\end{itemize}

 \[ \mathbf{g=7} \]
The general projective model is a complete intersection of $8$
hyperplanes in $\Sigma ^{10} _{12}$, as described in the beginning of
Section \ref{bncl}.

The other projective models are as follows:
\vspace{.4cm}

\begin{tabular}{|c|c|c|c|c|c|} \hline
$c$ & $D^2$ & scroll type & $\#$ mod. & type of $L$ & sing.  \\ \hline
$1$ & $0$   & $(2,2,1)$   & $18$ & $h^0(L-2D)=2$, $h^0(L-3D)=0$ & sm. \\ \hline
$1$ & $0$   & $(3,1,1)$   & $16$ & $h^0(L-2D)=2$, $h^0(L-3D)=1$ & sm. \\ \hline
$1$ & $0$   & $(3,2,0)$   & $17$ & $\{1,0\}^a$, $A^2=2$,  
                           $A >D$,     $A \not >2D$  & $A_1$ \\  \hline
$1$ & $0$   & $(4,1,0)$   & $16$ & $\{1,0\}^a$, $A^2=2$, 
                                    $A >2D$ $^{(i)}$ & $A_3$ \\  \hline
$2$ & $0$   & $(1,1,1,1)$ & $18$ &  $h^0(L-2D)=0$ & sm. \\ \hline
$2$ & $0$   & $(2,1,1,0)$ & $17$ & $\{2,0\}^a$, $A^2=-2$ & $A_1$\\ \hline
$2$ & $0$   & $(2,2,0,0)$ & $16$ & $\{2,0\}^b$ or $\{2,0\}^c$, $A \not
>D$ $^{(ii)}$ & $2A_1$\\ \hline          
$2$ & $0$   & $(3,1,0,0)$ & $15$ &     
               $\{2,0\}^b$ or $\{2,0\}^c$, $A >D$ $^{(iii)}$ & $2A_2$  \\ \hline
$2$ & $2$   & $(1,1,1,0,0)$ & $18$ &  $h^0(L-2D)=0$ & sm. \\ \hline
$2$ & $2$   & $(2,1,0,0,0)$ & $17$ & (E1) or (E2) $^{(iv)}$ & $2A_1$    \\ \hline

\end{tabular}

\vspace{.4cm}
In the cases $(c,D^2)=(1,0)$ $\O_{S'}$ has the following 
$\O_{\T}$-resolution:
\[
  0 \hpil \O_{\T}(-3\H+3\F) \hpil \O_{\T} \hpil \O_{S'} \hpil 0.
\]

In the cases $(c,D^2)=(2,0)$ $\O_{S'}$ has the following 
$\O_{\T}$-resolution:
\begin{eqnarray*}
  0 \hpil \O_{\T}(-4\H+(g-1)\F)  & \hpil &  \\
\O_{\T}(-2\H+b_1\F) \+ \O_{\T}(-2\H+b_2\F) & \hpil & \O_{\T} \hpil \O_{S'} 
\hpil 0,
\end{eqnarray*} 
with $(b_1,b_2)=(1,1)$ or $(2,0)$ for the scroll types $(1,1,1,1)$ and 
$(2,1,1,0)$ and $(b_1,b_2)=(2,0)$ for the scroll types $(2,2,0,0)$ and $(3,1,0,0)$.

In the cases $(c,D^2)=(2,2)$ then $S''$ is cut out in $\T_0$ by a section (which is constant on the fibers of $i$) of
\[ \oplus_{i=1}^4 \O_{\T_0}(2\H-b_i\F), \]
where $(b_1,b_2,b_3,b_4)=(1,1,1,0)$ or $(2,1,0,0)$ for the type $(1,1,1,0,0)$,
and $(2,1,0,0)$ for the type $(2,1,0,0,0)$.

\vspace{.4cm}
{\bf Comments on the types of $L$:}
\begin{itemize}
 \item[(i)] \hs $L \sim 4D + 3\Gamma + 2\Gamma_1 + \Gamma_2$, with the following 
configuration:
\[ \xymatrix{
{D} \ar@{-}[r] & {\Gamma} \ar@{-}[r]  & {\Gamma_1} \ar@{-}[r]  & {\Gamma_2} 
} \]
\item[(ii)] \hs The number of moduli of the case $\{2,0\}^c$ is $15$, with
  mildest singularity $A_3$.

\item[(iii)] \hs In the case $\{2,0\}^b$ we have $L \sim 3D + 2\Gamma_1 +2\Gamma_2+ 
   \Gamma_1' + \Gamma_2'$, with the following 
configuration:
\[ \xymatrix{
{D} \ar@{-}[dr] \ar@{-}[r] & {\Gamma_1} \ar@{-}[r]  & {\Gamma_1'} \\
& {\Gamma_2}  \ar@{-}[r] & {\Gamma_2',}  
} \]
and in the case $\{2,0\}^c$ we have $L \sim 3D + 4\Gamma_0 + 3\Gamma_1 +2\Gamma_2+ 
   2\Gamma_3 + \Gamma_4$, with the following 
configuration: 
\[ \xymatrix{
{D} \ar@{-}[r] & {\Gamma_0} \ar@{-}[d] \ar@{-}[r] 
               & {\Gamma_1} \ar@{-}[r] & {\Gamma_3} \ar@{-}[d] \\
& {\Gamma_2} & & {\Gamma_4.}  
} \]
The mildest singularity of this latter case is $A_5$.
\item[(iv)] \hs The number of moduli of the case (E2) is 16, with mildest
  singularity $A_3$.
\end{itemize}

 \[ \mathbf{g=8} \]
The general projective model is a complete intersection of $5$ hyperplanes in 
$\Grass (V^6, 2) \sub \PP ^{14}$.

The others are as follows:
\vspace{.4cm}

\begin{tabular}{|c|c|c|c|c|c|} \hline
$c$ & $D^2$ & scroll type & $\#$ mod. & type of $L$ & sing.  \\ \hline
$1$ & $0$ & $(2,2,2)$ & $18$ & $h^0(L-2D)=3$, $h^0(L-3D)=0$ & sm.  \\ \hline
$1$ & $0$ & $(3,2,1)$ & $17$ & $h^0(L-2D)=3$, $h^0(L-3D)=1$ & sm.  \\ \hline
$1$ & $0$ & $(4,2,0)$ & $17$ & $\{1,0\}^a$, $A^2=4$ & $A_1$ \\  \hline
$2$ & $0$ & $(2,1,1,1)$ & $18$ & $h^0(L-2D)=1$ & sm. \\ \hline
$2$ & $0$ & $(2,2,1,0)$ & $17$ & $\{2,0\}^a$, $A^2=0$, $A \not >D$ & $A_1$    
                                                                   \\ \hline
$2$ & $0$ & $(3,1,1,0)$ & $15$ & $\{2,0\}^a$, $A^2=0$, $A >D$ $^{(i)}$ & $A_2$
                                          \\ \hline
$2$ & $0$ & $(3,2,0,0)$ & $16$ & $\{2,0\}^b$ or $\{2,0\}^c$ , $A^2=2$, $A >D$  
                                            $^{(ii)}$        & $2A_1$ \\ \hline
$2$ & $2$ & $(2,1,1,0,0)$ & $18$ & $h^0(L-2D)=1$ & sm.   \\ \hline
$3$ & $2$ & $(1,1,1,0,0,0)$ & $18$ & (CG1) or (CG2) $^{(iii)}$ & sm.   \\ \hline
$3$ & $2$ & $(2,1,0,0,0,0)$ & $16$ & (CG3)-(CG7) $^{(iv)}$ &
$3A_1$   \\ \hline 
\end{tabular}

\vspace{.4cm}
In the cases $(c,D^2)=(1,0)$ $\O_{S'}$ has the following 
$\O_{\T}$-resolution:
\[
  0 \hpil \O_{\T}(-3\H+4\F) \hpil \O_{\T} \hpil \O_{S'} \hpil 0.
\]

In the cases $(c,D^2)=(2,0)$ $\O_{S'}$ has the following 
$\O_{\T}$-resolution:
\begin{eqnarray*}
  0 \hpil \O_{\T}(-4\H+(g-1)\F)  & \hpil &  \\
\O_{\T}(-2\H+b_1\F) \+ \O_{\T}(-2\H+b_2\F) & \hpil & \O_{\T} \hpil \O_{S'} 
\hpil 0,
\end{eqnarray*} 
with $(b_1,b_2)=(2,1)$, except for the type $(3,2,0,0)$, where 
$(b_1,b_2)=(3,0)$. 
In this latter case, $S$ also contains a different perfect Clifford divisor 
(by the footnote on page \pageref{footpage}), so $S'$ can also be described as 
for the case $(c,D^2)=(2,2)$.

In the cases $(c,D^2)=(2,2)$ then $S''$ is cut out in $\T_0$ by a section (which is constant on the fibers of $i$) of
\[ \oplus_{i=1}^4 \O_{\T_0}(2\H-b_i\F), \]
where $(b_1,b_2,b_3,b_4)=(3,2,0,0), (3,1,1,0)$ or $(2,2,1,0)$.

In the cases (CG1) and (CG2) then $S''$ is cut out in $\T_0$ by a section of
\[
 \begin{array}{lllllllr}
 \O_{\T_0}(2\H-2\F)  & \oplus & \O_{\T_0}(2\H-\F)^2 & \oplus & 
 \O_{\T_0}(2\H)^5 & \mbox{or} & & \\
 \O_{\T_0}(2\H-\F)^4 & \oplus & \O_{\T_0}(2\H)^4  & \mbox{or} & & & & \\
 \O_{\T_0}(2\H-2\F) & \oplus & \O_{\T_0}(2\H-\F)^3 & \oplus & 
 \O_{\T_0}(2\H)^3 & \oplus & \O_{\T_0}(2\H+\F) & 
\end{array}
\]
(which is constant on the fibers of $i$ in the first two cases).
In the cases (CG3)-(CG7) then $S''$ is cut out in $\T_0$ by a section (which is constant on the fibers of $i$) of
\[
 \begin{array}{lllllr}
 \O_{\T_0}(2\H-2\F)  & \oplus & \O_{\T_0}(2\H-\F)^2 & \oplus & 
 \O_{\T_0}(2\H)^5 &  
\end{array}
\]

\vspace{.4cm}
{\bf Comments on the types of $L$:}
\begin{itemize}
 \item[(i)] \hs Here there are two subcases, one of them is: $L \sim 3D + 2\Gamma + \Gamma' + \Gamma_1 + \Gamma_2$, with the following configuration:
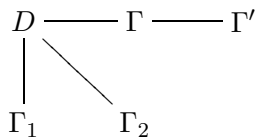
\[ \xymatrix{
{D} \ar@{-}[d] \ar@{-}[dr] \ar@{-}[r] & {\Gamma} \ar@{-}[r]  & {\Gamma'} \\
{\Gamma_1}   & {\Gamma_2} &
} \]
The number of moduli in this subcase is $15$, with
  mildest singularity $A_2$.

In the other subcase $L \sim 3D + 2\Gamma + \Gamma'+ 
   2\Gamma_0 + 2\Gamma_1 +  \cdots  + 2\Gamma_N + \Gamma_{N+1} + \Gamma_{N+2}$, for $N \geq 0$ (in general $N=0$) with the following 
configuration:
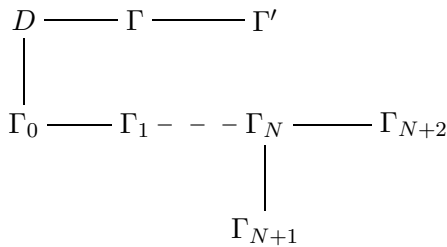
\[ \xymatrix{
{D} \ar@{-}[d] \ar@{-}[r] & {\Gamma} \ar@{-}[r]  & {\Gamma'} & \\
{\Gamma_0} \ar@{-}[r]   & {\Gamma_1} \ar@{--}[r] & 
            {\Gamma_N} \ar@{-}[d] \ar@{-}[r] &  {\Gamma_{N+2}}    \\
                      &  & {\Gamma_{N+1}} &
} \]
The number of moduli in this subcase is $14$, with
  mildest singularity $A_2+2A_1$. 
\item[(ii)] \hs The number of moduli of the case $\{2,0\}^c$ is $15$, with
  mildest singularity $A_3$.
\item[(iii)] \hs The number of moduli of the case (CG2) is $17$, with
  mildest singularity $A_1$.
\item[(iv)] \hs The number of moduli of the cases (CG4), (CG5), (CG6) and (CG7)
  are $15$, $14$, $13$ and $12$ respectively, with
  mildest singularities $A_1+A_3$, $A_5$, $D_6$ and $E_7$ respectively.
\end{itemize}

 \[ \mathbf{g=9} \]
The general projective model is a complete intersection of $4$ hyperplanes in 
$\Sigma ^6 _{16}$, as described in the beginning of Section \ref{bncl}.

The others are as follows: \label{tablepage}
\vspace{.4cm}

\begin{tabular}{|c|c|c|c|c|c|} \hline
$c$ & $D^2$ & scroll type & $\#$ mod. & type of $L$ & sing.  \\ \hline
$1$ & $0$ & $(3,2,2)$ & $18$ & $h^0(L-2D)=4$, $h^0(L-3D)=1$ & sm. \\  \hline
$1$ & $0$ & $(3,3,1)$ & $17$ & 
                   $h^0(L-iD)=4,2,0$, for $i=2,3,4$ & sm.\\  \hline
$1$ & $0$ & $(4,2,1)$   & $16$ & $h^0(L-iD)=4,2,1$, for $i=2,3,4$$^{(i)}$  & $A_1$   \\  \hline
$1$ & $0$ & $(5,2,0)$   & $17$ & $\{1,0\}^a$, $A^2=6$  & $A_1$ \\ \hline
$2$ & $0$ & $(2,2,1,1)$   & $18$ &  $h^0(L-2D)=2$, $h^0(L-3D)=0$  & sm.
                 \\ \hline
$2$ & $0$ & $(3,1,1,1)$   & $15$ &  $h^0(L-2D)=2$, $h^0(L-3D)=1$  & sm.
                                       \\ \hline
$2$ & $0$ & $(2,2,2,0)$   & $17$ & $\{2,0\}^a$, $A^2=2$, $A \not >D$ 
             &  $A_1$   \\ \hline
$2$ & $0$ & $(3,2,1,0)$   & $16$ & $\{2,0\}^a$, $A^2=2$, $A  >D$ $^{(ii)}$ 
            & $A_1$  \\ \hline
$2$ & $0$ & $(3,2,1,0)$   & $16$ & $\{2,0\}^a$, $A^2=2$, $A >D$ $^{(iii)}$ 
              & $A_2$  \\ \hline
$2$ & $0$ & $(4,2,0,0)$   & $17$ & $\{2,0\}^b$ ($L \sim 4D +2\Gamma_1
+ 2\Gamma_2$) or $\{2,0\}^c$$^{(iv)}$ 
              & $2A_1$  \\ \hline
$2$ & $2$ & $(2,2,1,0,0)$  & $18$ &  $h^0(L-2D)=2$, 
            $h^0(L-3D)=0$ &  sm. \\ \hline
$2$ & $4$ & $(2,1,1,0,0,0)$ & $19$ & $L \sim 2D$ & sm.   \\ \hline
$3$ & $0$ & $(1,1,1,1,1)$ & $18$ & $h^0(L-2D)=0$ & sm.   \\ \hline
$3$ & $0$ & $(2,1,1,1,0)$ & $17$ & $\{3,0\}^a$, $A^2=-2$  & $A_1$ \\ \hline
$3$ & $0$ & $(2,2,1,0,0)$ & $16$ & $\{3,0\}^b$ or $\{3,0\}^c$,
$A^2=0$, $A \not >D$ $^{(v)}$ & $2A_1$
                          \\ \hline
$3$ & $0$   & $(3,1,1,0,0)$   & $14$ & $\{3,0\}^b$ or $\{3,0\}^c$,
$A^2=0$, $A  >D$ $^{(vi)}$ & $2A_2$      \\ \hline
$3$ & $2$   & $(1,1,1,1,0,0)$   & $18$ &  $h^0(L-2D)=0$ & sm. \\ \hline
$3$ & $2$   & $(2,1,1,0,0,0)$   & $17$ & $\{3,2\}^a$  & $A_1$   \\ \hline
\end{tabular}

\vspace{.4cm}
In the cases $(c,D^2)=(1,0)$ $\O_{S'}$ has the following 
$\O_{\T}$-resolution:
\[
  0 \hpil \O_{\T}(-3\H+5\F) \hpil \O_{\T} \hpil \O_{S'} \hpil 0.
\]

In the cases $(c,D^2)=(2,0)$ $\O_{S'}$ has the following 
$\O_{\T}$-resolution:
\begin{eqnarray*}
  0 \hpil \O_{\T}(-4\H+(g-1)\F)  & \hpil &  \\
\O_{\T}(-2\H+b_1\F) \+ \O_{\T}(-2\H+b_2\F) & \hpil & \O_{\T} \hpil \O_{S'} 
\hpil 0,
\end{eqnarray*}
with $(b_1,b_2)=(2,2)$ or $(3,1)$ for the scroll type $(2,2,1,1)$;
$(b_1,b_2)=(2,2)$ for the scroll types $(3,1,1,1)$, $(2,2,2,0)$ and 
 $(3,2,1,0)$ ($A_2$-sing.); 
$(b_1,b_2)=(3,1)$ for the scroll type $(3,2,1,0)$ ($A_1$-sing.);
and $(b_1,b_2)=(4,0)$ for the scroll type $(4,2,0,0)$.
 In this latter case $S$ also contains a different perfect Clifford divisor 
(by the footnote on page \pageref{footpage}), so $S'$ can also be
 described as in the case $(c,D^2)=(2,2)$ (with perfect Clifford
 divisor $2D+\Gamma_1 + \Gamma_2$ or $D +\Delta_0$).
(The $\{2,0\}^c$ case of the table corresponds to ($L \sim 4D +2\Delta_0)$).
In the case $(c,D^2)=(2,2)$ then $S''$ is cut out in $\T_0$ by a  section (which is constant on the fibers of $i$) of
\[ \O_{\T_0}(2\H-3\F) \oplus \O_{\T_0}(2\H-2\F)^2 \oplus \O_{\T_0}(2\H). \]

In the case $(c,D^2)=(2,4)$ then $S''$ is cut out in $\T_0$ by a section (which is constant on the fibers of $i$) of:
\[ \O_{\T_0}(2\H-2\F)^3   \oplus  \O_{\T_0}(2\H-\F)^2  \oplus \O_{\T}(2\H)^2. \] 

We also have that $S'$ is the $2$-uple embedding of the quartic $\varphi_D(S)$ if and only if $D$ is not hyperelliptic. If $D$ is hyperelliptic,  then there is an elliptic pencil $|E|$ such that $E.D=2$. Then $E$ is also a free Clifford divisor for $L$ and defines a scroll $\T(2, E)$ containing $S'$. The $\O_{\T(2, E)}$-resolution of $\O_{S'}$ is given in Proposition \ref{non2-uple}. 

In the cases $(c,D^2)=(3,0)$ we have an $\O_{\T}$-resolution of $\O_{S'}$ of the following type:
\[
 \begin{array}{ccccccccc}
0 & \hpil & \O_{\T}(-5\H+8\F)  & \hpil & \+ _{i=1}^5 
\O_{\T}(-3\H+a_i\F)   \\
  & \hpil & \+ _{i=1}^5 \O_{\T}(-2\H + b_i\F) & \hpil & 
\O_{\T}  \hs \hpil  \hs \O_{S'}  \hs \hpil  \hs 0, 
 \end{array}
\]
with $a_i=3-b_i$, for all $i$.
For the smooth scroll type $(1,1,1,1,1)$ we have $(b_1,b_2,b_3,b_4)=(2,1,1,1,1)$ or $(2,2,2,0,0)$.
For the scroll type $(2,1,1,1,0)$ we have  $(b_1,b_2,b_3,b_4)=
 (2,1,1,1,1)$, $(2,2,1,1,0)$ or $(2,2,2,0,0)$.
For the remaining two singular scroll types we have
$(b_1,b_2,b_3,b_4)=(2,2,1,1,0)$ or $ (2,2,2,0,0)$.
 
In the case $(c,D^2)=(3,2)$ with scroll type $(1,1,1,1,0,0)$ then
$S''$ is cut out in $\T_0$ by a section (which is
constant on the fibers of $i$) of:
\[
 \begin{array}{lllllr}
 \O_{\T_0}(2\H-2\F)  & \oplus & \O_{\T_0}(2\H-\F)^5 & \oplus & 
              \O_{\T_0}(2\H)^2 &  \mbox{or} \\
 \O_{\T_0}(2\H-\F)^7  & \oplus & \O_{\T_0}(2\H) & & & \\
 \end{array}
\]

In the case $(c,D^2)=(3,2)$ with scroll type $(2,1,1,0,0,0)$ then
$S''$ is cut out in $\T_0$  by a  section (which is
constant on the fibers of $i$) of:
\[
 \begin{array}{lllllllr}
 \O_{\T_0}(2\H-2\F)^3  & \oplus & \O_{\T_0}(2\H-\F) & \oplus & \O_{\T_0}(2\H)^4 
                        & & & \mbox{or}\\
 \O_{\T_0}(2\H-2\F)^2  & \oplus & \O_{\T_0}(2\H-\F)^3 & \oplus & 
 \O_{\T_0}(2\H)^3. & & &   
 \end{array}
\]

\vspace{.4cm}
{\bf Comments on the types of $L$:}
\begin{itemize}
 \item[(i)] \hs $L \sim 4D + 2\Gamma_1 + \Gamma_2 + \Gamma_3$, with 
  the following configuration:
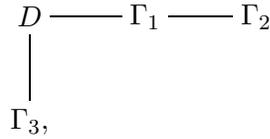
\[ \xymatrix{
{D} \ar@{-}[d] \ar@{-}[r] & {\Gamma_1} \ar@{-}[r] 
                                               & {\Gamma_2} \\
{\Gamma_3,}   & & &
} \]
or 
$L \sim 4D + 3\Gamma_1 +  \cdots  + 3\Gamma_N + 2\Gamma_{N+1} +
\Gamma_{N+2} + \Gamma_{N+3} $, for $N \geq 1$ (in general $N=0$) with 
  the following configuration:
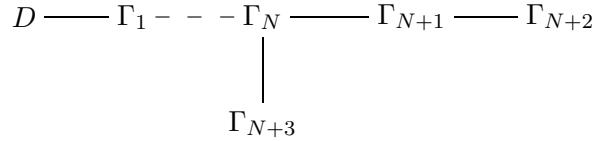
\[ \xymatrix{
{D} \ar@{-}[r] & {\Gamma_1} \ar@{--}[r] & {\Gamma_N} \ar@{-}[d]
             \ar@{-}[r] & {\Gamma_{N+1}} \ar@{-}[r] & {\Gamma_{N+2}}
             \\ 
             &  & {\Gamma_{N+3}}   & & &
} \]

\item[(ii)] \hs $L \sim 3D + 2\Gamma + \Gamma_1 + \Gamma_2$, with 
  the following configuration:
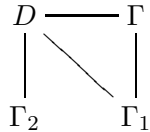
\[ \xymatrix{
{D} \ar@{-}[d] \ar@{-}[dr] \ar@{-}[r] & {\Gamma}  \ar@{-}[d] \\
{\Gamma_2}   & {\Gamma_1} 
} \]

\item[(iii)] \hs $L \sim 3D + 2\Gamma + \Gamma_1 + \Gamma_2$, with 
  the following configuration:
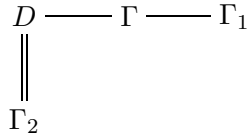
\[ \xymatrix{
{D} \ar@{=}[d] \ar@{-}[r] & {\Gamma}  \ar@{-}[r] & {\Gamma_1} \\
{\Gamma_2}   & & 
} \]
or 
$L \sim 3D + 2\Gamma + \Gamma_1 + \Gamma_2 +  \cdots  + \Gamma_{N+3}$, for 
$N \geq 0$ (in general $N=0$) with 
  the following configuration:
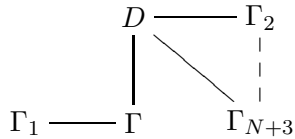
\[ \xymatrix{
          & {D} \ar@{-}[d] \ar@{-}[dr] \ar@{-}[r] & {\Gamma_2} \ar@{--}[d] \\
{\Gamma_1} \ar@{-}[r] & {\Gamma}  & {\Gamma_{N+3}}    
} \]

\item[(iv)] \hs The number of moduli of the case $\{2,0\}^c$ is $16$, with
  mildest singularity $A_3$.

\item[(v)] \hs  The number of moduli of the case $\{3,0\}^c$ is $15$, with
  mildest singularity $A_3$.

\item[(vi)] \hs In the case $\{3,0\}^b$ we have $L \sim 3D + 2\Gamma_1 + 2\Gamma_2 + \Gamma_3  + \Gamma_4 + \Gamma_5$, with 
  the following configuration:
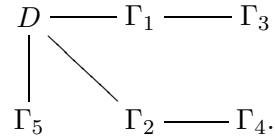
\[ \xymatrix{
{D} \ar@{-}[d] \ar@{-}[dr] \ar@{-}[r] & {\Gamma_1} \ar@{-}[r] & {\Gamma_3} \\
{\Gamma_5}                            & {\Gamma_2} \ar@{-}[r] & {\Gamma_4}. 
} \]
or $L \sim 3D + 2\Gamma_1 + \Gamma_1' + 3\Gamma_2 +  \cdots  + 3\Gamma_N  + 
2\Gamma_{N+1} + \Gamma_{N+2} + \Gamma_{N+3}$, with 
  the following configuration:
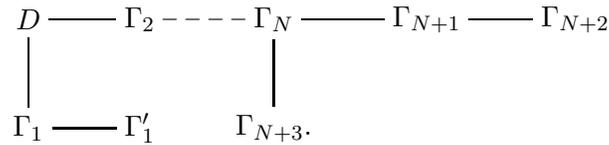
\[ \xymatrix{
{D} \ar@{-}[d] \ar@{-}[r] & {\Gamma_2} \ar@{--}[r] & {\Gamma_N} \ar@{-}[d] 
             \ar@{-}[r] & {\Gamma_{N+1}} \ar@{-}[r] & {\Gamma_{N+2}} \\
{\Gamma_1 } \ar@{-}[r] & {\Gamma_1'} & {\Gamma_{N+3}.}   & &                
} \]

In the case $\{3,0\}^c$ we have $L \sim 3D + \Gamma + 4\Gamma_0 + 3\Gamma_1  + 2\Gamma_2 + 2\Gamma_3
+ \Gamma_4$, 
  with the following configuration: 
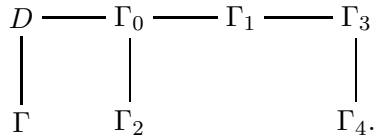
\[ \xymatrix{
{D} \ar@{-}[d] \ar@{-}[r] & {\Gamma_0} \ar@{-}[d] \ar@{-}[r] & {\Gamma_1} 
                                       \ar@{-}[r] & {\Gamma_3} \ar@{-}[d] \\
{\Gamma}                  & {\Gamma_2} & & {\Gamma_4}. 
} \]
The mildest singularity of this case is $A_5$.
\end{itemize}

 \[ \mathbf{g=10} \]

The general projective model is a complete intersection of $2$ hyperplanes in 
the homogeneous variety $\Sigma ^5 _{18}$, as described in Section \ref{bncl}.

The others are as follows:
\vspace{.4cm}

\begin{tabular}{|c|c|c|c|c|c|} \hline
$c$ & $D^2$ & scroll type & $\#$ mod. & type of $L$ & sing.  \\ \hline
$1$ & $0$   & $(3,3,2)$   & $18$ & $h^0(L-2D)=5$, $h^0(L-3D)=2$,  
             $h^0(L-4D)=0$  & sm. \\  \hline
$1$ & $0$   & $(4,2,2)$   & $16$ & $h^0(L-2D)=5$, $h^0(L-3D)=2$, $h^0(L-4D)=1$
             & sm. \\  \hline
$1$ & $0$   & $(4,3,1)$   & $17$ & $h^0(L-2D)=5$, $h^0(L-3D)=3$, $h^0(L-4D)=1$
            & sm.  \\  \hline
$1$ & $0$   & $(5,2,1)$   & $16$ & $h^0(L-5D)=1$ $^{(i)}$ & $A_2$ \\  \hline
$1$ & $0$   & $(6,2,0)$   & $18$ & $\{1,0\}^a$,
              $L \sim 6D+3\Gamma$   & $A_1$ \\  \hline

$2$ & $0$   & $(2,2,2,1)$   & $18$ &  $h^0(L-2D)=3$, $h^0(L-3D)=0$ &  sm.
\\ \hline
$2$ & $0$   & $(3,2,1,1)$   & $17$ &  $h^0(L-2D)=3$, $h^0(L-3D)=1$ $^{(ii)}$
      & $A_1$  \\ \hline
$2$ & $0$   & $(3,2,1,1)$   & $16$ &  $h^0(L-2D)=3$, $h^0(L-3D)=1$ $^{(iii)}$
 & sm.  \\ \hline
$2$ & $0$   & $(3,2,2,0)$  & $17$ & $\{2,0\}^a$, 
             $A^2=4$, $h^0(A-D)=1$, $A \not >2D$ $^{(iv)}$  & $A_1$  \\ \hline
$2$ & $0$   & $(3,3,1,0)$   & $16$ & $\{2,0\}^a$, 
          $A^2=4$, $h^0(A-D)=2$, $A \not >2D$ $^{(v)}$  & $A_2$   \\ \hline
$2$ & $0$   & $(4,2,1,0)$   & $16$ &   
               $\{2,0\}^a$, $A^2=4$, $A >2D$ $^{(vi)}$  &  $2A_1$ \\ \hline
$2$ & $2$   & $(3,2,1,0,0)$   & $19$ &  $L \sim 3D$    & sm.   \\ \hline

$3$ & $0$   & $(2,1,1,1,1)$   & $18$ &  $h^0(L-2D)=0$   & sm.  \\ \hline
$3$ & $0$   & $(2,2,1,1,0)$   & $17$ & $\{3,0\}^a$, $A^2=0$ & $A_1$  \\ \hline
$3$ & $0$   & $(2,2,2,0,0)$   & $16$ & $\{3,0\}^b$ or $\{3,0\}^c$,
$A^2=2$, 
$A \not >D$ $^{(vii)}$
              & $2A_1$    \\ \hline
$3$ & $0$   & $(3,2,1,0,0)$   & $15$ & $\{3,0\}^b$, $A^2=2$, $A  >D$ $^{(viii)}$
     & $2A_1$  \\ \hline
$3$ & $0$   & $(3,2,1,0,0)$   & $15$ & $\{3,0\}^b$, $A^2=2$, $A >D$ $^{(ix)}$ 
    & $A_1+A_2$ \\ \hline
$3$ & $0$   & $(3,2,1,0,0)$   & $14$ & $\{3,0\}^c$, $A^2=2$, $A  >D$ $^{(x)}$
& $A_3$ \\ \hline
$3$ & $0$   & $(3,2,1,0,0)$   & $14$ & $\{3,0\}^c$, $A^2=2$, $A  >D$ $^{(xi)}$
& $A_4$ \\ \hline

$3$ & $2$ & $(2,1,1,1,0,0)$   & $18$ & $h^0(L-2D)=1$    & sm. \\ \hline
$3$ & $4$ & $(2,1,1,0,0,0,0)$ & $18$ & (E0) & $A_1$ \\ \hline
$4$ & $2$ & $(1,1,1,1,0,0,0)$ & $18$ & (CG1)' or (CG2)' $^{(xii)}$ & sm.   \\ \hline
$4$ & $2$ & $(2,1,1,0,0,0,0)$ & $16$ & (CG3)' or (CG4)' $^{(xiii)}$ & $3A_1$   \\ \hline
\end{tabular}

\vspace{.4cm}
In the cases $(c,D^2)=(1,0)$ $\O_{S'}$ has the following 
$\O_{\T}$-resolution:
\[
  0 \hpil \O_{\T}(-3\H+6\F) \hpil \O_{\T} \hpil \O_{S'} \hpil 0.
\]

In the cases $(c,D^2)=(2,0)$ $\O_{S'}$ has the following 
$\O_{\T}$-resolution:
\begin{eqnarray*}
  0 \hpil \O_{\T}(-4\H+(g-1)\F)  & \hpil &  \\
\O_{\T}(-2\H+b_1\F) \+ \O_{\T}(-2\H+b_2\F) & \hpil & \O_{\T} \hpil \O_{S'} 
\hpil 0,
\end{eqnarray*}
with $(b_1,b_2)=(3,2)$ or $(4,1)$ for the scroll type $(4,2,1,0)$,
$(b_1,b_2)=(3,2)$ for the scroll types $(2,2,2,1)$, $(3,2,1,1)$ (smooth), 
$(3,2,2,0)$ and $(3,3,1,0)$, and 
$(b_1,b_2)=(4,1)$ for the scroll type $(3,2,1,1)$ ($A_1$-sing.).

In the case $(c,D^2)=(2,2)$ then $S''$ is cut out in $\T_0$ by a section (which is constant on the fibers of $i$) of
\[ \O_{\T_0}(2\H-4\F) \oplus \O_{\T_0}(2\H-3\F) \oplus \O_{\T_0}(2\H-2\F)
\oplus \O_{\T_0}(2\H). \] or of
\[\O_{\T_0}(2\H-3\F)^3  \oplus  \O_{\T_0}(2\H).\]

In the cases $(c,D^2)=(3,0)$ we have an $\O_{\T}$-resolution of $\O_{S'}$ of the following type:
\[
 \begin{array}{ccccccccc}
0 & \hpil & \O_{\T}(-5\H+9\F)  & \hpil & \+ _{i=1}^5 
\O_{\T}(-3\H+a_i\F)  \\
  & \hpil & \+ _{i=1}^5 \O_{\T}(-2\H + b_i\F) & \hpil & 
\O_{\T} \hs \hpil \hs  \O_{S'} \hs \hpil \hs  0, 
 \end{array}
\]
with $a_i=4-b_i$ for all $i$.

For the smooth scroll type $(2,1,1,1,1)$ we have $(b_1,b_2,b_3,b_4)=(2,2,2,1,1)$.
For the scroll type $(2,2,2,0,0)$ we have $(b_1,b_2,b_3,b_4)=(2,2,2,1,1)$ or \linebreak $(2,2,2,2,0)$.
For the remaining two singular scroll types we have \linebreak $(b_1,b_2,b_3,b_4)= (3,3,2,0,0)$, $(3,3,1,1,0)$, $(3,2,2,1,0)$, $(2,2,2,1,1)$  or \linebreak  $(2,2,2,2,0)$.

In the case $(c,D^2)=(3,2)$ then $S''$ is cut out in $\T_0$ by a section (which is constant on the fibers of $i$) of
\[ \O_{\T_0}(2\H-2\F)^3 \+ \O_{\T_0}(2\H-\F)^4 \+ 2\O_{\T_0}(\H). \]

In the case $(c,D^2)=(3,4)$ then $S''$ is cut out in $\T_0$ by a
section (which is constant on the fibers of $i$) of
\[ \O_{\T_0}(2\H-2\F)^3  \oplus  \O_{\T_0}(2\H-\F)^4  \oplus  \O_{\T}(2\H)^5. \]

In the cases (CG1)' and (CG2)' then $S''$ is cut out in $\T_0$ by a section of 
\[
 \begin{array}{lllllllr}
 \O_{\T_0}(2\H-2\F)  & \oplus & \O_{\T_0}(2\H-\F)^7 & \oplus & 
 \O_{\T_0}(2\H)^5 & \mbox{or} &  &   \\
 \O_{\T_0}(2\H-\F)^9 & \oplus & \O_{\T_0}(2\H)^4  & \mbox{or} & & & & \\
 \O_{\T_0}(2\H-2\F) & \oplus & \O_{\T_0}(2\H-\F)^8 & \oplus & 
 \O_{\T_0}(2\H)^3 & \oplus & \O_{\T_0}(2\H+\F), & 
\end{array}
\]
where the last option occurs only for the case (CG2)' (the section is
constant on the fibers of $i$ in the first two cases).

In the cases (CG3)', (CG4)' then $S''$ is cut out in $\T_0$ by a
section (which is constant on the fibers of $i$) of 
\[
 \begin{array}{llllllllr}
 \O_{\T_0}(2\H-2\F)^2  & \oplus & \O_{\T_0}(2\H-\F)^5 & \oplus & 
              \O_{\T_0}(2\H)^6 &  \mbox{or} & \\
 \O_{\T_0}(2\H-2\F)^3 & \oplus & \O_{\T_0}(2\H-\F)^3 & \oplus & 
              \O_{\T_0}(2\H)^7. &  & & & \\

\end{array}
\]

\vspace{.4cm}
{\bf Comments on the types of $L$:}
\begin{itemize}
 \item[(i)] \hs $L \sim 5D + 3\Gamma_1 + 2\Gamma_2 + \Gamma_3$, with 
           the following configuration:
\[ \xymatrix{
{D} \ar@{-}[r] & {\Gamma_1} \ar@{-}[r] & {\Gamma_2} \ar@{-}[r] & {\Gamma_3} 
} \]

 \item[(ii)] \hs $L$ is in general of the form $L \sim 3D + 2\Gamma_1 + \Gamma_2$, 
            with the following configuration:
\[ \xymatrix{
{D} \ar@{=}[r] & {\Gamma_1} \ar@{-}[r] & {\Gamma_2}  
} \]

 \item[(iii)] \hs $L$ is in general of the form $L \sim 3D + \Gamma_1 + \Gamma_2 + 
\Gamma_3$, with the following configuration:
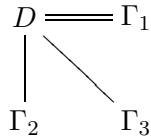
\[ \xymatrix{
{D} \ar@{-}[d] \ar@{-}[dr] \ar@{=}[r] & {\Gamma_1} \\
{\Gamma_2}                         & {\Gamma_3}
} \]
 
\item[(iv)] \hs $L$ is in general of the form $L \sim 3D + 2\Gamma + \Gamma_1$, 
with the following configuration:
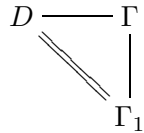
\[ \xymatrix{
{D} \ar@{=}[dr] \ar@{-}[r] & {\Gamma} \ar@{-}[d] \\
                            & {\Gamma_1}
} \]

\item[(v)] \hs $L \sim 3D + E+ 2\Gamma + \Gamma_1$, where $E$ is a smooth elliptic 
curve, with the following configuration:
\[ \xymatrix{
{D} \ar@{=}[d] \ar@{-}[r] & {\Gamma}  \ar@{-}[r] & {\Gamma_1} \\
 {E} & &
} \]

\item[(vi)] \hs $L$ is in general of the form $L \sim 4D + 2\Gamma + 2\Gamma_1 + 
\Gamma_2$, with the following configuration:
\[ \xymatrix{
{D} \ar@{-}[d] \ar@{-}[r] & {\Gamma_1}  \ar@{-}[r] & {\Gamma_2} \\
 {\Gamma} & &
} \]

\item[(vii)] \hs The number of moduli of the case $\{3,0\}^c$ is $15$, with
  mildest singularity $A_3$.

\item[(viii)] \hs $L$ is in general of the form $L \sim 3D + 2\Gamma_1 + 2\Gamma_2 + \Gamma_3 + \Gamma_4$, 
with the following configuration:
\[ \xymatrix{
             & {D} \ar@{-}[dl] \ar@{-}[d] \ar@{-}[r] & {\Gamma_1} \ar@{-}[d] \\
 {\Gamma_3} & {\Gamma_2} \ar@{-}[r] & {\Gamma_4}
} \]

\item[(ix)] \hs $L$ is in general of the form $L \sim 3D + 2\Gamma_1 + 2\Gamma_2 + \Gamma_3 + \Gamma_4$, 
with the following configuration:
\[ \xymatrix{
{D} \ar@{-}[d] \ar@{-}[dr] \ar@{-}[r] & {\Gamma_1} \ar@{-}[r] & {\Gamma_4} \\
{\Gamma_3} \ar@{-}[r] & {\Gamma_2}  &
} \]

\item[(x)] \hs $L$ is in general of the form $L \sim 3D + 4\Gamma_0 + 2\Gamma_1 + 
2\Gamma_2 + \Gamma_3 + \Gamma_4$, 
with the following configuration:
\[ \xymatrix{
{D} \ar@{-}[d] \ar@{-}[r] & {\Gamma_0} \ar@{-}[d] \ar@{-}[dr] \ar@{-}[r] 
       & {\Gamma_3} \\
{\Gamma_4} & {\Gamma_1} & {\Gamma_2}
} \]

\item[(xi)] \hs $L$ is in general of the form $L \sim 3D + 4\Gamma_0 + 3\Gamma_1 + 
2\Gamma_2 + 2\Gamma_3 + \Gamma_4$, 
with the following configuration:
\[ \xymatrix{
{D} \ar@{-}[d] \ar@{-}[r] & {\Gamma_0} \ar@{-}[dr] \ar@{-}[r] & {\Gamma_2} \\
{\Gamma_4} \ar@{-}[r] & {\Gamma_3} \ar@{-}[r] & {\Gamma_1}
} \]

\item[(xii)] \hs The number of moduli of the case (CG2)' is $17$, with mildest singularity $A_1$.
\item[(xiii)] \hs The number of moduli of the case (CG4)' is $15$, with mildest singularity 
$A_1 + A_3$.
\end{itemize}

\printindex

%%%%%%%%%%%%%%%%%%%%%%%%%%%%%%%%%%%%%%%%%%%%%%%%%%%%%%%%%%%%%%%%%%%%%%


\begin{thebibliography}{[A-C-G-H]}




\bibitem[A-C-G-H]{acgh} E.~Arbarello, M.~Cornalba, P.~A.~Griffiths,
  J.~Harris, \textit{Geometry of Algebraic Curves, Volume I}, number 267
  in \textit{Grundlehren der Mathematischen Wissenschaften}, Springer
  Verlag, Berlin, Heidelberg, New York, Tokyo (1985).

\bibitem[Ar]{Ar} M.~Artin, \textit{Some numerical criteria for contractability
of curves on algebraic surfaces}, American Journal of Math. \textbf{84}, 
485--496 (1962).

\bibitem[Ba]{bal} E.~Ballico,  \textit{A remark on linear series on general
$k$-gonal curves}, Boll. Un. Mat. Ital. A (7) \textbf{3}, no. 2, 195--197 (1989).

\bibitem[B-P-V]{BPV} W.~Barth, C.~Peters, A.~Van de Ven \textit{Compact Complex Surfaces}, Springer Verlag, Berlin, Heidelberg, New York, Tokyo (1984).

\bibitem[Be]{Bea} A.~Beauville, \textit{Determinantal hypersurfaces}, Michigan Math. Journal (Special volume dedicated to W. Fulton) \textbf{48} (1997).

\bibitem[B-S]{BS2} M.~Beltrametti, A.~J.~Sommese, \textit{Zero cycles and
    $k$-th order embeddings of smooth projective surfaces}, in
    \textit{Projective surfaces and their classification}, Cortona
    proceedings Symp. Math. INDAM 32, pp.~33-48, New York, London:
    Academic press 1988.

\bibitem[Bo]{bog} F.~Bogomolov, \textit{Holomorphic tensors and vector bundles
on projective varieties}, Math. USSR Izv. \textbf{13}, 499-555 (1979).

\bibitem[Br]{B} J.~Brawner, \textit{Tetragonal Curves, Scrolls and
$K3$-surfaces}, Trans. of the Am. Mat. Soc. \textbf{349}, No. 8, 3075--3091 (1997).



\bibitem[B-E]{B-E} D.~Buchsbaum, D.~Eisenbud \textit{Algebra structures for
finite free resolutions and some structure theorems for ideals in codimension 
3}, Am. J. Math.  \textbf{99}, 447--485 (1977).

\bibitem[C-F]{CF} F.~Catanese, M.~Franciosi, \textit{Divisors of small genus 
on algebraic surfaces and projective embeddings}, in \textit{Proceedings of 
the Hirzebruch 65 Conference on Algebraic Geometry (Ramat Gan, 1993)}, 
109--140, Israel Math. Conf. Proc. \textbf{9}, Bar-Ilan Univ., Ramat Gan 
(1996). 

\bibitem[C-P]{cp} C.~Ciliberto, G.~Pareschi, \textit{Pencils of minimal degree on curves on a $K3$ surface}, J. reine angew. Math. \textbf{460}, 15--36 (1995).

\bibitem[Cl]{C83} H.~Clemens,  \textit{Homological equivalence,
    modulo algebraic equivalence, is not finitely generated}, 
Publ. Math. IHES \textbf{58}, 19--38 (1983).

\bibitem[C-M]{cm} M.~Coppens, G.~Martens, \textit{Secant Spaces and Clifford's theorem}, Comp. Math. \textbf{68}, 337--341 (1991).

\bibitem[Co]{cos} F.~R.~ Cossec, \textit{Projective models of Enriques surfaces}, Math. Ann. \textbf{265}, 283--334 (1983).

\bibitem[CD]{CD} F.~R.~ Cossec, I. ~Dolgachev,  \textit{Enriques
    surfaces I}, Progress in Mathematics \textbf{76}, Birkhauser,  (1989).

\bibitem[D-M]{DM} R.~Donagi, D.~R.~Morrison, \textit{Linear systems
    on $K3$ sections}, J. Diff. Geom. \textbf{29}, 49--64 (1989).

%\bibitem[Ei]{eishyp} D.~Eisenbud, \textit{Transcanonical embeddings of hyperelliptic curves}, 
%J. Pure Appl. Alg.
%\textbf{19}, 77--83 (1980).

\bibitem[E-H]{E-H} D.~Eisenbud, J.~Harris, \textit{On Varieties of Minimal 
Degree (A Centennial Account)}, Proceedings of Symposia in Pure Mathematics 
\textbf{46}, 3--13 (1987).

\bibitem[E-L-M-S]{elms} D.~Eisenbud, H.~Lange, G.~Martens, F.-O.~Schreyer, 
\textit{The Clifford dimension of a projective curve}, Comp. Math. 
\textbf{72}, 173--204 (1989).

\bibitem[E-J-S]{EJS} T. Ekedahl, T. Johnsen, D.E. Sommervoll
\textit{Isolated rational curves on K3-fibered Calabi-Yau threefolds}, 
Manuscripta Math. \textbf{99}, 111--133 (1999).

\bibitem[Gr]{gre} M.~Green, \textit{Koszul cohomology and the geometry of 
projective varieties}, J. Differ. Geom. \textbf{19}, 125--171 (1984).

\bibitem[G-L1]{glroma} M.~Green, R.~Lazarsfeld, \textit{A simple proof of Petri's theorem on canonical curves}, In: geometry Today, progress in Math., birkh\"auser (1986).

\bibitem[G-L2]{glprojnorm} M.~Green, R.~Lazarsfeld, \textit{On the projective normality of complete linear series on an algebraic curve}, Invent. Math. \textbf{83}, 73-90 (1986).

\bibitem[G-L3]{gl2} M.~Green, R.~Lazarsfeld, \textit{Some results on the syzygies of finite sets and algebraic curves}, Comp. Math. \textbf{67}, 301--314 (1988).

\bibitem[G-L4]{gl} M.~Green, R.~Lazarsfeld, \textit{Special divisors on
    curves on a $K3$ surface}, Invent. Math. \textbf{89}, 357--370 (1987).

\bibitem[G-L5]{glapp} M.~Green, R.~Lazarsfeld, \textit{The nonvanishing of certain cohomology groups}, Appendix to \cite{gre}.

\bibitem[G-H]{gh}  P.~Griffiths, J.~Harris, \textit{Principles of algebraic geometry}, Wiley Classics Library, John Wiley \& Sons, Inc., New York, 1994.

\bibitem[G-L-P]{glp} L.~Gruson, R.~Lazarsfeld, C. ~Peskine, \textit{On a theorem
of Castelnuovo and the equations defining space curves}, Invent. Math. 
\textbf{72}, 491--506 (1983).

\bibitem[Han]{han} G.M.~Hana, \textit{Projective Models of Enriques
    Surfaces in Scrolls}, Preprint, 2003.

\bibitem[Har]{H} J.~Harris, \textit{A Bound on the Geometric Genus of 
Projective Varieties}, Ann. Scuola Norm. Sup. Pisa Cl. Sci, \textbf{8} (4),
35--68 (1981).

\bibitem[H-M]{hm} J.~Harris, D.~Mumford, \textit{On the Kodaira dimension of the moduli space of curves}, 
Invent. Math. \textbf{67}, no. 1, 23--86 (1982).

\bibitem[Hrts]{Ha} R.~Hartshorne, \textit{Algebraic Geometry}, Graduate Texts
in Mathematics, \textbf{52}, Springer Verlag (1977).


\bibitem[I-R1]{ir1} A.~Iliev, K.~Ranestad, \textit{Canonical curves and varieties of sums of powers of cubic polynomials},  J. Algebra \textbf{246}, no. 1, 385--393 (2001). 
 
\bibitem[I-R2]{ir2} A.~Iliev, K.~Ranestad, \textit{$K3$ surfaces of genus 8 and varieties of sums of powers of cubic fourfolds}, Trans. Amer. Math. Soc. \textbf{353}, no. 4, 1455--1468 (2001).

\bibitem[J-K]{JoKn} T. ~Johnsen, A.~L.~Knutsen, \textit{Rational
 Curves in Calabi-Yau threefolds}, Comm. in Algebra \textbf{31} (8), 3917--3953 (2003).


\bibitem[Kn1]{kn4} A.~L.~Knutsen, \textit{Exceptional curves on Del Pezzo surfaces},
Mathematische Nachrichten {\bf 256}, 1, 58--81 (2003).


\bibitem[Kn2]{kn3} A.~L.~Knutsen, \textit{Gonality and Clifford index
    of curves on $K3$ surfaces}, Archiv der Mathematik  \textbf{80} (3), 235--238 (2003).

\bibitem[Kn3]{kn5} A.~L.~Knutsen, \textit{Higher order birational embeddings of Del Pezzo
surfaces}, Mathematische Nachrichten {\bf 254-255}, 1, 183-196 (2003).

\bibitem[Kn4]{kn1} A.~L.~Knutsen, \textit{On $k$th order embeddings of $K3$ 
surfaces and Enriques surfaces}, Manuscripta Math. \textbf{104}, 211--237 (2001).

\bibitem[Kn5]{kn2} A.~L.~Knutsen, \textit{Smooth curves on projective $K3$ 
surfaces}, Math. Scand. \textbf{90}, No. 2, 215--231 (2002).

\bibitem[La1]{lazsamp} R.~Lazarsfeld, \textit{A sampling of vector bundle techniques in the study of linear series}, Sugaku \textbf{47}, 125-144 (1995).



\bibitem[La2]{la} R.~Lazarsfeld, \textit{Brill-Noether-Petri without degenerations}, J. Diff. Geom. \textbf{23}, 299--307 (1986).

\bibitem[GMa]{ma} G.~Martens, \textit{On curves on $K3$ surfaces},
  Springer LNM 1398,  174--182 (1989). 

\bibitem[HMa]{Mar} H.H.~Martens, \textit{On the varieties of special
divisors on a curve, II}, J. reine und angew. Math. \textbf{233},  
89--100 (1968).

\bibitem[Mo]{Morrison} D.~R.~Morrison, \textit{On $K3$ surfaces with large 
Picard number}, Invent. Math. \textbf{75}, 105--121 (1984).

\bibitem[Mu1]{Mu1} S.~Mukai, \textit{Curves, $K3$ surfaces and Fano 3-folds of
Genus $ \leq 10$}, Algebraic Geometry and Commutative Algebra
(In Honor of M. Nagata) \textbf{1},  357--377 (1988). 

\bibitem[Mu2]{Mu2} S.~Mukai, \textit{New development of theory of Fano 3-folds:
vector bundle method and moduli problem}, Sugaku 
\textbf{47},  125--144 (1995).

\bibitem[Na]{N} M.~Nagata, \textit{Local rings},
  Interscience Publishers, Wiley (1989).

\bibitem[Ni]{nik} V.~V.~Nikulin, \textit{Integer symmetric bilinear forms and some of their geometric applications}, Izv. Akad. Nauk SSSR Ser. Mat. \textbf{43}, no. 1, 
111--177, 238 (1979). 

\bibitem[No]{noe} M.~Noether, \textit{{\"U}ber die invariante Darstellung algebraischer Funktionen einer Variablen}, Math. Ann. \textbf{88}, 163--284 (1880).

\bibitem[Pa]{par} G.~Pareschi, \textit{Exceptional linear systems on curves on Del Pezzo surfaces}, Math. Ann. \textbf{291}, 17-38 (1991).



\bibitem[Pe]{pet} K.~Petri, \textit{{\"U}ber die invariante Darstellung algebraischer Funktionen}, Math. Ann. \textbf{17},  243--289 (1923).


\bibitem[P-S]{PS} R.~Piene, G. ~Sacchiero, \textit{Duality for rational
normal scrolls}, Comm. in Algebra \textbf{12} (9),
1041--1076 (1984).

\bibitem[R-S]{rs} K~.Ranestad, F.~O.~ Schreyer, \textit{Varieties of sums of powers}, 
J. Reine Angew. Math. \textbf{525}, 147--181 (2000).

\bibitem[R-R-W]{RRW} M.~S.~Ravi, J.~Rosenthal, X.~Wang, \textit{Dynamic pole
assignment and Schubert calculus}, SIAM J. of Control and Optimization, 
\textbf{34} (3), 813--832 (1996). 

\bibitem[Re1]{re} M.~Reid, \textit{``Bogomolov's theorem ${c_1}^2 \leq 4c_2$''},
Internat. Sympos. on Algebraic Geometry, Kyoto, 623-642 (1977). 



\bibitem[Re2]{Re} M.~Reid, \textit{Chapters on Algebraic Surfaces}, In: Complex
algebraic geometry (Park City 1993), IAS/Park City Math., Ser. \textbf{3}, 
3--159 (1997).

\bibitem[Re3]{reid} M.~Reid, \textit{Special linear systems on curves lying on a K3 surface}, J. London Math. Soc. (2) \textbf{13}, no. 3, 454-458 (1976).

\bibitem[Rdr]{reider} I.~Reider, \textit{Vector bundles of rank $2$ and linear systems on algebraic surfaces}, Annals of Math. \textbf{127}, 309--316 (1988). 

\bibitem[SD]{S-D} B.~Saint-Donat, \textit{Projective models of $K3$
  Surfaces}, Amer. J. Math. \textbf{96}, 602--639 (1974).

\bibitem[Sc]{Sc} F.~O.~ Schreyer, \textit{Syzygies of canonical curves and 
special linear series}, Math. Ann. \textbf{275}, 105--137 (1986).

\bibitem[Se1]{Se1} C.~Segre, \textit{Sulle rigate razionali in uno spazio 
lineare qualunque}, Atti della R. Acc. delle scienze di Torino, \textbf{XIX}, 
265--282 (1883-4).

\bibitem[Se2]{Se2} C.~Segre, \textit{Sulle varieta a tre dimensioni composte di serie semplici razionali di piani},  Atti della R. Acc. delle scienze di Torino \textbf{XXI}, 95--115 (1885-6).

\bibitem[Si]{Sim} A.~Simis, \textit{Multiplicities and Betti numbers of
homogeneous ideals}, In: Space Curves (Rocca di Papa 1985), Springer Lecture 
Notes in Mathematics, \textbf{1266}, 232--250 (1987).

\bibitem[Ste]{St} J.~Stevens, \textit{Rolling factors deformations and
extensions of canonical curves}, Doc. Math. \textbf{6}, 185--226 (2001).

\bibitem[Str]{Str} S.~A.~ Str{\o}mme, \textit{On parametrized rational curves 
in Grassmann varieties},  In: Space Curves (Rocca di Papa 1985), Springer 
LNM 1266, 251--272 (1987).

\bibitem[Ty]{tyu} A.~N.~Tyurin, \textit{Cycles, curves and vector bundles on an algebraic surface}, Duke Math. J. \textbf{54}, no. 1, 1--26 (1987). 


\bibitem[Vo1]{voi1} C.~Voisin, \textit{Green's canonical syzygy conjecture for generic curves of odd genus}, math.AG/0301359 (2003). 

\bibitem[Vo2]{voi2} C.~Voisin, \textit{Green's generic syzygy conjecture for curves of even genus lying on a $K3$ surface}, J. Eur. Math. Soc. (JEMS) \textbf{4}, no. 4, 363--404 (2002). 

\bibitem[Wa]{W} C.~Walter, \textit{Pfaffian Subschemes}, Jour. of Alg. Geom.
\textbf{5}, 671--704 (1996).

\bibitem[We]{weil} A.~Weil, \textit{\OE vres Scientifiques (Collected
    Papers)}, II,  Springer, Heidelberg (1980).

%
% Monographs
%\bibitem[KLR73]{monograph} Kagan, A.M., Linnik, Y.V., Rao, C.R.:
%Characterization Problems in Mathematical Statistics. Wiley, New York (1973)
%
% Contributed Works
%\bibitem[Mey89]{contribution} Meyer, P.A.: A short presentation of
%stochastic calculus. In: Emery, M. (ed) Stochastic Calculus in
%Manifolds. Springer, Berlin Heidelberg New York (1989)
%
% Journal
%\bibitem[MR97]{journal} Miller, B.M., Runggaldier, W.J.: Kalman
%filtering for linear systems with coefficients driven by a hidden Markov
%jump process. Syst. Control Lett., \textbf{31}, 93--102 (1997)
%
% Theses
%\bibitem[Ros77]{thesis} Ross, D.W.: Lysosomes and storage diseases. MA
%Thesis, Columbia University, New York (1977)

\end{thebibliography}
\end{document}